\documentclass[a4paper,12pt,amsart,frenchb]{article}

\usepackage{amsmath,amsbsy,amsfonts,amssymb}
\usepackage[french]{babel}
\newcommand{\ass}[2]{\vskip0.3cm\noindent
{\bf {#1}}. { \sl {#2}}\vskip0.3cm\noindent
}

\begin{document}

\title{  Une formule int\'egrale reli\'ee \`a la conjecture locale de Gross-Prasad, $2^{\grave{e}me}$ partie: extension aux repr\'esentations temp\'er\'ees}
\author{J.-L. Waldspurger}
\date{1er avril 2009}
\maketitle

\bigskip

{\bf Introduction}

\bigskip

Cet article fait suite \`a [W1]. Rappelons les d\'efinitions des principaux objets. Soit $F$ un corps local non archim\'edien de caract\'eristique nulle. Soit $(V,q_{V})$ un espace quadratique, c'est-\`a-dire que $V$ est un espace vectoriel de dimension finie sur $F$ et $q_{V}$ est une forme quadratique non d\'eg\'en\'er\'ee sur $V$. Soit $(W,q_{W})$ un autre espace quadratique. On suppose que  l'on a une d\'ecomposition orthogonale $V=W\oplus D_{0}\oplus Z$, o\`u $D_{0}$ est une droite et $Z$ est muni d'une base $\{v_{i}; i=\pm 1,...,\pm r\}$ telle que $q_{V}(v_{i},v_{j})=\delta_{i,-j}$ pour tous $i,j$. On note $G$ et $H$ les groupes sp\'eciaux orthogonaux de $V$ et $W$. Le groupe $H$ se plonge naturellement dans $G$. Introduisons le sous-groupe parabolique de $G$ form\'e des \'el\'ements qui conservent le drapeau
$$Fv_{r}\subset Fv_{r}\oplus Fv_{r-1}\subset....\subset Fv_{r}\oplus...\oplus Fv_{1}.$$
Notons $U$ son radical unipotent. Fixons un \'el\'ement non nul $v_{0}\in D_{0}$ et un caract\`ere continu non trivial $\psi$ de $F$. On d\'efinit un caract\`ere $\xi$ de $U(F)$ par l'\'egalit\'e
$$\xi(u)=\psi(\sum_{i=0,...,r-1}q_{V}(uv_{i},v_{-i-1})).$$
Soient $\pi$, resp. $\rho$, une repr\'esentation admissible irr\'eductible de $G(F)$, resp. $H(F)$, dans un espace complexe $E_{\pi}$, resp. $E_{\rho}$. On note $Hom_{H,\xi}(\pi,\rho)$ l'espace des applications lin\'eaires $\varphi:E_{\pi}\to E_{\rho}$ telles que
$$\varphi(\pi(hu)e)=\xi(u)\rho(h)\varphi(e)$$
pour tous $u\in U(F)$, $h\in H(F)$, $e\in E_{\pi}$. On note $m(\rho,\pi)$ la dimension de cet espace. D'apr\`es [AGRS] th\'eor\`eme 1' et [GGP] corollaire 20.4, ce nombre vaut $0$ ou $1$. Il est ind\'ependant des divers choix effectu\'es.

Supposons $G$ et $H$ quasi-d\'eploy\'es sur $F$ et affectons les notations d'un indice $i$: $V_{i}$, $G_{i}$ etc... Supposons pour cette introduction $dim(W_{i})\geq3$. A \'equivalence pr\`es, il y a un unique espace quadratique que nous notons $(V_{a},q_{V_{a}})$ tel que $dim(V_{a})=dim(V_{i})$, que les discriminants de $q_{V_{i}}$ et $q_{V_{a}}$ soient \'egaux mais leurs indices de Witt soient distincts. On introduit de m\^eme un espace quadratique $(W_{a},q_{W_{a}})$. Le couple $(V_{a},W_{a})$ v\'erifie les m\^emes propri\'et\'es que ci-dessus. Les groupes sp\'eciaux orthogonaux $G_{a}$, resp. $H_{a}$, de $V_{a}$, resp. $W_{a}$, sont des formes int\'erieures de $G_{i}$, resp. $H_{i}$. On note $Temp(G_{i})$, $Temp(G_{a})$ etc... les ensembles de repr\'esentations temp\'er\'ees et irr\'eductibles de $G_{i}(F)$, $G_{a}(F)$ etc... On admet que ces ensembles sont unions disjointes de $L$-paquets v\'erifiant certaines propri\'et\'es encore conjecturales. Pr\'ecis\'ement on admet les propri\'et\'es (1), (2) et (3) de [W1] 13.2. Soient $\Pi_{i}$, resp. $\Sigma_{i}$, un $L$-paquet dans $Temp(G_{i})$, resp. $Temp(H_{i})$. Il peut correspondre \`a $\Pi_{i}$ un $L$-paquet dans $Temp(G_{a})$, que l'on note $\Pi_{a}$. Ou bien il n'y a pas de tel $L$-paquet et on pose $\Pi_{a}=\emptyset$. On d\'efinit de m\^eme $\Sigma_{a}$.
La multiplicit\'e $m(\rho,\pi)$ est bien d\'efinie pour tout $(\rho,\pi)\in (\Sigma_{i}\times \Pi_{i})\cup(\Sigma_{a}\times \Pi_{a})$.

\ass{Th\'eor\`eme}{Sous ces hypoth\`eses, il existe un unique couple $(\rho,\pi)\in (\Sigma_{i}\times \Pi_{i})\cup(\Sigma_{a}\times \Pi_{a})$ tel que $m(\rho,\pi)=1$. }

C'est une partie de la conjecture 6.9 de [GP]. Ce th\'eor\`eme r\'esulte ais\'ement d'une formule qui calcule $m(\rho,\pi)$ comme une somme d'int\'egrales de fonctions qui se d\'eduisent des caract\`eres de $\rho$ et $\pi$. Plus pr\'ecis\'ement, revenons aux notations sans indices du d\'ebut de cette introduction. Soient $\pi$ et $\rho$ des repr\'esentations admissibles irr\'eductibles de $G(F)$ et $H(F)$. On introduit une expression $m_{geom}(\rho,\pi)$ pour la d\'efinition de laquelle on renvoie \`a l'introduction de [W1].

\ass{Th\'eor\`eme}{Supposons $\pi$ et $\rho$ temp\'er\'ees et irr\'eductibles. Alors on a l'\'egalit\'e $m(\rho,\pi)=m_{geom}(\rho,\pi)$.}

Dans [W1], on avait d\'emontr\'e cette \'egalit\'e sous les hypoth\`eses que $\pi$ \'etait cuspidale et $\rho$ admissible. Ici, on \'elargit l'hypoth\`ese sur $\pi$ qui n'est plus que temp\'er\'ee. Par contre, on impose une hypoth\`ese plus forte \`a $\rho$ qui est elle-aussi temp\'er\'ee. Comme dans [W1], le second th\'eor\`eme implique le premier. Evidemment, dans [W1], l'hypoth\`ese de cuspidalit\'e pr\'esente dans le second th\'eor\`eme se retrouvait dans le premier. C'est cette hypoth\`ese que nous faisons dispara\^{\i}tre dans le pr\'esent article.

D\'ecrivons l'id\'ee principale de la preuve du second th\'eor\`eme. Rappelons que, pour une fonction $f\in C_{c}^{\infty}(G(F))$, on dit que $f$ est tr\`es cuspidale si et seulement si, pour tout sous-groupe parabolique propre $P=MU$ de $G$ (avec une notation famili\`ere) et pour tout $m\in M(F)$, on a l'\'egalit\'e
$$\int_{U(F)}f(mu)du=0.$$
Soient $\rho\in Temp(H)$ et $f$ une fonction tr\`es cuspidale sur $G(F)$. On note $\theta_{\rho}$ le caract\`ere de $\rho$. Pour tout $N\in {\mathbb N}$, on introduit une fonction $\kappa_{N}$ sur $G(F)$ qui est la fonction caract\'eristique de l'image r\'eciproque d'un sous-ensemble compact de $H(F)U(F)\backslash G(F)$ qui est de plus en plus grand quand $N$ tend vers l'infini. Posons
$$I_{N}(\theta_{\rho},f)=\int_{H(F)U(F)\backslash G(F)}\int_{H(F)}\int_{U(F)}\theta_{\rho}(h)f(g^{-1}hug)\xi(u)\kappa_{N}(g)du\,dh\,dg.$$
On montre que, quand $N$ tend vers l'infini, cette expression a une limite. En fait, et c'est cela qui est fructueux, il y a deux fa\c{c}ons de calculer la limite. L'une, que l'on peut qualifier de g\'eom\'etrique, a \'et\'e d\'evelopp\'ee en [W1], et conduit \`a une \'egalit\'e
$$lim_{N\to \infty}I_{N}(\theta_{\rho},f)=I_{geom}(\theta_{\rho},f),$$
o\`u le membre de droite est une somme d'int\'egrales sur certains sous-tores de $H(F)$. Dans le pr\'esent article, on calcule la limite d'une autre fa\c{c}on, que l'on peut qualifier de spectrale. On obtient une \'egalit\'e (cf. th\'eor\`eme 6.1):
$$lim_{N\to \infty}I_{N}(\theta_{\rho},f)=I_{spec}(\theta_{\rho},f),$$
o\`u
$$I_{spec}(\theta_{\rho},f)=\sum_{L\in {\cal L}(M_{min})}\vert W^L\vert \vert W^G\vert ^{-1}(-1)^{a_{L}}\sum_{{\cal O}\in \{\Pi_{ell}(L)\}; m({\cal O},\rho)=1}$$
$$[i{\cal A}^{\vee}_{{\cal O}}:i{\cal A}_{L,F}^{\vee}]^{-1}t(\pi)^{-1}\int_{i{\cal A}^*_{L,F}}J_{L}^G(\pi_{\lambda},f)d\lambda.$$
Tous les termes de cette formule seront d\'efinis dans l'article. Disons simplement ici que, dans le cas o\`u $L=G$, les objets ${\cal O}$ sont simplement les repr\'esentations irr\'eductibles temp\'er\'ees et elliptiques de $G(F)$ et, si l'on pose plus simplement $\pi={\cal O}$, la condition $m({\cal O},\rho)=1$ n'est autre que $m(\rho,\pi)=1$ tandis que $J_{G}^G(\pi,f)=\theta_{\pi}(f)$. Pour $L$ quelconque, $J_{L}^G(\pi_{\lambda},f)$ est la valeur sur $f$ du caract\`ere pond\'er\'e associ\'e \`a $\pi_{\lambda}$.

On a donc l'\'egalit\'e
$$I_{geom}(\theta_{\rho},f)=I_{spec}(\theta_{\rho},f)$$
qui, bien s\^ur, rappelle fortement la formule des traces locale d'Arthur. De fait, la preuve reprend tr\`es largement celle de [A3].  Dans les deux membres de la formule apparaissent des distributions qui ne sont pas invariantes: int\'egrales orbitales pond\'er\'ees et caract\`eres pond\'er\'es. Le proc\'ed\'e mis au point par Arthur, appliqu\'e en particulier dans [A5] \`a la formule des traces locale, permet de transformer la formule ci-dessus en une autre o\`u n'apparaissent que des distributions invariantes. Le terme de droite de cette formule continue de distinguer les repr\'esentations $\pi$ de $G(F)$ telles que $m(\rho,\pi)=1$. Le second th\'eor\`eme r\'esulte facilement de cette formule "invariante".

Expliquons encore deux points. Dans la formule non invariante, la fonction $f$ est suppos\'ee tr\`es cuspidale, ce qui est assez restrictif. Cela parce que nous ne savons pas calculer la limite de $I_{N}(\theta_{\rho},f)$ pour une fonction qui ne v\'erifie pas cette hypoth\`ese (le r\'esultat rend d'ailleurs douteuse la possibilit\'e d'\'etendre nos calculs \`a des fonctions ne v\'erifiant pas cette hypoth\`ese). Mais, une fois la formule rendue invariante, on peut supposer $f$ seulement cuspidale (c'est-\`a-dire les int\'egrales orbitales $J^G(x,f)$ sont nulles pour tout \'el\'ement $x\in G(F)$ qui est semi-simple, fortement r\'egulier et non elliptique). Cela r\'esulte du lemme suivant (lemme 2.7).

\ass{Lemme}{Soit $f\in C_{c}^{\infty}(G(F))$ une fonction cuspidale. Alors il existe une fonction tr\`es cuspidale $f'\in C_{c}^{\infty}(G(F))$ telle que $D(f)=D(f')$ pour toute distribution $D$ sur $G(F)$ invariante par conjugaison.}

Cet affaiblissement de la condition sur $f$ est n\'ecessaire pour achever la preuve (on prend pour $f$ un pseudo-coefficient d'une repr\'esentation temp\'er\'ee et elliptique).

Le deuxi\`eme point est l'apparition de la condition $m(\rho,\pi)=1$ dans le terme $I_{spec}(\theta_{\rho},f)$. Fixons ici une repr\'esentation $\pi\in Temp(G)$. L'espace $Hom_{H,\xi}(\pi,\rho)$ dont $m(\rho,\pi)$ est la dimension est d\'efini de fa\c{c}on abstraite. Il ne peut pas intervenir directement dans $I_{spec}(\theta_{\rho},f)$ qui est une int\'egrale explicite. Ce qui intervient dans ce terme, c'est la forme sesquilin\'eaire ${\cal L}_{\pi,\rho}$ sur $E_{\rho}\otimes_{{\mathbb C}}E_{\pi}$ d\'efinie par
$${\cal L}_{\pi,\rho}(\epsilon'\otimes e',\epsilon\otimes e)=\int_{H(F)U(F)}(\rho(h)\epsilon',\epsilon)(e',\pi(hu)e)\bar{\xi}(u)du\,dh,$$
pour $\epsilon,\epsilon'\in E_{\rho}$ et $e,e'\in E_{\pi}$ (les produits $(.,)$ sont des produits hermitiens invariants sur $E_{\rho}$ et $E_{\pi}$). L'int\'egrale ci-dessus n'est pas absolument convergente, mais on peut la d\'efinir comme une limite d'int\'egrales absolument convergentes, cf. 5.1. N\'egligeons cette question de convergence.   Fixons $\epsilon$ et $e'$. D\'efinissons une application $l:E_{\pi}\to E_{\rho}$ par l'\'egalit\'e $(\epsilon',l(e))={\cal L}_{\pi,\rho}(\epsilon'\otimes e',\epsilon\otimes e)$ pour tous $\epsilon'\in E_{\rho}$ et $e\in E_{\pi}$. On v\'erifie que $l\in Hom_{H,\xi}(\pi,\rho)$. Si ${\cal L}_{\pi,\rho}$ n'est pas nulle, cet espace $Hom_{H,\xi}(\pi,\rho)$ ne l'est pas non plus et $m(\rho,\pi)=1$. On a besoin de la r\'eciproque, qui s'av\`ere vraie.

\ass{Proposition}{Soient $\pi\in Temp(G)$ et $\rho\in Temp(H)$. Alors $m(\rho,\pi)=1$ si et seulement si la forme sesquilin\'eaire ${\cal L}_{\pi,\rho}$ est non nulle.}

Cf. proposition 5.7. Signalons que cette fa\c{c}on concr\`ete de construire l'espace $Hom_{H,\xi}(\pi,\rho)$ se trouve d\'ej\`a dans l'article [II] de Ikeda et Ichino.

La premi\`ere section est consacr\'ee aux notations et \`a des rappels sur les op\'erateurs d'entrelacement et la formule de Plancherel. La deuxi\`eme l'est aux propri\'et\'es des fonctions cuspidales ou tr\`es cuspidales et aux quasi-caract\`eres qu'elles permettent de d\'efinir. Les sections 3 et 4 sont franchement p\'enibles. On y d\'emontre diverses majorations n\'ecessaires pour la suite (pour le groupe $GL_{k}$ dans la section 3, pour un groupe sp\'ecial orthogonal dans la section 4). On s'inspire ici plus que largement des travaux d'Harish-Chandra. Signalons \`a ce propos que l'on fait constamment r\'ef\'erence \`a l'article [W2]. Mais l'apparence est trompeuse puisque dans [W2], on s'\'etait content\'e de r\'ediger des r\'esultats non publi\'es d'Harish-Chandra. D'autre part, dans [W2], on avait cru judicieux de modifier la d\'efinition de l'homomorphisme habituel $H_{G}$ en y glissant un signe $-$. On persiste \`a penser que, sur un corps de base $p$-adique, c'est une meilleure d\'efinition. Mais, pour utiliser les r\'esultats d'Arthur, il vaut mieux reprendre ses d\'efinitions. C'est ce que l'on fait, mais cela induit des changements de signe dans les r\'ef\'erences que l'on fera \`a [W2]: cela \'echange une chambre positive avec son oppos\'ee. La section 5 est consacr\'ee \`a la d\'efinition et l'\'etude des formes sesquilin\'eaires ${\cal L}_{\pi,\rho}$ \'evoqu\'ees ci-dessus. La preuve de l'\'egalit\'e $lim_{N\to\infty}I_{N}(\theta_{\rho},f)=I_{spec}(\theta_{\rho},f)$ se trouve dans la section 6. Il s'agit pour l'essentiel de recopier [A3]. On en d\'eduit dans la section 7 les deux th\'eor\`emes \'enonc\'es ci-dessus.

 \bigskip

\section{Notations et rappels}

\bigskip

\subsection{Notations g\'en\'erales}

On utilise les notations introduites dans [W1], qui sont la plupart du temps celles d'Arthur et d'Harish-Chandra. Soit $F$ un corps local non archim\'edien de caract\'eristique nulle. On note $\mathfrak{o}_{F}$ son anneau d'entiers, $\mathfrak{p}_{F}$ l'id\'eal maximal de $\mathfrak{o}_{F}$, $q$ le nombre d'\'el\'ements du corps r\'esiduel, $val_{F}$ et $\vert .\vert _{F}$ les valuation et valeur absolue usuelles et on fixe une uniformisante $\varpi_{F}$. Soit $G$ un groupe r\'eductif connexe d\'efini sur $F$. On note $\mathfrak{g}$ l'alg\`ebre de Lie de $G$. On note $A_{G}$ le plus grand tore d\'eploy\'e central dans $G$, $X(G)$ le groupe des caract\`eres de $G$ d\'efinis sur $F$, ${\cal A}_{G}=Hom(X(G),{\mathbb R})$ et ${\cal A}_{G}^*=X(G)\otimes_{{\mathbb Z}}{\mathbb R}$ le dual de ${\cal A}_{G}$. On d\'efinit l'homomorphisme habituel $H_{G}:G(F)\to {\cal A}_{G}$. On note ${\cal A}_{G,F}$, resp. $\tilde{{\cal A}}_{G,F}$, l'image de $G(F)$, resp. $A_{G}(F)$, par cet homomorphisme. On note ${\cal A}_{G,F}^{\vee}$, resp. $\tilde{{\cal A}}_{G,F}^{\vee}$, le sous-groupe des $\lambda\in {\cal A}_{G}^*$ tels que $\lambda(\zeta)\in 2\pi{\mathbb Z}$ pour tout $\zeta\in {\cal A}_{G,F}$, resp. $\zeta\in \tilde{{\cal A}}_{G,F}$. On note $a_{G}$ la dimension de ${\cal A}_{G}$. 

Soit $K$ un sous-groupe compact sp\'ecial de $G(F)$. Soit $P=MU$ un sous-groupe parabolique de $G$. Rappelons nos conventions: $P$ est implicitement suppos\'e d\'efini sur $F$ et la notation $P=MU$ signifie que $M$ est une composante de L\'evi de $P$, d\'efinie sur $F$, et $U$ est le radical unipotent de $P$. Supposons que $K$ soit en bonne position relativement \`a $M$. Pr\'ecis\'ement, il existe un sous-tore d\'eploy\'e maximal $A_{0}$ de $M$ tel que $K$ fixe un point de l'appartement associ\'e \`a $A_{0}$ dans l'immeuble de $G$. On a l'\'egalit\'e $G(F)=M(F)U(F)K$. Pour tout $g\in G(F)$, on fixe des \'el\'ements $m_{P}(g)\in M(F)$, $u_{P}(g)\in U(F)$, $k_{P}(g)\in K$ tels que $g=m_{P}(g)u_{P}(g)k_{P}(g)$. On  prolonge l'application  $H_{M}:M(F)\to {\cal A}_{M}$ en une fonction $H_{P}:G(F)\to {\cal A}_{M}$ par $H_{P}(g)=H_{M}(m_{P}(g))$.

Supposons fix\'e un L\'evi minimal $M_{min}$ de $G$. On pose $W^G=Norm_{G(F)}(M_{min})/M_{min}(F)$.

On note $\Xi^G$ la fonction d'Harish-Chandra ([W2] II.1). Elle d\'epend de $K$. Mais elle ne nous sert qu'\`a r\'esoudre des questions de majorations. Or changer de groupe $K$ remplace $\Xi^G$ par une fonction \'equivalente. Il est donc loisible d'utiliser cette fonction sans pr\'eciser le groupe $K$ qui permet de la d\'efinir. On utilise aussi la fonction $\sigma$. Rappelons que, dans [W1], on a l\'eg\`erement modifi\'e la d\'efinition d'Harish-Chandra en posant $\sigma(g)=sup(1,log(\vert \vert g\vert \vert ))$.  On a la relation $\sigma(gg')\leq \sigma(g)+\sigma(g')\leq 2\sigma(g)\sigma(g')$ pour tous $g,g'\in G(F)$. Pour tout r\'eel $b\geq0$, on note ${\bf 1}_{\sigma< b}$, resp. ${\bf 1}_{\sigma\geq b}$, la fonction caract\'eristique de l'ensemble des $g\in G(F)$ tels que $\sigma(g)< b$, resp. $\sigma(g)\geq b$.

Quand deux nombres r\'eels positifs ou nuls $a$ et $b$ d\'ependent d'un certain nombre de variables $x_{1},...,x_{n}$, on dit que $a$ est essentiellement major\'e par $b$, ce que l'on note $a<<b$, s'il existe un r\'eel $c>0$ tel que $a\leq cb$ pour tous $x_{1},...,x_{n}$. Cette notation est quelque peu impr\'ecise mais nous \'evite d'introduire une kyrielle de constantes superflues.

On introduit l'espace ${\cal S}(G(F))$ des fonctions de Schwartz-Harish-Chandra sur $G(F)$. C'est l'ensemble des fonctions $f:G(F)\to {\mathbb C}$ qui sont biinvariantes par un sous-groupe ouvert compact et telles que, pour tout r\'eel $R\geq0$, on ait une majoration
$$\vert  f(g)\vert << \Xi^G(g)\sigma(g)^{-R}$$
pour tout $g\in G(F)$.  L'espace ${\cal S}(G(F))$ contient l'espace $C_{c}^{\infty}(G(F))$ des fonctions localement constantes \`a support compact.

Soit $\pi$ une repr\'esentation admissible de $G(F)$. On note sans plus de commentaire $E_{\pi}$ un espace complexe dans lequel elle se r\'ealise. Si $K$ est un sous-groupe de $G(F)$, on note $E_{\pi}^K$ le sous-espace des \'el\'ements de $E_{\pi}$ invariants par $K$. Supposons $\pi$ unitaire. On fixe  une forme hermitienne d\'efinie positive $(.,.)$ sur $E_{\pi}$ invariante par l'action de $G(F)$. On appelle une telle forme un produit scalaire invariant. Pr\'ecisons notre convention sur les formes sesquilin\'eaires: la forme $(.,.)$ v\'erifie la relation $(\lambda'e',\lambda e)=\bar{\lambda}'\lambda(e',e)$ pour tous $\lambda,\lambda'\in {\mathbb C}$, $e,e'\in E_{\pi}$.  Nous dirons que $\pi$ est temp\'er\'ee si elle est unitaire, de longueur finie, et  qu'il existe un entier $D$ tel que, pour tous $e,e'\in E_{\pi}$, on ait une majoration
$$\vert (e',\pi(g)e)\vert <<\Xi^G(g)\sigma(g)^D.$$
En fait, l'entier $D$ ne sert \`a rien: on peut prendre $D=0$, cf. [W2] lemme VI.2.2. Supposons $G(F)$ muni d'une mesure de Haar. Si $\pi$ est temp\'er\'ee, l'action de $C_{c}^{\infty}(G(F))$ dans $E_{\pi}$ se prolonge en une action de ${\cal S}(G(F))$. Pour $f\in {\cal S}(G(F))$ et $e,e'\in E_{\pi}$, on a l'\'egalit\'e
$$(e',\pi(f)e)=\int_{G(F)}f(g)(e',\pi(g)e)dg.$$
Cette int\'egrale est absolument convergente. On note $Temp(G)$ l'ensemble des classes d'\'equivalence de repr\'esentations temp\'er\'ees irr\'eductibles de $G(F)$. 

On fixe  un caract\`ere $\psi:F\to {\mathbb C}^{\times}$, continu et non trivial. On note $c_{\psi}$ le plus petit entier relatif $c$ tel que $\psi$ soit trivial sur $ \mathfrak{p}_{F}^c$.

 \bigskip

\subsection{Mesures}

Dans la  suite de l'article, la situation sera la suivante. Le groupe $G$ est fix\'e,  ainsi qu'un L\'evi minimal $M_{min}$ de $G$ et un sous-groupe compact sp\'ecial $K$ de $G(F)$, en bonne position relativement \`a $M_{min}$.   

On munit $K$ de la mesure de Haar de masse totale $1$. On munit $G(F)$ d'une mesure de Haar (pour laquelle $mes(K)$ n'est pas forc\'ement \'egale \`a $1$). Soit $P=MU\in {\cal F}(M_{min})$ (les notations ${\cal F}(L)$, ${\cal P}(L)$, ${\cal L}(L)$ sont celles d'Arthur, cf. [W1] 1.1). On munit $U(F)$ de l'unique mesure de Haar telle que
$$\int_{U(F)}\delta_{\bar{P}}(m_{\bar{P}}(u))du=1,$$
o\`u $\bar{P}=M\bar{U}$ est le parabolique oppos\'e \`a $P$ et $\delta_{\bar{P}}$ est le module usuel.
On munit $M(F)$ de l'unique mesure de Haar telle que, pour toute $f\in C_{c}^{\infty}(G(F))$, on ait l'\'egalit\'e
$$\int_{G(F)}f(g)dg=\int_{K}\int_{U(F)}\int_{M(F)}f(muk)dm\,du\,dk.$$
Le point est que cette mesure sur $M(F)$ ne d\'epend pas du sous-groupe parabolique $P\in {\cal P}(M)$ utilis\'e pour la d\'efinir ([A3] 1.2).  

 On munit l'espace $i{\cal A}^*_{M}\subset {\cal A}^*_{M}\otimes_{{\mathbb R}}{\mathbb C}$ de la mesure de Haar telle que le quotient $i{\cal A}^*_{M}/i\tilde{{\cal A}}^{\vee}_{M,F}$ soit de mesure $1$. On pose $i{\cal A}_{M,F}^*=i{\cal A}_{M}/i{\cal A}_{M,F}^{\vee}$ et on le munit de la mesure telle que l'application naturelle de $i{\cal A}_{M}^*$ dans $i{\cal A}_{M,F}^*$ pr\'eserve localement les mesures.

 Soit $T$ un tore. Si $T$ est d\'eploy\'e, on  munit $T(F)$ de la mesure de Haar telle que le sous-groupe compact maximal de $T(F)$ soit de mesure $1$. En g\'en\'eral, on munit $A_{T}(F)$ de la mesure que l'on vient de d\'efinir et $T(F)$ de la mesure telle que $T(F)/A_{T}(F)$ soit de mesure $1$ pour la mesure quotient. 
 
 \ass{ Remarques}{ 1. Dans le cas o\`u $M_{min}$ est un tore, les d\'efinitions pr\'ec\'edentes peuvent entrer en conflit. On croit qu'en pratique, il n'y aura pas d'ambig\" uit\'e.
 
 2. Dans les sections 3, 4 et 5, on se pr\'eoccupera de questions de convergence pour lesquelles les choix de mesures sont sans importance. On ne tiendra pas compte des normalisations ci-dessus. On supposera au contraire que les mesures sont choisies de telle sorte que toutes les constantes qui apparaissent \`a cause d'elles dans les calculs soient \'egales \`a $1$.}

\bigskip

\subsection{Repr\'esentations induites, op\'erateurs d'entrelacement}

  Soit $P=MU$ un sous-groupe parabolique de $G$ et $\tau$ une repr\'esentation admissible de $M(F)$. On d\'efinit la repr\'esentation induite $Ind_{P}^G(\tau)$.  On note $E^G_{P,\tau}$ son espace. C'est celui des fonctions $e:G(F)\to E_{\tau}$ qui sont invariantes \`a droite par un sous-groupe ouvert compact de $G(F)$ et v\'erifient
$$e(mug)=\delta_{P}(m)^{1/2}\tau(m)e(g)$$
pour tous $m\in M(F)$, $u\in U(F)$ et $g\in G(F)$. Pour $g\in G(F)$, ou $f\in C_{c}^{\infty}(G(F))$, on note $Ind_{P}^G(\tau,g)$, ou $Ind_{P}^G(\tau,f)$, l'action de $g$, ou $f$, dans $E^G_{P,\tau}$.    Pour $\lambda\in {\cal A}_{M}^*\otimes_{{\mathbb R}}{\mathbb C}$, on d\'efinit la repr\'esentation $\tau_{\lambda}$ de $M(F)$ par $\tau_{\lambda}(m)=exp(\lambda(H_{M}(m)))\tau(m)$ et la repr\'esentation induite $Ind_{P}^G(\tau_{\lambda})$. Remarquons que ces repr\'esentations ne d\'ependent que de l'image de $\lambda$ dans  $({\cal A}^*_{M}\otimes_{{\mathbb R}}{\mathbb C})/i{\cal A}^{\vee}_{M,F}$.  Supposons $M_{min}\subset M$. Notons ${\cal K}^G_{P,\tau}$ l'espace des fonctions $e:K\to E_{\tau}$  qui sont invariantes \`a droite par un sous-groupe ouvert compact de $K$ et v\'erifient la m\^eme relation que ci-dessus, pour $m\in K\cap M(F)$, $u\in K\cap U(F)$ et $g\in K$. Par restriction \`a $K$, $ E_{P,\tau_{\lambda}}^G$ s'identifie \`a ${\cal K}^G_{P,\tau}$, ce dernier espace est donc un mod\`ele commun \`a toutes les repr\'esentations $Ind_{P}^G(\tau_{\lambda})$. Supposons $\tau$ unitaire. On d\'efinit un produit hermitien sur ${\cal K}^G_{P,\tau}$ par
$$(e',e)=\int_{K}(e'(k),e(k))dk.$$
C'est un produit scalaire invariant pour la repr\'esentation $Ind_{P}^G(\tau_{\lambda})$ pour tout $\lambda\in i{\cal A}^*_{M,F}$.

  Laissons $M$ fix\'e mais faisons varier $P$ parmi les \'el\'ements de ${\cal P}(M)$.   Pour $P=MU,P'=MU'\in {\cal P}(M)$ et $\lambda\in {\cal A}^*_{M}\otimes_{{\mathbb R}}{\mathbb C}$, on d\'efinit l'op\'erateur d'entrelacement
$$J_{P'\vert P}(\tau_{\lambda}):E_{P,\tau_{\lambda}}^G\to E^G_{P',\tau_{\lambda}}$$
Quand la partie r\'eelle de $\lambda$ est dans un certain c\^one, il est d\'efini par la formule
$$(J_{P'\vert P}(\tau_{\lambda})e)(g)=\int_{(U(F)\cap U'(F))\backslash U'(F)}e(ug)du.$$
En g\'en\'eral, il est d\'efini par prolongement m\'eromorphe ( il est m\^eme rationnel, si l'on consid\`ere $({\cal A}^*_{M}\otimes_{{\mathbb R}}{\mathbb C})/i{\cal A}^{\vee}_{M,F}$ comme un tore alg\'ebrique complexe). Par restriction \`a $K$, on peut consid\'erer $J_{P'\vert P}(\tau_{\lambda})$ comme un homomorphisme de ${\cal K}^G_{P,\tau}$ dans ${\cal K}^G_{P',\tau}$. C'est ce point de vue que l'on adopte dans la suite. 

Supposons $\tau $ irr\'eductible. L'op\'erateur $J_{P\vert \bar{P}}(\tau_{\lambda})J_{\bar{P}\vert P}(\tau_{\lambda})$ est une homoth\'etie. Notons $j(\tau_{\lambda})$ le rapport d'homoth\'etie. Il ne d\'epend pas de $P$. On peut normaliser l'op\'erateur d'entrelacement. On introduit une fonction $r_{P'\vert P}(\tau_{\lambda})$ \`a valeurs complexes, qui est m\'eromorphe et m\^eme rationnelle, de sorte qu'en posant 
$$R_{P'\vert P}(\tau_{\lambda})=r_{P'\vert P}(\tau_{\lambda})^{-1}J_{P'\vert P}(\tau_{\lambda}),$$
cet op\'erateur v\'erifie les conditions du th\'eor\`eme 2.1 de [A4]. Les principales conditions sont

- pour $P,P',P''\in {\cal P}(M)$, $R_{P''\vert P'}(\tau_{\lambda})R_{P'\vert P}(\tau_{\lambda})=R_{P''\vert P}(\tau_{\lambda})$;

- supposons $\tau $ temp\'er\'ee;   pour $\lambda\in i{\cal A}^*_{M,F}$, $R_{P'\vert P}(\tau_{\lambda})$ est holomorphe et son adjoint pour le produit scalaire est $R_{P\vert P'}(\tau_{\lambda})$.

La d\'efinition des op\'erateurs normalis\'es s'\'etend au cas o\`u $\tau$ est semi-simple. En particulier, soit $P^M=M_{0}U^M_{0}$ un sous-groupe parabolique de $M$ tel que $M_{min}\subset M_{0}$, soit $\tau_{0}$ une repr\'esentation temp\'er\'ee irr\'eductible de $M_{0}(F)$, supposons que $\tau=Ind_{P^M}^M(\tau_{0})$. Pour $P=MU\in {\cal P}(M)$, introduisons le groupe $P_{0}=P^MU\in{\cal P}(M_{0})$. L'espace ${\cal K}_{P,\tau}^G$ s'identifie \`a ${\cal K}_{P_{0},\tau_{0}}^G$. Pour $P,P'\in {\cal P}(M)$ et $\lambda\in i{\cal A}_{M,F}^*$, l'op\'erateur $R_{P'\vert P}(\tau_{\lambda})$ s'identifie \`a $R_{P'_{0}\vert P_{0}}(\tau_{0,\lambda})$.

\bigskip

\subsection{Caract\`eres pond\'er\'es}

On conserve les donn\'ees $M$ et $\tau$ du paragraphe pr\'ec\'edent. On suppose que $\tau$ est   temp\'er\'ee, donc semi-simple d'apr\`es la d\'efinition que l'on a adopt\'ee. Pour tous $P,P'\in {\cal P}(M)$, l'op\'erateur $R_{P'\vert P}(\tau)$ est bien d\'efini et inversible. Plus g\'en\'eralement, pour tout $\lambda\in i{\cal A}_{M}^*$, l'op\'erateur $R_{P'\vert P}(\tau_{\lambda})$ est bien d\'efini et inversible. Fixons $P$. Pour tout $P'\in {\cal P}(M)$, consid\'erons la fonction ${\cal R}_{P'}(\tau)$ sur $i{\cal A}_{M}^*$ d\'efinie par
$$ {\cal R}_{P'}(\tau,\lambda)=R_{P'\vert P}(\tau)^{-1}R_{P'\vert P}(\tau_{\lambda}).$$
Elle prend ses valeurs dans l'espace d'endomorphismes de ${\cal K}_{P,\tau}^G$. La famille $({\cal R}_{P'}(\tau))_{P'\in {\cal P}(M)}$ est une $(G,M)$-famille \`a valeurs op\'erateurs ([A1] paragraphe 7).  Cela entra\^{\i}ne que la fonction 
$$\lambda\mapsto \sum_{P'\in {\cal P}(M)}{\cal R}_{P'}(\tau,\lambda)\theta_{P'}(\lambda)^{-1}$$
sur $i{\cal A}_{M}^*$ est $C^{\infty}$ (la fonction $\theta_{P'}$ est d\'efinie en [A1] p.15). On note ${\cal R}_{M}(\tau)$ la valeur de cette fonction en $\lambda=0$. C'est un endomorphisme de ${\cal K}_{P,\tau}$. Plus g\'en\'eralement, soient $\tilde{M}\in {\cal L}(M)$ et $Q=LU\in {\cal F}(\tilde{M})$. On d\'efinit une $(L,\tilde{M})$-famille $({\cal R}^Q_{\tilde{P}^L}(\tau))_{\tilde{P}^L\in {\cal P}^L(\tilde{M})}$ de la fa\c{c}on suivante: ${\cal R}^Q_{\tilde{P}^L}(\tau)$ est la restriction \`a $i{\cal A}_{\tilde{M}}^*$ de la fonction ${\cal R}_{P'}(\tau)$, o\`u $P'$ est un \'el\'ement quelconque de ${\cal P}(M)$ tel que $P'\subset Q$ et $P'\cap L\subset \tilde{P}^L$. Comme ci-dessus, on associe \`a cette $(L,\tilde{M})$-famille un op\'erateur ${\cal R}_{\tilde{M}}^Q(\tau)$.

Le caract\`ere pond\'er\'e de la repr\'esentation $\tau$ est la distribution $f\mapsto J_{M}^G(\tau,f) $d\'efinie par
$$J_{M}^G(\tau,f)=trace({\cal R}_{M}(\tau)Ind_{P}^G(f))$$
pour toute $f\in C_{c}^{\infty}(G(F))$. Cette distribution est d\'efinie \`a l'aide du sous-groupe parabolique $P$ que nous avons fix\'e, mais on montre qu'elle ne d\'epend pas de ce choix. Plus g\'en\'eralement, pour $\tilde{M}$ et $Q$ comme ci-dessus, on d\'efinit une distribution $f\mapsto J_{\tilde{M}}^Q(\tau,f)$. 

Dans le cas o\`u $M=G$,  on pose simplement $\theta_{\tau}(f)=J_{G}^G(\tau,f)$.  La distribution $f\mapsto \theta_{\tau}(f)$ est le caract\`ere usuel de $\tau$.

\bigskip

\subsection{Le $R$-groupe}

  Soient $M $ un L\'evi de $G$ et $\tau$ une repr\'esentation admissible  de $M(F)$. Soit $g\in G(F)$. On d\'efinit la repr\'esentation $g\tau$ de $gM(F)g^{-1}$ par $(g\tau)(gmg^{-1})=\tau(m)$. Son espace $E_{g\tau}$ est \'egal \`a $E_{\tau}$. Sa classe d'isomorphie ne d\'epend que de l'image  de $g$ dans l'ensemble de classes $G(F)/M(F)$. La conjugaison par $g$ induit un isomorphisme de ${\cal A}_{M}^*\otimes_{{\mathbb R}}{\mathbb C}$ sur ${\cal A}_{gMg^{-1}}^*\otimes_{{\mathbb R}}{\mathbb C}$ que l'on note $\lambda\mapsto g\lambda$. On a l'\'egalit\'e $g(\tau_{\lambda})=(g\tau)_{g\lambda}$ pour tout $\lambda\in {\cal A}_{M}^*\otimes_{{\mathbb R}}{\mathbb C}$.
  
   Supposons $M_{min}\subset M$ et $\tau$ irr\'eductible et de la s\'erie discr\`ete. Notons $Norm_{G(F)}(\tau)$ le sous-groupe des $g\in Norm_{G(F)}(M)$  tels que $g\tau\simeq \tau$. Posons $W(\tau)=Norm_{G(F)}(\tau)/M(F)$. Simplifions la th\'eorie en supposant v\'erifi\'ees les conditions suivantes:

- la repr\'esentation $\tau$ se prolonge en une repr\'esentation $\tau^N$ de $Norm_{G(F)}(\tau)$;

- l'homomorphisme naturel de $K\cap Norm_{G(F)}(\tau)$ dans $W(\tau)$ admet une section $\iota$ qui est un homomorphisme de groupes.

 Fixons $\tau^N$ et $\iota$. Pour simplifier la notation, on identifie tout \'el\'ement $w$ de $W(\tau)$ \`a son image $\iota(w)$. Remarquons que $\tau^N$ est n\'ecessairement unitaire. Pour $P\in {\cal P}(M)$ et $w\in W(\tau)$, on d\'efinit un homomorphisme
$$A_{P}(w):{\cal K}^G_{w^{-1}Pw,w^{-1}\tau}\to {\cal K}^G_{P,\tau}$$
par $(A_{P}(w)e)(k)=\tau^N(w)e(w^{-1}k)$. Pour $\lambda\in i{\cal A}^*_{M,F}$, on d\'efinit un endomomorphisme $R_{P}(w,\tau_{\lambda})$ de ${\cal K}^G_{P,\tau}$ par
$$R_{P}(w,\tau_{\lambda})=A_{P}(w)R_{w^{-1}Pw\vert P}(\tau_{\lambda})=R_{P\vert  wPw^{-1}}((w\tau)_{w\lambda})A_{wPw^{-1}}(w).$$
C'est un op\'erateur unitaire. On a la relation d'entrelacement
$$Ind_{P}^G((w\tau)_{w\lambda},g)R_{P}(w,\tau_{\lambda})=R_{P}(w,\tau_{\lambda})Ind_{P}(\tau_{\lambda},g)$$
et la relation de composition
$$R_{P}(w_{1}w_{2},\tau_{\lambda})=R_{P}(w_{1},\tau_{w_{2}\lambda})R_{P}(w_{2},\lambda).$$
Apppliquons ceci pour $\lambda=0$.  Notons $W'(\tau)$ le sous-groupe des $w\in W(\tau)$ tels que $R_{P}(w,\tau)$ soit une homoth\'etie.  C'est le groupe de Weyl d'un syst\`eme de racines $\Sigma'$ dont tout \'el\'ement est proportionnel \`a une racine de $A_{M}$ dans $\mathfrak{g}$. Ce syst\`eme est conserv\'e par l'action de $W(\tau)$. Fixons un sous-ensemble $\Sigma^{_{'}+}$ de racines positives et notons $R(\tau)$ le sous-groupe des \'el\'ements de $W(\tau)$ qui conservent $\Sigma^{_{'}+}$. On a la d\'ecomposition $W(\tau)=W'(\tau)\rtimes R(\tau)$. L'application $w\mapsto R_{P}(w,\tau)$ se prolonge en un isomorphisme de l'alg\`ebre de groupe ${\mathbb C}[R(\tau)]$ sur l'alg\`ebre commutante de la repr\'esentation $Ind_{P}^G(\tau)$. Ces propri\'et\'es forment la th\'eorie du $R$-groupe, qui est due \`a Silberger dans le cas $p$-adique.

 Simplifions encore en supposant le groupe $R(\tau)$ ab\'elien. On note $R(\tau) ^{\vee}$ le groupe dual de $R(\tau)$. Pour tout caract\`ere $\zeta\in R(\tau)^{\vee}$, notons ${\cal K}^G_{P,\tau,\zeta}$ le sous-espace des \'el\'ements $e\in {\cal K}^G_{P,\tau}$ tels que $R_{P}(w,\tau)e=\zeta(w)e$ pour tout $w\in R(\tau)$. Alors ${\cal K}^G_{P,\tau,\zeta}$ est stable par la repr\'esentation $Ind^G_{P}(\tau)$. Notons $Ind^G_{P}(\tau,\zeta)$ la restriction de $Ind^G_{P}(\tau)$ \`a ce sous-espace. Alors $Ind^G_{P}(\tau,\zeta)$ est irr\'eductible, sa classe ne d\'epend pas de $P$ et $Ind^G_{P}(\tau,\zeta)$ est isomorphe \`a $Ind^G_{P}(\tau,\zeta')$ si et seulement si $\zeta=\zeta'$.
 
 Soit $\pi$ une repr\'esentation admissible   de $G(F)$. On dit qu'elle est proprement induite s'il existe un \'el\'ement $Q=LU\in{\cal F}(M_{min})$ et une repr\'esentation admissible irr\'eductible $\sigma$ de $L(F)$ tels que $Q\not=G$ et $\pi\simeq Ind_{Q}^G(\sigma)$. Soit $\pi\in Temp(G)$. 	Nous dirons que $\pi$ est elliptique si elle n'est pas proprement induite. Revenons \`a la situation pr\'ec\'edente. Notons $W(M)=Norm_{G(F)}(M)/M(F)$ et $W(M)_{reg}$ le sous-ensemble des \'el\'ements de $W(M)$ qui op\`erent sans points fixes non nuls sur ${\cal A}_{M}/{\cal A}_{G}$. Une repr\'esentation $Ind_{P}^G(\tau,\zeta)$ comme ci-dessus est elliptique si et seulement si  $R(\tau)\cap W(M)_{reg}\not=\emptyset$. Si cette condition est v\'erifi\'ee, on a $W'(\tau)=\{1\}$. Inversement, pour toute repr\'esentation elliptique $\pi\in Temp(G)$, il existe $M$, $\tau$ et $\zeta$ v\'erifiant toutes les conditions ci-dessus de sorte que $\pi\simeq Ind_{P}^G(\tau,\zeta)$. La classe de conjugaison par $W^G$ du triplet $(M,\tau,\zeta)$ est bien d\'etermin\'ee.

\bigskip

\subsection{La formule de Plancherel-Harish-Chandra}

    Pour tout $M\in {\cal L}(M_{min})$, fixons un \'el\'ement $P\in {\cal P}(M)$. Notons $\Pi_{2}(M)$ l'ensemble des classes d'isomorphie de repr\'esentations irr\'eductibles et de la s\'erie discr\`ete de $M(F)$. Cet ensemble se d\'ecompose en orbites pour l'action $\tau\mapsto \tau_{\lambda}$ de $i{\cal A}^*_{M,F}$. Notons $\{\Pi_{2}(M)\}$ l'ensemble des orbites. Pour chaque orbite ${\cal O}$, fixons un \'el\'ement $\tau$ de cette orbite. Notons $i{\cal A}^{\vee}_{{\cal O}}$ le groupe des $\lambda\in i{\cal A}_{M}^*$ tels que  les repr\'esentations $\tau $ et $\tau_{\lambda}$ soient \'equivalentes.  Pour tout $\lambda\in i{\cal A}_{M}^*$, on  d\'efinit la mesure de Plancherel
$$m(\tau_{\lambda})=j(\tau_{\lambda})^{-1}d(\tau),$$
o\`u $d(\tau)$ est le degr\'e formel de $\tau$. Soit $f\in {\cal S}(G(F))$.   La formule de Plancherel-Harish-Chandra affirme l'\'egalit\'e
$$f(g)=\sum_{M\in {\cal L}(M_{min})}\vert W^M\vert \vert W^G\vert ^{-1}\sum_{{\cal O}\in \{\Pi_2(M)\}}[i{\cal A}^{\vee}_{{\cal O}}:i{\cal A}_{M,F}^{\vee}]^{-1}$$
$$\int_{i{\cal A}_{M,F}^*}m(\tau_{\lambda})trace(Ind^G_{P}(\tau_{\lambda},g^{-1})Ind_{P}^G(\tau_{\lambda},f))d\lambda$$
pour tout $g\in G(F)$.     Seules interviennent de fa\c{c}on non nulle les orbites ${\cal O}$ pour lesquelles une repr\'esentation $Ind_{P}^G(\tau_{\lambda})$ admet des vecteurs non nuls invariants par un sous-groupe ouvert compact de $G(F)$ tel que $f$ soit biinvariante par ce sous-groupe. Ces orbites sont en nombre fini. La formule ci-dessus est d\'emontr\'ee dans [W2] th\'eor\`eme VIII.1.1. Dans cette r\'ef\'erence, il y a quelques constantes suppl\'ementaires dues aux normalisations diff\'erentes des mesures. Arthur a introduit les normalisations que nous utilisons et qui font dispara\^{\i}tre ces constantes.

Nous aurons aussi besoin d'une autre formule. Fixons $P=MU\in {\cal F}(M_{min})$ et une repr\'esentation admissible irr\'eductible $\tau$ de $M(F)$, de la s\'erie discr\`ete. Soient $e,e'\in {\cal K}^G_{P,\tau}$ et $\varphi$ une fonction $C^{\infty}$ sur $i{\cal A}_{M,F}^*$. D\'efinissons une fonction $f_{e,e',\varphi}$ sur $G(F)$ par
$$f(g)=\int_{i{\cal A}_{M,F}^*}\varphi(\lambda)(Ind_{P}^G(\tau_{\lambda},g)e',e)m(\tau_{\lambda})d\lambda.$$
  Cette fonction appartient \`a ${\cal S}(G(F))$.  Identifions tout \'el\'ement de $W(M)$ \`a un repr\'esentant dans $K\cap Norm_{G(F)}(M)$. Notons ${\cal E}(\tau)$ l'ensemble des couples $(w,\mu)\in W(M)\times i{\cal A}_{M,F}^*$ tels que $w^{-1}\tau\simeq \tau_{\mu}$. Pour $(w,\mu)\in {\cal E}(\tau)$, fixons un automorphisme unitaire $\tau(w,\mu)$ de ${\cal K}^G_{P,\tau}$ tel que 
$$\tau(w,\mu)\tau_{\mu}(m)=(w^{-1}\tau)(m)\tau(w,\mu)$$
pour tout $m\in M(F)$. D\'efinissons l'homomorphisme $A(w,\mu):{\cal K}^G_{w^{-1}Pw,\tau}\to {\cal K}_{P,\tau}^G$ par 
$$(A(w,\mu)e)(g)=\tau(w,\mu)e(w^{-1}g).$$
Pour $\lambda\in i{\cal A}_{M,F}^*$, d\'efinissons l'endomorphisme $R(w,\mu,\lambda)$ de ${\cal K}^G_{P,\tau}$ par
$$R(w,\mu,\lambda)=A(w,\mu)R_{w^{-1}Pw\vert P}(\tau_{\lambda+\mu}).$$
Il v\'erifie la relation d'entrelacement 
$$R(w,\mu,\lambda)Ind_{P}^G(\tau_{\lambda+\mu},g)=Ind_{P}^G(\tau_{w\lambda},g)R(w,\mu,\lambda).$$
Posons simplement   $R(w,\mu)=R(w,\mu,0)$. Soient $e_{0},e'_{0}\in {\cal K}^G_{P,\tau}$. Alors on a l'\'egalit\'e
$$\int_{G(F)}f_{e,e',\varphi}(g)(e'_{0},Ind_{P}^G(\tau,g)e_{0})dg=\sum_{(w,\mu)\in {\cal E}(\tau)}\varphi(\mu)(R(w,\mu)e',e_{0})(e'_{0},R(w,\mu)e).$$
C'est une autre fa\c{c}on d'\'ecrire la proposition VII.2 de [W2]. Ici encore, les constantes disparaissent gr\^ace aux normalisations d'Arthur.

 L'ensemble ${\cal E}(\tau)$ est fini, on peut donc choisir un voisinage $\omega$ de $0$ dans $i{\cal A}_{M,F}^*$ tel que $(w,\mu)\in {\cal E}(\tau)$ et $\mu\in \omega$ entra\^{\i}nent $\mu=0$. Evidemment l'application $w\mapsto (w,0)$ est un isomorphisme de $W(\tau)$ sur le sous-ensemble des \'el\'ements de ${\cal E}(\tau)$ de la forme $(w,0)$. Supposons valides les hypoth\`eses du paragraphe pr\'ec\'edent. Pour $w\in W(\tau)$, l'op\'erateur $R(w,0)$ est \'egal au $R_{P}(w,\tau)$ du paragraphe pr\'ec\'edent (du moins, on peut effectuer les divers choix de sorte qu'il en soit ainsi). Pour $w\in W'(\tau)$, c'est une homoth\'etie dont le rapport est de module $1$ par un argument d'unitarit\'e. Seuls les \'el\'ements du $R$-groupe interviennent de fa\c{c}on non triviale dans la somme ci-dessus. On obtient alors
 
 (1) supposons le support de $\varphi$ contenu dans $\omega$; alors on a l'\'egalit\'e
$$\int_{G(F)}f_{e,e',\varphi}(g)(e'_{0},Ind_{P}^G(\tau,g)e_{0})dg=\vert W'(\tau)\vert \varphi(0)\sum_{w\in R(\tau)}(R_{P}(w,\tau)e',e_{0})(e'_{0},R_{P}(w,\tau)e).$$

\bigskip

\section{Fonctions tr\`es cuspidales}

\bigskip

\subsection{Un lemme d'annulation}

 Soit $\pi$ une repr\'esentation admissible de $G(F)$. Introduisons sa contragr\'ediente $\check{\pi}$. Soient $B$ une forme bilin\'eaire sur $E_{\check{\pi}}\times E_{\pi}$ et $f\in C_{c}^{\infty}(G(F))$. Fixons un sous-groupe ouvert compact $K_{f}$ de $G(F)$ tel que $f$ soit biinvariante par $K_{f}$. Fixons une base ${\cal B}^{K_{f}}$ du sous-espace $E_{\pi}^{K_{f}}$  et introduisons la base duale $\{\check{e}; e\in {\cal B}^{K_{f}} \}$ de $E_{\check{\pi}}^{K_{f}}$. Posons
 $$trace_{B}(\pi(f))=\sum_{e\in {\cal B}^{K_{f}}}B(\check{e},\pi(f)e).$$
 On v\'erifie que ce terme ne d\'epend ni du choix de $K_{f}$, ni de celui de la base. Remarquons que la trace usuelle $\theta_{\pi}(f))$ s'obtient comme cas particulier en prenant pour $B$ l'accouplement naturel sur  $E_{\check{\pi}}\times E_{\pi}$. On note $<.,.>$ cet accouplement.
 
 Soient $P=MU\in {\cal F}(M_{min})$ et $\tau$ une repr\'esentation admissible de $M(F)$. Soient $B$ une forme bilin\'eaire sur $E_{P,\check{\tau}}^G\times E_{P,\tau}^G$ et $f\in C_{c}^{\infty}(G(F))$. On impose les hypoth\`eses suivantes
 
 (1) $P\not=G$;
 
 (2) $f$ est tr\`es cuspidale (cf. [W1] 5.1);
 
 (3) soient $e\in E_{P,\tau}^G$ et $e'\in E_{P,\check{\tau}}^G$ tels que $e'(g)\otimes e(g)=0$ pour tout $g\in G(F)$; alors $B(e',e)=0$. 
 
 {\bf Remarque.} $e'(g)\otimes e(g)$ est un \'el\'ement de $E_{\check{\tau}}\otimes_{{\mathbb C}}E_{\tau}$.
 
 \ass{Lemme}{Sous ces hypoth\`eses, on a $trace_{B}(Ind_{P}^G(\tau,f))=0$.}
 
 Preuve.  On fixe   un sous-groupe ouvert compact $K_{f}$ de $K$ tel que $f$ soit biinvariante par $K_{f}$. Fixons un ensemble de repr\'esentants $\Gamma$ de l'ensemble de doubles classes $P(F)\backslash G(F)/K_{f}$. On peut choisir une base ${\cal B}^{K_{f}}$ de $(E_{P,\tau}^G)^{K_{f}}$ telle que, pour tout $e\in {\cal B}^{K_{f}}$, il existe $\gamma\in \Gamma$ de sorte que le support de $e$ soit contenu dans $P(F)\gamma K_{f}$. L'\'el\'ement correspondant $\check{e}$ de la base duale v\'erifie la m\^eme propri\'et\'e, avec le m\^eme $\gamma$. Pour tout $e\in {\cal B}^{K_{f}}$, on a
 $$Ind_{P}^G(\tau,f)e=\sum_{e'\in {\cal B}^{K_{f}}}<\check{e}',Ind_{P}^G(\tau,f)e>e',$$
 d'o\`u
 $$trace_{B}(Ind_{P}^G(\tau,f))=\sum_{e,e'\in {\cal B}^{K_{f}}}B(\check{e},e')<\check{e}',Ind_{P}^G(\tau,f)e>.$$
   Fixons $e,e'\in {\cal B}^{K_{f}}$. Il suffit de prouver que le terme index\'e par $e,e'$ dans cette somme est nul. Soient $\gamma,\gamma'\in \Gamma$ tels que le support de $e$, resp. $e'$, soit contenu dans $P(F)\gamma K_{f}$, resp. $P(F)\gamma'K_{f}$. Si $\gamma\not=\gamma'$, on a $\check{e}(g)\otimes e'(g)=0$ pour tout $g\in G(F)$, donc $B(\check{e},e')=0$. Supposons $\gamma=\gamma'$.  On a
 $$<\check{e}',Ind_{P}^G(\tau,f)e>=\int_{K}\int_{G(F)}f(g)<\check{e}'(x),e(xg)>dg\,dx,$$
 o\`u l'accouplement int\'erieur est celui sur $E_{\check{\tau}}\times E_{\tau}$. On effectue le changement de variable $g\mapsto x^{-1}g$ puis on d\'ecompose $g$ en $g=muk$, avec $m\in M(F)$, $u\in U(F)$, $k\in K$. D'o\`u
$$<\check{e}',Ind_{P}^G(\tau,f)e>=\int_{K}\int_{K}\int_{M(F)}\int_{U(F)}f(x^{-1}muk)$$
$$<\check{e}'(x),\tau(m)e(k)>\delta_{P}(m)^{1/2}du\,dm\,dk\,dx.$$
Fixons $x,k,m$ et supposons $<\check{e}'(x),\tau(m)e(k)>\not=0$. Cela entra\^{\i}ne $x\in P(F)\gamma K_{f}$ et $k\in P(F)\gamma K_{f}$. Donc $k\in P(F)xK_{f}$. Ecrivons $k=m'u'xk'$, avec $m'\in M(F)$, $u'\in U(F)$, $k'\in K_{f}$ et consid\'erons l'int\'egrale int\'erieure de la formule ci-dessus. Puisque $f$ est invariante \`a droite par $K_{f}$, le $k'$ dispara\^{\i}t. Par le changement de variable $u\mapsto m'uu^{_{'}-1}m^{_{'}-1}$,  cette int\'egrale  devient
$$\delta_{P}(m')^{1/2}\int_{U(F)}f(x^{-1}mm'ux)du.$$
Elle est nulle puisque $f$ est tr\`es cuspidale. Donc $<\check{e}',Ind_{P}^G(\tau,f)e>=0$, ce qui ach\`eve la preuve. $\square$

Ce lemme admet plusieurs variantes.  Supposons $M_{min}\subset M$. Au lieu des mod\`eles $E_{P,\tau}^G$ et $E_{P,\check{\tau}}^G$, on peut aussi consid\'erer les mod\`eles ${\cal K}_{P,\tau}^G$ et ${\cal K}_{P,\check{\tau}}^G$ et une forme bilin\'eaire $B$ sur ${\cal K}_{P,\check{\tau}}^G\times {\cal K}_{P,\tau}^G$. Le lemme reste valide si l'on remplace l'hypoth\`ese (3) par

(3')  soient $e\in {\cal K}_{P,\tau}^G$ et $e'\in {\cal K}_{P,\check{\tau}}^G$ tels que $e'(k)\otimes e(k)=0$ pour tout $k\in K$; alors $B(e',e)=0$.

Dans le cas o\`u $\tau$ unitaire, on peut aussi consid\'erer une forme sesquilin\'eaire $B$ sur ${\cal K}_{P,\tau}\times {\cal K}_{P,\tau}$ v\'erifiant la m\^eme condition (3'). Le lemme reste valide. 

\bigskip

\subsection{Caract\`eres pond\'er\'es et fonctions tr\`es cuspidales}

Soient $M\in {\cal L}(M_{min})$, $\tau$ une repr\'esentation   temp\'er\'ee de $M(F)$ et $f\in C_{c}^{\infty}(G(F))$.

\ass{Lemme}{Supposons $f$ tr\`es cuspidale.

(i) Soient $\tilde{M}\in {\cal L}(M)$ et $Q\in {\cal F}(\tilde{M})$. Si $\tilde{M}\not=M$ ou si $Q\not=G$, on a $J_{\tilde{M}}^Q(\tau,f)=0$.

(ii) Si $\tau$ est  proprement induite, on a $J_{M}^G(\tau,f)=0$.}

Preuve. Soit $Q=LU\in {\cal F}(M)$ tel que $Q\not=G$. Fixons $P\in{\cal P}(M)$ tel que $P\subset Q$. Pour $P'\in {\cal P}(M)$ tel que $P'\subset Q$, d\'efinissons une fonction $c_{P'}$ sur $i{\cal A}_{M}^*$ par
$$c_{P'}(\lambda)=trace(R_{P'\vert P}(\tau)^{-1}R_{P'\vert P}(\tau_{\lambda})Ind_{P}^G(\tau,f)).$$
Il s'agit de la trace d'un endomorphisme de ${\cal K}_{P,\tau}^G$. Par d\'efinition, $J_{M}^Q(\tau,f)$ est la valeur en $\lambda=0$ de la fonction
$$\sum_{P'\in {\cal P}(M); P'\subset Q}c_{P'}(\lambda)\theta^L_{P'\cap L}(\lambda)^{-1}.$$
Pour d\'emontrer que $J_{M}^Q(\tau,f)=0$, il suffit de prouver que $c_{P'}(\lambda)=0$ pour tous $P'$ et $\lambda$. Fixons $P'$ et $\lambda$. Introduisons les repr\'esentations $\pi=Ind_{P\cap L}^L(\tau)$ et $\pi'=Ind_{P'\cap L}^L(\tau)$, que l'on r\'ealise dans les espaces ${\cal K}_{P\cap L,\tau}^L$ et ${\cal K}_{P'\cap L,\tau}^L$.  On peut identifier ${\cal K}_{P,\tau}^G$, resp. ${\cal K}_{P',\tau}^G$, \`a ${\cal K}_{Q,\pi}^G$, resp. ${\cal K}_{Q,\pi'}^G$. On dispose de l'op\'erateur $R_{P'\cap L\vert P\cap L}^L(\tau_{\lambda}):{\cal K}^L_{P\cap L,\tau}\to {\cal K}^L_{P'\cap L,\tau}$. Modulo les identifications pr\'ec\'edentes, $R_{P'\vert P}(\tau_{\lambda})$ s'identifie \`a l'op\'erateur $e\mapsto R_{P'\cap L\vert P\cap L}^L(\tau_{\lambda})\circ e$ de ${\cal K}_{Q,\pi}^G$ dans ${\cal K}_{Q,\pi'}^G$. Il en est de m\^eme pour l'op\'erateur $R_{P'\vert P}(\tau)$. Introduisons la forme sesquilin\'eaire $B$ sur ${\cal K}_{Q,\pi'}^G\times {\cal K}_{Q,\pi}^G$ d\'efinie par
$$B(e',e)=(e',R_{P'\cap L\vert P\cap L}^L(\tau)^{-1}R_{P'\cap L\vert P\cap L}^L(\tau_{\lambda})\circ e).$$
On a alors $c_{P'}(\lambda)=trace_{B}(\pi(f))$. Les conditions du lemme pr\'ec\'edent sont v\'erifi\'ees, si l'on remplace dans ce lemme $P$ et $\tau$ par $Q$ et $\pi$. Le lemme entra\^{\i}ne $c_{P'}(\lambda)=0$ comme on le voulait.

Soient $\tilde{M}$ et $Q=LU$ comme en (i). On peut appliquer \`a la $(G,M)$-famille $({\cal R}^G_{P^L}(\tau))_{P^L\in {\cal P}^L(M)}$ les formules de descente d'Arthur, en particulier le corollaire 7.2 de [A2]. On en d\'eduit l'\'egalit\'e
$$J_{\tilde{M}}^Q(\tau,f)=\sum_{L'\in {\cal L}^L(M)}d^L_{M}(\tilde{M},L')J_{M}^{Q'}(\tau,f).$$
Le sous-groupe parabolique $Q'$ appartient \`a ${\cal P}(L')$ et est contenu dans $Q$. Si $Q\not=G$, tous ces $Q'$ sont aussi diff\'erents de $G$. Il en est de m\^eme si $\tilde{M}\not=M$ car dans ce cas, la condition $d_{M}^L(\tilde{M},L')\not=0$ implique que $L'\subsetneq L$. Alors le r\'esultat pr\'ec\'edent entra\^{\i}ne $J_{\tilde{M}}^Q(\tau,f)=0$, ce qui prouve (i).

Supposons $\tau $ proprement induite, fixons $P'=M'U'\in {\cal F}^M(M_{min})$ et une repr\'esentation irr\'eductible $\tau'$ de $M'(F)$ tels que $P'\not=M$ et $\tau=Ind_{P'}^M(\tau')$. La repr\'esentation $\tau'$ est temp\'er\'ee. Il r\'esulte des d\'efinitions que l'on a l'\'egalit\'e
$$J_{M}^G(\tau,f)=J_{M}^G(\tau',f).$$
On applique le (i) en rempla\c{c}ant $M$ et $\tilde{M}$ par $M'$ et $M$. On obtient la nullit\'e du terme de droite ci-dessus, d'o\`u le (ii) de l'\'enonc\'e. $\square$

Le terme $J_{M}^G(\tau,f)$ d\'epend a priori des facteurs $r_{P'\vert P}(\tau,\lambda)$ utilis\'es pour d\'efinir les op\'erateurs d'entrelacement normalis\'es. En fait

(1) pour $f$ tr\`es cuspidale, $J_{M}^G(\tau,f)$ ne d\'epend pas des facteurs de normalisation.

 En effet, consid\'erons deux familles de facteurs, que l'on indexe par les nombres $1$ et $2$. On en d\'eduit deux $(G,M)$-familles $(c_{P',1})_{P'\in {\cal P}(M)}$ et $(c_{P',2})_{P'\in {\cal P}(M)}$ comme dans la preuve ci-dessus et il suffit de prouver que $c_{M,1}=c_{M,2}$. Or il existe une $(G,M)$-famille $(d_{P'})_{P'\in {\cal P}(M)}$, construite \`a l'aide des rapports $r_{P'\vert P,1}(\tau,\lambda)r_{P'\vert P,2}(\tau,\lambda)^{-1}$, telle que $c_{P',1}=c_{P',2}d_{P'}$ pour tout $P'$.  On a alors la formule de descente
$$c_{M,1}=\sum_{L',L''\in {\cal L}(M)}d_{M}^G(L',L'')c_{M,2}^{Q'}d_{M}^{Q''},$$
cf. [A2] corollaire 7.4. Le terme $Q'$ est un \'el\'ement de ${\cal P}(M')$ et $c_{M,2}^{Q'}=J_{M}^{Q'}(\tau,f)$, ce terme \'etant calcul\'e \`a l'aide de la seconde famille de facteurs. Si $L'\not=G$, il est nul d'apr\`es le lemme ci-dessus. Dans la somme ci-dessus, il ne reste que la contribution du couple $(L',L'')=(G,M)$. Pour ce couple, $c_{M,2}^{Q'}d_{M}^{Q''}=c_{M,2}$, d'o\`u l'\'egalit\'e cherch\'ee.

\bigskip

\subsection{Induction de quasi-caract\`eres}

Pour ce paragraphe, oublions les  choix de mesures de Haar et de sous-groupe compact sp\'ecial que l'on a effectu\'es. Soit $M$ un L\'evi de $G$. On munit $G(F)$ et $M(F)$ de mesures de Haar. Soit $D^M$ une distribution sur $M(F)$ invariante par conjugaison. On sait d\'efinir la distribution induite $D=Ind_{M}^G(D^M)$, qui est une distribution invariante sur $G(F)$. Rappelons sa d\'efinition. Fixons un \'el\'ement $P=MU\in {\cal P}(M)$ et un sous-groupe compact sp\'ecial $K$ de $G(F)$, en bonne position relativement \`a $M$. Munissons $K$ et $U(F)$ de mesures de Haar, compatibles au sens habituel avec les mesures sur $M(F)$ et $G(F)$. Pour $f\in C_{c}^{\infty}(G(F))$, on d\'efinit $f_{P}\in C_{c}^{\infty}(M(F))$ par
$$f_{P}(m)=\delta_{P}(m)^{1/2}\int_{K}\int_{U(F)}f(k^{-1}muk)du\,dk.$$
On pose $D(f)=D^M(f_{P})$. Cela ne d\'epend pas des choix de $P$ et $K$. Soit maintenant $\theta^M$ une fonction d\'efinie presque partout sur $M(F)$, localement int\'egrable et invariante par conjugaison. Soit $D^M$ la distribution associ\'ee, c'est-\`a-dire que 
$$D^M(\varphi)=\int_{M(F)}\varphi(m)\theta^M(m)dm.$$
A l'aide de la formule d'int\'egration de Weyl, on v\'erifie que $D$ est elle-aussi associ\'ee \`a une fonction $\theta$ sur $G(F)$, localement int\'egrable et invariante par conjugaison. Pour tout $x\in G(F)$, fixons un ensemble ${\cal X}^M(x)$ de repr\'esentants des classes de conjugaison par $M(F)$ dans l'ensemble des \'el\'ements de $M(F)$ qui sont conjugu\'es \`a $x$ par un \'el\'ement de $G(F)$. Pour $x\in G_{reg}(F)$, on a l'\'egalit\'e
$$(1) \qquad \theta(x)=\sum_{x'\in {\cal X}^M(x)}D^G(x)^{-1/2}D^{M}(x')^{1/2}\theta^M(x').$$
Cette formule montre que $\theta$ est ind\'ependante des choix de mesures sur $G(F)$ et $M(F)$. On note $Ind_{M}^G(\theta^M)=\theta$. 

Rappelons que l'on a d\'efini en [W1] 4.1 la notion de quasi-caract\`ere sur $G(F)$. Soit $\theta$ une fonction d\'efinie presque partout sur $G(F)$ et invariante par conjugaison. On dit que c'est un quasi-caract\`ere si et seulement si, pour tout \'el\'ement semi-simple $x$ de $G(F)$, il existe un bon voisinage $\omega$ de $0$ dans $\mathfrak{g}_{x}(F)$ et, pour tout ${\cal O}\in Nil(\mathfrak{g}_{x})$, il existe $c_{\theta,{\cal O}}(x)\in {\mathbb C}$ de sorte que l'on ait l'\'egalit\'e
$$(2) \qquad \theta(xexp(X))=\sum_{{\cal O}\in Nil(\mathfrak{g}_{x})}c_{\theta,{\cal O}}(x)\hat{j}({\cal O},X)$$
presque partout pour $x\in \omega$. Donnons quelques explications. On note $G_{x}$ la composante neutre du centralisateur $Z_{G}(x)$ de $x$ dans $G$ et, comme toujours $\mathfrak{g}_{x}$ son alg\`ebre de Lie. On renvoie \`a [W1] 3.1 pour la notion de bon voisinage. On note $Nil(\mathfrak{g}_{x})$ l'ensemble des orbites nilpotentes dans $\mathfrak{g}_{x}$. La fonction $X\mapsto\hat{j}({\cal O},X)$ est la fonction associ\'ee \`a la distribution transform\'ee de Fourier de l'int\'egrale orbitale $J_{{\cal O}}$ associ\'ee \`a ${\cal O}$, normalis\'ee comme en [W1] 1.2.

D'autre part, soit ${\cal O}^M\in Nil(\mathfrak{m})$. On sait d\'efinir "l'orbite induite" de ${\cal O}^M$. Plus exactement, c'est  la r\'eunion d'un certain nombre d'\'el\'ements de $Nil(\mathfrak{g})$. Fixons $P=MU\in {\cal P}(M)$. Une orbite ${\cal O}\in Nil(\mathfrak{g})$ est incluse dans cette orbite induite si et seulement si l'intersection ${\cal O}\cap({\cal O}^M+\mathfrak{u}(F))$ contient un ouvert non vide de ${\cal O}^M+\mathfrak{u}(F)$. On pose $[{\cal O}:{\cal O}^M]=1$ si ${\cal O}$  est incluse dans l'orbite induite de ${\cal O}^M$, $[{\cal O}:{\cal O}^M]=0$ sinon. Remarquons que si $[{\cal O}:{\cal O}^M]=1$ et si l'une des deux orbites est r\'eguli\`ere, l'autre l'est aussi.

Pour un \'el\'ement semi-simple $x\in G(F)$, on a fix\'e ci-dessus un ensemble ${\cal X}^M(x)$. Pour tout \'el\'ement $x'$ de cet ensemble, on  note $\Gamma_{x'}$ l'ensemble des $g\in G(F)$ tels que $gxg^{-1}=x'$. C'est un espace principal homog\`ene pour l'action \`a droite de $Z_{G}(x)(F)$. Pour tout $g\in \Gamma_{x'}$, la conjugaison par $g$ envoie $Nil(\mathfrak{g}_{x})$ sur $Nil(\mathfrak{g}_{x'})$. On note ${\cal O}\mapsto g{\cal O}$ cette application.

\ass{Lemme}{Soit $\theta^M$ un quasi-caract\`ere de $M(F)$ et $\theta=Ind_{M}^G(\theta^M)$. Alors

(i) $\theta$ est un quasi-caract\`ere de $G(F)$;

(ii) soient $x$ un \'el\'ement semi-simple de $G(F)$ et ${\cal O}\in Nil(\mathfrak{g}_{x})$ une orbite r\'eguli\`ere; on a l'\'egalit\'e
$$c_{\theta,{\cal O}}(x)=\sum_{x'\in {\cal X}^M(x)}\sum_{g\in \Gamma_{x'}/G_{x}(F)}\sum_{{\cal O}'\in Nil(\mathfrak{m}_{x'})}D^G(x)^{-1/2}D^M(x')^{1/2}$$
$$[Z_{M}(x')(F):M_{x'}(F)]^{-1}[g{\cal O}:{\cal O}']c_{\theta^M,{\cal O}'}(x').$$}

Preuve. Soit $x$ un \'el\'ement semi-simple de $G(F)$. Consid\'erons un bon voisinage $\omega$ de $0$ dans $\mathfrak{g}_{x}$. Pour $x'\in {\cal X}^M(x)$, posons $\omega_{x'}=g\omega g^{-1}$, o\`u $g$ est un \'el\'ement quelconque de $\Gamma_{x'}$. C'est un bon voisinage de $0$ dans $\mathfrak{g}_{x'}$. En prenant $\omega$ assez petit, on peut supposer que $\omega_{x'}^M=\omega_{x'}\cap \mathfrak{m}_{x'}(F)$ est un bon voisinage de $0$ dans $\mathfrak{m}_{x'}(F)$ et que le quasi-caract\`ere $\theta^M$ admet un d\'eveloppement de la forme (2) dans $x'exp(\omega_{x'}^M)$. On d\'efinit $\theta^M_{x',\omega_{x'}^M}$: c'est la fonction sur $\mathfrak{m}_{x'}(F)$, \`a support dans $\omega_{x'}^M$ et telle que $\theta^M_{x',\omega_{x'}^M}(Y)=\theta^M(x'exp(Y))$ pour tout $Y\in \omega_{x'}^M$. C'est un quasi-caract\`ere sur $\mathfrak{m}_{x'}(F)$. En adaptant les d\'efinitions ci-dessus aux alg\`ebres de Lie, on d\'efinit la fonction localement int\'egrable $\phi_{x',\omega_{x'}}=Ind_{M_{x'}}^{G_{x'}}(\theta^M_{x',\omega^M_{x'}})$ sur $\mathfrak{g}_{x'}(F)$. On va prouver

(3) pour tout $X\in \omega\cap \mathfrak{g}_{x,reg}(F)$, on a l'\'egalit\'e
$$\theta(xexp(X))=\sum_{x'\in {\cal X}^M(x)}\sum_{g\in \Gamma_{x'}/G_{x}(F)}D^{G}(x)^{-1/2}D^{M}(x')^{1/2}$$
$$[Z_{M}(x')(F):M_{x'}(F)]^{-1}\phi_{x',\omega_{x'}}(gXg^{-1}).$$

Fixons $X$. Pour tout $x'\in {\cal X}^M(x)$ et tout $g\in \Gamma_{x'}$, fixons un ensemble ${\cal X}^{M_{x'}}(gXg^{-1})$ de repr\'esentants des classes de conjugaison par $M_{x'}(F)$ dans l'ensemble des \'el\'ements de $M_{x'}(F)$ qui sont conjugu\'es \`a $gXg^{-1}$ par un \'el\'ement de $G_{x'}(F)$. Il est inclus dans $\omega_{x'}^M$ et on peut supposer qu'il ne d\'epend que de l'image de $g$ dans $\Gamma_{x'}/G_{x}(F)$. En appliquant la formule (1) aux fonctions $\phi_{x',\omega_{x'}}$, le membre de droite de la formule (3) est \'egal \`a
$$\sum_{x'\in {\cal X}^M(x)}\sum_{g\in \Gamma_{x'}/G_{x}(F)}D^{G}(x)^{-1/2}D^{M}(x')^{1/2}[Z_{M}(x')(F):M_{x'}(F)]^{-1}$$
$$\sum_{Y\in {\cal X}^{M_{x'}}(gXg^{-1})}D^{G_{x'}}(Y)^{-1/2}D^{M_{x'}}(Y)^{1/2}\theta^M(x'exp(Y)).$$
Soient $x'$, $g$ et $Y$ apparaissant dans cette somme. On a
$$D^{G}(x)D^{G_{x'}}(Y)=D^G(x')D^{G_{x'}}(Y)=D^G(x'exp(Y)),$$
$$D^M(x')D^{M_{x'}}(Y)=D^M(x'exp(Y)).$$
D'autre part, $x'exp(Y)$ est un \'el\'ement de $M(F)$ conjugu\'e \`a $xexp(X)$ par un \'el\'ement de $G(F)$. Il est donc conjugu\'e par un \'el\'ement de $M(F)$ \`a un unique \'el\'ement de ${\cal X}^M(xexp(X))$. Notons cet \'el\'ement $y(x',g,Y)$. On a
$$D^G(x'exp(Y))^{-1/2}D^M(x'exp(Y))^{1/2}\theta^M(x'exp(Y))=$$
$$D^G(xexp(X))^{-1/2}D^M(y(x',g,Y))^{1/2}\theta^M(y(x',g,Y)).$$
La formule ci-dessus s'\'ecrit donc
$$\sum_{y\in {\cal X}^M(xexp(X))}c(y)D^G(xexp(X))^{-1/2}D^M(y)^{1/2}\theta^M(y),$$
o\`u 
$$c(y)=\sum_{x',g,Y; y(x',g,Y)=y}[Z_{M}(x')(F):M_{x'}(F)]^{-1}.$$
En comparant avec la formule (1), on voit que, pour d\'emontrer (3), il suffit de prouver que

(4) $c(y)=1$ pour tout $y\in {\cal X}^M(xexp(X))$.

Soit $y\in {\cal X}^M(xexp(X))$. Fixons $\gamma\in G(F)$ tel que $\gamma xexp(X)\gamma^{-1}=y$. Puisque $y\in M(F)\cap G_{reg}(F)$, le centralisateur $G_{y}$ de $y$ est contenu dans $M$. Mais $\gamma x\gamma^{-1}$ appartient \`a $G_{y}(F)$. Il appartient donc \`a $M(F)$. Il existe donc $x'\in {\cal X}^M(x)$ et $m\in M(F)$ tel que $\gamma x\gamma^{-1}=mx'm^{-1}$. Posons $g=m^{-1}\gamma$. Alors $g\in \Gamma_{x'}$ et $\gamma=mg$. On a $y=mx'exp(gXg^{-1})m^{-1}$. Puisque $y\in M(F)$, on a aussi $x'exp(gXg^{-1})\in M(F)$, donc $gXg^{-1}\in \mathfrak{m}_{x'}(F)$. Alors $gXg^{-1}$ est conjugu\'e par un \'el\'ement de $M_{x'}(F)$ \`a un \'el\'ement $Y\in {\cal X}^{M_{x'}}(gXg^{-1})$ et $y$ est conjugu\'e par un \'el\'ement de $M(F)$ \`a $x'exp(Y)$. Pour ces choix de $x',g,Y$, on a $y=y(x',g,Y)$. 
Soit $(x'_{1},g_{1},Y_{1})$ un autre triplet, supposons $y=y(x'_{1},g_{1},Y_{1})$.  Quitte \`a multiplier $g$ \`a gauche par un \'el\'ement de $G_{x'}(F)$ (ce qui revient \`a le multiplier \`a droite par un \'el\'ement de $G_{x}(F)$), on peut supposer $Y=gXg^{-1}$. De m\^eme, on peut supposer $Y_{1}=g_{1}Xg_{1}^{-1}$. Soit $\mu\in M(F)$ tel que $\mu x'exp(Y)\mu^{-1}=x'_{1}exp(Y_{1})$. Alors $\mu gxexp(X)g^{-1}\mu^{-1}=g_{1}xexp(X)g_{1}^{-1}$. En posant $h=g_{1}^{-1}\mu g$, on a $hxexp(X)h^{-1}=xexp(X)$. Puisque $xexp(X)$ est r\'egulier, cela entra\^{\i}ne $h\in G_{xexp(X)}(F)$. Cet ensemble est contenu dans $G_{x}(F)$ d'apr\`es les propri\'et\'es des bons voisinages. Donc $h\in G_{x}(F)$. On a alors
$$\mu x'\mu^{-1}=\mu gxg^{-1}\mu^{-1}=g_{1}hxh^{-1}g_{1}^{-1}=x'_{1}.$$
Par d\'efinition de l'ensemble ${\cal X}^M(x)$, cela entra\^{\i}ne $x'_{1}=x'$. A fortiori, la constante $[Z_{M}(x_{1}')(F):M_{x_{1}'}(F)]^{-1}$ qui intervient dans la d\'efinition de $c(y)$ est \'egale \`a $[Z_{M}(x')(F):M_{x'}(F)]^{-1}$ et ne d\'epend pas du triplet. Puisque $x'_{1}=x'$, la relation ci-dessus entra\^{\i}ne $\mu\in Z_{M}(x')(F)$. En revenant \`a la d\'efinition de $\mu$, on a $\mu Y\mu^{-1}=Y_{1}$. Le couple $(g_{1}G_{x}(F),Y_{1})$  appartient donc \`a l'orbite de $(gG_{x}(F),Y)$ pour l'action de $Z_{M}(x')(F)$ ainsi d\'efinie: l'action de $\mu\in Z_{M}(x')(F)$ envoie $(gG_{x}(F),Y)$ sur le couple $(g_{1}G_{x}(F),Y_{1})$ tel que $g_{1}G_{x}(F)=\mu gG_{x}(F)$ et que $Y_{1}$ soit l'unique \'el\'ement de $ {\cal X}^{M_{x'}}(g_{1}Xg_{1}^{-1})$ conjugu\'e \`a $\mu Y\mu^{-1}$ par un \'el\'ement de $M_{x'}(F)$. Inversement, on v\'erifie que tout couple ainsi obtenu convient. L'action de $Z_{M}(x')(F)$ que l'on vient de d\'efinir se quotiente en une action de $Z_{M}(x')(F)/M_{x'}(F)$. Remarquons que $Z_{M}(x')\cap G_{x'}=M_{x'}$ car ces deux ensembles sont \'egaux au commutant de $A_{M}$ dans $G_{x'}$. Il en r\'esulte que l'action de $Z_{M}(x')(F)/M_{x'}(F)$ est libre: son action sur la premi\`ere composante l'est. Le nombre de triplet est donc \'egal au nombre d'\'el\'ements de ce groupe, ce qui entra\^{\i}ne (4) et ach\`eve la preuve de (3).

La formule (3) nous ram\`ene au probl\`eme suivant. Soit maintenant $\theta^M$ un quasi-caract\`ere sur $\mathfrak{m}(F)$, dont on \'ecrit le d\'eveloppement \`a l'origine
$$\theta^M(Y)=\sum_{{\cal O}^M\in Nil(\mathfrak{m})}c_{\theta^M,{\cal O}^M}\hat{j}^M({\cal O}^M,Y).$$
Soit $\theta=Ind_{M}^G(\theta^M)$. On doit prouver que $\theta$ poss\`ede un d\'eveloppement \`a l'origine de la forme
$$(5) \qquad \theta(X)=\sum_{{\cal O}\in Nil(\mathfrak{g})}c_{\theta,{\cal O}}\hat{j}^G({\cal O},X),$$
et prouver que, pour ${\cal O}$ r\'eguli\`ere, on a l'\'egalit\'e
$$(6) \qquad c_{\theta,{\cal O}}=\sum_{{\cal O}^M\in Nil(\mathfrak{m})}[{\cal O}:{\cal O}^M].$$
Comme on l'a expliqu\'e, l'induction ne d\'epend pas des mesures de Haar, si on la consid\`ere comme une application portant sur des fonctions localement int\'egrables. On peut donc supposer que les mesures sont normalis\'ees comme en [W1] 1.2. L'analogue pour les alg\`ebres de Lie de l'application $f\mapsto f_{P}$  "commute" \`a la transformation de Fourier. On en d\'eduit que l'induite d'une fonction $Y\mapsto \hat{j}^M({\cal O}^M,Y)$ est la fonction associ\'ee \`a la transform\'ee de Fourier de la distribution induite de l'int\'egrale orbitale  $J^M_{{\cal O}^M}$. Il est bien connu que cette distribution induite est combinaison lin\'eaire des int\'egrales orbitales $J^G_{{\cal O}}$ pour des \'el\'ements ${\cal O}\in Nil(\mathfrak{g})$ inclus dans l'orbite induite de ${\cal O}^M$. En tout cas, l'induite d'une fonction $Y\mapsto \hat{j}^M({\cal O}^M,Y)$ est combinaison lin\'eaire de fonctions $X\mapsto j^G({\cal O},X)$, ce qui prouve l'existence du d\'eveloppement (5). Pour prouver (6),  on voit qu'il suffit de prouver que, pour ${\cal O}^M$ r\'eguli\`ere, la distribution induite de $J^M_{{\cal O}^M}$ est \'egale \`a
$$\sum_{{\cal O}\in Nil(\mathfrak{g})}[{\cal O}:{\cal O}^M]J^G_{{\cal O}}.$$
On peut supposer $M$ et $G$ quasi-d\'eploy\'es, sinon il n'y a pas d'orbites nilpotentes r\'eguli\`eres et la question est vide. Toute orbite nilpotente r\'eguli\`ere ${\cal O}$ de $\mathfrak{g}(F)$ appara\^{\i}t dans  l'orbite induite  d'une unique orbite nilpotente r\'eguli\`ere ${\cal O}^M$ de $\mathfrak{m}(F)$. En effet, fixons $P\in {\cal P}(M)$ et un sous-groupe de Borel $B$ de $G$ tel que $B\subset P$. Soient $Y,Y'\in \mathfrak{m}(F)$, $N,N'\in \mathfrak{u}(F)$, supposons que $Y+N$ et $Y'+N'$ appartiennent \`a ${\cal O}$. Quitte \`a effectuer des conjugaisons par des \'el\'ements de $M(F)$, on peut supposer $Y,Y'\in \mathfrak{b}(F)\cap \mathfrak{m}(F)$. Soit $g\in G(F)$ tel que $g(Y+N)g^{-1}=Y'+N'$. L'\'el\'ement $Y'+N'$ appartient aux deux sous-alg\`ebres de Borel $\mathfrak{b}$ et $g\mathfrak{b}g^{-1}$. Mais $Y'+N'$ est r\'egulier donc n'appartient qu'\`a une seule telle sous-alg\`ebre. Donc $g\mathfrak{b}g^{-1}=\mathfrak{b}$ et $g$ appartient \`a $B(F)$. En \'ecrivant $g=mu$, avec $m\in M(F)$ et $u\in U(F)$, on a alors $Y'=mYm^{-1}$, c'est-\`a-dire que $Y$ et $Y'$ sont dans la m\^eme orbite. Cette unicit\'e nous permet de transformer notre probl\`eme en le suivant: prouver que la distribution induite de
$$\sum_{{\cal O}^M \text{ r\'eguli\`ere}}J_{{\cal O}^M}^M$$
est \'egale \`a
$$\sum_{{\cal O}\text{ r\'eguli\`ere}}J_{{\cal O}}^G.$$
Introduisons un sous-groupe de Borel $B$ comme ci-dessus et un sous-tore maximal $T\subset B\cap M$. Fixons un \'el\'ement $X\in \mathfrak{t}(F)\cap \mathfrak{g}_{reg}(F)$. En utilisant un r\'esultat de Shelstad, on a prouv\'e en [W1] lemme 11.4 que la premi\`ere distribution ci-dessus \'etait la limite simple des distributions $f\mapsto J^M_{M}(zX,f)$ sur $M(F)$ quand $z\in F^{\times}$ tend vers $0$ (dans [W1], notre groupe \'etait un groupe sp\'ecial orthogonal, mais la d\'emonstration de cette propri\'et\'e n'utilisait pas cette particularit\'e). De m\^eme, la seconde distribution est la limite simple des distributions $f\mapsto J_{G}(zX,f)$ sur $G(F)$. Mais il r\'esulte des d\'efinitions que la distribution $f\mapsto J_{G}(zX,f)$ est l'induite de la distribution $f\mapsto J^M_{M}(zX,f)$. La conclusion s'ensuit. $\square$

\bigskip

\subsection{Int\'egrales orbitales pond\'er\'ees invariantes}

Soient $f,f'\in C_{c}^{\infty}(G(F))$. Nous dirons que $f$ et $f'$ sont \'equivalentes si et seulement si $D(f)=D(f')$  pour toute distribution $D$ sur $G(F)$ invariante par conjugaison. Comme on le sait, cette condition est \'equivalente \`a l'une ou l'autre des deux conditions suivantes

(1) $J_{G}(x,f)=J_{G}(x,f')$ pour tout $x\in G(F)$;

(2) $\theta_{\pi}(f)=\theta_{\pi}(f')$ pour toute repr\'esentation  $\pi\in Temp(G)$. 

Soient $M\in {\cal L}(M_{min})$, $x\in M(F)\cap G_{reg}(F)$ et $f\in C_{c}^{\infty}(G(F))$. Arthur a d\'efini l'int\'egrale orbitale pond\'er\'ee $J_{M}(x,f)$. On a rappel\'e la d\'efinition en [W1] 2.3. Il a aussi d\'efini l'int\'egrale pond\'er\'ee invariante $I_{M}(x,f)$. Rappelons la d\'efinition. Pour $Z\in {\cal A}_{G,F}$, notons ${\bf 1}_{H_{G}=Z}$ la fonction caract\'eristique de l'ensemble des $x\in G(F)$ tels que $H_{G}(x)=Z$. Notons ${\cal H}_{ac}(G(F))$ l'ensemble des fonctions $f:G(F)\to {\mathbb C}$ qui v\'erifient les deux conditions suivantes

(3) $f$ est biinvariante par un sous-groupe ouvert compact  de $G(F)$;

(4) pour tout $Z\in {\cal A}_{G,F}$, la fonction $f{\bf 1}_{H_{G}=Z}$ appartient \`a $C_{c}^{\infty}(G(F))$.

Remarquons que plusieurs d\'efinitions pos\'ees pour les fonctions appartenant \`a $C_{c}^{\infty}(G(F))$ se g\'en\'eralisent aux \'el\'ements de ${\cal H}_{ac}(G(F))$. Par exemple les int\'egrales orbitales pond\'er\'ees (on pose $J_{M}(x,f)=J_{M}(x,f{\bf 1}_{H_{G}=H_{G}(x)})$) ou la notion d'\'equivalence introduite ci-dessus.

Soient $L\in {\cal L}(M_{min})$ et $f\in C_{c}^{\infty}(G(F))$. Arthur montre qu'il existe une fonction $\phi_{L}(f)\in {\cal H}_{ac}(L(F))$ telle que, pour toute repr\'esentation $\pi\in Temp(L)$ et tout $Z\in {\cal A}_{L,F}$, on ait l'\'egalit\'e
$$(5) \qquad \int_{i{\cal A}_{L,F}^*}J_{L}(\pi_{\lambda},f)exp(-\lambda(Z))d\lambda=\theta_{\pi}(\phi_{L}(f){\bf 1}_{H_{L}=Z}).$$
La fonction $\phi_{L}(f)$ est bien d\'efinie \`a \'equivalence pr\`es. On d\'efinit $I_{M}(x,f)$ par r\'ecurrence sur $a_{M}-a_{G}$ par la formule
$$J_{M}(x,f)=\sum_{L\in {\cal L}(M)}I_{M}^L(x,\phi_{L}(f){\bf 1}_{H_{L}=H_{L}(x)}).$$
Bien s\^ur, $I_{M}(x,f)$ ne d\'epend que de la classe de conjugaison par $M(F)$ de $x$. La distribution $f\mapsto I_{M}(x,f)$ est invariante par conjugaison par $G(F)$ et ne d\'epend pas du choix du groupe $K$. La propri\'et\'e suivante en r\'esulte, par simple transport de structure. Soit $g\in G(F)$ tel que $gMg^{-1}\in {\cal L}(M_{min})$. Alors on a l'\'egalit\'e $I_{gMg^{-1}}(gxg^{-1},f)=I_{M}(x,f)$.

Pour $f\in C_{c}^{\infty}(G(F))$, on d\'efinit une fonction $I\theta_{f}$ sur $G_{reg}(F)$ de la fa\c{c}on suivante. Soit $x\in G_{reg}(F)$.   Notons $M(x)$ le commutant de $A_{G_{x}}$ dans $G$. C'est un L\'evi de $G$ et c'est le plus petit L\'evi contenant $x$. Choisissons $g\in G(F)$ tel que $gM(x)g^{-1}\in {\cal L}(M_{min})$. On pose
$$I\theta_{f}(x)=(-1)^{a_{M(x)}-a_{G}}D^G(x)^{-1/2}I_{gM(x)g^{-1}}(gxg^{-1},f).$$
Cela ne d\'epend pas du choix de $g$. La fonction $I\theta_{f}$ est invariante par conjugaison et localement constante sur $G_{reg}(F)$.  Remarquons que $I\theta_{f}=I\theta_{f'}$ si $f$ et $f'$ sont \'equivalentes, puisque les distributions $f\mapsto I_{M}(x,f)$ sont invariantes par conjugaison.

\bigskip

\subsection{Fonctions cuspidales et quasi-caract\`eres invariants}

  Soit $f\in C_{c}^{\infty}(G(F))$. On dit que $f$ est cuspidale si et seulement si, pour tout groupe de L\'evi $M\subsetneq G$ et pour tout $x\in G_{reg}(F)\cap M(F)$, on a $J_{G}(x,f)=0$.   Cette condition est \'equivalente \`a ce que $\theta_{\pi}(f)=0$ pour toute repr\'esentation $\pi$ de $G(F)$ qui est temp\'er\'ee et proprement induite. Une fonction tr\`es cuspidale est cuspidale.
  
  \ass{Lemme}{Soit $f\in C_{c}^{\infty}(G(F))$, supposons $f$ cuspidale. Alors $I\theta_{f}$ est un quasi-caract\`ere de $G(F)$.}
  
  Preuve. Arthur d\'efinit un ensemble de repr\'esentations virtuelles $T_{ell}(G)$. Tout \'el\'ement $\pi$ de $T_{ell}(G)$ est une combinaison lin\'eaire \`a coefficients complexes de repr\'esentations elliptiques. Par lin\'earit\'e, on d\'efinit la contragr\'ediente $\check{\pi}$, le caract\`ere $\theta_{\pi}$ et, pour $\lambda\in i{\cal A}_{G}^*$, la repr\'esentation virtuelle $\pi_{\lambda}$ qui appartient aussi \`a $T_{ell}(G)$. On note $\{T_{ell}(G)\}$ l'ensemble des orbites dans $T_{ell}(G)$ pour l'action $\lambda\mapsto \pi_{\lambda}$.  Si $K'$ est un sous-groupe ouvert compact de $G(F)$, il n'y a qu'un nombre fini d'orbites ${\cal O}\in T_{ell}(G)$ pour lesquelles il existe $\pi\in {\cal O}$ et une fonction $f'\in C_{c}^{\infty}(G(F))$, biinvariante par $K'$, de sorte que $\theta_{\pi}(f')\not=0$. Pour toute orbite ${\cal O}$, on fixe $\pi\in {\cal O}$ et on d\'efinit un certain coefficient $c({\cal O})>0$. Cela \'etant, en [A5] th\'eor\`eme 5.1, Arthur d\'emontre que, pour toute fonction cuspidale $f\in C_{c}^{\infty}(G(F))$, pour tout $M\in {\cal L}(M_{min})$ et pour tout \'el\'ement $y\in M(F)\cap G_{reg}(F)$ qui est elliptique dans $M(F)$, on a l'\'egalit\'e
  $$D^G(y)^{-1/2}(-1)^{a_{M}-a_{G}}I_{M}(y,f)=\sum_{{\cal O}\in \{T_{ell}(G)\}}c({\cal O})\int_{i{\cal A}_{G,F}^*}\theta_{\pi_{\lambda}}(y)\theta_{(\pi_{\lambda})\check{}}(f)d\lambda.$$
  Elle \'equivaut \`a
  $$D^G(y)^{-1/2}(-1)^{a_{M}-a_{G}}I_{M}(y,f)=\sum_{{\cal O}\in \{T_{ell}(G)\}}c({\cal O})\theta_{\pi}(y)\theta_{\check{\pi}}(f{\bf 1}_{H_{G}=H_{G}(y)}).$$
  Soit $x$ un \'el\'ement semi-simple de $G(F)$. Pour $y$ dans un certain voisinage de $x$, on a $H_{G}(y)=H_{G}(x)$. Nos d\'efinitions et la formule ci-dessus entra\^{\i}nent que, pour $y\in G_{reg}(F)$ dans ce voisinage, on a l'\'egalit\'e
  $$(1) \qquad I\theta_{f}(y)=\sum_{{\cal O}\in \{T_{ell}(G)\}}c({\cal O})\theta_{\pi}(y)\theta_{\check{\pi}}(f{\bf 1}_{H_{G}=H_{G}(x)}).$$
  Comme on l'a dit, la somme est en fait finie. Donc $I\theta_{f}$ co\"{\i}ncide dans ce voisinage de $x$ avec une combinaison lin\'eaire finie de caract\`eres de repr\'esentations admissibles irr\'eductibles. D'apr\`es Harish-Chandra ([HCDeBS] th\'eor\`eme 16.2), tout tel caract\`ere est un quasi-caract\`ere. La notion de quasi-caract\`ere \'etant de nature locale, la conclusion s'ensuit. $\square$
  
  On appelle $I\theta_{f}$ le quasi-caract\`ere invariant associ\'e \`a $f$. 
    
  La notion de cuspidalit\'e   se g\'en\'eralise aux \'el\'ements de ${\cal H}_{ac}(G(F))$: $f\in {\cal H}_{ac}(G(F))$ est cuspidale si et seulement si $f{\bf 1}_{H_{G}=Z}$ l'est pour tout $Z\in {\cal A}_{G,F}$. La d\'efinition de $I\theta_{f}$ aussi: $I\theta_{f}=\sum_{Z\in {\cal A}_{G,F}}I\theta_{f{\bf 1}_{H_{G}=Z}}$, cette somme \'etant localement finie.
  
  \bigskip
  
  \subsection{Quasi-caract\`ere et quasi-caract\`ere invariant}
  
  Soit $f\in C_{c}^{\infty}(G(F))$ une fonction tr\`es cuspidale. On vient de lui associer un quasi-caract\`ere $I\theta_{f}$ sur $G(F)$. Dans [W1] 5.6 et 5.9, on lui a aussi associ\'e un quasi-caract\`ere $\theta_{f}$. En fait, cette d\'efinition d\'epend des choix de mesures. Nous modifions la d\'efinition de [W1] 5.6 en utilisant nos pr\'esentes mesures plut\^ot que celles de cette r\'ef\'erence. Il convient de comparer $\theta_{f}$ et $I\theta_{f}$.
  
  \ass{Lemme}{Soit $f\in C_{c}^{\infty}(G(F))$ une fonction tr\`es cuspidale. Alors
  
  (i) pour tout $L\in {\cal L}(M_{min})$, la fonction $\phi_{L}(f)$ est cuspidale;
  
  (ii) on a l'\'egalit\'e 
  $$\theta_{f}=\sum_{L\in {\cal L}(M_{min})}\vert W^L\vert \vert W^G\vert ^{-1}(-1)^{a_{L}-a_{G}}Ind_{L}^G(I\theta^L_{\phi_{L}(f)}).$$}
  
  Preuve. Pour prouver (i), on doit montrer que, pour tout $Z\in {\cal A}_{L,F}$ et toute repr\'esentation  temp\'er\'ee proprement induite $\pi$ de $L(F)$, on a $\theta_{\pi}(\phi_{L}(f){\bf 1}_{H_{L}=Z})=0$. Cela r\'esulte de l'\'egalit\'e 2.4(5) et du lemme 2.2(ii) qui affirme que $J_{L}^G(\pi_{\lambda},f)=0$ pour tout $\lambda\in i{\cal A}_{L}^*$.
  
  Soit $\varphi\in C_{c}^{\infty}(G(F))$. On peut \'ecrire la formule d'int\'egration de Weyl sous la forme
  $$(1) \qquad \int_{G(F)}\theta_{f}(g)\varphi(g)dg=\sum_{M\in {\cal L}(M_{min})}\vert W^M\vert \vert W^G\vert ^{-1}$$
  $$\sum_{T\in {\cal T}_{ell}(M)}\vert W(M,T)\vert ^{-1}\int_{T(F)}\theta_{f}(t)J_{G}(t,\varphi)D^G(t)^{1/2}dt,$$
  avec les notations d'Arthur que l'on a rappel\'ees en [W1] 2.4. Pour tout $L\in {\cal L}(M_{min})$, fixons $Q_{L}\in {\cal P}(L)$. On a de m\^eme
  $$\int_{G(F)}Ind_{L}^G(I\theta^L_{\phi_{L}(f)})(g)\varphi(g)dg=\int_{L(F)}I\theta^L_{\phi_{L}(f)}(l)\varphi_{Q_{L}}(l)dl$$
  $$=\sum_{M\in {\cal L}^L(M_{min})}\vert W^M\vert \vert W^L\vert ^{-1}\sum_{T\in {\cal T}_{ell}(M)}\vert W(M,T)\vert ^{-1}\int_{T(F)}I\theta^L_{\phi_{L}(f)}(t)J_{L}(t,\varphi_{Q_{L}})D^L(t)^{1/2}dt.$$
  Le terme $J_{L}(t,\varphi_{Q_{L}})$ intervenant ci-dessus est \'egal \`a $J_{G}(t,\varphi)$. Notons $\theta'_{f}$ la fonction figurant dans le membre de droite du (ii) de l'\'enonc\'e. En sommant les \'egalit\'es ci-dessus, on obtient
  $$(2) \qquad \int_{G(F)}\theta'_{f}(g)\varphi(g)dg$$
  $$=\sum_{M\in {\cal L}(M_{min})}\vert W^M\vert \vert W^G\vert ^{-1}\sum_{T\in {\cal T}_{ell}(M)}\vert W(M,T)\vert ^{-1}\int_{T(F)}\gamma_{M,T}(t)J_{G}(t,\varphi)dt,$$
  o\`u  on a pos\'e
  $$\gamma_{M,T}(t)=\sum_{L\in {\cal L}(M)}(-1)^{a_{L}-a_{G}}I\theta^L_{\phi_{L}(f)}(t)D^L(t)^{1/2}.$$
  Soient $M\in {\cal L}(M_{min})$, $T\in {\cal T}_{ell}(M)$ et $t\in T(F)\cap G_{reg}(F)$. Pour $L\in {\cal L}(M)$, appliquons la d\'efinition de $I\theta^L_{\phi_{L}(f)}(t)$ donn\'ee en 2.4. Puisque $T$ est elliptique dans $M$, le L\'evi $M(t)$  est \'egal \`a $M$ (que le groupe ambiant soit $G$ ou $L$).  Donc 
  $$I\theta^L_{\phi_{L}(f)}(t)= (-1)^{a_{M}-a_{L}}D^L(t)^{-1/2}I_{M}^L(t,\phi_{L}(f)).$$
  On peut aussi bien remplacer $\phi_{L}(f)$ par $\phi_{L}(f){\bf 1}_{H_{L}=H_{L}(t)}$. Alors
  $$\gamma_{M,T}(t)=(-1)^{a_{M}-a_{G}}\sum_{L\in {\cal L}(M)}I_{M}^L(t,\phi_{L}(f){\bf 1}_{H_{L}=H_{L}(t)})$$
  $$=(-1)^{a_{M}-a_{G}}J_{M}^G(t,f).$$
  En se reportant \`a la d\'efinition de [W1] 5.3, on obtient $\gamma_{M,T}(t)=D^G(t)^{1/2}\theta_{f}(t)$. On conclut en comparant les \'egalit\'es (1) et (2). $\square$
  
  \bigskip
  
  \subsection{Fonctions cuspidales et fonctions tr\`es cuspidales}
  
  \ass{Lemme}{Soit $f\in C_{c}^{\infty}(G(F))$ une fonction cuspidale. Alors il existe une fonction tr\`es cuspidale $f'\in C_{c}^{\infty}(G(F))$  qui est \'equivalente \`a $f$.}
  
  Preuve. Par un proc\'ed\'e de partition de l'unit\'e tel que celui de la preuve de [W1] proposition 6.4, il suffit de prouver l'assertion suivante
  
  (1) soit $x\in G(F)$ un \'el\'ement semi-simple; alors il existe un $G$-domaine $\Omega$ dans $G(F)$ et une fonction tr\`es cuspidale $f'\in C_{c}^{\infty}(G(F))$ tels que $x\in \Omega$ et $J_{G}(y,f')=J_{G}(y,f)$ pour tout $y\in \Omega\cap G_{reg}(F)$.
  
  Supposons $A_{G_{x}}\not=A_{G}$. Il existe un $G$-domaine $\Omega$ contenant $x$ tel que, pour $y\in \Omega\cap G_{reg}(F)$, on ait $A_{G_{y}}\not=A_{G}$, autrement dit $y$ n'est pas elliptique dans $G(F)$. Alors $J_{G}(y,f)=0$ et il suffit de prendre $f'=0$ pour v\'erifier l'assertion. Supposons maintenant $A_{G_{x}}=A_{G}$. Fixons un bon voisinage $\omega$ de $0$ dans $\mathfrak{g}_{x}(F)$. Le quasi-caract\`ere $I\theta_{f}$ se descend en un quasi-caract\`ere $I\theta_{f,x,\omega}$ sur $\mathfrak{g}_{x}(F)$, cf. [W1] 4.3, qui est \'evidemment \`a support compact modulo conjugaison et invariant par l'action de $Z_{G}(x)(F)$. En combinant la proposition 6.4 et le lemme 6.2 de [W1], on voit qu'il existe une fonction tr\`es cuspidale $f'\in C_{c}^{\infty}(G(F))$ telle que $\theta_{f',x,\omega}=I\theta_{f,x,\omega}$. Posons $\Omega=\{g^{-1}xexp(X)g; X\in \omega, g\in G(F)\}$ et soit $y\in \Omega\cap G_{reg}(F)$. Si $y$ n'est pas elliptique dans $G(F)$, on a $J_{G}(y,f')=0=J_{G}(y,f)$. Supposons $y$ elliptique, \'ecrivons $y=g^{-1}xexp(X)g$ avec $g\in G(F)$ et $X\in \omega$. D'apr\`es les d\'efinitions, on a $J_{G}(y,f')=\theta_{f',x,\omega}(X)$ et $J_{G}(y,f)=I\theta_{f,x,\omega}(X)$. D'o\`u l'\'egalit\'e $J_{G}(y,f')=J_{G}(y,f)$. Cela prouve (1) et le lemme. $\square$

\bigskip

\section{Majorations pour le groupe lin\'eaire $GL_{k}$}

\bigskip

\subsection{Le groupe lin\'eaire}

Soient $k\geq1$ un entier, $V$ un espace vectoriel sur $F$ de dimension $k$ et $(v_{i})_{i=1,...,k}$ une base de $V$. On note simplement $GL_{k}$ le groupe (alg\'ebrique) des automorphismes lin\'eaires de $V$. Pour $g\in GL_{k}(F)$, on note $(g_{i,j})_{i,j=1,...,k}$ sa matrice dans la base fix\'ee. On note $B_{k}$ le sous-groupe de Borel triangulaire sup\'erieur de $GL_{k}$, $U_{k}$ son radical unipotent et $A_{k}$ le sous-tore diagonal. Pour $a\in A_{k}(F)$, on note simplement $a_{i}=a_{i,i}$ son coefficient diagonal, pour $i=1,...,k$. On note $K_{k}$ le sous-groupe compact sp\'ecial de $GL_{k}(F)$ form\'e des \'el\'ements \`a coefficients entiers.

La th\'eorie du $R$-groupe est "triviale" pour le groupe lin\'eaire. C'est-\`a-dire que les repr\'esentations temp\'er\'ees irr\'eductibles et elliptiques de $GL_{k}(F)$ sont de la s\'erie discr\`ete.
 
Ces notations seront utilis\'ees pour divers espaces, parfois sans que l'on pr\'ecise leur base. Ou bien le choix de cette base sera implicite, ou bien il n'aura pas d'importance.

Dans la suite de cette section, on fixe un entier $k\geq1$, on pose $G=GL_{k}$, et on utilise les notations ci-dessus dont on supprime l'indice $k$.

\bigskip

\subsection{Une majoration}

  Pour tout $g\in G(F)$, on note $g=a_{B}(g)u_{B}(g)k_{B}(g)$ une d\'ecomposition de $g$ telle que $a_{B}(g)\in A(F)$, $u_{B}(g)\in U(F)$, $k_{B}(g)\in K$.  Pour un entier $c\geq 1$, on note $U(F)_{c}$ le sous-groupe des \'el\'ements $u\in U(F)$ tels que $val_{F}(u_{i,i+1})\geq-c$ pour tout $i=1,...,k-1$. Soient $D\in {\mathbb R}$ et $g\in G(F)$, posons
$$I(c,D,g)=\int_{U(F)_{c}}\Xi^G(ug)\sigma(ug)^Ddu.$$

\ass{Proposition}{Cette int\'egrale est convergente. Pour tout $D$, il existe un r\'eel $R$ tel que
$$I(c,D,g)<<c^R\sigma(g)^R\delta_{B}(a_{B}(g))^{1/2}$$
pour tous $c\geq1$ et tout $g\in G(F)$.}

La preuve sera donn\'ee en 3.4

\bigskip

\subsection{Un lemme auxiliaire}

 Supposons $k\geq2$, notons $P=MU_{P}$ le sous-groupe parabolique des \'el\'ements de $G$ qui conservent la droite $Fv_{1}$. Soient $c\geq1$ un entier, $D$ un r\'eel et $m\in M(F)$. Posons $U_{P}(F)_{c}=U_{P}(F)\cap U(F)_{c}$ et
$$I_{P}(c,D,m)=\int_{U_{P}(F)_{c}}\Xi^G(um)\sigma(um)^Ddu.$$

\ass{Lemme}{Cette int\'egrale est convergente. Pour tout $D$, il existe un r\'eel $R$ tel que
$$I_{P}(c,D,m)<<c^R\sigma(m)^R\delta_{P}(m)^{1/2}\Xi^M(m)$$
pour tout $c\geq1$ et tout $m\in M(F)$.}

Preuve.  Le r\'eel $D$ est fix\'e. Soit $b\geq0$ un r\'eel. On a introduit en 1.1 la fonction ${\bf 1}_{\sigma\geq b}$. Notons $U_{\natural}$ le sous-groupe des \'el\'ements $u\in U_{1}$ tels que $u_{1,2}=0$. Posons
$$I_{\natural}(b,D)=\int_{U_{\natural}(F)}{\bf 1}_{\sigma\geq b}(u)\delta_{\bar{B}}(a_{\bar{B}}(u))^{1/2}\sigma(u)^Ddu.$$
Montrons que

(1) cette int\'egrale est convergente et il existe $\epsilon>0$ tel que 
$$I_{\natural}(b,D)<<exp(-\epsilon b)$$
pour tout $b\geq0$.

Consid\'erons le sous-espace $V''$ de $V$ engendr\'e par les \'el\'ements $v_{1},v_{3},...,v_{k}$. Soit $G''=GL_{k-1}$ son groupe d'automorphismes, qui s'identifie \`a un sous-groupe de $G$. Le groupe $B''=G''\cap B$ est le sous-groupe de Borel standard de $G''$, $P''=G''\cap P$ est un sous-groupe parabolique de $G''$  et $U_{\natural}$ n'est autre que le radical unipotent de $P''$. D'apr\`es [W2] lemme II.4.2, il existe un entier $d\geq0$ tel que l'int\'egrale
$$(2) \qquad \int_{U_{\natural}(F)}\delta_{\bar{B}''}(a_{\bar{B}''}(u))^{1/2}\sigma(u)^{-d}du$$
soit convergente. Soit $u\in U_{\natural}(F)$. On peut supposer $a_{\bar{B}}(u)=a_{\bar{B}''}(u)$. Notons $a_{1},...,a_{k}$ les coefficients diagonaux de cet \'el\'ement. On peut supposer  $a_{2}=1$ et  $\prod_{i=1,...,k}a_{i}=1$. On a
$$\delta_{\bar{B}}(a_{\bar{B}}(u))^{1/2}=\prod_{i=1,...,k}\vert a_{i}\vert_{F} ^{i-(k+1)/2},$$
tandis que
$$\delta_{\bar{B}''}(a_{\bar{B}''}(u))^{1/2}=\vert a_{1}\vert_{F} ^{1-k/2}\prod_{i=3,...,k}\vert a_{i}\vert_{F} ^{i-1-k/2}.$$
D'o\`u
$$\delta_{\bar{B}}(a_{\bar{B}}(u))^{1/2}=\vert a_{1}\vert _{F}^{-1/2}\left(\prod_{i=3,...,k}\vert a_{i}\vert _{F}^{1/2}\right)\delta_{\bar{B}''}(a_{\bar{B}''}(u))^{1/2}=\vert a_{1}\vert _{F}^{-1}\delta_{\bar{B}''}(a_{\bar{B}''}(u))^{1/2}.$$
L'\'egalit\'e $u=a_{\bar{B}''}(u)u_{\bar{B}''}(u)k_{\bar{B}''}(u)$ et un calcul matriciel entra\^{\i}nent  que 
$$val_{F}(a_{1})=inf\{val_{F}(u_{1,j}); j=1,3,...,k\},$$
d'o\`u $-val_{F}(a_{1})>> \sigma(u)$. Il existe donc $\epsilon>0$ tel que $\vert a_{1}\vert_{F} ^{-1}<<exp(-2\epsilon\sigma(u))$. Si ${\bf 1}_{\sigma\geq b}(u)=1$, on transforme cette relation en $\vert a_{1}\vert_{F} ^{-1}<<exp(-\epsilon b)exp(-\epsilon\sigma(u))$. On obtient
$$I_{\natural}(b,D)<<exp(-\epsilon b)\int_{U_{\natural}(F)}\delta_{\bar{B}''}(a_{\bar{B}''}(u))^{1/2}exp(-\epsilon\sigma(u))du.$$
La derni\`ere int\'egrale est convergente d'apr\`es la convergence de (2). D'o\`u (1).

Introduisons un r\'eel $b\geq0$, que nous pr\'eciserons par la suite. On peut \'ecrire
$$(3) \qquad I_{P}(c,D,m)=I_{<b}(c,D,m)+I_{\geq b}(c,D,m),$$
o\`u 
$$I_{<b}(c,D,m)=\int_{U_{P}(F)_{c}}\Xi^G(um)\sigma(um)^D{\bf 1}_{\sigma< b}(um)du,$$
et
$$I_{\geq b}(c,D,m)=\int_{U_{P}(F)_{c}}\Xi^G(um)\sigma(um)^D{\bf 1}_{\sigma\geq b}(um)du.$$
Dans la premi\`ere int\'egrale, on a $\sigma(um)<b$ et l'int\'egrale est donc major\'ee par
$$b^{D'}\int_{U_{P}(F)}\Xi^G(um)\sigma(um)^{D-D'}du$$
pour tout r\'eel $D'\geq0$. D'apr\`es [W2] proposition II.4.5, on peut choisir $D'$ tel que cette derni\`ere int\'egrale soit convergente et essentiellement born\'ee par $\delta_{P}(m)^{1/2}\Xi^M(m)$. Pour un tel $D'$, on a donc
$$(4) \qquad I_{< b}(c,D,m)<<b^{D'}\delta_{P}(m)^{1/2}\Xi^M(m).$$
Introduisons le sous-groupe $U_{1,2}$ de $U_{P}$ form\'e des \'el\'ements $u$ dont la seule coordonn\'ee non diagonale non nulle est $u_{1,2}$ et posons $U_{1,2}(F)_{c}=U_{1,2}(F)\cap U(F)_{c}$. On a $U_{P}(F)_{c}=U_{\natural}(F)U_{1,2}(F)_{c}$, d'o\`u 
$$I_{\geq b}(c,D,m)=\int_{U_{1,2}(F)_{c}}\int_{U_{\natural}(F)}\Xi^G(u_{\natural}vm)\sigma(u_{\natural}vm)^D{\bf 1}_{\sigma\geq b}(u_{\natural}vm)du_{\natural}dv.$$
Fixons $v\in U_{1,2}(F)_{c}$ et consid\'erons l'int\'egrale int\'erieure, que l'on note $I_{\geq b}(c,D,m,v)$. On a

(5) il existe $\alpha>0$ tel que
$$\Xi^G(gg')<<\Xi^G(g)exp(\alpha\sigma(g'))$$
pour tous $g,g'\in G(F)$.

En effet, par d\'efinition,
$$\Xi^G(g)=\int_{K}\delta_{B}(a_{B}(kg))^{1/2}dk.$$
On a $a_{B}(kgg')=a_{B}(kg)a_{B}(k_{B}(kg)g')$ et il existe $\alpha>0$ tel que $\delta_{B}(a_{B}(k_{B}(kg)g'))^{1/2}<<exp(\alpha\sigma(g'))$. D'o\`u le r\'esultat.

{\bf Remarque.} La propri\'et\'e (5) est vraie pour tout groupe r\'eductif connexe.
\bigskip

 On a la majoration $\sigma(gg')\leq \sigma(g)+\sigma(g')$ pour tous $g,g'\in G(F)$.  On a aussi une majoration $\sigma(v)<<c$. Il existe donc $c_{1}>0$ tel que $\sigma(u_{\natural})\geq c_{1}\sigma(u_{\natural}vm)-c-\sigma(m)$. Quand ${\bf 1}_{\sigma\geq b}(u_{\natural}vm)=1$, on a $\sigma(u_{\natural})\geq c_{1}b-c-\sigma(m)$. Imposons \`a $b$ la condition $c_{1}b-c-\sigma(m)>0$. D'autre part, d'apr\`es [W2] lemmes II.1.1 et II.3.2, il existe un r\'eel $D''$ tel que l'on ait une majoration
$$\Xi^G(g)<<\delta_{\bar{B}}(a_{\bar{B}}(g))^{1/2}\sigma(g)^{D''}.$$
  Ces relations entra\^{\i}nent la majoration
$$I_{\geq b}(c,D,m,v)<<exp(\alpha\sigma(vm))\sigma(vm)^{D+D''}\int_{U_{\natural}(F)}\delta_{\bar{B}}(a_{\bar{B}}(u))^{1/2}\sigma(u)^{D+D''}{\bf 1}_{\sigma\geq c_{1}b-c-\sigma(m)}(u)du$$
$$<< exp(c_{2}c)exp(c_{2}\sigma(m))I_{\natural}(c_{1}b-c-\sigma(m),D+D'')$$
pour un $c_{2}>0$ convenable. Cette expression ne d\'epend plus de $v$. Le terme $I_{\geq b}(c,D,m)$ est major\'e par la m\^eme expression, multipli\'ee par $mes(U_{1,2}(F)_{c})$. Cette mesure est major\'ee par $exp(c_{3}c)$ pour un $c_{3}>0$ convenable. D'o\`u
$$I_{\geq b}(c,D,m)<<exp(c_{4}c)exp(c_{2}\sigma(m))I_{\natural}(c_{1}b-c-\sigma(m),D+D''),$$
o\`u $c_{4}=c_{2}+c_{3}$. Il existe aussi $c_{5}>0$ tel que l'on ait la minoration
$$exp(-c_{5}\sigma(m))<<\delta_{P}(m)^{1/2}\Xi^M(m).$$
D'o\`u
$$I_{\geq b}(c,D,m)<<\delta_{P}(m)^{1/2}\Xi^M(m)exp(c_{4}c)exp(c_{6}\sigma(m))I_{\natural}(c_{1}b-c-\sigma(m),D+D''),$$
o\`u $c_{6}=c_{2}+c_{5}$. Utilisons la relation (1). On voit qu'il existe $c_{7}>0$ tel que, pour $b=c_{7}(c+\sigma(m))$, le produit des trois derniers termes ci-dessus est born\'e. Choisissons $b$ ainsi. Alors 
$$I_{\geq b}(c,D,m)<<\delta_{P}(m)^{1/2}\Xi^M(m).$$
Cette majoration et les relations (3) et (4) entra\^{\i}nent celle de l'\'enonc\'e. $\square$

\bigskip

\subsection{Preuve de la proposition 3.2}

On d\'emontre la proposition par r\'ecurrence sur $k$. Le cas $k=1$ est \'evident. Supposons $k\geq2$. 
Remarquons tout d'abord que l'on peut se limiter \`a d\'emontrer la majoration de l'\'enonc\'e pour $g=a\in A(F)$. En effet, pour $g$ quelconque, on \'ecrit $g=vak$, avec $v\in U(F)$, $a\in A(F)$, $k\in K$. Effectuons le changement de variable $u\mapsto uv^{-1}$.  Le nouveau domaine d'int\'egration est $U(F)_{c}v$. Mais il existe $c_{1}>0$ tel que cet ensemble soit inclus dans $U(F)_{c+c_{1}\sigma(g)}$. Alors
$$I(c,D,g)\leq I(c+c_{1}\sigma(g),D,a).$$
Si le deuxi\`eme terme v\'erifie une majoration comme dans l'\'enonc\'e, le premier terme aussi. Supposons donc $g=a\in A(F)$. Avec les notations du paragraphe pr\'ec\'edent, on a
$$I(c,D,a)=\int_{M(F)\cap U(F)_{c}}\int_{U_{P}(F)_{c}}\Xi^G(uva)\sigma(uva)^Ddu\,dv$$
$$=\int_{M(F)\cap U(F)_{c}}I_{P}(c,D,va)dv.$$
En appliquant le lemme 3.3, on a
$$I(c,D,a)<<c^R\int_{M(F)\cap U(F)_{c}}\delta_{P}(va)^{1/2}\Xi^M(va)\sigma(va)^Rdv.$$
Ecrivons $M=GL_{1}\times G'$, o\`u $G'=GL_{k-1}$ et affectons d'un $'$ les objets relatifs \`a $G'$. Ecrivons aussi $a=(a_{1},a')$, avec $a_{1}\in F^{\times}$ et $a'\in A'(F)$. On a $\delta_{P}(va)^{1/2}=\delta_{B}(a)^{1/2}\delta_{B'}(a')^{-1/2}$, $\Xi^M(va)=\Xi^{G'}(va')$ et $\sigma(va)^R<<\sigma(a)^R\sigma(va')^R$. On  obtient
$$I(c,D,a)<<c^R\sigma(a)^R\delta_{B}(a)^{1/2}\delta_{B'}(a')^{-1/2}I'(c,R,a').$$
 En utilisant la majoration du dernier terme fournie par l'hypoth\`ese de r\'ecurrence, on obtient celle cherch\'ee. $\square$
 
 \bigskip
 
 \subsection{Mod\`eles de Whittaker et int\'egrales de coefficients}

  On d\'efinit un caract\`ere $\xi$ de $U(F)$ par la formule
$$\xi(u)=\psi(\sum_{i=1,...,k-1}u_{i,i+1}).$$

 Soit $\mu\in Temp(G)$. On appelle fonctionnelle de Whittaker sur $E_{\mu}$ une application lin\'eaire $\phi:E_{\mu}\to {\mathbb C}$ telle que $\phi(\mu(u)e)=\xi(u)\phi(e)$ pour tous $u\in U(F)$ et $e\in E_{\mu}$. Comme on le sait, l'espace des fonctionnelles de Whittaker sur $E_{\mu}$ est une droite. Soit $c\geq1$ un entier. D\'efinissons une forme sesquilin\'eaire ${\cal L}_{\mu,c}$ sur $ E_{\mu}$ (ce qui est un raccourci pour dire qu'il s'agit d'une forme sur $E_{\mu}\times E_{\mu}$) par
 $${\cal L}_{\mu,c}(e',e)=\int_{U(F)_{c}}(e',\mu(u)e)\bar{\xi}(u)du.$$
 Cette int\'egrale est absolument convergente d'apr\`es la proposition 3.2. Notons $\omega_{[1,k-1]}(c)$ le sous-groupe des $a\in A(F)$ tels que $a_{k}=1$ et $val_{F}(1-a_{i})\geq c$ pour tout $i=1,..,k-1$.
 
 \ass{Lemme}{Pour tous $e,e'\in E_{\mu}$, il existe un entier $c_{0}\geq1$ tel que ${\cal L}_{\mu,c}(e',e)$ soit ind\'ependant de $c$ pour $c\geq c_{0}$. Plus pr\'ecis\'ement, soit $c'\geq1$ un entier. Il existe $c_{0}$ tel que cette conclusion  soit v\'erifi\'ee pour tous $e,e'\in E_{\mu}^{\omega_{[1,k-1]}(c')}$.}
 
 Preuve. Soient $e,e'\in E_{\mu}^{\omega_{[1,k-1]}(c')}$. Choisissons $c_{0}$ tel que $c_{0}\geq1$ et $-c_{0}+c'\leq c_{\psi}$. Pour $c\geq c_{0}$, notons $U(F)_{c}-U(F)_{c_{0}}$ le compl\'ementaire de $U(F)_{c_{0}}$ dans $U(F)_{c}$. Il suffit de prouver que
 $$\int_{U(F)_{c}-U(F)_{c_{0}}}(e',\mu(u)e)\bar{\xi}(u)du=0.$$
 Soit $a\in \omega_{[1,k-1]}(c')$. Dans l'int\'egrale pr\'ec\'edente, on peut remplacer $e'$ et $e$ par $\mu(a)e'$ et $\mu(a)e$. Par le changement de variables $u\mapsto aua^{-1}$, l'int\'egrale devient
 $$\int_{U(F)_{c}-U(F)_{c_{0}}}(e',\mu(u)e)\bar{\xi}(aua^{-1})du.$$
 Elle est donc aussi \'egale \`a
 $$mes(\omega_{[1,k-1]}(c'))^{-1}\int_{\omega_{[1,k-1]}(c')}\int_{U(F)_{c}-U(F)_{c_{0}}}(e',\mu(u)e)\bar{\xi}(aua^{-1})du\,da.$$
 Cette expression est absolument convergente et on peut permuter les int\'egrales. Mais
 $$\int_{\omega_{[1,k-1]}(c')}\bar{\xi}(aua^{-1})da=0$$
 pour tout $u\in U(F)_{c}-U(F)_{c_{0}}$. Cela prouve la nullit\'e cherch\'ee et le lemme. $\square$
 
 On d\'efinit une forme sesquilin\'eaire ${\cal L}_{\mu}$ sur $E_{\mu}$ par
 $${\cal L}_{\mu}(e',e)=lim_{c\to \infty}{\cal L}_{\mu,c}(e',e).$$
 Cette forme v\'erifie les relations
 $${\cal L}_{\mu}(\mu(u')e',\mu(u)e)=\xi(u)\bar{\xi}(u'){\cal L}_{\mu}(e',e).$$
 Fixons une fonctionnelle de Whittaker $\phi$ sur $E_{\mu}$, non nulle. Alors il existe $C_{\mu}\in {\mathbb C}$ tel que
 $${\cal L}_{\mu}(e',e)=C_{\mu}\overline{\phi(e')}\phi(e)$$
 pour tous $e,e'$. On montrera plus loin que $C_{\mu}\not=0$.
 
 Pour un entier $c'\in {\mathbb N}$, notons $\iota_{c'}$ la fonction caract\'eristique du sous-ensemble des $a\in A(F)$ tels que $val_{F}(a_{i})\geq val_{F}(a_{i+1})-c'$ pour tout $i=1,...,k-1$. On a
 
 (1) il existe un r\'eel $R$ et, pour tous $e,e'\in E_{\mu}$, il existe un entier $c'\in {\mathbb N}$ tel que
 $$\vert {\cal L}_{\mu}(\mu(a')e',\mu(a)e)\vert<< \iota_{c'}(a')\delta_{B}(a')^{1/2}\sigma(a')^R\iota_{c'}(a)\delta_{B}(a)^{1/2}$$
 pour tous $a,a'\in A(F)$.
 
 Si $C_{\mu}=0$, c'est \'evident: ${\cal L}_{\mu}=0$. Sinon, fixons $e_{0}\in E_{\mu}$ tel que $\phi(e_{0})=1$. On a 
 $${\cal L}_{\mu}(\mu(a')e',\mu(a)e)=C_{\mu}\overline{\phi(\mu(a')e')}\phi(\mu(a)e)=C_{\mu}\overline{\phi(\mu(a')e')}\phi(e_{0})\overline{\phi(e_{0})}\phi(\mu(a)e)$$
 $$=\overline{C_{\mu}}^{-1}\overline{{\cal L}_{\mu}(e_{0},\mu(a')e')}{\cal L}_{\mu}(e_{0},\mu(a)e).$$
 Cela nous ram\`ene \`a majorer $\vert {\cal L}_{\mu}(e_{0},\mu(a)e)\vert $. D'apr\`es le lemme ci-dessus, on peut remplacer ${\cal L}_{\mu}$ par ${\cal L}_{\mu,c}$ pour $c$ assez grand. La majoration voulue r\'esulte de la proposition 3.2  et du fait bien connu  qu'il existe $c'$ tel que le support de la fonction $a\mapsto \phi(\mu(a)e)$ soit contenu dans celui de $\iota_{c'}$ 
  
 \bigskip
 
 \subsection{Quelques \'egalit\'es d'int\'egrales}
 
 Soit $h$ un entier tel que $1\leq h\leq k$. Notons $P^h$ le sous-groupe parabolique de $G$ form\'e des \'el\'ements   qui conservent le drapeau 
 $$  Fv_{1}\oplus...\oplus Fv_{h}\subset Fv_{1}\oplus...\oplus Fv_{h+1}\subset...\subset Fv_{1}\oplus...\oplus Fv_{k}.$$
 Ecrivons $P^h=M^hU^h$, o\`u $M^h$ est la composante de L\'evi qui contient $A$. On a $M^h=GL_{h}\times GL_{1}\times ...\times GL_{1} $, avec $k-h$ termes $GL_{1}$. Identifions $GL_{h-1}$  au groupe des automorphismes lin\'eaires du sous-espace de $V$ engendr\'e par $v_{1},...,v_{h-1}$. Pour un entier $c\geq1$, notons $\omega_{[h,k-1]}(c)$ le sous-groupe des $\gamma\in A(F)$ tels que $\gamma_{1}=...=\gamma_{h-1}=1$, $\gamma_{k}=1$ et $val_{F}(1-\gamma_{i})\geq c$ pour $i=h,...,k-1$.  Posons $U^h(F)_{c}=U(F)_{c}\cap U^h(F)$. Supposons $c+c_{\psi}\geq1$. Pour $\mu\in Temp(G)$ et $e,e'\in E_{\mu}$, posons
 $$I_{c}^h(e',e)=\int_{\omega_{[h,k-1]}(c+c_{\psi})}\int_{U_{h-1}(F)\backslash GL_{h-1}(F)}{\cal L}_{\mu}(\mu(\gamma g)e',\mu(\gamma g)e)\vert det(g)\vert _{F}^{h-k}dg \,d\gamma,$$
 $$J_{c}^h(e',e)=\int_{U^h(F)_{c}}(e',\mu(u)e)\bar{\xi}(u)du.$$
 
 \ass{Lemme}{Soient $h$ un entier tel que $1\leq h\leq k$ et $c$ un entier tel que $c\geq1$ et $c+c_{\psi}\geq1$. Les int\'egrales ci-dessus sont absolument convergentes. Il existe $C>0$ tel que
 $$J_{c}^h(e',e)=CI_{c}^h(e',e)$$
 pour tous $e,e'\in E_{\mu}$.}
 
 Preuve. Prouvons la convergence de $I_{c}^h(e',e)$. L'int\'egrale sur le groupe compact $\omega_{[h,k-1]}(c+c_{\psi})$ est insignifiante, on peut l'oublier. On d\'ecompose $g\in U_{h-1}(F)\backslash GL_{h-1}(F)$ en $g=ak$, avec $a\in A_{h-1}(F)$, $k\in K$. La mesure devient $\delta_{B_{h-1}}(a)^{-1}dadk$. De nouveau, on peut oublier l'int\'egrale sur $K$ et on doit majorer
 $$\int_{A_{h-1}(F)}\vert {\cal L}_{\mu}(\mu(a)e',\mu(a)e)\vert  \vert det(a)\vert _{F}^{h-k}\delta_{B_{h-1}}^{-1}(a)da.$$
 Gr\^ace \`a 3.5(1), c'est major\'e par
 $$\int_{A_{h-1}(F)}\iota_{c'}(a)\delta_{B}(a)\vert det(a)\vert _{F}^{h-k}\delta_{B_{h-1}}^{-1}(a)da$$
 pour un entier $c'$ convenable. Pour $a\in A_{h-1}(F)$, on a
 $$\delta_{B}(a)\vert det(a)\vert _{F}^{h-k}\delta_{B_{h-1}}^{-1}(a)= \vert det(a)\vert _{F}$$
 et il est imm\'ediat que l'int\'egrale
 $$\int_{A_{h-1}(F)}\iota_{c'}(a)\vert det(a)\vert _{F}da$$
 est convergente.
 
 Prouvons la convergence de $J_{c}^h(e',e)$. Il suffit de prouver que
 $$\int_{U^h(F)_{c}}\Xi^G(u)du$$
 est convergente. Puisqu'on en aura besoin plus loin, prouvons la propri\'et\'e plus forte suivante. On s'autorise pour un instant \`a faire varier l'entier $c$. Alors
 
 (1) il existe un r\'eel $R$ tel que
 $$\int_{U^h(F)_{c}}\Xi^G(m^{-1}um)du <<c^R\sigma(m)^R\delta_{P^h}(m)$$
 pour tous $c\geq1$ et tout $m\in M^h(F)$
 
 Notons $X(m)$ l'int\'egrale ci-dessus. On peut \'ecrire $m=vak$, avec $v\in U(F)\cap M^h(F)$, $a\in A(F)$ et $k\in K\cap M^h(F)$. La conjugaison par $v$ conserve $U^h(F)_{c}$, ce qui permet de faire dispara\^{\i}tre $v$. Le terme $k$ dispara\^{\i}t \'egalement. Puisque $\sigma(a)<<\sigma(m)$ et $\delta_{P^h}(a)=\delta_{P^h}(m)$, on est ramen\'e au cas o\`u $m=a$. Effectuons le changement de variable $u\mapsto aua^{-1}$. Cela remplace $du$ par $\delta_{P^h}(a)du$ et le domaine d'int\'egration $U^h(F)_{c}$ par l'ensemble des $u\in U^h(F)$ tels que
 $$val_{F}(u_{i,i+1})\geq -c+val_{F}(a_{i+1})-val_{F}(a_{i})$$
 pour tout $i=h,...,k-1$. Il existe $c_{1}>0$ tel que $val_{F}(a_{i+1})-val_{F}(a_{i})\geq -c_{1}\sigma(a)$. L'ensemble ci-dessus est donc contenu dans $U^h(F)_{c+c_{1}\sigma(a)}$  on obtient
  $$X(a)\leq \delta_{P^h}(a)\int_{U^h(F)_{c+c_{1}\sigma(a)}}\Xi^G(u)du.$$
 Pour $u'\in M^h(F)\cap U(F)\cap K$, on a $\Xi^G(u)=\Xi^G(uu')$. On peut donc remplacer ci-dessus $u$ par $uu'$,  puis int\'egrer en $u'$.  D'o\`u
 $$X(a)<<\delta_{P^h}(a) \int_{M^h(F)\cap U(F)\cap K}\int_{U^h(F)_{c+c_{1}\sigma(a)}}\Xi^G(uv)du\,dv.$$
 L'ensemble d'int\'egration est contenu dans $U(F)_{c+c_{1}\sigma(a)}$. Donc
 $$X(a)<<\delta_{P^h}(a)\int_{U(F)_{c+c_{1}\sigma(a)}}\Xi^G(u)du.$$
 Il ne reste plus qu'\`a faire appel \`a la proposition 3.2 pour obtenir (1).

 Le groupe $A_{h-1}(F)$ agit \`a gauche sur $U_{h-1}(F)\backslash G_{h-1}(F)$. Introduisons le sous-groupe $\omega_{[1,h-1]}(c+c_{\psi})$ de  $A_{h-1}(F)$. Dans la d\'efinition de $I_{c}^h(e',e)$, on peut remplacer $g$ par $\gamma'g$ pour $\gamma'\in \omega_{[1,h-1]}(c+c_{\psi})$. On peut ensuite int\'egrer en $\gamma'$, \`a condition de diviser le tout par $mes(\omega_{[1,h-1]}(c+c_{\psi}))$. Puisque $\omega_{[1,h-1]}(c+c_{\psi})\omega_{[h,k-1]}(c+c_{\psi})=\omega_{[1,k-1]}(c+c_{\psi})$, on obtient
 $$I_{c}^h(e',e)=mes(\omega_{[1,h-1]}(c+c_{\psi}))^{-1}\int_{\omega_{[1,k-1]}(c+c_{\psi})}\int_{U_{h-1}(F)\backslash GL_{h-1}(F)}$$
 $${\cal L}_{\mu}(\mu(\gamma g)e',\mu(\gamma g)e)\vert det(g)\vert _{F}^{h-k}dg \,d\gamma.$$
 Le m\^eme calcul qu'en 3.5 montre que
 $$\int_{\omega_{[1,k-1]}(c+c_{\psi})}{\cal L}_{\mu}(\mu(\gamma g)e',\mu(\gamma g)e)d\gamma=mes(\omega_{[1,k-1]}(c+c_{\psi})){\cal L}_{\mu,c}(\mu(g)e',\mu(g)e).$$
 D'o\`u
 $$(2)\qquad I_{c}^h(e',e)=mes(\omega_{[h,k-1]}(c+c_{\psi}))\int_{U_{h-1}(F)\backslash GL_{h-1}(F)}{\cal L}_{\mu,c}(\mu(g)e',\mu(g)e)\vert det(g)\vert _{F}^{h-k}dg.$$
 
 On d\'emontre maintenant le lemme par r\'ecurrence sur $h$. Pour $h=1$, l'int\'egrale ci-dessus dispara\^{\i}t et $I_{c}^1(e',e)=C{\cal L}_{\mu,c}(e',e)$, o\`u $C>0$. Mais $U^h=U$ et $J_{c}^1(e',e)={\cal L}_{\mu,c}(e',e)$ par d\'efinition. Supposons maintenant $h\geq2$ et le lemme vrai pour $h-1$. Notons $\tilde{Y}$ le sous-groupe des \'el\'ements $y\in GL_{h-1}(F)$ qui v\'erifient

- pour $i=1,...,h-2$, $y_{i,i}=1$;

- pour $i,j=1,...,h-1$, avec $i\not=j$ et $i\not=h-1$, $y_{i,j}=0$.

Par l'application $y\mapsto (y_{h-1,1},...,y_{h-1,h-1})$, $\tilde{Y}$ s'identifie au compl\'ementaire d'un hyperplan (l'hyperplan $y_{h-1,h-1}=0$) dans un espace vectoriel $Y$ de dimension $h-1$ sur $F$. On note $dy$ la  mesure de Haar sur $Y$ et sa restriction \`a $\tilde{Y}$ (ce n'est pas une mesure de Haar sur cet ensemble). On v\'erifie qu'il existe $C_{0}>0$ tel que, pour toute fonction int\'egrable $\varphi$ sur $U_{h-1}(F)\backslash GL_{h-1}(F)$, on ait l'\'egalit\'e
$$\int_{U_{h-1}(F)\backslash GL_{h-1}(F)}\varphi(g)dg=C_{0}\int_{\tilde{Y}}\int_{U_{h-2}(F)\backslash GL_{h-2}(F)}\varphi(g'y)\vert det(g')\vert ^{-1}dg'\vert y_{h-1,h-1}\vert _{F}^{-1}dy.$$
On d\'eduit de cette \'egalit\'e et de (2) la relation
$$I_{c}^h(e',e)=C_{1}\int_{\tilde{Y}}I_{c}^{h-1}(\mu(y)e',\mu(y)e)\vert y_{h-1,h-1}\vert _{F}^{h-k-1}dy$$
pour un $C_{1}>0$ convenable. Gr\^ace \`a l'hypoth\`ese de r\'ecurrence, on en d\'eduit qu'il existe $C_{2}>0$ tel que 
$$I_{c}^h(e',e)=C_{2}\int_{\tilde{Y}}J_{c}^{h-1}(\mu(y)e',\mu(y)e)\vert y_{h-1,h-1}\vert _{F}^{h-k-1}dy$$
$$=C_{2}\int_{\tilde{Y}} \int_{U^{h-1}(F)_{c}}(\mu(y)e',\mu(uy)e)\bar{\xi}(u)du\vert y_{h-1,h-1}\vert _{F}^{h-k-1}dy.$$
Pour $n\in {\mathbb N}$, notons $Y_{n}$ le sous-ensemble des $y\in Y$ tels que $val_{F}(y_{h-1,j})\geq -n$ pour tout $j=1,...,h-1$. On a
$$(3) \qquad I_{c}^h(e',e)=C_{2}lim_{n\to \infty}X_{n},$$
o\`u
$$X_{n}=\int_{Y_{n}\cap \tilde{Y}}\int_{U^{h-1}(F)_{c}}(\mu(y)e',\mu(uy)e)\bar{\xi}(u)du\vert y_{h-1,h-1}\vert _{F}^{h-k-1}dy.$$
Cette derni\`ere expression est absolument convergente. En effet, rempla\c{c}ons tous les termes par leurs valeurs absolues. D'apr\`es (1), l'int\'egrale int\'erieure est major\'ee par $\sigma(y)^R\delta_{P^{h-1}}(y)=\sigma(y)^R\vert y_{h-1,h-1}\vert _{F}^{k+1-h}$. L'expression totale est donc major\'ee par
$$\int_{Y_{n}}\sigma(y)^Rdy$$
qui est convergente. Posons $Z=U^{h-1}(F)\cap M^h(F)$. Pour $y\in Y$ et $z\in Z$, posons $x(y,z)=\sum_{i=1,...,h-1}y_{h-1,i}z_{i,h}$. Pour $y\in Y$, notons $Z(y)$ le sous-ensemble des $z\in Z$ tels que $val_{F}(x(y,z))\geq-c$ et , pour $z\in Z$, notons $Y(z)$ le sous-ensemble des $y\in Y$ v\'erifiant la m\^eme condition. Dans $X_{n}$, effectuons le changement de variable $u\mapsto yuy^{-1}$. Cela remplace le domaine d'int\'egration $U^{h-1}(F)_{c}$ par $U^h(F)_{c}Z(y)$, donc la variable $u$ par $vz$, avec $v\in U^h(F)_{c}$ et $z\in Z(y)$, la mesure $du $ par $\vert y_{h-1,h-1}\vert _{F}^{k+1-h}dz\,dv$ et  $\xi(u)$ par  $\xi(v)\psi(x(y,z))$.  On obtient
  $$X_{n}=\int_{Y_{n}\cap \tilde{Y}}\int_{U^{h}(F)_{c}}\int_{Z(y)}(e',\mu(vz)e)\bar{\xi}(v)\bar{\psi}(x(y,z))dz\,dv\,dy.$$
  D'apr\`es l'absolue convergence de cette expression, on peut permuter les int\'egrales et on obtient
  $$X_{n}=\int_{U^h(F)_{c}}X_{n}(v)\bar{\xi}(v)dv,$$
  o\`u
  $$X_{n}(v)=\int_{Z}(e',\mu(vz)e)\int_{Y_{n}\cap Y(z)\cap \tilde{Y}} \bar{\psi}(x(y,z))dy\,dz$$
  $$=\int_{Z}(e',\mu(vz)e)\int_{Y_{n}\cap Y(z)} \bar{\psi}(x(y,z))dy\,dz.$$
  Fixons $z\in Z$.  L'ensemble $Y_{n}\cap Y(z)$ est un $\mathfrak{o}_{F}$-r\'eseau dans  $Y$ et l'application $y\mapsto \bar{\psi}(x(y,z))$ est un caract\`ere de $Y$. Son int\'egrale sur le r\'eseau est nulle si le caract\`ere y est non trivial et vaut la mesure du r\'eseau si le caract\`ere y est trivial. On a
  
  (4) le caract\`ere $y\mapsto \bar{\psi}(x(y,z))$ est trivial sur $Y_{n}\cap Y(z)$ si et seulement si $val_{F}(z_{i,h})\geq n+c_{\psi}$ pour tout $i=1,...,h-1$.
  
  En effet, ce caract\`ere est trivial si et seulement si $val_{F}(x(y,z))\geq c_{\psi}$ pour tout $y\in Y_{n}\cap Y(z)$. Cette condition est satisfaite si $z$ v\'erifie les conditions de (3). Inversement, s'il existe $i$ tel que $val_{F}(z_{i,h})<n+c_{\psi}$, soit $y\in Y$ dont la seule coordonn\'ee non nulle soit $y_{h-1,i}$ de valuation $c_{\psi}-val_{F}(z_{i,h)}-1$. Ce nombre est $\geq -n$, donc $y\in Y_{n}$. On a $val_{F}(x(y,z))=c_{\psi}-1$. On a suppos\'e $c+c_{\psi}\geq1$, donc $val_{F}(x(y,z))\geq-c$, ce qui entra\^{\i}ne $y\in Y(z)$. La condition $val_{F}(x(y,z))\geq c_{\psi}$ n'est pas v\'erifi\'ee pour cet $y$. D'o\`u (4).
  
 Notons $Z^n$ l'ensemble des $z\in Z$ v\'erifiant les conditions de (4). Alors
 $$X_{n}(v)=\int_{Z^n}(e',\mu(vz)e)mes(Y_{n}\cap Y(z))dz.$$
 L'in\'egalit\'e $c+c_{\psi}\geq1$ entra\^{\i}ne que, pour $z\in Z^n$, on a $Y_{n}\cap Y(z)=Y_{n}$. D'autre part, si $n$ est assez grand, on a $\mu(z)e=e$ pour tout $z\in Z^n$. Alors
 $$X_{n}(v)=mes(Y_{n})mes(Z^n)(e',\mu(v)e).$$
 Le produit des mesures est une constante positive, disons $C_{3}$. Pour $n$ assez grand, on a donc
 $$X_{n}=C_{3}\int_{U^h(F)_{c}}(e',\mu(v)e)\bar{\xi}(v)dv=C_{3}J_{c}^h(e',e).$$
 En reportant cette \'egalit\'e dans (3), on obtient
 $$I_{c}^h(e',e)=C_{2}C_{3}J_{c}^h(e',e),$$
 ce qui ach\`eve la preuve. $\square$
 
 \bigskip
 
 \subsection{Propri\'et\'es des fonctionnelles de Whittaker}
 
 Soient $\mu\in Temp(G)$. Appliquons le lemme pr\'ec\'edent pour $h=k$. On obtient une \'egalit\'e
 $$\int_{U_{k-1}(F)\backslash GL_{k-1}(F)}{\cal L}_{\mu}(\mu(g)e',\mu(g)e)dg=C(e',e)$$
 pour tous $e,e'\in E_{\mu}$, o\`u $C$ est une constante positive. Il en r\'esulte que ${\cal L}_{\mu}$ est non nulle, autrement dit que la constante $C_{\mu}$ du paragraphe 3.5 est non nulle. On peut alors r\'ecrire la relation 3.5(1) et le lemme 3.6 sous la forme suivante.
 
 \ass{Lemme}{Soient $\mu\in Temp(G)$ et $\phi$ une fonctionnelle de Whittaker non nulle sur $E_{\mu}$. 
 
(i) Il existe un r\'eel $R$ et, pour tout $e\in E_{\mu}$, il existe $c'\in {\mathbb N}$ tel que
 $$\vert \phi(\mu(a)e)\vert << \iota_{c'}(a)\delta_{B}(a)^{1/2}\sigma(a)^R$$
 pour tout $a\in A(F)$.
 
 (ii) Pour tout $h=1,...,k$ et tout entier $c$ tel que $c\geq1$ et $c+c_{\psi}\geq1$, il existe $C>0$ tel que l'on ait l'\'egalit\'e
   $$ \int_{\omega_{[h,k-1]}(c+c_{\psi})}\int_{U_{h-1}(F)\backslash GL_{h-1}(F)}\overline{\phi(\mu(ag)e')}\phi(\mu(ag)e)\vert det(g)\vert _{F}^{h-k}dg \,da$$
 $$=C\int_{U^h(F)_{c}}(e',\mu(u)e)\bar{\xi}(u)du$$
pour tous $e,e'\in E_{\mu}$.}

\bigskip

\section{Majorations pour un groupe sp\'ecial orthogonal}

\bigskip

\subsection{Les groupes sp\'eciaux orthogonaux}

Soit $(V,q_{V})$ un espace quadratique sur $F$, c'est-\`a-dire que $V$ est un espace vectoriel de dimension finie sur $F$ et $q_{V}$ est une forme bilin\'eaire sym\'etrique non d\'eg\'en\'er\'ee sur $V$ (on dira souvent que $V$ est un espace quadratique, la forme $q_{V}$ \'etant sous-entendue).  On note aussi $q_{V}$ la forme quadratique d\'efinie par  $q_{V}(v)=q_{V}(v,v)/2$. On note $d_{V}$ la dimension de $V$ et $G$ le groupe sp\'ecial orthogonal de $V$. Consid\'erons un syst\`eme hyperbolique maximal $(v_{\pm i})_{i=1,...,l}$ dans $V$ ("syst\`eme hyperbolique" signifie que $q_{V}(v_{i},v_{j})=\delta_{i,-j}$ pour tous $i,j$, o\`u $\delta_{i,-j}$ est le symbole de Kronecker). Notons $Z$ le sous-espace de $V$ engendr\'e par ce syst\`eme et $V_{an}$ l'orthogonal de $Z$ dans $V$. La restriction $q_{V_{an}}$ de $q_{V}$ \`a $V_{an}$ est anisotrope. On note $d_{an,V}$ la dimension de $V_{an}$. Fixons un r\'eseau sp\'ecial $R_{an}\subset V_{an}$ ([W1] 7.1). On peut choisir un r\'eseau $R_{Z}$ de $Z$ ayant une base form\'ee de vecteurs proportionnels aux $v_{i}$, de sorte que $R=R_{Z}\oplus R_{an}$ soit sp\'ecial.   On note $K$ le stabilisateur de $R$ dans $G(F)$. C'est un sous-groupe compact sp\'ecial de $G(F)$. Consid\'erons une suite d'entiers $(k_{1},...,k_{s})$ telle que $k_{j}\geq1$ pour tout $j$ et $\sum_{j=1,...,s}k_{j}\leq l$. Pour tout $j$, notons $Z_{j}$, resp. $Z_{-j}$, le sous-espace de $V$ engendr\'e par les $v_{i}$, resp. $v_{-i}$, pour $i=k_{1}+...+k_{j-1}+1,...,k_{1}+...+k_{j}$. Notons $\tilde{V}$ l'orthogonal dans $V$ de la somme des $Z_{\pm j}$ et notons $\tilde{G}$ le groupe sp\'ecial orthogonal de $\tilde{V}$. Notons $P$ le sous-groupe parabolique de $G$ form\'e des \'el\'ements qui conservent le drapeau de sous-espaces
$$Z_{1}\subset Z_{1}\oplus Z_{2}\subset...\subset Z_{1}\oplus...\oplus Z_{s}.$$
Notons $M$ la composante de L\'evi de $P$ form\'ee des \'el\'ements qui conservent chaque sous-espace $Z_{\pm j}$. On a
$$(1) \qquad M\simeq GL_{k_{1}}\times...\times GL_{k_{s}}\times \tilde{G}.$$
On sait que $K$ est en bonne position relativement \`a $M$. Inversement, si $P=MU$ est un sous-groupe parabolique de $G$, si $K$ est un sous-groupe compact sp\'ecial de $G(F)$ en bonne position relativement \`a $M$, on peut trouver un syst\`eme hyperbolique, un r\'eseau sp\'ecial et une suite d'entiers de sorte que $P$, $M$ et $K$ soient d\'etermin\'es comme ci-dessus (ces donn\'ees ne sont pas uniques). Si  $s=l$ et $k_{j}=1$ pour tout $j$, $M$ est un L\'evi minimal et inversement, si $M$ est un L\'evi minimal, on peut supposer ces \'egalit\'es v\'erifi\'ees. 

Dans la situation ci-dessus, supposons $M$ minimal et notons-le plut\^ot $M_{min}$.  L'application naturelle $K\cap Norm_{G(F)}(M_{min})\to W^G$ est surjective et il est utile de remarquer qu'elle poss\`ede une section $\iota:W^G\to K\cap Norm_{G(F)}(M_{min})$ qui est un homomorphisme de groupes. En effet, pour tout $i=\pm 1,...,\pm l$, fixons un \'el\'ement $v'_{i}\in Fv_{i}$ de sorte que $(v'_{i})_{i=\pm 1,...,\pm l}$ soit une base sur $\mathfrak{o}_{F}$ de $R_{Z}$. Parce que $R$ est un r\'eseau sp\'ecial, on peut supposer que $q_{V}(v'_{i},v'_{-i})=q_{V}(v'_{i'},v'_{-i'})$ pour tous $i,i'=1,...,l$. Si $V_{an}\not=\{0\}$, fixons un \'el\'ement $g_{an}$ du groupe orthogonal de $V_{an}$  tel que $det(g_{an})=-1$ et $g_{an}^2=1$. L'action de tout \'el\'ement de ce groupe orthogonal, en particulier l'action de $g_{an}$, conserve le r\'eseau $R_{an}$ . Notons $W^K$ le sous-ensemble des \'el\'ements $g\in G(F)$ qui agissent par permutation sur l'ensemble $\{v_{\pm i}; i=1,...,l\}$ et agissent sur $V_{an}$,  soit par l'identit\'e, soit comme $g_{an}$. On v\'erifie que $W^K$ est un sous-groupe de $K\cap Norm_{G(F)}(M_{min})$ et que l'application naturelle de $W^K$ dans $W^G$ est un isomorphisme.

Les hypoth\`eses de 1.5 sont v\'erifi\'ees pour le groupe $G$. Soient $M$ un L\'evi que l'on \'ecrit sous la forme (1) et $\tau$ une repr\'esentation admissible irr\'eductible et de la s\'erie discr\`ete de $M(F)$. On a
 $$\tau\simeq \mu_{1}\otimes...\otimes\mu_{s}\otimes\tilde{\tau},$$
 o\`u $\mu_{j}$, resp. $\tilde{\tau}$, est une repr\'esentation de la s\'erie discr\`ete de $GL_{k_{j}}(F)$, resp. $\tilde{G}(F)$. L'espace ${\cal A}_{M}$ s'identifie naturellement \`a ${\mathbb R}^s$. Supposons  $R(\tau)\cap W(M)_{reg}\not=\emptyset$. Alors le groupe $R(\tau)$ s'identifie \`a un sous-groupe de $\{\pm 1\}^s$. Un \'el\'ement $\epsilon=(\epsilon_{1},...,\epsilon_{s})$ de ce groupe agit sur ${\mathbb R}^s$ par
 $$(x_{1},...,x_{s})\mapsto (\epsilon_{1}x_{1},...,\epsilon_{s}x_{s}).$$
 Le groupe $R(\tau)$ contient l'\'el\'ement $t=(-1,...,-1)$ de $ \{\pm 1\}^s$, qui est l'unique \'el\'ement de $R(\tau)\cap W(M)_{reg}$. On a $\vert det(t-1)_{\vert {\cal A}_{M}}\vert =2^{a_{M}}$. 
 
 Soit $\pi$ une repr\'esentation temp\'er\'ee irr\'eductible et elliptique de $G(F)$. On peut trouver $M$, $\tau$ comme ci-dessus, et $\zeta\in R(\tau)^{\vee}$ de sorte que $\pi=Ind_{P}^G(\tau,\zeta)$, o\`u $P$ est un \'el\'ement de ${\cal P}(M)$. Puisque la classe de conjugaison du couple $(M,\tau)$ est bien d\'etermin\'ee, on peut poser $r(\pi)=\vert  R(\tau)\vert $ et $t(\pi)=2^{a_{M}}$.

\bigskip

\subsection{Espaces quadratiques compatibles}

 Soient $(V,q_{V})$ et $(W,q_{W})$ deux espaces quadratiques. Notons $G$ et $H$ leurs groupes sp\'eciaux orthogonaux,  $d_{V}$ et $d_{W}$ les dimensions de $V$ et $W$. Supposons par exemple $d_{W}\leq d_{V}$. On dit que les deux espaces quadratiques sont compatibles si $d_{V}$ et $d_{W}$ sont de parit\'es distinctes et si $W$ est isomorphe (comme espace quadratique) \`a un sous-espace de $V$ dont l'orthogonal est d\'eploy\'e, c'est-\`a-dire une somme orthogonale $D_{0}\oplus Z$, o\`u $D_{0}$ est une droite et $Z$ est hyperbolique. On peut alors identifier $W$ \`a un sous-espace de $V$ et $H$ \`a un sous-groupe de $G$. D'apr\`es le th\'eor\`eme de Witt, cette identification est unique \`a conjugaison pr\`es par un \'el\'ement de $G(F)$.

Supposons que $V$ soit la somme directe orthogonale de deux sous-espaces $V'$ et $Z'$, avec $Z'$ hyperbolique. Alors $V$ et $W$ sont compatibles si et seulement si $V'$ et $W$ le sont.

Soient $V$ et $W$ deux espaces quadratiques compatibles, avec $d_{W}<d_{V}$. On fixe un isomorphisme $V=W\oplus D_{0}\oplus Z$ avec les propri\'et\'es ci-dessus, une base hyperbolique $(v_{\pm i})_{i=1,...,r}$ de $Z$ et un \'el\'ement non nul $v_{0}\in D_{0}$. On pose $V_{0}=W\oplus D_{0}$ et on note $G_{0}$ son groupe sp\'ecial orthogonal. On note $A$ le sous-tore maximal du groupe sp\'ecial orthogonal de $Z$ qui conserve chaque droite $Fv_{\pm i}$. Pour $a\in A(F)$ et $i=\pm 1,...,\pm r$, on note $a_{i}$ la valeur propre de $a$ sur le vecteur $v_{i}$. On note $P$ le sous-groupe parabolique de $G$ form\'e des \'el\'ements qui conservent le drapeau
$$Fv_{r}\subset Fv_{r}\oplus Fv_{r-1}\subset...\subset Fv_{r}\oplus...\oplus Fv_{1}$$
de $V$. On note $U$ le radical unipotent de $P$ et $M$ sa composante de L\'evi qui contient $A$. On a l'\'egalit\'e $M=AG_{0}$. On d\'efinit un caract\`ere $\xi$ de $U(F)$ par
$$\xi(u)=\psi(\sum_{i=0,...,r-1}q_{V}(uv_{i},v_{-i-1})).$$

On fixe un r\'eseau sp\'ecial $R_{0}\subset V_{0}$, cf. [W1] 7.1. On choisit, ainsi qu'il est loisible, un r\'eseau $R_{Z}\subset Z$ poss\'edant une base sur $\mathfrak{o}_{F}$ form\'ee de vecteurs proportionnels aux $v_{\pm i}$ et tel que le r\'eseau $R=R_{0}\oplus R_{Z}$ soit sp\'ecial. On note $K$ le stabilisateur de $R$ dans $G(F)$. Pour un entier $N\geq1$, on d\'efinit une fonction $\kappa_{N}$ sur $G(F)$ de la fa\c{c}on suivante. Elle est invariante \`a gauche par $U(F)$ et \`a droite par $K$. Sa restriction \`a $M(F)$ est la fonction caract\'eristique des \'el\'ements $ag_{0}$, avec $a\in A(F)$ et $g_{0}\in G_{0}(F)$, qui v\'erifient les conditions $\vert val_{F}(a_{i})\vert \leq N$ pour tout $i=1,...,r$ et $g_{0}^{-1}e_{0}\in \mathfrak{p}_{F}^{-N}R_{0}$.

{\bf Remarque.} Les constructions et notations ci-dessus seront utilis\'ees sans plus de commentaires chaque fois que l'on se donnera des espaces quadratiques compatibles $V$ et $W$ avec $d_{W}<d_{V}$.

\bigskip

\subsection{Les r\'esultats}

On \'enonce ici toutes les majorations que la section est destin\'ee \`a prouver. On fixe pour toute cette section deux espaces quadratiques compatibles $(V,q_{V})$ et $(W,q_{W})$ tels que $d_{V}>d_{W}$.

(1) Il existe un r\'eel $R$ tel que
$$\int_{G(F)}{\bf 1}_{\sigma< b}(g)dg<<exp(Rb)$$
pour tout r\'eel $b\geq0$.

{\bf Remarque.} Cette assertion vaut en fait pour tout groupe r\'eductif connexe $G$.

Pour un entier $N\geq1$ et un r\'eel $D$, posons
$$I(N,D)=\int_{G(F)}\Xi^G(g)^2\kappa_{N}(g)\sigma(g)^Ddg.$$

(2) Cette int\'egrale est convergente; le r\'eel $D$ \'etant fix\'e, il existe un r\'eel $R$ tel que
$$I(N,D)<<N^R$$
pour tout entier $N\geq1$.  

Pour $u\in U(F)$ et tout $i=1,...,r$, on note $u_{i,i-1}$ la coordonn\'ee $u_{i,i-1}=q_{V}(uv_{i-1},v_{-i})$. Pour tout entier $c\geq1$, on note $U(F)_{c}$ l'ensemble des $u\in U(F)$ tels que $val_{F}(u_{i,i-1})\geq-c$ pour tout $i=1,...,r$. C'est un sous-groupe de $U(F)$ conserv\'e par conjugaison par $H(F)$.
 Pour un r\'eel $D$ et un \'el\'ement $m\in M(F)$, on pose
 $$X(c,D,m)=\int_{U(F)_{c}}\Xi^G(um)\sigma(um)^Ddu.$$
 
 (3) Cette expression est convergente. Pour $D$ fix\'e, il existe un r\'eel $R$ tel que
 $$X(c,D,m)<<c^R\sigma(m)^R\delta_{P}(m)^{1/2}\Xi^M(m)$$
 pour tous $c\geq1$ et tout $m\in M(F)$.

(4) Pour tout r\'eel $D$ et tout entier $c\geq1$, l'int\'egrale
$$\int_{H(F)U(F)_{c}}\Xi^H(h)\Xi^G(hu)\sigma(hu)^Ddu\,dh$$
est convergente.

(5) Pour tout r\'eel $D$ et tout entier $c\geq1$, l'int\'egrale
$$\int_{H(F)U(F)_{c}}\int_{H(F)U(F)_{c}}\Xi^G(hu)\Xi^H(h'h)\Xi^G(h'u')\sigma(hu)^D\sigma(h'u')^Ddu'\,dh'\,du\,dh$$
est convergente.

Soient $D$ et $C$ deux r\'eels, $c$,  $c'$ et $N$  trois entiers. On suppose $C,c,c',N\geq1$. Pour $m\in M(F)$, $h\in H(F)$, $u,u'\in U(F)$, posons
$$\phi(m,h,u,u';D)=\Xi^H(h)\Xi^G(u'm)\Xi^G(u^{-1}h^{-1}u'm)\kappa_{N}(m)\sigma(u')^D\sigma(u)^D\sigma(h)^D\sigma(m)^D\delta_{P}(m)^{-1}.$$
Posons
  $$I(c,N,D)=\int_{M(F)}\int_{H(F)U(F)_{c}}\int_{U(F)} \phi(m,h,u,u';D)du'\,du\,dh\,dm,$$
 $$I(c,c',N,D)=\int_{M(F)}\int_{H(F)U(F)_{c}}\int_{U(F)-U(F)_{c'}} \phi(m,h,u,u';D)du'\,du\,dh\,dm,$$
$$I(c,c',N,C,D)=\int_{M(F)}\int_{H(F)U(F)_{c}}\int_{U(F)_{c'}}{\bf 1}_{\sigma\geq C}(hu) \phi(m,h,u,u';D)du'\,du\,dh\,dm.$$ 

(6) L'int\'egrale $I(c,N,D)$ est convergente; les termes $c$ et $D$ \'etant fix\'es, il existe un r\'eel $R$ tel que
$$I(c,N,D)<<N^R$$
pour tout $N\geq1$.

(7) L'int\'egrale $I(c,c',N,D)$ est convergente; les termes $c$ et $D$ \'etant fix\'es, pour tout r\'eel $R$, il existe $\alpha>0$  tel que
$$I(c,c',N,D)<<N^{-R}$$
pour tout $N\geq2$ et tout $c'\geq \alpha log(N)$.

(8) l'int\'egrale $I(c,c',N,C,D)$ est convergente; les termes $c$ et $D$ \'etant fix\'es, pour tout r\'eel $R$, il existe $\alpha>0$ tel que
$$I(c,c',N,C,D)<<N^{-R}$$
pour tout $N\geq1$, tout $c'\geq1$ et tout $C\geq \alpha(log(N)+c')$. 

\bigskip

\subsection{Preuve de la majoration 4.3(1)}

Fixons un L\'evi minimal $M_{min}$ de $G$ tel que $K$ soit en bonne position relativement \`a $M_{min}$. Soit $P_{min}=M_{min}U_{min}\in {\cal P}(M_{min})$. On a
$$\int_{G(F)}{\bf 1}_{\sigma<b}(g)dg=\int_{K}\int_{U_{min}(F)}\int_{M_{min}(F)}{\bf 1}_{\sigma<b}(muk)dm\,du\,dk.$$
D'apr\`es [W2] lemme II.3.1, il existe $c_{1}>0$ tel que la relation $\sigma(muk)<b$ entra\^{\i}ne $\sigma(m)<c_{1}b$, $\sigma(u)<c_{1}b$. Alors
$$(1) \qquad \int_{G(F)}{\bf 1}_{\sigma<b}(g)dg<<\int_{U_{min}(F)}{\bf 1}_{\sigma< c_{1}b}(u)du\int_{M_{min}(F)}{\bf 1}_{\sigma< c_{1}b}(m)dm.$$
Montrons que la premi\`ere int\'egrale v\'erifie la condition requise. On peut fixer des coordonn\'ees $(u_{j})_{j\in J}$ sur $U_{min}$ de sorte que $du=\prod_{j\in J}du_{j}$ et que la condition ${\bf 1}_{\sigma< c_{1}b}(u)$ entra\^{\i}ne $val_{F}(u_{j})>-c_{2}b$ pour tout $j$, pour une constante $c_{2}>0$ convenable. Or il existe $c_{3}$ tel que la mesure de l'ensemble $\{x\in F; val_{F}>-c_{2}b\}$ soit born\'ee par $exp(c_{3}b)$, d'o\`u le r\'esultat. On doit borner la seconde int\'egrale de  (1), ce qui nous ram\`ene au cas o\`u $G=M_{min}$. Dans ce cas, $G(F)$ est le produit de $A_{G}(F)$ et d'un sous-ensemble compact. Cela nous ram\`ene au cas o\`u $G=A_{G}$ est un tore d\'eploy\'e.  On se ram\`ene imm\'ediatement au cas o\`u ce tore est de dimension $1$. Alors il existe une constante $c_{4}>0$ telle que notre int\'egrale soit la mesure de l'ensemble $\{x\in F^{\times}; \vert val_{F}(x)\vert <c_{4}b\}$ (pour la mesure de Haar multiplicative). Cette mesure est essentiellement born\'ee par $b$. $\square$

\bigskip

\subsection{Majoration d'une int\'egrale unipotente, cas $r\geq2$}

 Supposons $r\geq1$. Pour tout $i=1,...,r$, on note $P_{i}$ le sous-groupe parabolique
des \'el\'ements de $G$ qui conservent le drapeau
 $$Fv_{r}\subset Fv_{r}\oplus Fv_{r-1}\subset... \subset Fv_{r}\oplus...\oplus Fv_{i}.$$
 On note $U_{i}$ le radical unipotent de $P_{i}$ et $M_{i}$ la composante de L\'evi qui contient $M$.
  Notons $U_{r,\natural}$ le sous-groupe des \'el\'ements $u\in U_{r}$ tels que $u_{r,r-1}=0$.  Pour un r\'eel $D$, posons
$$I_{r,\natural}(b,D)=\int_{U_{r,\natural}(F)}{\bf 1}_{\sigma\geq b}(u)\delta_{\bar{P}}(m_{\bar{P}}(u))^{1/2}\Xi^M(m_{\bar{P}}(u))\sigma(u)^Ddu.$$

\ass{Lemme}{Supposons $r\geq2$.  L'int\'egrale ci-dessus est convergente. Pour $D$ fix\'e, il existe $\epsilon>0$ tel que
$$I_{r,\natural}(b,D)<<exp(-\epsilon b)$$
pour tout $b\geq0$.}

Preuve. Notons $V_{\flat}$ l'orthogonal dans $V$ de l'espace de dimension $4$ engendr\'e par $v_{r}$, $v_{r-1}$, $v_{1-r}$ et $v_{-r}$.  Pour $x\in V_{\flat}$ et $y\in F$, notons $u(x,y)$ l'unique \'el\'ement de $U_{\natural}(F)$ tel que $u(x,y)v_{-r}=v_{-r}+x+yv_{r-1}-q_{V}(x)v_{r}$. Alors $(x,y)\mapsto u(x,y)$ est un isomorphisme de $V_{\flat}\times F$ sur $U_{\natural}(F)$. On pose simplement $u(x)=u(x,0)$ et $\underline{y}=u(0,y)$. Une description analogue vaut sur le corps de base $\bar{F}$. On note $U_{r,\sharp}$, resp. $Y$, le sous-groupe de $U_{r,\natural}$ form\'e des $u(x)$, resp. $\underline{y}$.
 On a $U_{r,\natural}=U_{r,\sharp}Y$. Notons $V_{\sharp}$ l'orthogonal dans $V$ du plan engendr\'e par $v_{r-1}$ et $v_{1-r}$. Notons $G_{\sharp}$ le groupe sp\'ecial orthogonal de $V_{\sharp}$ et affectons d'un indice $\sharp$ les intersections avec $G_{\sharp}$ des groupes que l'on a introduits. En particulier $P_{r,\sharp}=G_{\sharp}\cap P_{r}$. On voit que $U_{r,\sharp}$ n'est autre que le radical unipotent de $P_{r,\sharp}$.  Soit $x\in V_{\flat}$, introduisons l'\'el\'ement $m_{\bar{P}_{\sharp}}(u(x))$, que l'on \'ecrit  
 $a(x)g_{0}(x)$, avec $a(x)\in A(F)$, $g_{0}(x)\in G_{0}(F)$. On a $a(x)_{r-1}=1$ puisque  $a(x)\in G_{\sharp}(F)$. Pour tout $v\in V$, notons $val_{R}(v)$ le plus grand entier $n\in {\mathbb Z}$ tel que $v\in \mathfrak{p}_{F}^nR$. On a
 
  (1) il existe $c_{1},c_{2}\in {\mathbb Z}$ tel que $val_{F}(a(x)_{r})=inf(0,c_{1}+val_{R}(x),c_{2}+val_{F}(q_{V}(x)))$; il existe $\epsilon_{1}>0$ tel que $\vert a(x)_{r}\vert _{F}^{-1}<< exp(-\epsilon_{1}\sigma(u(x)))$.
 
 En effet, posons $k=u(x)^{-1}m_{\bar{P}_{\sharp}}(u(x))u_{\bar{P}_{\sharp}}(u(x))$. On a $k=k_{\bar{P}_{\sharp}}(u(x))^{-1}\in K$.  Donc $val_{R}(kv_{-r})=val_{R}(v_{-r})$. Mais $kv_{-r}=a(x)_{r}^{-1}(v_{-r}-x-q_{V}(x)v_{r})$, d'o\`u
 $$val_{R}(kv_{-r})=-val_{F}(a(x))+val_{R}(v_{-r}-x-q_{V}(x)v_{r}).$$
 D'apr\`es la d\'efinition de $R$, on a
 $$val_{R}(v_{-r}-x-q_{V}(x)v_{r})=inf(val_{R}(v_{-r}),val_{R}(x),val_{F}(q_{V}(x)+val_{R}(v_{r})).$$
 En utilisant toutes ces \'egalit\'es, on obtient
 $$val_{F}(a(x)_{r})=inf(0,-val_{R}(v_{-r})+val_{R}(x),val_{R}(v_{r})-val_{R}(v_{-r})+val_{F}(q_{V}(x))).$$
D'o\`u la premi\`ere assertion de (1). La seconde s'en d\'eduit imm\'ediatement.

On a l'\'egalit\'e $m_{\bar{P}}(u(x))=m_{\bar{P}_{\sharp}}(u(x))$. On calcule
$$\delta_{\bar{P}}(m_{\bar{P}}(u(x)))^{1/2}\Xi^M(m_{\bar{P}}(u(x)))=\Xi^{G_{0}}(g_{0}(x))\prod_{i=1,...,r}\vert a(x)_{i}\vert _{F}^{1-i-d_{V_{0}}/2},$$
$$\delta_{\bar{P}_{\sharp}}(m_{{\bar{P}_{\sharp}}}(u(x)))^{1/2}\Xi^{M_{\sharp}}(m_{\bar{P}_{\sharp}}(u(x)))=\Xi^{G_{0}}(g_{0}(x))\vert a(x)_{r}\vert_{F}^{2-r-d_{V_{0}}/2}\prod_{i=1,...,r-2}\vert a(x)_{i}\vert _{F} ^{1-i-d_{V_{0}}/2},$$
d'o\`u
$$(2) \qquad \delta_{\bar{P}}(m_{\bar{P}}(u(x)))^{1/2}\Xi^M(m_{\bar{P}}(u(x)))=\vert a(x)_{r}\vert _{F}^{-1}\delta_{\bar{P}_{\sharp}}(m_{{\bar{P}_{\sharp}}}(u(x)))^{1/2}\Xi^{M_{\sharp}}(m_{\bar{P}_{\sharp}}(u(x))).$$

Soient $x\in V_{\flat}$ et $y\in F$. On a
$$m_{\bar{P}}(u(x)\underline{y})=m_{\bar{P}}(u(x))m_{\bar{P}}(k_{\bar{P}}(u(x)\underline{y}).$$
Il existe donc $R_{1}>0$ tel que
$$(3) \qquad \delta_{\bar{P}}(m_{\bar{P}}(u(x)\underline{y}))^{1/2}\Xi^M(m_{\bar{P}}(u(x)\underline{y}))<<\delta_{\bar{P}}(m_{\bar{P}}(u(x)))^{1/2}\Xi^M(m_{\bar{P}}(u(x))exp(R_{1} \sigma(y)),$$
o\`u $\sigma(y)=sup(1,-val_{F}(y))$. Introduisons un r\'eel  $\mu>0$  que nous fixerons plus tard et posons
$$I^1_{r,\natural}(b,D)=\int_{y\in F; val_{F}(y)\geq -\mu b}\int_{V_{\flat}}{\bf 1}_{\sigma\geq b}(u(x)\underline{y})\delta_{\bar{P}}(m_{\bar{P}}(u(x)\underline{y}))^{1/2}$$
$$\Xi^M(m_{\bar{P}}(u(x)\underline{y}))\sigma(u(x)\underline{y})^D dx \,dy.$$
Pour $y$ tel que $val_{F}(y)\geq -\mu b$, on a $\vert \sigma(u(x)\underline{y})-\sigma(u(x))\vert <<\mu b$. Si $\mu$ est assez petit, la condition ${\bf 1}_{\sigma\geq b}(u(x)\underline{y})=1$ entra\^{\i}ne ${\bf 1}_{\sigma\geq b/2}(u(x))=1$. Gr\^ace \`a (3), on obtient
$$I^1_{r,\natural}(b,D)<<\int_{V_{\flat}}{\bf 1}_{\sigma\geq b/2}(u(x))\delta_{\bar{P}}(m_{\bar{P}}(u(x)))^{1/2}\Xi^M(m_{\bar{P}}(u(x)))\sigma(u(x))^Ddx$$
$$ \int_{y\in F; val_{F}(y)\geq -\mu b}(\mu b)^Dexp(R_{1} \sigma(y))dy.$$
Il existe $R_{2}$ tel la derni\`ere int\'egrale soit born\'ee par $exp(R_{2}\mu b)$. Dans la premi\`ere int\'egrale, on utilise (1) et (2). Pour ${\bf 1}_{\sigma\geq b/2}(u(x))=1$, on a
 $$\delta_{\bar{P}}(m_{\bar{P}}(u(x)))^{1/2}\Xi^M(m_{\bar{P}}(u(x))<<$$
 $$exp(-\epsilon_{1}b/4-\epsilon_{1}\sigma(u(x))/2)\delta_{\bar{P}}(m_{{\bar{P}_{\sharp}}}(u(x)))^{1/2}\Xi^{M_{\sharp}}(m_{\bar{P}_{\sharp}}(u(x))).$$
Alors
$$I^1_{r,\natural}(b,D)<<exp(R_{2}\mu b-\epsilon_{1}b/4)\int_{V_{\flat}}exp(-\epsilon_{1}\sigma(u(x))/2)\delta_{\bar{P}}(m_{{\bar{P}_{\sharp}}}(u(x)))^{1/2}$$
$$\Xi^{M_{\sharp}}(m_{\bar{P}_{\sharp}}(u(x)))\sigma(u(x))^Ddx.$$
D'apr\`es [W2] lemme II.4.3, cette int\'egrale est convergente. Choisissons $\mu$ tel que $\epsilon_{1}/4-R_{2}\mu>0$ et notons $\epsilon_{2}$ ce dernier terme. On obtient alors
$$I^1_{r,\natural}(b,D)<<exp(-\epsilon_{2} b).$$
Le r\'eel $\mu$ \'etant maintenant fix\'e, posons
$$I^2_{r,\natural}(b,D)=\int_{y\in F; val_{F}(y)< -\mu b}\int_{V_{\flat}}{\bf 1}_{\sigma\geq b}(u(x)\underline{y})\delta_{\bar{P}}(m_{\bar{P}}(u(x)\underline{y}))^{1/2}$$
$$\Xi^M(m_{\bar{P}}(u(x)\underline{y}))\sigma(u(x)\underline{y})^D dx \,dy.$$
Il nous reste \`a majorer cette expression.  Soit $y\in F$ tel que $val_{F}(y)<\mu b$. Le groupe $Y$ est en fait le sous-groupe unipotent $U_{2}$ du groupe $GL_{2}$ agissant dans le plan engendr\'e par $v_{r-1}$ et $v_{-r}$. Dans $GL_{2}$, on a l'\'egalit\'e habituelle
$$(4) \qquad \left(\begin{array}{cc}1&y\\ 0&1\\ \end{array}\right)=\left(\begin{array}{cc}y&0\\ 0&y^{-1}\\ \end{array}\right)\left(\begin{array}{cc}1&0\\ y&1\\ \end{array}\right)\left(\begin{array}{cc}y^{-1}&1\\ -1&0\\ \end{array}\right).$$
On en d\'eduit une \'egalit\'e $\underline{y}=a(y)\bar{n}(y)k(y)$ o\`u 

 $a(y)$ est l'\'el\'ement de $A(F)$ tel que $a(y)_{r}=a(y)_{r-1}=y$ et $a(y)_{i}=1$ pour $i=1,...,r-2$;
 
 $\bar{n}(y)$ est l'\'el\'ement de $U_{\bar{P}}(F)$ qui envoie $v_{r}$ sur $v_{r}-yv_{1-r}$, $v_{r-1}$ sur $v_{r-1}+yv_{-r}$ et fixe les autres vecteurs $v_{\pm i}$ ainsi que $V_{0}$;
 
 $k(y)$ est un \'el\'ement qui reste born\'e.
 
 Soit $x\in V_{\flat}$. On a $u(x)a(y)=a(y)u(y^{-1}x)$. Posons $x'=y^{-1}x$. Un calcul matriciel montre que $u(x')\bar{n}(y)u(x')^{-1}=\bar{n}(x',y)n(x',y)$, o\`u $\bar{n}(x',y)\in U_{\bar{P}}(F)$ et $n(x',y)$ est l'\'el\'ement de $U_{P}(F)$ qui envoie $v_{r-1}$ sur $v_{r-1}-yq_{V}(x')v_{r}$, $v_{-r}$ sur $v_{-r}+yq_{V}(x')v_{1-r}$ et qui fixe les autres $v_{\pm i}$ ainsi que $V_{0}$. D'o\`u
  $$u(x)\underline{y}=a(y)\bar{n}(x',y)n(x',y)u(x')k(y)$$
  $$=a(y)\bar{n}(x',y)n(x',y)a(x')g_{0}(x')u_{\bar{P}_{\sharp}}(u(x'))k_{\bar{P}_{\sharp}}(u(x'))k(y).$$
  Posons $n'(x',y)=a(x')^{-1}n(x',y)a(x')$. C'est l'\'el\'ement de $U_{P}(F)$ qui envoie $v_{r-1}$ sur $v_{r-1}-ya(x')_{r}^{-1}q_{V}(x')v_{r}$, $v_{-r}$ sur $v_{-r}+ya(x')_{r}^{-1}q_{V}(x')v_{1-r}$ et qui fixe les autres $v_{\pm i}$ ainsi que $V_{0}$. Cet \'el\'ement $n'(x',y)$ appartient au sous-groupe unipotent $U_{2}$ du groupe $GL_{2}$ agissant dans le plan engendr\'e par $v_{1-r}$ et $v_{-r}$. En utilisant (4), on a $n'(x',y)=a(x',y)\bar{n}'(x',y)k(x',y)$ o\`u $\bar{n}'(x',y)\in U_{\bar{P}}(F)$, $k(x',y)$ reste born\'e et $a(x',y)$ est l'\'el\'ement de $A(F)$ tel que
  
- si $val_{F}(q_{V}(x'))+val_{F}(y)-val_{F}(a(x')_{r})\geq0$, $a(x',y)=1$;
 
 - si $val_{F}(q_{V}(x'))+val_{F}(y)-val_{F}(a(x')_{r})<0$, $a(x',y)_{r}=a(x')_{r}^{-1}yq_{V}(x)$, $a(x',y)_{r-1}=a(x')_{r}y^{-1}q_{V}(x')^{-1}$, et $a(x',y)_{i}=1$ pour tout $i=1,...,r-2$.
 
 D'o\`u 
 $$u(x)\underline{y}=a(y)\bar{n}(x',y)a(x')a(x',y)\bar{n}'(x',y)k(x',y)g_{0}(x')u_{\bar{P}_{\sharp}}(u(x'))k_{\bar{P}_{\sharp}}(u(x'))k(y).$$
 L'\'el\'ement $k(x',y)$ appartient au m\^eme groupe $GL_{2}$ que ci-dessus. Ce groupe commute \`a $G_{0}$. Il conserve le radical unipotent du sous-groupe parabolique de $G$ form\'e des \'el\'ements qui conservent le drapeau
 $$Fv_{-r}\oplus Fv_{1-r}\subset Fv_{-r}\oplus Fv_{1-r}\oplus Fv_{2-r}\subset...\subset Fv_{-r}\oplus...\oplus Fv_{-1}.$$
 Ce radical unipotent est contenu dans $U_{\bar{P}}$ et contient $U_{\bar{P}_{\sharp}}$. Donc  $k(x',y)u_{\bar{P}_{\sharp}}(u(x'))k(x',y)^{-1}\in U_{\bar{P}}(F)$. Finalement
 $$u(x)\underline{y}\in U_{\bar{P}}(F)a(y)a(x')a(x',y)g_{0}(x')\Gamma,$$
 o\`u $\Gamma$ est un ensemble born\'e. On en d\'eduit, gr\^ace \`a (2)
 $$\delta_{\bar{P}}(m_{\bar{P}}(u(x)\underline{y}))^{1/2}\Xi^M(m_{\bar{P}}(u(x)\underline{y}))<<\delta_{\bar{P}}(a(y)a(x',y))^{1/2}\delta_{\bar{P}}(m_{\bar{P}}(u(x')))^{1/2}\Xi^M(m_{\bar{P}}(u(x')))$$
 $$<<\delta_{\bar{P}}(a(y)a(x',y))^{1/2}\vert a(x')_{r}\vert _{F}^{-1}\delta_{\bar{P}_{\sharp}}(m_{\bar{P}_{\sharp}}(u(x')))^{1/2}\Xi^{M_{\sharp}}(m_{\bar{P}_{\sharp}}(u(x'))).$$
 On calcule
 $$\delta_{\bar{P}}(a(y))^{1/2}=\vert y\vert _{F}^{3-d_{V}},$$
 $$\delta_{\bar{P}}(a(x',y))^{1/2}=\vert a(x',y)_{r-1}\vert _{F}.$$
 Dans l'int\'egrale $I^2_{r,\natural}(b,D)$, effectuons le changement de variable $x\mapsto yx$, autrement dit $x'\mapsto x$. Cela remplace $dx$ par $\vert y\vert _{F}^{d_{V}-4}dx$. En majorant la fonction ${\bf 1}_{\sigma\geq b}(u(x)\underline{y})$ par $1$, on obtient
 $$(5) \qquad I^2_{r,\natural}(b,D)<<\int_{y\in F; val_{F}(y)<-\mu b}\int_{V_{\flat}}\vert a(x) _{r}\vert ^{-1}\delta_{\bar{P}_{\sharp}}(m_{\bar{P}_{\sharp}}(u(x)))^{1/2}$$
 $$\Xi^{M_{\sharp}}(m_{\bar{P}_{\sharp}}(u(x)))\vert a(x,y)_{r-1}\vert _{F}\sigma(y)^D\sigma(u(x))^D\vert y\vert _{F}^{-1}dx\,dy.$$
 Posons $R_{\flat}=R\cap V_{\flat}$ et notons $R_{\flat}^{\vee}$ son dual, c'est-\`a-dire $R_{\flat}^{\vee}=\{v\in V_{\flat}; \forall v'\in R_{\flat}, q_{V}(v,v')\in \mathfrak{o}_{F}\}$. D\'ecomposons le membre de droite de (5) en deux int\'egrales, $I^3_{r,\natural}(b,D)+I^4_{r,\natural}(b,D)$, o\`u, dans $I^3_{r,\natural}(b,D)$, resp. $I^4_{r,\natural}(b,D)$, on restreint l'int\'egration aux $x\in \mathfrak{p}_{F}^{-c_{1}}R_{\flat}$, resp. $x\not\in \mathfrak{p}_{F}^{-c_{1}}R_{\flat}$. On doit d\'emontrer que chacune des deux int\'egrales v\'erifie la majoration de l'\'enonc\'e.
 
  Dans $I^3_{r,\natural}(b,D)$, $u(x)$ reste dans un compact et tous les termes d\'ependant de $x$ sont born\'es. D'o\`u
 $$I^3_{r,\natural}(b,D)<<\int_{y\in F; val_{F}(y)< -\mu b}X(y) \sigma(y)^D\vert y\vert _{F}^{-1}dy,$$
   o\`u
 $$X(y)=\int_{\mathfrak{p}_{F}^{-c_{1}}R_{\flat}}\vert a(x,y)_{r-1}\vert _{F}dx.$$
 Il existe $c_{3}>0$ tel que $\vert val_{F}(a(x)_{r})\vert \leq c_{3}$ pour tout $x\in \mathfrak{p}_{F}^{-c_{1}}R_{\flat}$.  On a alors 
 $$\vert a(x,y)_{r-1}\vert _{F}\leq \left\lbrace\begin{array}{cc}1,&\text{ si }val_{F}(q_{V}(x))+val_{F}(y)\geq -c_{3},\\ q^{c_{3}}\vert yq_{V}(x)\vert _{F}^{-1},&\text{ si }val_{F}(q_{V}(x))+val_{F}(y)<-c_{3}.\\ \end{array}\right.$$
 Pour tout entier $k\in {\mathbb Z}$, notons $m(k)$ la mesure de l'ensemble des $x\in \mathfrak{p}_{F}^{-c_{1}}R_{\flat}$ tels que $val_{F}(q_{V}(x))=k$. Certainement, cette valuation est born\'ee inf\'erieurement, donc il existe $c_{4}$ tel que $m(k)=0$ si $k<c_{4}$. On obtient
 $$X(y)<<\sum_{k; c_{4}\leq k< -c_{3}-val_{F}(y)}m(k)\vert y\vert _{F}^{-1}q^k+\sum_{-c_{3}-val_{F}(y)\leq k}m(k).$$
 On v\'erifie que l'on a une majoration $m(k)<<q^{-k/2}$. Alors
 $$X(y)<<\vert y\vert _{F}^{-1/2}\sigma(y).$$
 En reportant cette majoration dans $I^3_{r,\natural}(b,D)$, on obtient
 $$I^3_{r,\natural}(b,D)<<\int_{y\in F; val_{F}(y)< -\mu b} \sigma(y)^{D+1}\vert y\vert _{F}^{-3/2}dy.$$
 Mais il existe $\epsilon_{3}>0$ tel que cette int\'egrale soit born\'ee par $exp(-\epsilon_{3}b)$. C'est ce qu'on voulait.
 
 Traitons maintenant l'int\'egrale $I^4_{r,\natural}(b,D)$. Choisissons un entier $c_{4}$ tel que 
 
 (6)(i) $\mathfrak{p}_{F}^{c_{4}}R_{\flat}^{\vee}\subset \mathfrak{p}_{F}^{-c_{1}}R_{\flat}$;
 
 (6)(ii) $val_{F}(q_{V}(v))\geq sup(0,-c_{1}+c_{4})$ pour  tout $v\in \mathfrak{p}_{F}^{c_{4}}R_{\flat}^{\vee}$;
 
 (6)(iii) $u(v)\in K$ pour  tout $v\in \mathfrak{p}_{F}^{c_{4}}R_{\flat}^{\vee}$;
 
 (6)(iv) $c_{4}+c_{2}-c_{1}>0$.
 
 Pour $y\in F$ tel que $val_{F}(y)<-\mu b$, posons $n(y)=c_{4}+[\frac{-val_{F}(y)}{2}]$. D'apr\`es (6)(i), l'ensemble $V_{\flat}-\mathfrak{p}_{F}^{-c_{1}}R_{\flat}$ est stable par translation par $\mathfrak{p}_{F}^{n(y)}R_{\flat}^{\vee}$. Dans l'int\'egrale int\'erieure de $I^4_{r,\natural}(b,D)$, on peut remplacer $x$ par $x+v$, avec $v\in \mathfrak{p}_{F}^{n(y)}R_{\flat}^{\vee}$, puis int\'egrer sur $v$, tout en divisant par $mes(\mathfrak{p}_{F}^{n(y)}R_{\flat}^{\vee})$. On a $u(x+v)=u(x)u(v)$ et $u(v)\in K$ d'apr\`es (6)(iii). Donc $\sigma(u(x+v))=\sigma(u(x))$, $m_{\bar{P}_{\sharp}}(u(x+v))=m_{\bar{P}_{\sharp}}(u(x))$ et $a(x+v)=a(x)$. On obtient
 $$I^4_{r,\natural}(b,D)=\int_{y\in F; val_{F}(y)<-\mu b}\int_{V_{\flat}-\mathfrak{p}_{F}^{-c_{1}}R_{\flat}}\vert a(x) _{r}\vert ^{-1}\delta_{\bar{P}_{\sharp}}(m_{\bar{P}_{\sharp}}(u(x)))^{1/2}\Xi^{M_{\sharp}}(m_{\bar{P}_{\sharp}}(u(x))) \alpha(x,y)$$
 $$\sigma(y)^D\sigma(u(x))^D\vert y\vert _{F}^{-1}dx\,dy,$$
  o\`u
 $$\alpha(x,y)=mes(\mathfrak{p}_{F}^{n(y)}R_{\flat}^{\vee})^{-1}\int_{\mathfrak{p}_{F}^{n(y)}R_{\flat}^{\vee}} \vert a(x+v,y)_{r-1}\vert _{F}dv.$$
 Fixons $x,y$ intervenant dans l'int\'egrale ci-dessus.  Montrons que
 
 (7) on a $\alpha(x,y)<< \vert y\vert _{F}^{-1/2}\sigma(y)$.
 
 On a $val_{R}(x)<-c_{1}$ donc, d'apr\`es (1)
 $$val_{F}(a(x)_{r})=inf(val_{R}(x)+c_{1},val_{F}(q_{V}(x))+c_{2})<0.$$
 Pour $v\in \mathfrak{p}_{F}^{n(y)}R_{\flat}^{\vee}$, on a $a(x+v)=a(x)$ comme on vient de le dire et la valeur de $a(x+v,y)_{r-1}$ est fonction de $val_{F}(q_{V}(x+v))+val_{F}(y)-val_{F}(a(x)_{r})$. Notons $N(v)$ cet entier. On a $q_{V}(x+v)=q_{V}(x)+q_{V}(x,v)+q_{V}(v)$ et $val_{F}(q_{V}(x,v))\geq val_{R}(x)+n(y)$. Supposons d'abord $val_{F}(q_{V}(x))< val_{R}(x)+n(y)$. Gr\^ace \`a (6)(ii), on a 
 $$val_{R}(x)+n(y)<-c_{1}+c_{4}+[\frac{-val_{F}(y)}{2}]\leq -c_{1}+c_{4}+2[\frac{-val_{F}(y)}{2}]\leq val_{F}(q_{V}(v)).$$
 Donc $val_{F}(q_{V}(x+v))=val_{F}(q_{V}(x))$. Si $val_{F}(a(x)_{r})=val_{F}(q_{V}(x))+c_{2}$, alors
 $ N(v)=val_{F}(y)-c_{2}\leq -c_{2}$. D'apr\`es la d\'efinition de $a(x+v,y)_{r-1}$, on a 
 $$\vert a(x+v)_{r-1}\vert _{F}\leq q^{N(v)}<<\vert y\vert _{F}^{-1},$$
 et $\alpha(x,y)<<\vert y\vert _{F}^{-1}$. Si $val_{F}(a(x)_{r})=val_{R}(x)+c_{1}$, on a 
 $$N(v)<val_{R}(x)+n(y)+val_{F}(y)-val_{F}(a(x)_{r})=n(y)+val_{F}(y)-c_{1}\leq c_{4}-c_{2}+val_{F}(y)/2,$$
 d'o\`u $\vert a(x+v)_{r-1}\vert _{F}\leq q^{N(v)}<<\vert y\vert _{F}^{-1/2}$, puis $\alpha(x,y)<<\vert y\vert _{F}^{-1/2}$.
 Supposons maintenant $val_{F}(q_{V}(x))\geq val_{R}(x)+n(y)$. Gr\^ace \`a (6)(iv), $val_{F}(q_{V}(x))+c_{2}>val_{R}(x)+c_{1}$, donc $val_{F}(a(x)_{r})=val_{R}(x)+c_{1}$. Choisissons une base $(x_{j})_{j=1,...,m}$ de $R_{\flat}$ sur $\mathfrak{o}_{F}$ de sorte que $x=\varpi_{F}^{val_{R}(x)}x_{1}$. Introduisons la base duale $(x_{j}^{\vee})_{j=1,...,m}$ de $R_{\flat}^{\vee}$. Introduisons des coordonn\'ees sur $\mathfrak{p}_{F}^{n(y)}R^{\vee}$ en posant $v=\varpi_{F}^{n(y)}((z_{1}-\varpi_{F}^{-val_{R}(x)-n(y)}q_{V}(x))x_{1}^{\vee}+\sum_{j=2,...,m}z_{j}x_{j}^{\vee})$. Les $z_{j}$ parcourent $\mathfrak{o}_{F}$. On peut supposer $dv=\prod_{j=1,...,m}dz_{j}$ et alors $mes(\mathfrak{p}_{F}^{n(y)}R_{\flat}^{\vee})$ est ind\'ependant de $y$. On a $q_{V}(x)+q_{V}(x,v)=\varpi_{F}^{val_{R}(x)+n(y)}z_{1}$. On a
 $$val_{F}(q_{V}(x)+q_{V}(x,v))+val_{F}(y)-val_{F}(a(x)_{r})=val_{F}(z_{1})-n'(y), $$
 o\`u $n'(y)=c_{1}-val_{F}(y)-n(y)=c_{1}-c_{4}+[\frac{1-val_{F}(y)}{2}]$.
 Remarquons que, gr\^ace \`a (6)(ii), on a 
 $$val_{F}(q_{V}(v))+val_{F}(y)-val_{F}(a(x)_{r})\geq 2[\frac{-val_{F}(y)}{2}]+val_{F}(y)+1\geq 0.$$
 Si $val_{F}(z_{1})< n'(y)$, alors $N(v)=val_{F}(z_{1})-n'(y)$ et $\vert a(x+v)_{r-1}\vert_{F}=q^ {-n'(y)}\vert z_{1}\vert _{F}^{-1}$. Si $val_{F}(z_{1})\geq n'(y)$, alors $N(v)\geq0$ et $\vert a(x+v)_{r-1}\vert _{F}=1$. D'o\`u
 $$\alpha(x,y)<<\int_{z_{1}\in \mathfrak{o}_{F}; val_{F}(z_{1})<n'(y)}q^{-n'(y)}\vert z_{1}\vert _{F}^{-1}dz_{1}+\int_{z_{1}\in \mathfrak{o}_{F}; val_{F}(z_{1})\geq n'(y)}dz_{1}.$$ 
 On voit que $\alpha(x,y)<<(1+\vert n'(y)\vert )q^{-n'(y)}<<\sigma(y)\vert y\vert ^{-1/2}$. On a obtenu la majoration (7) dans tous les cas.
 
Alors
$$I^4_{r,\natural}(b,D)<<\int_{y\in F; val_{F}(y)<-\mu b}\sigma(y)^{D+1}\vert y\vert _{F}^{-3/2}dy$$
$$\int_{V_{\flat}}\vert a(x) _{r}\vert ^{-1}\delta_{\bar{P}_{\sharp}}(m_{\bar{P}_{\sharp}}(u(x)))^{1/2}\Xi^{M_{\sharp}}(m_{\bar{P}_{\sharp}}(u(x)))  \sigma(u(x))^D dx.$$
Il existe $\epsilon_{4}>0$ tel que la premi\`ere int\'egrale soit major\'ee par $exp(-\epsilon_{4}b)$. Lorsque l'on a trait\'e l'int\'egrale $I^1_{r,\natural}(b,D)$, on a montr\'e que la seconde int\'egrale \'etait convergente. D'o\`u la majoration cherch\'ee pour l'int\'egrale $I^4_{r,\natural}(b,D)$, ce qui ach\`eve la preuve. $\square$
 
 \bigskip
 
 \subsection{Majoration d'une int\'egrale unipotente, cas $r=1$}
 
 \ass{Lemme}{Supposons $r=1$. L'int\'egrale d\'efinissant $I_{1,\natural}(b,D)$ est convergente et, $D$ \'etant fix\'e, il existe $\epsilon>0$ tel que
 $$I_{1,\natural}(b,D)<<exp(-\epsilon b)$$
 pour tout $b\geq0$. }
 
 Preuve. Puisque $r=1$, on a $P_{1}=P$. Notons $V_{\sharp}$ le sous-espace de $V$ orthogonal \`a la droite $D_{0}$ port\'ee par $v_{0}$. Avec les m\^emes notations que dans le paragraphe pr\'ec\'edent, on a $U_{1,\natural}=U_{\sharp}$. Fixons un sous-groupe compact sp\'ecial $K_{\sharp}$ de $G_{\sharp}(F)$, en bonne position relativement \`a $M_{\sharp}$. Il n'est pas forc\'ement inclus dans $K$. Pour $u\in U_{1,\natural}(F)$, on a  $u=m_{\bar{P}_{\sharp}}(u)u_{\bar{P}_{\sharp}}(u)k_{\bar{P}_{\sharp}}(u) $ et on en d\'eduit que $m_{\bar{P}}(u)m_{\bar{P}_{\sharp}}(u)^{-1}$ reste dans un compact.
 Ecrivons $m_{\bar{P}_{\sharp}}(u)=a(u)g_{0}(u)$, avec $a(u)\in A(F)$ et $g_{0}(u)\in G_{0}(F)$. Comme dans le paragraphe pr\'ec\'edent, on voit qu'il existe $\epsilon_{1}>0$ tel que $\vert a(u)_{1}\vert _{F}^{-1}<<exp(-\epsilon_{1} \sigma(u))$ et que l'on a l'\'egalit\'e
 $$\delta_{\bar{P}}(m_{\bar{P}}(u))^{1/2}\Xi^M(m_{\bar{P}}(u))<<\vert a(u)_{1}\vert_{F} ^{-1/2}\delta_{\bar{P}_{\sharp}}(m_{\bar{P}_{\sharp}}(u))^{1/2}\Xi^{M_{\sharp}}(m_{\bar{P}_{\sharp}}(u)).$$
 Alors
 $$I_{1,\natural}(b,D)<<exp(-\epsilon_{1}b/4)\int_{U_{\sharp}(F)}\delta_{\bar{P}_{\sharp}}(m_{\bar{P}_{\sharp}}(u))^{1/2}\Xi^{M_{\sharp}}(m_{\bar{P}_{\sharp}}(u))exp(-\epsilon_{1}\sigma(u)/4)\sigma(u)^D du.$$
 L'int\'egrale est convergente d'apr\`es le lemme II.4.3 de [W2]. D'o\`u le r\'esultat. $\square$

\bigskip

\subsection{Preuve de 4.3(3)}

Si $r=0$, $U=\{1\}$ et l'assertion 4.3(3) est \'evidente. Supposons $r\geq1$. On a introduit le sous-groupe parabolique $P_{r}=M_{r}U_{r}$ en 4.5. Pour $x\in M_{r}(F)$, posons
$$X_{r}(c,D,x)=\int_{U_{r}(F)\cap U(F)_{c}}\Xi^G(ux)\sigma(ux)^Ddu.$$
On va montrer que

(1) cette int\'egrale est convergente; le r\'eel $D$ \'etant fix\'e, il existe un r\'eel $D'$ tel que
$$X_{r}(c,D,x)<<c^{D'}\sigma(x)^{D'}\delta_{P_{r}}(x)^{1/2}\Xi^{M_{r}}(x)$$
pour tous $c,x$.

Auparavant, montrons que cette assertion entra\^{\i}ne  4.3(3). On a
$$X(c,D,m)=\int_{M_{r}(F)\cap U(F)_{c}}\int_{U_{r}(F)\cap U(F)_{c}}\Xi^G(uvm)\sigma(uvm)^Ddu\,dv$$
$$=\int_{M_{r}(F)\cap U(F)_{c}}X_{r}(c,D,vm)dv.$$
Gr\^ace \`a (1), on a
$$X(c,D,m)<<c^{D'}\int_{U_{r}(F)\cap U(F)_{c}}\delta_{P_{r}}(vm)^{1/2}\Xi^{M_{r}}(vm)\sigma(vm)^{D'}dv.$$
On a $\delta_{P_{r}}(vm)=\delta_{P_{r}}(m)$ pour tout $v$. On a $M_{r}=GL_{1}\times \tilde{G}$, o\`u $\tilde{G}$ est le groupe sp\'ecial orthogonal de l'orthogonal dans $V$ du plan engendr\'e par $v_{r}$ et $v_{-r}$. On a $M_{r}\cap U=\tilde{G}\cap U$. Ecrivons $m=(a,\tilde{m})$, avec $a\in GL_{1}(F)$, $\tilde{m}\in \tilde{G}(F)$. Alors $\Xi^{M_{r}}(vm)=\Xi^{\tilde{G}}(v\tilde{m})$. D'o\`u
$$X(c,D,m)<<c^{D'}\sigma(a)^{D'}\delta_{P_{r}}(m)^{1/2}\int_{\tilde{G}(F)\cap U(F)_{c}}\Xi^{\tilde{G}}(v\tilde{m})\sigma(v\tilde{m})^{D'}dv.$$
Mais cette int\'egrale est analogue \`a celle de d\'epart: on a remplac\'e $G$ par $\tilde{G}$, $D$ par $D'$ et $m$ par $\tilde{m}$. En passant de $G$ \`a $\tilde{G}$, on remplace $r$ par $r-1$. En raisonnant par r\'ecurrence sur $r$, on peut supposer qu'il existe $R'$ tel que l'int\'egrale soit essentiellement major\'ee par $c^{R'}\sigma(\tilde{m})^{R'}\delta_{P\cap \tilde{G}}(\tilde{m})^{1/2}\Xi^{M}(\tilde{m})$. Alors
$$X(c,D,m)<<c^{D'+R'}\sigma(a)^{D'}\sigma(\tilde{m})^{R'}\delta_{P_{r}}(m)^{1/2}\delta_{P\cap \tilde{G}}(\tilde{m})^{1/2}\Xi^M(\tilde{m}).$$
Mais $\delta_{P_{r}}(m)\delta_{P\cap \tilde{G}}(\tilde{m})=\delta_{P}(m)$ et $\Xi^M(\tilde{m})=\Xi^M(m)$, d'o\`u la majoration cherch\'ee.

La preuve de (1) est analogue \`a celle du lemme 3.3. Introduisons un r\'eel $b>0$ que nous pr\'eciserons par la suite. On d\'ecompose $X_{r}(c,D,x)$ en $X_{r,<b}(c,D,x)+X_{r,\geq b}(c,D,x)$, o\`u
$$X_{r,<b}(c,D,x)=\int_{U_{r}(F)\cap U(F)_{c}}{\bf 1}_{\sigma<b}(ux)\Xi^G(ux)\sigma(ux)^Ddu,$$
$$X_{r,\geq b}(c,D,x)=\int_{U_{r}(F)\cap U(F)_{c}}{\bf 1}_{\sigma\geq b}(ux)\Xi^G(ux)\sigma(ux)^Ddu.$$
 Dans la premi\`ere, on a $\sigma(ux)<b$, donc
 $$X_{r,<b}(c,D,x)<<b^{D'}\int_{U_{r}(F)}\Xi^G(ux)\sigma(ux)^{D-D'}du$$
 pour tout r\'eel $D'>0$. D'apr\`es [W2] proposition II.4.5, si $D'$ est assez grand, cette int\'egrale est convergente et essentiellement born\'ee par $\delta_{P_{r}}(x)^{1/2}\Xi^{M_{r}}(x)$. Pour un tel $D'$, on a donc
 $$(2) \qquad X_{r,<b}(c,D,x)<<b^{D'}\delta_{P_{r}}(x)^{1/2}\Xi^{M_{r}}(x).$$
 
 Introduisons le "sous-groupe radiciel" \'evident $U_{r,r-1}$ de $U_{r}$. On a la d\'ecomposition $U_{r}=U_{r,\natural}U_{r,r-1}$, avec la notation de 4.5 et $U_{r}(F)\cap U(F)_{c}=U_{r,\natural}(F)U_{r,r-1}(F)_{c}$, o\`u $U_{r,r-1}(F)_{c}=U_{r,r-1}(F)\cap U(F)_{c}$. D'o\`u
 $$X_{r,\geq b}(c,D,x)=\int_{U_{r,r-1}(F)_{c}}\int_{U_{r,\natural}(F)}{\bf 1}_{\sigma\geq b}(uv)\Xi^G(uvx)\sigma(uvx)^Ddu\,dv.$$
 Il existe $c_{1}>0$ tel que $\sigma(v)<c_{1 }c$ pour tout $v\in U_{r,r-1}(F)_{c}$.  Si $b>2c_{1}c$, la condition $\sigma(uv)\geq b$ entra\^{\i}ne $\sigma(u)\geq b/2$. On suppose $b>2c_{1}c$. On a encore la majoration 3.3(5), d'o\`u, comme en 3.3,
 $$X_{r,\geq b}(c,D,x)<<exp(\alpha c_{1}c+\alpha\sigma(x))\int_{U_{r,\natural}(F)}{\bf 1}_{\sigma\geq b/2}\Xi^G(u)\sigma(u)^Ddu,$$
 pour un $\alpha>0$ convenable.  D'apr\`es [W2], lemmes II.1.1 et II.3.2, il existe un r\'eel $D''$ tel que
 $$\Xi^G(g)<<\delta_{\bar{P}}(m_{\bar{P}}(g))^{1/2}\Xi^M(m_{\bar{P}}(g))\sigma(g)^{D''}$$
 pour tout $g\in G(F)$. Alors l'int\'egrale ci-dessus est born\'ee par l'int\'egrale $I_{r,\natural}(b/2,D+D'')$ de 4.5. En utilisant les lemmes 4.5 ou 4.6 selon la valeur de $r$, elle est essentiellement born\'ee par $exp(-\epsilon b) $ pour un $\epsilon>0$ convenable. Enfin, on v\'erifie qu'il existe $c_{2}>0$ tel que l'on ait la majoration
 $$exp(-c_{2}\sigma(x))<<\delta_{P_{r}}(x)^{1/2}\Xi^{M_{r}}(x).$$
 Alors 
 $$X_{r,\geq b}(c,D,x)<<exp(\alpha c_{1}c+(\alpha+c_{2})\sigma(x)-\epsilon b)\delta_{P_{r}}(x)^{1/2}\Xi^{M_{r}}(x).$$
 On choisit maintenant $b=2c_{1}c+\epsilon^{-1}\alpha c_{1}c+\epsilon^{-1}(\alpha+c_{2})\sigma(x)$. Alors
 $$X_{r,\geq b}(c,D,x)<<\delta_{P_{r}}(x)^{1/2}\Xi^{M_{r}}(x).$$
 Cette majoration et (2) entra\^{\i}nent (1). $\square$
 
 \bigskip
 
 \subsection{Majoration d'int\'egrales doubles sur $U(F)$}
 
 Soient $D$ un r\'eel, $c$ et $c'$ deux entiers. Supposons $c'\geq c>0$. Remarquons que l'ensemble $U(F)-U(F)_{c'}$ est invariant par translation par $U(F)_{c}$. Pour $m,m'\in M(F)$, posons
 $$X(c,D,m,m')=\int_{U(F)/U(F)_{c}}\int_{U(F)_{c}}\int_{U(F)_{c}}\Xi^G(uvm)\Xi^G(uv'm')\sigma(uvm)^D\sigma(uv'm')^Ddv'\,dv\,du,$$
 $$X(c,c',D,m,m')=\int_{(U(F)-U(F)_{c'})/U(F)_{c}}\int_{U(F)_{c}}\int_{U(F)_{c}}\Xi^G(uvm)$$
 $$\Xi^G(uv'm')\sigma(uvm)^D\sigma(uv'm')^Ddv'\,dv\,du,$$
 $$X(D,m,m')=\int_{U(F)}\Xi^G(um)\Xi^G(um')\sigma(um)^D\sigma(um')^Ddu.$$
 
 \ass{Lemme}{(i) Ces int\'egrales sont convergentes.
 
 (ii) Le r\'eel $D$ \'etant fix\'e, il existe un r\'eel $R$ tel que
 $$X(c,D,m,m')<<c^R\sigma(m)^R \delta_{P}(m)^{1/2}\Xi^M(m)\sigma(m')^R\delta_{P}(m')^{1/2}\Xi^M(m')$$
 pour tous $m,m'\in M(F)$ et tout $c\geq1$.
 
 (iii) Les termes $c$ et $D$ \'etant fix\'es, il existe un r\'eel $R$ et un r\'eel $\epsilon>0$ tels que
 $$X(c,c',D,m,m')<<exp(-\epsilon c')\sigma(m)^R \delta_{P}(m)^{1/2}\Xi^M(m)\sigma(m')^R\delta_{P}(m')^{1/2}\Xi^M(m')$$
pour tous $m,m'\in M(F)$ et tout $c'\geq c$.

(iv) Le r\'eel $D$ \'etant fix\'e, il existe un r\'eel $R$ tel que
 $$X(D,m,m')<<\sigma(m)^R \delta_{P}(m)^{1/2}\Xi^M(m)\sigma(m')^R\delta_{P}(m')^{1/2}\Xi^M(m')$$
 pour tous $m,m'\in M(F)$.}
 
 Preuve.  L'application
 $$\begin{array}{ccc}U(F)/U(F)_{c}&\to& (F/\mathfrak{p}_{F}^{-c})^r\\ u&\mapsto& (u_{r,r-1}+\mathfrak{p}_{F}^{-c},...,u_{1,0}+\mathfrak{p}_{F}^{-c})\\ \end{array}$$
 est un isomorphisme. D\'efinissons une fonction $val_{F}^c$ sur $F $ par $val_{F}^c(x)=0$ si $x\in \mathfrak{p}_{F}^{-c}$, $val_{F}^c(x)=val_{F}(x)+c$ si $x\not\in \mathfrak{p}_{F}^{-c}$. Elle se quotiente en une fonction sur $F/\mathfrak{p}_{F}^{-c}$. Pour $u\in U(F)/U(F)_{c}$, posons $val_{F}^c(u)=\sum_{i=1,...,r}val_{F}^c(u_{i,i-1})$. Montrons que:
 
 (1) il existe un r\'eel $D_{1}$ tel que l'on ait une majoration
 $$\int_{U(F)_{c}}\Xi^G(uvm)\sigma(uvm)^Ddv<<(c- val_{F}^c(u)) ^{D_{1}}q^{val_{F}^c(u)}\sigma(m)^{D_{1}}\delta_{P}(m)^{1/2}\Xi^M(m)$$
 pour tout $m\in M(F)$, tout $c\geq1$ et tout $u\in U(F)/U(F)_{c}$.
 
 Soit $a\in A(F)\cap K$. On peut remplacer $\Xi^G(uvm)\sigma(uvm)^D$ par $\Xi^G(auvma^{-1})\sigma(auvma^{-1})^D$. On peut ensuite int\'egrer en $a$. Puisque $a$ commute \`a $m$ et normalise $U(F)_{c}$, on obtient
$$ \int_{U(F)_{c}}\Xi^G(uvm)\sigma(uvm)^Ddv<<\int_{A(F)\cap K}\int_{U(F)_{c}}\Xi^G(aua^{-1}vm)\sigma(aua^{-1}vm)^Ddv\,da.$$
Consid\'erons l'application
$$\begin{array}{ccc}A(F)\cap K&\to& U(F)/U(F)_{c}\\ a&\mapsto&aua^{-1}U(F)_{c}.\\ \end{array}$$
On v\'erifie que son jacobien est born\'e par $q^{val_{F}^c(u)}$. Son image est contenue dans $U(F)_{-val_{F}^c(u)+c}/U(F)_{c}$.
D'o\`u
 $$ \int_{U(F)_{c}}\Xi^G(uvm)\sigma(uvm)^Ddv<<q^{val_{F}^c(u)}\int_{U(F)_{-val_{F}^c(u)+c}/U(F)_{c}}\int_{U(F)_{c}}\Xi^G(uvm)\sigma(uvm)^Ddv\,du$$
 $$<< q^{val_{F}^c(u)}\int_{U(F)_{-val_{F}^c(u)+c}}\Xi^G(vm)\sigma(vm)^Ddv.$$
 Il reste \`a faire appel \`a 4.3(3) pour obtenir l'assertion (1).
 
Gr\^ace \`a (1), on a
$$X(c,D,m,m')<<\sigma(m)^{D_{1}}\delta_{P}(m)^{1/2}\Xi^M(m)\sigma(m')^{D_{1}}\delta_{P}(m')\Xi^M(m')$$
$$\int_{U(F)/U(F)_{c}}(c- val_{F}^c(u)) ^{2D_{1}}q^{2val_{F}^c(u)}du.$$ 
Cette derni\`ere int\'egrale est produit d'int\'egrales du type
$$(2)\qquad\int_{F/\mathfrak{p}_{F}^{-c}}(c-val_{F}^c(x))^{2D_{1}}q^{2val_{F}^c(x)}dx.$$
On v\'erifie qu'il existe un r\'eel $D_{2}$ tel que cette expression soit essentiellement born\'ee par $c^{D_{2}}$.  On en d\'eduit une majoration similaire pour l'int\'egrale intervenant dans la formule ci-dessus. Alors $X(c,D,m,m')$ v\'erifie la majoration  du (ii) de l'\'enonc\'e.

Gr\^ace \`a (1), on a
$$X(c,D,m,m')<<\sigma(m)^{D_{1}}\delta_{P}(m)^{1/2}\Xi^M(m)\sigma(m')^{D_{1}}\delta_{P}(m')\Xi^M(m')$$
$$\int_{(U(F)-U(F)_{c'})/U(F)_{c}}(c- val_{F}^c(u)) ^{2D_{1}}q^{2val_{F}^c(u)}du.$$ 
Cette derni\`ere int\'egrale est combinaison lin\'eaire de termes qui sont des produits d'int\'egrales du type (2) et d'au moins une int\'egrale du type
$$(3) \qquad \int_{(F-\mathfrak{p}_{F}^{-c'})/\mathfrak{p}_{F}^{-c}}(c-val_{F}^c(x))^{2D_{1}}q^{2val_{F}^c(x)}dx.$$
On v\'erifie qu'il existe $\epsilon>0$ tel que  cette expression soit essentiellement born\'ee par $exp(-\epsilon c')$. On en d\'eduit le (iii) de l'\'enonc\'e.

  Soit $v\in U(F)\cap K$. Dans l'int\'egrale $X(D,m,m')$, on peut remplacer $\Xi^G(um)$ par $\Xi^G(vum)$, puis int\'egrer sur $v$. Donc
$$X(D,m,m')<<\int_{U(F)\cap K}\int_{U(F)}\Xi^G(vum)\Xi^G(um')\sigma(vum)^D\sigma(um')^Ddu\, dv.$$ 
 Choisissons $c$ tel que $U(F)\cap K\subset U(F)_{c}$. On peut remplacer l'int\'egrale sur $U(F)\cap K$ par l'int\'egrale sur $U(F)_{c}$. Ce groupe \'etant distingu\'e dans $U(F)$, on peut remplacer $vu$ par $uv$. On peut ensuite d\'ecomposer l'int\'egrale sur $U(F)$ en la compos\'ee d'une int\'egrale sur $U(F)_{c}$ et d'une int\'egrale sur $U(F)/U(F)_{c}$. On obtient alors
 $$X(D,m,m')<<X(c,D,m,m'),$$
 et  l'assertion (iv) r\'esulte de (ii). $\square$

 \bigskip
 
 \subsection{Comparaison de $\Xi^G$ et $\Xi^H$}
 
 \ass{Lemme}{Supposons $r=0$. Il existe $\epsilon>0$ tel que $\Xi^G(h)<<exp(-\epsilon \sigma(h))\Xi^H(h)$ pour tout $h\in H(F)$.}
 
 Preuve.  On a n\'ecessairement $d_{an,V}=d_{an,W}\pm 1$. Si $d_{an,V}=d_{an,W}+1$, fixons un syst\`eme hyperbolique maximal $(e_{\pm j})_{j=1,...,n}$ de $W$. C'est aussi un syst\`eme hyperbolique maximal de $V$. Si $d_{an,V}=d_{an,W}-1$, fixons un syst\`eme hyperbolique maximal $(e_{\pm j})_{j=2,...,n}$ de $W$. Notons $W_{an}$ l'orthogonal dans $W$ du sous-espace engendr\'e par ces vecteurs. Fixons un syst\`eme hyperbolique maximal $(e_{\pm 1})$ de $W_{an}\oplus D_{0}$. Alors $(e_{\pm j})_{j=1,...,n}$ est un syst\`eme hyperbolique maximal de $V$. On pose $\iota=1$ dans le premier cas, $\iota=2$ dans le second. Notons $A_{min}$ le sous-tore d\'eploy\'e maximal de $G$ form\'e des \'el\'ements qui conservent chaque droite $Fe_{\pm j}$ pour $j=1,...,n$ et qui agissent trivialement sur l'orthogonal du sous-espace engendr\'e par ces vecteurs. Pour $a\in A_{min}(F)$, et $j=1,...,n$, on note $a_{j}$ la valeur propre de $a$ sur le vecteur $e_{j}$. Posons $A_{min}^H=A_{min}\cap H$. C'est un sous-tore d\'eploy\'e maximal de $H$.  On sait qu'il existe un sous-ensemble compact $\Gamma$ de $H(F)$ tel que $H(F)=\Gamma A_{min}^H(F)\Gamma$.  Il suffit donc de d\'emontrer le lemme pour les \'el\'ements de $A_{min}^H(F)$.  On va d\'emontrer 

(1) il existe un r\'eel $D$ tel que 
$$\Xi^G(h)<<\Xi^H(h)\sigma(h)^D q^{-\sum_{j=\iota,...,n}\vert val_{F}(h_{j})\vert/2}$$
pour tout $h\in A_{min}^H(F)$. 

Cette relation implique la majoration de l'\'enonc\'e. Notons $\mathfrak{S}$ le groupe des permutations $s$ de $\{\pm 1,...,\pm n\}$ telles que $s(-j)=-s(j)$ pour tout $j\in \{\pm 1,...,\pm n\}$. Si $d_{an,V}=d_{an,W}+1$, on pose $\mathfrak{S}^H=\mathfrak{S}$; si $d_{an,V}=d_{an,W}-1$, on note $\mathfrak{S}^H$ le sous-groupe des $s\in \mathfrak{S}$ qui fixent $1$ et $-1$. Pour tout $s\in \mathfrak{S}$, notons $A_{min}(F)_{s}^-$ l'ensemble des $a\in A_{min}(F)$ tels que
$$val_{F}(a_{sn})\geq...\geq val_{F}(a_{s1})\geq0.$$
Le groupe $A_{min}^H(F)$ est contenu dans la r\'eunion des $A_{min}(F)_{s}^-$ quand $s$ d\'ecrit $\mathfrak{S}^H$. On peut fixer $s$ et se limiter \`a prouver (1) pour $h\in A_{min}^H(F)\cap A_{min}(F)_{s}^-$. Fixons donc $s$. Quitte \`a r\'eindexer notre syst\`eme hyperbolique, on peut supposer $s=1$. On abandonne les indices $s$, en conservant les exposants $-$. Notons $P_{min}$ le sous-groupe parabolique de $G$ form\'e des \'el\'ements qui conservent le drapeau
 $$Fe_{n}\subset Fe_{n}\oplus Fe_{n-1}\subset...\subset Fe_{n}\oplus...\oplus Fe_{1}.$$
 Posons $P_{min}^H=P_{min}\cap H$. D'apr\`es [W2] lemme II.1.1, pour $h\in A_{min}^H(F)\cap A_{min}(F)^-$, on a
 $$\Xi^H(h)>> \delta_{P_{min}^H}(h)^{1/2}.$$
 Ce dernier terme est \'egal \`a
 $$\prod_{j=\iota,...,n}\vert h_{j}\vert_{F} ^{j-\iota+d_{an,W}/2}$$
 D'apr\`es la m\^eme r\'ef\'erence, il existe un r\'eel $D$ tel que
 $$\Xi^G(h)<<\sigma(h)^D\delta_{P_{min}^G}(h)^{1/2}.$$
 Ce dernier terme s'\'ecrit
 $$\prod_{j=1,...,n}\vert h_{j}\vert_{F} ^{j-1+d_{an,V}/2}.$$
 Remarquons que $d_{an,V}/2-1=1/2+ d_{an,W}/2-\iota$ et $h_{1}=1$ si $\iota=2$.  Ces relations entra\^{\i}nent (1). $\square$
   
   \bigskip
   
   \subsection{Preuve des  relations 4.3(4) et 4.3(5)}
   
   On veut prouver la convergence de 
   $$(1) \qquad \int_{H(F)U(F)_{c}}\Xi^H(h)\Xi^G(hu)\sigma(hu)^Ddu\,dh.$$
   On utilise 4.3(3) pour majorer l'int\'egrale sur $U(F)_{c}$. En remarquant que $\delta_{P}(h)=1$ et $\Xi^M(h)=\Xi^{G_{0}}(h)$, on obtient qu'il existe un r\'eel $D'$ tel que l'int\'egrale ci-dessus soit essentiellement born\'ee par
   $$\int_{H(F)}\Xi^H(h)\Xi^{G_{0}}(h)\sigma(h)^{D'}dh.$$
   En appliquant le lemme 4.9 \`a $G_{0}$, il existe $\epsilon>0$ tel que cette expression soit essentiellement born\'ee par
   $$\int_{H(F)}\Xi^H(h)^2 exp(-\epsilon\sigma(h))dh.$$
   Or cette int\'egrale est convergente d'apr\`es [W2] lemme II.1.5.  D'o\`u la relation 4.3(4).
   
   On veut prouver la convergence de
     $$(2) \qquad \int_{H(F)U(F)_{c}}\int_{H(F)U(F)_{c}}\Xi^G(hu)\Xi^H(h'h)\Xi^G(h'u')\sigma(hu)^D\sigma(h'u')^Ddu'\,dh'\,du\,dh.$$
     Soit $K^H$ le sous-groupe compact sp\'ecial de $H(F)$ sous-jacent \`a la d\'efinition de $\Xi^H$. On peut remplacer $h$ par $kh$, avec $k\in K^H$, puis int\'egrer sur $k$, tout en divisant par $mes(K^H)$. Or $\Xi^G(khu)<<\Xi^G(hu)$ et $\sigma(khu)<<\sigma(hu)$. Le proc\'ed\'e ci-dessus revient donc \`a remplacer le terme $\Xi^H(h'h)$ par
     $$mes(K^H)^{-1}\int_{K^H }\Xi^H(h'kh)dk.$$
     D'apr\`es [W2] lemme II.1.3, ceci n'est autre que $\Xi^H(h')\Xi^H(h)$. Alors l'expression (2) appara\^{\i}t comme le carr\'e de l'expression (1). Elle est convergente puisque (1) l'est. Cela prouve la relation 4.3(5).

\bigskip

\subsection{Majoration d'une int\'egrale de fonctions d'Harish-Chandra, cas $r=0$}

On suppose dans ce paragraphe $r=0$. Soit $D$ un r\'eel. Pour $h\in H(F)$ et $N\geq1$ un entier, posons
$$\chi(h,N,D)=\int_{G(F)}\Xi^G(hx)\Xi^G(x)\kappa_{N}(x)\sigma(x)^Ddx.$$

\ass{Lemme}{Cette int\'egrale est convergente. Le r\'eel $D$ \'etant fix\'e, il existe    un r\'eel $R$ tel que
$$\chi(h,N,D)<< \Xi^G(h) N^R\sigma(h)^R$$
pour tout $h\in H(F)$ et tout entier $N\geq1$.}

Preuve. Si $V$ est anisotrope, le groupe $G(F)$ est compact et l'assertion est \'evidente. On suppose que $V$ n'est pas anisotrope. Comme dans la preuve de 4.9, dont on reprend les notations, on introduit un syst\`eme hyperbolique maximal $(e_{\pm j})_{j=1,...,n}$ de $V$. On note $V_{an}$ l'orthogonal dans $V$ du sous-espace engendr\'e par ces vecteurs. Si $d_{an,V}=d_{an,W}+1$, $v_{0}$ appartient \`a $V_{an}$. Si $d_{an,V}=d_{an,W}-1$, on peut choisir $e_{1}$ et $e_{-1}$ de sorte que $v_{0}=e_{1}+\nu_{0}e_{-1}$, o\`u $\nu_{0}=q_{V}(v_{0})$.  Il suffit de d\'emontrer le lemme pour $h\in A_{min}^H(F)$. En effet, il existe un sous-ensemble compact $\Gamma^H$ de $H(F)$ tel que $H(F)=\Gamma^H A_{min}^H(F)\Gamma^H$. Ecrivons $h=\gamma'a\gamma$, avec $\gamma,\gamma'\in \Gamma^H$ et $a\in A_{min}^H(F)$. On effectue le changement de variable $x\mapsto \gamma^{-1} x$. Puisque la fonction $\kappa_{N}$ est invariante \`a gauche par $H(F)$, on obtient
$$\chi(h,N,D)=\int_{G(F)}\Xi^G(\gamma'ax)\Xi^G(\gamma^{-1}x)\kappa_{N}(x)\sigma(\gamma^{-1}x)^Ddx.$$
Mais cette expression est essentiellement major\'ee par $\chi(a,N,D)$, d'o\`u l'assertion.

On suppose donc $h\in A_{min}^H(F)$. Fixons un  sous-groupe d'Iwahori $I$ de $G(F)$ en bonne position relativement \`a $A_{min}$.  Il est loisible de  supposer que $\Xi^G$ est biinvariante par $I$ (par contre, on ne suppose pas que $I$ soit inclus dans le sous-groupe compact sp\'ecial $K$ que l'on a fix\'e, c'est-\`a-dire dans le fixateur du r\'eseau $R$).
D'apr\`es Bruhat-Tits, il existe un sous-ensemble ouvert compact $\Gamma$ de $G(F)$ tel que $G(F)=IA_{min}(F)\Gamma$. On fixe $\Gamma$ et on suppose $I\subset \Gamma$. Le groupe $I\cap A_{min}(F)$ est le sous-groupe compact maximal de $A_{min}(F)$. L'application
$$\begin{array}{ccc}A_{min}(F) &\to&{\mathbb Z}^n\\ a&\mapsto& (val_{F}(a_{1}),...,val_{F}(a_{n}))\\ \end{array}$$
se quotiente en un isomorphisme de $A_{min}(F)/(I\cap A_{min}(F))$ sur ${\mathbb Z}^n$. Fixons un sous-ensemble $\Lambda\subset A_{min}(F)$ qui s'envoie bijectivement sur ${\mathbb Z}^n$. Alors
$$\chi(h,N,D)\leq \sum_{a\in \Lambda}\chi(h,N,D,a),$$
o\`u
$$\chi(h,N,D,a)=\int_{Ia\Gamma}\Xi^G(hx)\Xi^G(x)\kappa_{N}(x)\sigma(x)^Ddx.$$
Comme en 4.9, introduisons le groupe de permutations $\mathfrak{S}$ et,  pour tout $s\in \mathfrak{S}$, le sous-ensemble $A_{min}(F)_{s}^-$. Posons $\Lambda_{s}^-=\Lambda\cap A_{min}(F)_{s}^-$.  Alors $\Lambda$ est r\'eunion des $\Lambda_{s}^-$. Il nous suffit de fixer $s$ et de majorer $\chi(h,N,D)_{s}^-$, o\`u
$$\chi(h,N,D)_{s}^-=\sum_{a\in \Lambda_{s}^-}\chi(h,N,D,a).$$
Fixons donc $s$. Quitte \`a r\'eindexer notre syst\`eme hyperbolique, on peut supposer $s=1$. On abandonne les indices $s$ dans les notations pr\'ec\'edentes, en conservant seulement les exposants $-$. Il faut prendre garde au fait que, dans le cas o\`u $d_{an,V}=d_{an,W}-1$, la r\'eindexation n'a pas de raison de conserver les vecteurs $e_{1}$ et $e_{-1}$.  Ceux-ci se transforment en des vecteurs que nous notons $e_{\pm t}$, avec $t\in \{1,...,n\}$. On a $e_{0}=e_{t}+\nu_{0}e_{-t}$ ou $e_{0}=\nu_{0}e_{t}+e_{-t}$. On introduit le sous-groupe parabolique  $P_{min}=M_{min}U_{min}$   de $G$ comme en 4.9. Soit $a\in \Lambda^-$. On a l'\'egalit\'e $I=(I\cap U_{min}(F))(I\cap \bar{P}_{min}(F))$ et $a^{-1}(I\cap \bar{P}_{min}(F))a\subset I\subset \Gamma$. Donc $Ia\Gamma=(I\cap U_{min}(F))a\Gamma$. On v\'erifie que la mesure de cet ensemble est essentiellement born\'ee par $\delta_{P_{min}}(a)^{-1}$. Donc
$$\chi(h,N,D,a)<<\delta_{P_{min}}(a)^{-1}\int_{I\cap U_{min}(F)}\int_{\Gamma}\Xi^G(hya\gamma)\Xi^G(ya\gamma)\kappa_{N}(ya\gamma)\sigma(ya\gamma)^Dd\gamma\,dy.$$
Il existe un entier $c_{1}\geq0$ tel que $\Gamma R\subset \mathfrak{p}_{F}^{-c_{1}}R$. De la d\'efinition de la fonction $\kappa_{N}$ r\'esulte que $\kappa_{N}(ya\gamma)\leq \kappa_{N+c_{1}}(ya)$. Alors
$$\chi(h,N,D,a)<<\delta_{P_{min}}(a)^{-1}\Xi^G(a)\sigma(a)^D\int_{I\cap U_{min}(F)}\Xi^G(hya)\kappa_{N+c_{1}}(ya)dy.$$
Gr\^ace au lemme II.1.1 de [W2], il existe un r\'eel $D_{1}$ tel que
$$\delta_{P_{min}}(a)^{-1}\Xi^G(a)\sigma(a)^D<<\delta_{P_{min}}(a)^{-1/2}\sigma(a)^{D_{1}}.$$
D'autre part, pour le r\'esultat que l'on veut obtenir, le changement de $N$ en $N+c_{1}$ est insignifiant.  Il nous suffit de majorer
$$(1) \qquad \sum_{a\in \Lambda^-}\delta_{P_{min}}(a)^{-1/2}\sigma(a)^{D_{1}}Y(h,N,a),$$
o\`u
$$Y(h,N,a)=\int_{I\cap U_{min}(F)}\Xi^G(hya)\kappa_{N}(ya)dy.$$
D\'efinissons des entiers $N(h)$ et $b(h,N,a)$ par $N(h)=sup(N,N-val_{F}(h_{n}))$, $b(h,N,a)=sup(0,val_{F}(a_{n})-N(h))$. Montrons que

(2) il existe $\epsilon'>0$ tel que
$$Y(h,N;a)<<exp(-\epsilon' b(h,N,a))\int_{I\cap U_{min}(F)}\Xi^G(hya)dy$$
pour tous $h\in A_{min}^H(F)$, $N\geq1$, $a\in \Lambda^-$.

Si $val_{F}(a_{n})\leq N$, il suffit de majorer $\kappa_{N}(ya)$ par $1$. Supposons $val_{F}(a_{n})>N$.
Supposons d'abord que $v_{0}$ est orthogonal \`a $e_{n}$ et $e_{-n}$, c'est-\`a-dire que  $d_{an,V}=d_{an,W}+1$ ou $d_{an,V}=d_{an,W}-1$ et $t\not=n$. Introduisons le sous-groupe de $G(F)$ form\'e des \'el\'ements $u(z)$, pour $z\in F$, d\'efinis ainsi: $u(z)$ envoie $v_{0}$ sur $v_{0}+ze_{n}$, $e_{-n}$ sur $e_{-n}-\frac{z}{2\nu_{0}}v_{0}-\frac{z^2}{2\nu_{0}}e_{n}$, et fixe $e_{n}$ ainsi que l'orthogonal du sous-espace engendr\'e par $e_{-n}$, $v_{0}$ et $e_{n}$.  Il existe un entier $c_{2}$ tel que $u(z)\in I\cap U_{min}(F)$ pour $val_{F}(z)\geq c_{2}$. On a $h^{-1}u(z)h=u(h_{n}^{-1}z)$. Pour $ z\in \mathfrak{p}_{F}^{ sup(c_{2},val_{F}(h_{n})+c_{2})}$, les deux \'el\'ements $u(z)$ et $h^{-1}u(z)h$ appartiennent \`a $I\cap U_{min}(F)$. Pour un tel $z$, on  ne modifie pas $Y(h,N,a)$ en rempla\c{c}ant $h$ par $u(z)h$. On peut ensuite int\'egrer en $z$, tout en divisant par $mes(\mathfrak{p}_{F}^{ sup(c_{2},val_{F}(h_{n})+c_{2})})$.
On effectue le changement de variable $y\mapsto h^{-1}u(-z)hy$ et on obtient
$$(3) \qquad Y(h,N,a)= \int_{I\cap U_{min}(F)}\Xi^G(hya)\kappa_{N,h}(ya)dy,$$
o\`u
$$\kappa_{N,h}(ya)=mes(\mathfrak{p}_{F}^{ sup(c_{2},val_{F}(h_{n})+c_{2})})^{-1}\int_{\mathfrak{p}_{F}^{ sup(c_{2},val_{F}(h_{n})+c_{2})}}\kappa_{N}(u(-h_{n}^{-1}z)ya)dz.$$
Supposons $\kappa_{N}(u(-h_{n}^{-1}z)ya)=1$. Alors $a^{-1}y^{-1}u(h_{n}^{-1}z)v_{0}\in \mathfrak{p}^{-N}R$, c'est-\`a-dire $a^{-1}y^{-1}(v_{0}+h_{n}^{-1}ze_{n})\in \mathfrak{p}^{-N}R$. Il existe $c_{3}$ tel que $val_{F}(q_{V}(e_{-n},v))\geq c_{3}$ pour tout $v\in R$. Donc 
$$val_{F}(q_{V}(e_{-n},a^{-1}y^{-1}(v_{0}+h_{n}^{-1}ze_{n})))\geq c_{3}-N.$$
Puisque $y\in U_{min}(F)$, on a $y^{-1}e_{n}=e_{n}$, donc
$$q_{V}(e_{-n},a^{-1}y^{-1}(v_{0}+h_{n}^{-1}ze_{n}))=q_{V}(ae_{-n},y^{-1}v_{0}+h_{n}^{-1}ze_{n})=a_{n}^{-1}q_{V}(e_{-n},y^{-1}v_{0})+a_{n}^{-1}h_{n}^{-1}z.$$
En posant $z(h,y)=-h_{n}q_{V}(e_{-n},y^{-1}v_{0})$ et $c(h,N,a)=val_{F}(a_{n})+val_{F}(h_{n})+c_{3}-N$,  on a donc $z\in z(h,y)+\mathfrak{p}_{F}^{ c(h,N,a)}$. Alors
$$\kappa_{N,h}(ya)\leq mes(\mathfrak{p}_{F}^{ sup(c_{2},val_{F}(h_{n})+c_{2})})^{-1}mes(\mathfrak{p}_{F}^{sup(c_{2},val_{F}(h_{n})+c_{2})}\cap (z(h,y)+\mathfrak{p}_{F}^{c(h,N,a)}))$$
$$\leq inf(1,mes(\mathfrak{p}_{F}^{ sup(c_{2},val_{F}(h_{n})+c_{2})})^{-1}mes(\mathfrak{p}_{F}^{c(h,N,a)})).$$
On v\'erifie qu'il existe $\epsilon'>0$ tel que cette derni\`ere expression soit essentiellement born\'ee par $exp(-\epsilon'b(h,N,a))$. Alors (3) entra\^{\i}ne la majoration (2). Supposons maintenant que $v_{0}$ n'est pas orthogonal aux deux vecteurs $e_{n}$ et $e_{-n}$. C'est-\`a-dire que $d_{an,V}=d_{an,W}-1$ et $t=n$. On peut \'ecrire $v_{0}=\nu_{n}e_{n}+\nu_{-n}e_{-n}$. Introduisons le sous-groupe de $A_{min}(F)$ form\'e des \'el\'ements $a(z)$, pour $z\in 1+\mathfrak{p}_{F}$, d\'efinis ainsi: $a(z)_{n}=z$ et $a(z)_{j}=1$ pour tout $j=1,...,n-1$. Il existe un entier $c_{4}>0$ tel que $a(z)\in I\cap K$ pour tout $z\in 1+\mathfrak{p}_{F}^{c_{4}}$. Un tel $a(z)$ normalise $I\cap U_{min}(F)$, on peut donc effectuer dans $Y(h,N,a)$ le changement de variable $y\mapsto a(z)^{-1}ya(z)$, puis int\'egrer en $z$, \`a condition de diviser par $mes(1+\mathfrak{p}_{F}^{c_{4}})$. Puisque $a(z)$ commute \`a $a$ et $h$, on obtient
$$(4)\qquad Y(h,N,a)= \int_{I\cap U_{min}(F)}\Xi^G(hya)\kappa_{N}^*(ya)dy,$$
o\`u
$$\kappa_{N}^*(ya)=mes(1+\mathfrak{p}_{F}^{c_{4}})^{-1}\int_{1+\mathfrak{p}_{F}^{c_{4}}}\kappa_{N}(a(z)^{-1}ya)dz.$$
Supposons $\kappa_{N}(a(z)^{-1}ya)=1$. Alors $a^{-1}y^{-1}a(z)v_{0}\in \mathfrak{p}_{F}^{-N}R$ et, comme ci-dessus, 
$$val_{F}(q_{V}(e_{-n},a^{-1}y^{-1}a(z)v_{0}))\geq c_{3}-N.$$
 On a
$$q_{V}(e_{-n},a^{-1}y^{-1}a(z)v_{0})=q_{V}(ae_{-n},y^{-1}a(z)(\nu_{-n}e_{-n}+\nu_{n}e_{n}))$$
$$=a_{n}^{-1}q_{V}(e_{-n},y^{-1}(z^{-1}\nu_{-n}e_{-n}+z\nu_{n}e_{n}))=a_{n}^{-1}z^{-1}\nu_{-n}q_{V}(e_{-n},y^{-1}e_{-n})+a_{n}^{-1}z\nu_{n}.$$
Posons $z(y)=\nu_{-n}\nu_{n}^{-1}q_{V}(e_{-n},y^{-1}e_{-n})$ et $c(N,a)=val_{F}(a_{n})+c_{3}-N-val_{F}(\nu_{n})$. Alors $val_{F}(z(y)+z^2)\geq c(N,a)$. Notons ${\cal Z}(y,N,a)$ l'ensemble des $z$ qui v\'erifient cette condition. Alors, comme ci-dessus,
$$\kappa_{N}^*(ya)\leq inf(1,mes(1+\mathfrak{p}_{F}^{c_{4}})^{-1}mes({\cal Z}(y,N,a))).$$
On a $mes({\cal Z}(y,N,a)<< mes(\mathfrak{p}_{F}^{c(N,a)})$. Puisque $t=n$, on a $h_{n}=1$ par d\'efinition de notre syst\`eme hyperbolique. On v\'erifie qu'il existe $\epsilon''>0$ tel que l'expression ci-dessus soit essentiellement born\'ee par $exp(-\epsilon''b(h,N,a))$. Alors (4) entra\^{\i}ne (2), ce qui ach\`eve la preuve de cette relation. 

Montrons

(5) il existe des r\'eels $D_{2}$ et $D_{3}$ tels que 
$$\int_{I\cap U_{min}(F)}\Xi^G(hya)dy <<\sigma(h)^{D_{2}}\Xi^G(h)\sigma(a)^{D_{3}}\delta_{P_{min}}(a)^{1/2}$$
pour tout $h\in A_{min}(F)$ et tout $a\in\Lambda^-$.

On peut fixer $s\in \mathfrak{S}$ et supposer $h\in A_{min}(F)_{s}^-$. L'\'el\'ement $s$ d\'etermine un sous-groupe parabolique minimal $P_{min,s}=M_{min}U_{min,s}$ form\'e des \'el\'ements de $G$ qui conservent le drapeau
$$Fe_{sn}\subset Fe_{sn}\oplus Fe_{s(n-1)}\subset...\subset Fe_{sn}\oplus...\oplus F_{s1}.$$
On note $\bar{U}_{min,s}$ le radical unipotent de $\bar{P}_{min,s}$. On a l'\'egalit\'e $I\cap U_{min}(F)=(I\cap U_{min}(F)\cap U_{min,s}(F))(I\cap U_{min}(F)\cap \bar{U}_{min,s}(F))$. Pour $y\in I\cap U_{min}(F)\cap U_{min,s}(F)$, on a $hyh^{-1}\in I$. Donc
 $$\int_{I\cap U_{min}(F)}\Xi^G(hya)dy <<\int_{I\cap U_{min}(F)\cap \bar{U}_{min,s}(F)}\Xi^G(hya)dy$$
 $$<<\delta_{0}(h)\int_{h(I\cap U_{min}(F)\cap \bar{U}_{min,s}(F))h^{-1}}\Xi^G(yha)dy,$$
 o\`u $\delta_{0}(h)$ est la valeur absolue du d\'eterminant de $ad(h^{-1})$ agissant sur $\mathfrak{u}_{min}(F)\cap \bar{\mathfrak{u}}_{min,s}(F)$.  Pour $v\in I\cap U_{min}(F)\cap U_{min,s}(F)$, on a $\Xi^G(vyha)=\Xi^G(yha)$. Donc
 $$\int_{I\cap U_{min}(F)}\Xi^G(hya)dy << \delta_{0}(h)\int_{I\cap U_{min}(F)\cap U_{min,s}(F)}\int_{h(I\cap U_{min}(F)\cap \bar{U}_{min,s}(F))h^{-1}}\Xi^G(vyha)dy\, dv.$$
 Dans le domaine d'int\'egration, on  a $\sigma(vy)<<\sigma(h)$. Pour tout r\'eel $D_{3}>0$, on a donc
 $$\int_{I\cap U_{min}(F)}\Xi^G(hya)dy << \delta_{0}(h)\int_{I\cap U_{min}(F)\cap U_{min,s}(F)}$$
 $$\int_{h(I\cap U_{min}(F)\cap \bar{U}_{min,s}(F))h^{-1}}\Xi^G(vyha)\sigma(h)^{D_{3}}\sigma(vy)^{-D_{3}}dy\, dv$$
 $$<<\delta_{0}(h)\sigma(h)^{D_{3}}\int_{U_{min}(F)}\Xi^G(uha) \sigma(u)^{-D_{3}}du.$$
  D'apr\`es [W2], proposition II.4.5, il existe un r\'eel $D_{3}\geq0$ tel que la derni\`ere int\'egrale  soit convergente et essentiellement born\'ee par $\sigma(ha)^{D_{3}}\delta_{P_{min}}(ha)^{1/2}$ pour tout $x\in A_{min}(F)$. Fixons un tel $D_{3}$. On obtient
   $$(6) \qquad \int_{I\cap U_{min}(F)}\Xi^G(hya)dy <<\sigma(h)^{2D_{3}}\sigma(a)^{D_{3}}\delta_{0}(h)\delta_{P_{min}}(ha)^{1/2}.$$
 On calcule $\delta_{0}(h)\delta_{P_{min}}(h)^{1/2}=\delta_{P_{min,s}}(h)^{1/2}$. D'apr\`es [W2], lemme II.1.1, on a
 $$\delta_{P_{min,s}}(h)^{1/2}<<\Xi^G(h)$$
 puisque $h\in A_{min}(F)_{s}^-$. Alors (6) entra\^{\i}ne (5).
 
 Gr\^ace \`a (2) et (5), l'expression (1) est born\'ee par
 $$\sigma(h)^{D_{2}}\Xi^G(h)\sum_{a\in \Lambda^-}\sigma(a)^{D_{1}+D_{3}}exp(-\epsilon'b(h,N,a)).$$
 On peut identifier   $\Lambda^-$ \`a  l'ensemble ${\cal M}$ des  familles $\underline{m}=(m_{n},...,m_{1})$ d'entiers telles que $m_{n}\geq m_{n-1}\geq...\geq m_{1}\geq0$.   Pour une telle famille, posons $b(h,N,\underline{m})=b(h,N,a)$ o\`u $a\in \Lambda^-$ correspond \`a la famille $\underline{m}$. C'est-\`a-dire que $b(h,N,\underline{m})=sup(0,m_{n}-N(h))$. On v\'erifie que
  $$\sum_{\underline{m}\in{\cal M}} exp(-\epsilon'b(h,N,\underline{m}))$$
  est convergente et qu'il existe un r\'eel $D_{4}$ tel que cette expression soit essentiellement major\'ee par $N(h)^{D_{4}}$. De plus $N(h)^{D_{4}}<<N^{D_{4}}(1+\vert val_{F}(h_{n})\vert )^{D_{4}}<<N^{D_{4}}\sigma(h)^{D_{4}}$. Alors l'expression (1) est born\'ee par
  $$N^{D_{4}}\sigma(h)^{D_{2}+D_{4}}\Xi^G(h).$$
 C'est ce qu'on voulait d\'emontrer. $\square$
  
   \bigskip
  
  \subsection{Majoration d'une int\'egrale de fonctions d'Harish-Chandra, cas $r>0$}
  
  Soient $D$ un r\'eel, $C>0$ un r\'eel et $c\geq1$, $N\geq1$ deux entiers. Posons
  $$\chi(c,C,N,D)=\int_{M(F)}\int_{U(F)_{c}}{\bf 1}_{\sigma\geq C}(u)\Xi^M(m)\Xi^G(um)\kappa_{N}(m)\delta_{P}(m)^{-1/2}\sigma(u)^D\sigma(m)^Ddu\,dm.$$
  
  \ass{Lemme}{Cette expression est convergente. Le r\'eel $D$ \'etant fix\'e, pour tout r\'eel $R$, il existe $\alpha>0$ tel que
  $$\chi(c,C,N,D)<<exp(-cR)N^{-R}$$
  pour tout $c\geq1$, $N\geq1$ et tout $C$ tel que $C\geq \alpha(log(N)+c)$.}
  
  Preuve. Pour $i=1,...,r$, on a introduit en 4.5 le sous-groupe parabolique $P_{i}=M_{i}U_{i}$. Posons $U'_{i}=M_{i+1}\cap U_{i}$, avec la convention $M_{r+1}=G$. Le groupe $U$ est produit de ces groupes $U'_{i}$. D'o\`u
  $$\chi(c,C,N,D)=\int_{M(F)}\int_{U'_{1}(F)\cap U(F)_{c}}...\int_{U'_{r}(F)\cap U(F)_{c}}{\bf 1}_{\sigma\geq C}(u_{r}...u_{1})\Xi^M(m)\Xi^G(u_{r}...u_{1}m)$$
  $$\kappa_{N}(m)\delta_{P}(m)^{-1/2}\sigma(u_{r}...u_{1})^D\sigma(m)^Ddu_{r}...du_{1}\,dm.$$
  Introduisons une fonction $b$ sur $\{1,...,r\}$, \`a valeurs r\'eelles strictement positives, que nous pr\'eciserons par la suite.  Si nous supposons
  $$(1) \qquad C\geq  \sum_{i=1,...,r}b(i),$$
  la condition ${\bf 1}_{\sigma\geq C}(u_{r}...u_{1})=1$ entra\^{\i}ne qu'il existe $i$ tel que ${\bf 1}_{\sigma\geq b(i)}(u_{i})=1$. Donc $\chi(c,C,N,D)$ est major\'ee par la somme sur les sous-ensembles non vides $J$ de $\{1,...,r\}$ des $\chi(c,C,N,D;J)$, o\`u, dans ce dernier terme, on restreint l'int\'egration aux $u_{i}$ v\'erifiant les conditions
  
  - si $i\in J$, $\sigma(u_{i})\geq b(i)$;
  
  - si $i\not\in J$, $\sigma(u_{i})< b(i)$.
  
  On peut fixer $J$ et majorer $\chi(c,C,N,D;J)$. Notons $j$ le plus petit \'el\'ement de $J$. On a
  $$(2) \quad \chi(c,C,N,D;J)<<\int_{M(F)}\int_{U'_{1}(F) }...\int_{U'_{j-1}(F)}\int_{U'_{j}(F)\cap U(F)_{c}}\int_{U_{j+1}(F)\cap U(F)_{c}}(\prod_{i=1,...,j-1}{\bf 1}_{\sigma<b(i)}(u_{i}))$$
  $${\bf 1}_{\sigma\geq b(j)}(u_{j})\Xi^M(m)\Xi^G(u_{j+1}u_{j}...u_{1}m)\kappa_{N}(m)\delta_{P}(m)^{-1/2}\sigma(u_{j+1}u_{j}...u_{1})^D\sigma(m)^Ddu_{j+1}...du_{1}\,dm.$$
  La m\^eme preuve qu'en 4.7 permet de majorer 
  $$\int_{U_{j+1}(F)\cap U(F)_{c}}\Xi^G(u_{j+1}u_{j}...u_{1}m)\sigma(u_{j+1}...u_{1})^Ddu_{j+1}$$
  par 
  $$c^{D_{1}}\sigma(u_{j}...u_{1})^{D_{1}}\sigma(m)^{D_{1}}\delta_{P_{j+1}}(u_{j}...u_{1}m)^{1/2}\Xi^{M_{j+1}}(u_{j}....u_{1}m)$$
  pour un r\'eel $D_{1}$ convenable.  On $\delta_{P_{j+1}}(u_{j}...u_{1}m)=\delta_{P_{j+1}}(m)$. Le groupe $M_{j+1}$ est de la forme $A_{j+1}\tilde{G}$, o\`u $\tilde{G}$ est un groupe sp\'ecial orthogonal et $A_{j+1}$ est le plus grand sous-tore central et d\'eploy\'e de $M_{j+1}$. On a $M=A_{j+1}\tilde{M}$, o\`u $\tilde{M}=M\cap \tilde{G}$. On remarque que, pour $a\in A_{j+1}(F)$ et $\tilde{m}\in \tilde{M}(F)$, on a $\Xi^{M_{j+1}}(u_{j}...u_{1}a\tilde{m})=\Xi^{\tilde{G}}(u_{j}...u_{1}\tilde{m})$, $\Xi^M(a\tilde{m})=\Xi^{M\cap \tilde{G}}(\tilde{m})$,  $\delta_{P}(a\tilde{m})^{-1/2}\delta_{P_{j+1}}(u_{j}...u_{1}m)^{1/2}=\delta_{P\cap \tilde{G}}(\tilde{m})^{-1/2}$ et $\kappa_{N}(a\tilde{m})=\kappa_{N}(a)\kappa^{\tilde{G}}_{N}(\tilde{m})$, o\`u $\kappa_{N}^{\tilde{G}}$ est l'analogue de $\kappa_{N}$ pour le groupe $\tilde{G}$. Alors $\chi(c,C,N,D;J)$ est born\'e par le produit de $c^{D_{1}}$, de l'int\'egrale
  $$\int_{A_{j+1}(F)}\kappa_{N}(a)\sigma(a)^{D+D_{1}}da$$
  et de l'int\'egrale
  $$\int_{\tilde{M}(F)}\int_{U'_{1}(F) }...\int_{U'_{j-1}(F)}\int_{U'_{j}(F)\cap U(F)_{c}}(\prod_{i=1,...,j-1}{\bf 1}_{\sigma<b(i)}(u_{i})){\bf 1}_{\sigma\geq b(j)}(u_{j})\Xi^{M\cap \tilde{G}}(\tilde{m})\Xi^{\tilde{G}}(u_{j}...u_{1}\tilde{m})$$
  $$\kappa^{\tilde{G}}_{N}(\tilde{m})\delta_{P\cap \tilde{G}}(\tilde{m})^{-1/2}\sigma(u_{j}...u_{1})^{D+D_{1}}\sigma(\tilde{m})^{D+D_{1}}du_{j}...du_{1}\,d\tilde{m}.$$
 L'int\'egrale sur $A_{j+1}(F)$ est produit d'int\'egrales
  $$\int_{x\in F; \vert val_{F}(x)\vert \leq N}(1+\vert val_{F}(x)\vert )^{D'}\vert x\vert _{F}^{-1}dx$$
  qui sont convergentes et born\'ees par une puissance de $N$. La deuxi\`eme int\'egrale est exactement le membre de droite de (2) quand on remplace $G$ par $\tilde{G}$, $j$ par le nombre $\tilde{r}$ analogue de $r$ pour $\tilde{G}$, et $D$ par $D+D_{1}$. Cela nous ram\`ene au cas o\`u $j=r$, ce que nous supposons d\'esormais. On doit se rappeler qu'\`a la fin, il faudra multiplier la majoration obtenue par $c^{D_{1}}N^{D_{2}}$, pour un certain r\'eel $D_{2}$.
  
Supposons donc $j=r$, c'est-\`a-dire $J=\{r\}$. On effectue le changement de variable $u_{r}\mapsto u_{r-1}...u_{1}u_{r}u_{1}^{-1}...u_{r-1}^{-1}$. Cela ne change pas le domaine d'int\'egration, mais change ${\bf 1}_{\sigma\geq b(r)}(u_{r})$ en ${\bf 1}_{\sigma\geq b(r)}(u_{r-1}...u_{1}u_{r}u_{1}^{-1}...u_{r-1}^{-1})$. D\'efinissons une fonction $b_{1}$ sur $\{1,...,r\}$ par
$b_{1}(i)=b(i)-2\sum_{i'<i}b(i')$.  Supposons $b_{1}(i)>0$ pour tout $i$. Si $\sigma(u_{i})<b(i)$ pour tout $i=1,...,r-1$ et  ${\bf 1}_{\sigma\geq b(r)}(u_{r-1}...u_{1}u_{r}u_{1}^{-1}...u_{r-1}^{-1})=1$, on a ${\bf 1}_{\sigma\geq b_{1}(r)}(u_{r})=1$. D'o\`u
  $$ \chi(c,C,N,D;\{r\})<<\int_{M(F)}\int_{U'_{1}(F) }...\int_{U'_{r-1}(F)}\int_{U'_{r}(F)\cap U(F)_{c}} (\prod_{i=1,...,j-1}{\bf 1}_{\sigma<b(i)}(u_{i})){\bf 1}_{\sigma\geq b_{1}(r)}(u_{r})$$
  $$\Xi^M(m)\Xi^G(u_{r-1}...u_{1}u_{r}m)\kappa_{N}(m)\delta_{P}(m)^{-1/2}\sigma(u_{r-1}...u_{1}u_{r})^D\sigma(m)^Ddu_{r}...du_{1}\,dm.$$
  On a $U'_{r}=U_{r}$ et on peut d\'ecomposer ce groupe en $U_{r,r-1}U_{r,\natural}$. Il existe $c_{1}>0$ tel que, pour $v\in U_{r,r-1}(F)$ et $u\in U_{r,\natural}(F)$, les conditions $vu\in U(F)_{c}$ et ${\bf 1}_{\sigma\geq b_{1}(r)}(vu)$ entra\^{\i}nent $\sigma(u)\geq b_{1}(r)-c_{1}c$. Renfor\c{c}ons la condition $b_{1}(i)>0$ en
  
  (3) $b_{1}(i)>c_{1}c$ pour tout $i$,
  
  \noindent et posons $b_{2}(i)=b_{1}(i)-c_{1}c$. Alors
  $$ \chi(c,C,N,D;\{r\})<<\int_{M(F)}\int_{U'_{1}(F) }...\int_{U'_{r-1}(F)}\int_{U_{r,r-1}(F)\cap U(F)_{c}}\int_{U_{r,\natural}(F)} (\prod_{i=1,...,j-1}{\bf 1}_{\sigma<b(i)}(u_{i}))$$
  $${\bf 1}_{\sigma\geq b_{2}(r)}(u)\Xi^M(m)\Xi^G(u_{r-1}...u_{1}vum)\kappa_{N}(m)\delta_{P}(m)^{-1/2}\sigma(u_{r-1}...u_{1}vu)^D\sigma(m)^Ddu\,dv\,du_{r-1}...du_{1}\,dm.$$
  Gr\^ace \`a 3.3(5), il existe $c_{2}>0$ tel que
  $$\Xi^G(u_{r-1}...u_{1}vum)<<exp(c_{2}\sigma(u_{r-1}...u_{1}v))\Xi^G(um).$$
  Dans le domaine d'int\'egration, on a $\sigma(u_{r-1}...u_{1}v)<<c+\sum_{i=1,...,r-1}b(i)$. Il existe donc $c_{3}>0$ tel que
  $$\Xi^G(u_{r-1}...u_{1}vum)<<exp(c_{3}(c+\sum_{i=1,...,r-1}b(i)))\Xi^G(um).$$
  Alors $\chi(c,C,N,D;\{r\})$ est born\'e par le produit de $exp(c_{3}(c+\sum_{i=1,...,r-1}b(i)))$, de l'int\'egrale
  $$\int_{U'_{1}(F) }...\int_{U'_{r-1}(F)}\int_{U_{r,r-1}(F)\cap U(F)_{c}}(\prod_{i=1,...,j-1}{\bf 1}_{\sigma<b(i)}(u_{i}))\sigma(u_{r-1}...u_{1}v)^Ddv\,du_{r-1}...du_{1},$$
  et de l'int\'egrale
  $$Z(b_{2}(r),N,D)=\int_{M(F)}\int_{U_{r,\natural}(F)}{\bf 1}_{\sigma\geq b_{2}(r)}(u)\Xi^M(m)\Xi^G(um)\kappa_{N}(m)\delta_{P}(m)^{-1/2} \sigma(u)^D\sigma(m)^Ddu\,dm.$$
  La m\^eme preuve qu'en 4.4 montre que la premi\`ere int\'egrale est born\'ee par une exponentielle de $c+\sum_{i=1,...,r-1}b(i)$. Il existe donc $c_{4}>0$ tel que
  $$(4) \qquad \chi(c,C,N,D;\{r\})<<exp(c_{4}(c+\sum_{i=1,...,r-1}b(i)))Z(b_{2}(r),N,D).$$
 On doit majorer $Z(b_{2}(r),N,D)$. On commence par changer la variable $u$ en $u^{-1}$. D'apr\`es [W2], lemme II.1.1 et II.3..2, il existe un r\'eel $D_{3}$ tel que, pour tout $g\in G(F)$, on ait 
 $$\Xi^G(g)=\Xi^G(g^{-1})<<\delta_{\bar{P}}(m_{\bar{P}}(g^{-1}))^{1/2}\Xi^M(m_{\bar{P}}(g^{-1}))\sigma(g)^{D_{3}}.$$
  On applique cette relation \`a $g=u^{-1}m$. On a $m_{\bar{P}} (g^{-1})=m^{-1}m_{\bar{P}}(u)$ et $\delta_{\bar{P}}(m^{-1})=\delta_{P}(m)$, d'o\`u
  $$Z(b_{2}(r),N,D)<<\int_{M(F)}\int_{U_{r,\natural}(F)}{\bf 1}_{\sigma\geq b_{2}(r)}(u)\Xi^M(m)\delta_{\bar{P}}(m_{\bar{P}}(u))^{1/2}\Xi^M(m^{-1}m_{\bar{P}}(u))\kappa_{N}(m)$$
  $$\sigma(u)^{D+D_{3}}\sigma(m)^{D+D_{3}}du\,dm.$$
  On d\'ecompose $M$ en $AG_{0}$.   Comme plus haut, on peut majorer l'int\'egrale sur $A(F)$ par une puissance de $N$ et on obtient qu'il existe $D_{4}$ tel que
    $$Z(b_{2}(r),N,D)<<N^{D_{4}}\int_{G_{0}(F)}\int_{U_{r,\natural}(F)}{\bf 1}_{\sigma\geq b_{2}(r)}(u)\Xi^{G_{0}}(x)\delta_{\bar{P}}(a(u))^{1/2}\Xi^{G_{0}}(x^{-1} g_{0}(u))\kappa_{N}(x)$$
    $$\sigma(u)^{D+D_{3}}\sigma(x)^{D+D_{3}}du\,dx,$$
    o\`u, pour tout $u\in U_{r,\natural}(F)$, on  a \'ecrit $m_{\bar{P}}(u)=a(u)g_{0}(u)$, avec $a(u)\in A(F)$ et $g_{0}(u)\in G_{0}(F)$. Supposons d'abord $r\geq2$. On remarque qu'alors $U_{r,\natural}(F)$ est invariant par conjugaison par $G_{0}(F)$, a fortiori par $K\cap G_{0}(F)$. Soit $k\in K\cap G_{0}(F)$. On peut remplacer la variable $u$ par $kuk^{-1}$, puis int\'egrer en $k$. Changer $u$ en $kuk^{-1}$ ne modifie qu'une seule des fonctions que l'on int\`egre, \`a savoir $\Xi^{G_{0}}(x^{-1}g_{0}(u))$. On a $g_{0}(kuk^{-1})=kg_{0}(u)k^{-1}$ d'o\`u, puisque $\Xi^{G_{0}}$ est invariante par $K\cap G_{0}(F)$,  $\Xi^{G_{0}}(x^{-1}g_{0}(kuk^{-1}))=\Xi^{G_{0}}(x^{-1}kg_{0}(u))$. Or, d'apr\`es [W2], lemme II.1.3, l'int\'egrale de ce terme sur $k\in K\cap G_{0}(F)$ est \'egale \`a $\Xi^{G_{0}}(x^{-1})\Xi^{G_{0}}(g_{0}(u))$, ou encore \`a $\Xi^{G_{0}}(x)\Xi^{G_{0}}(g_{0}(u))$.   On voit alors que
    $$Z(b_{2}(r),N,D)<<N^{D_{4}}\chi^{G_{0}}(1,N,D+D_{3})I_{r,\natural}(b_{2}(r),D+D_{3}),$$
    avec les notations introduites en 4.5 et 4.11 (l'exposant $G_{0}$ indiquant que le groupe ambiant est $G_{0}$ au lieu de $G$). D'apr\`es les lemmes de ces paragraphes, il y a un r\'eel $D_{5}$ et un r\'eel $\epsilon>0$ tel que
    $$(5) \qquad Z(b_{2}(r),N,D)<<N^{D_{5}}exp(-\epsilon b_{2}(r)).$$
    Supposons maintenant $r=1$. Dans ce cas, on introduit le sous-espace $V_{\sharp}$ de $V$ orthogonal \`a $D_{0}$ et son groupe sp\'ecial orthogonal $G_{\sharp}$. On fixe un sous-groupe compact sp\'ecial $K_{\sharp}$ de $G_{\sharp}(F)$. Le groupe $U_{1,\natural}$ est contenu dans le groupe $G_{\sharp}$.  Ecrivons $m_{\bar{P}_{\sharp}}(u)=a_{\sharp}(u)g_{0,\sharp}(u)$, avec $a_{\sharp}(u)\in A(F)$ et $g_{0,\sharp}(u)\in G_{0}(F)\cap G_{\sharp}(F)=H(F)$. Comme dans la preuve du lemme 4.6, les \'el\'ements $a(u)a_{\sharp}(u)^{-1}$ et $g_{0}(u)g_{0,\sharp}(u)^{-1}$ restent dans des compacts. En utilisant les relations $\Xi^{G_{0}}(x^{-1}g_{0}(u))<<\Xi^{G_{0}}(x^{-1}g_{0,\sharp}(u))=\Xi^{G_{0}}(g_{0,\sharp}(u)^{-1}x)$, on voit que  l'on a
      $$Z(b_{2}(r),N,D)<<N^{D_{4}}\int_{U_{r,\natural}(F)}{\bf 1}_{\sigma\geq b_{2}(r)}(u) \delta_{\bar{P}}(a(u))^{1/2} \chi^{G_{0}}(g_{0,\sharp}(u)^{-1},N,D+D_{3})\sigma(u)^{D+D_{3}} du.$$
      En utilisant les lemmes 4.11 puis 4.6,   on obtient encore une majoration de la forme (5).
      
       Les formules (4) et (5) fournissent une majoration de $\chi(c,C,N,D;\{r\})$. Revenons au terme plus g\'en\'eral $\chi(c,C,N,D;J)$. On doit remplacer $r$ par $j$. On doit aussi multiplier la majoration issue de (4) et (5) par $c^{D_{1}}N^{D_{2}}$. Mais le terme $c^{D_{1}}$ est absorb\'e par le facteur $exp(c_{4}c)$ qui figure dans (4), quitte \`a augmenter $c_{4}$. On obtient qu'il existe des r\'eels $c_{5}$, $D_{5}$, $\epsilon$, tous strictement positifs, tels que
       $$(6)\qquad \chi(c,C,N,D;J)<<N^{D_{5}}exp(c_{5}(c+\sum_{i=1,...,j-1}b(i))-\epsilon b_{2}(j)).$$
    Soit $R$ un r\'eel.   Fixons, ind\'ependamment de $c$, $C$ et $N$,  une fonction $b^*$ sur $\{1,...,r\}$, \`a valeurs r\'eelles strictement positives. On en d\'eduit comme ci-dessus une fonction $b^*_{1}$. Supposons que $b^*_{1}$  v\'erifie
    $$b^*_{1}(i)> c_{1},$$
           $$\epsilon b^*_{1}(i) -c_{5}\sum_{i'=1,...,i-1}b^*(j)>sup(R+c_{5}+\epsilon c_{1},R+D_{5})$$
           pour tout $i$.  Une telle fonction existe, ces conditions pouvant \^etre impos\'ees par r\'ecurrence sur $i$. Prenons pour fonction $b$ la fonction $b(i)=(c+log(N))b^*(i)$. Cette fonction v\'erifie (3). On a
           $$c_{5}(c+\sum_{i=1,...,j-1}b(i))-\epsilon b_{2}(j)<-(R+D_{5})log(N)-Rc.$$
           Alors la majoration (6) devient
           $$\chi(c,C,N,D;J)<<N^{-R}exp(-Rc),$$
           c'est-\`a-dire celle que l'on voulait. La seule condition que l'on a impos\'ee \`a $C$ est la condition (1), qui s'\'ecrit $C\geq \alpha(c+log(N))$, o\`u $\alpha=\sum_{i=1,...,r}b^*(i)$. $\square$

  \bigskip
  
  \subsection{Preuve de la relation 4.3(2)}
  
  On veut majorer
  $$I(N,D)=\int_{G(F)}\Xi^G(g)^2\kappa_{N}(g)\sigma(g)^Ddg.$$
  Par la d\'ecomposition usuelle de la mesure $dg$, on a
  $$I(N,D)=\int_{K}\int_{M(F)}\int_{U(F)}\Xi^G(umk)^2\kappa_{N}(umk)\sigma(umk)^D\delta_{P}(m)^{-1}du\,dm\,dk.$$
  Le $k$ dispara\^{\i}t et l'int\'egrale sur $K$ \'egalement. Puisque $\kappa_{N}$ est invariante \`a gauche par $U(F)$, l'int\'egrale en $u$ est celle not\'ee $X(D,m,m)$ en 4.8. D'apr\`es le (iv) du lemme de ce paragraphe, on obtient
  $$I(N,D)<<\int_{M(F)}\Xi^M(m)^2\kappa_{N}(m)\sigma(m)^{D'}dm,$$
  pour un r\'eel $D'$ convenable. On d\'ecompose $M$ en $AG_{0}$.  Comme dans la preuve pr\'ec\'edente, on obtient 
  $$I(N,D)<<\int_{A(F)}\sigma(a)^{D'}\kappa_{N}(a)da\int_{G_{0}(F)}\sigma(g_{0})^{D'}\Xi^{G_{0}}(g_{0})^2\kappa_{N}(g_{0})dg_{0}.$$
   La premi\`ere int\'egrale est convergente et born\'ee par une puissance de $N$. La seconde int\'egrale n'est autre que $\chi(1,N,D'/2)$, avec la notation de 4.11 appliqu\'ee au groupe $G_{0}$. Par le lemme de ce paragraphe, elle est convergente et born\'ee par une puissance de $N$. D'o\`u la relation 4.3(2).
  
  \bigskip
  
  \subsection{Preuve de la relation 4.3(7)} 
  
   La conclusion que l'on veut obtenir nous autorise \`a supposer $c'\geq c$. Alors  l'ensemble $U(F)-U(F)_{c'}$ est invariant par translation par $U(F)_{c}$ et on peut d\'ecomposer l'int\'egrale en $u'\in U(F)-U(F)_{c'}$ en compos\'ee d'une int\'egrale sur $U(F)_{c}$ et d'une int\'egrale sur $(U(F)-U(F)_{c'})/U(F)_{c}$. C'est-\`a-dire
   $$I(c,c',N,D)= \int_{M(F)}\int_{H(F)U(F)_{c}}\int_{(U(F)-U(F)_{c'})/U(F)_{c}}\int_{U(F)_{c}}\phi(m,h,u,u'v_{1};D)dv_{1}\,du'\,du\,dh\,dm.$$
   Posons $v_{2}=u^{_{'}-1}u^{-1}h^{-1}u'v_{1}h$. Le groupe $H(F)U(F)_{c}$ normalise $U(F)_{c}$ et agit trivialement sur $U(F)/U(F)_{c}$. Quand $u$ d\'ecrit $U(F)_{c}$, $v_{2}$ d\'ecrit le m\^eme ensemble. On remplace la variable $u$ par $v_{2}$. Le jacobien de cette transformation est $1$ et on obtient
   $$I(c,c',N,D)<<\int_{M(F)}\int_{H(F)}\Xi^H(h)\sigma(h)^{D_{1}}\kappa_{N}(m)\sigma(m)^D\delta_{P}(m)^{-1}
   \int_{(U(F)-U(F)_{c'})/U(F)_{c}}$$
   $$\int_{U(F)_{c}}\int_{U(F)_{c}}\Xi^G(u'v_{1}m)\Xi^G(u'v_{2}h^{-1}m)\sigma(u'v_{1})^D\sigma(u'v_{2})^{D_{1}}dv_{1}\,dv_{2}\,du'\,dh\,dm$$
   pour un r\'eel $D_{1}$ convenable. La triple int\'egrale int\'erieure est essentiellement $X(c,c',D,m,h^{-1}m)$, avec la notation de 4.8. En appliquant le (iii) du lemme de ce paragraphe, et en remarquant que $\delta_{P}(h)=1$, on obtient
   $$I(c,c',N,D)<<exp(-\epsilon c')\int_{M(F)}\int_{H(F)}\Xi^H(h)\Xi^M(m)\Xi^M(h^{-1}m)\kappa_{N}(m)\sigma(h)^{D_{2}}\sigma(m)^{D_{2}}$$
   pour des r\'eels $D_{2}$ et $\epsilon>0$ convenables. Comme dans le paragraphe pr\'ec\'edent, on d\'ecompose l'int\'egrale sur $M(F)$ en produit d'int\'egrales sur $A(F)$ et $G_{0}(F)$. L'int\'egrale sur $A(F)$ est born\'ee par une puissance de $N$. D'apr\`es le lemme 4.11, l'int\'egrale sur $G_{0}(F)$ est born\'ee par $exp(-\epsilon' \sigma(h))\Xi^H(h)N^{D_{3}}$ pour des r\'eels $D_{3}$ et $\epsilon'>0$ convenables. Il reste une int\'egrale 
   $$\int_{H(F)}\Xi^H(h)^2\sigma(h)^{D_{2}}exp(-\epsilon'\sigma(h))dh,$$
   qui est convergente. Finalement, on a une majoration
   $$I(c,c',N,D)<<exp(-\epsilon c')N^{D_{4}},$$
   pour un r\'eel $D_{4}$ convenable. Soit $R$ un r\'eel. Il existe $\alpha>0$ tel que, si $c'\geq \alpha log(N)$, l'expression ci-dessus est major\'ee par $N^{-R}$. C'est ce qu'il fallait d\'emontrer. $\square$

 \bigskip
 
 \subsection{Preuve de la relation 4.3(8)}
 
 On a
 $$I(c,c',N,C,D)<<I(sup(c,c'),sup(c,c'),N,C,D).$$
 Puisque $c$ est fix\'e,  on peut aussi bien majorer le membre de droite, autrement dit supposer $c=c'$.  Introduisons un r\'eel $b>0$ que nous pr\'eciserons plus tard. On a
  $$I(c',c',N,C,D)=I_{\geq  b}(c',c',N,C,D)+I_{<b}(c',c',N,C,D),$$
 o\`u
 $$I_{\geq  b}(c',c',N,C,D)=\int_{M(F)}\int_{H(F)U(F)_{c'}}\int_{U(F)_{c'}}{\bf 1}_{\sigma\geq b}(h){\bf 1}_{\sigma\geq C}(hu)\phi(m,h,u,u';D)du'\,du\,dh\,dm,$$
 et
 $$I_{<b}(c',c',N,C,D)=\int_{M(F)}\int_{H(F)U(F)_{c'}}\int_{U(F)_{c'}}{\bf 1}_{\sigma< b}(h){\bf 1}_{\sigma\geq C}(hu)\phi(m,h,u,u';D)du'\,du\,dh\,dm.$$
 Dans la premi\`ere int\'egrale, on majore ${\bf 1}_{\sigma\geq C}(hu)$ par $1$. On effectue le changement de variable $u\mapsto h^{-1}u'hu^{-1}$ et on obtient
 $$I_{\geq  b}(c',c',N,C,D)<<\int_{M(F)}\int_{H(F)}{\bf 1}_{\sigma\geq b}(h)\Xi^H(h)\kappa_{N}(m)\sigma(h)^{D_{1}}\sigma(m)^D\delta_{P}(m)^{-1}$$
 $$\int_{U(F)_{c'}}\int_{U(F)_{c'}}\Xi^G(u'm)\Xi^G(uh^{-1}m)\sigma(u)^{D}\sigma(u')^{D_{1}}du'\,du\,dh\,dm,$$
 pour un r\'eel $D_{1}$ convenable. Les int\'egrales sur $U(F)_{c'}$ sont major\'ees par 4.3(3). On en d\'eduit une majoration
 $$I_{\geq  b}(c',c',N,C,D)<<c^{_{'}D_{2}}\int_{M(F)}\int_{H(F)}{\bf 1}_{\sigma\geq b}(h)\Xi^H(h)\Xi^M(m)\Xi^M(h^{-1}m)\kappa_{N}(m)$$
 $$\sigma(h)^{D_{2}}\sigma(m)^{D_{2}}dh\,dm,$$
 pour un r\'eel $D_{2}$ convenable. Ainsi qu'on l'a d\'ej\`a fait plusieurs fois, on d\'ecompose l'int\'egrale sur $M(F)$ en le produit d'int\'egrales sur $A(F)$ et $G_{0}(F)$. L'int\'egrale sur $A(F)$ est born\'ee par une puissance de $N$. L'int\'egrale sur $G_{0}(F)$ est l'expression $\chi^{G_{0}}(h^{-1},N,D_{2})$ de 4.11. En utilisant le lemme de ce paragraphe, on obtient
 $$I_{\geq  b}(c',c',N,C,D)<<c^{_{'}D_{3}}N^{D_{3}}\int_{H(F)}{\bf 1}_{\sigma\geq b}(h)\Xi^H(h)\Xi^{G_{0}}(h)\sigma(h)^{D_{3}}dh,$$
 pour un r\'eel $D_{3}$ convenable. D'apr\`es le lemme 4.9, il existe $\epsilon>0$ tel que ceci soit born\'e par
 $$c^{_{'}D_{3}}N^{D_{3}}\int_{H(F)}{\bf 1}_{\sigma\geq b}(h)\Xi^H(h)^2exp(-\epsilon \sigma(h))dh.$$
 Pour $\sigma(h)\geq b$, on a $exp(-\epsilon\sigma(h))\leq exp(-\epsilon b/2)exp(-\epsilon\sigma(h)/2)$, d'o\`u
 $$I_{\geq  b}(c',c',N,C,D)<<c^{_{'}D_{3}}N^{D_{3}}exp(-\epsilon b/2)\int_{H(F)}\Xi^H(h)^2exp(-\epsilon\sigma(h)/2)dh.$$
 La derni\`ere int\'egrale est convergente, d'o\`u
 $$(1) \qquad I_{\geq  b}(c',c',N,C,D)<<c^{_{'}D_{3}}N^{D_{3}}exp(-\epsilon b/2).$$
 
 Consid\'erons maintenant l'int\'egrale $I_{<b}(c',c',N,C,D)$.  On effectue encore le changement de variable $u\mapsto h^{-1}u'hu^{-1}$. La condition $\sigma(hu)\geq C$ devient $\sigma(u'hu^{-1})\geq C$. Jointe \`a la condition $\sigma(h)<b$, elle entra\^{\i}ne $\sigma(u)+\sigma(u')\geq C-b$, a fortiori $sup(\sigma(u),\sigma(u'))\geq(C-b)/2$. Supposons $  C-b>0$. Alors
 $$I_{<b}(c',c',N,C,D)<<I'_{<}(c',c',N,C,D)+I''_{<}(c',c',N,C,D),$$
 o\`u, pour un r\'eel $D_{4}$ convenable,
 $$I'_{<}(c',c',N,C,D)=\int_{M(F)}\int_{H(F)}\int_{U(F)_{c'}}\int_{U(F)_{c'}}{\bf 1}_{\sigma< b}(h){\bf 1}_{\sigma\geq ((C-b)/2)}(u') \Xi^H(h)\Xi^G(u'm)$$
 $$\Xi^G(uh^{-1}m)\kappa_{N}(m)\sigma(u')^{D_{4}}\sigma(u)^{D_{4}}\sigma(h)^{D_{4}}\sigma(m)^{D_{4}}\delta_{P}(m)^{-1}du'\,du\,dh\,dm,$$
 et $I''_{<}(c',c',N,D)$ est l'expression analogue o\`u ${\bf 1}_{\sigma\geq (C-b)/2}(u')$ est remplac\'e par ${\bf 1}_{\sigma\geq (C-b)/2}(u)$. En fait, le changement de variables $(h,m)\mapsto (h^{-1}, h^{-1}m)$ rend les deux expressions similaires. Il suffit donc de borner $I'_{<}(c',c',N,C,D)$. L'int\'egrale en $u$ se majore gr\^ace \`a 4.3(3). On obtient
 $$I'_{<}(c',c',N,C,D)<<c^{_{'}D_{5}}\int_{M(F)}\int_{H(F)}\int_{U(F)_{c'}}{\bf 1}_{\sigma< b}(h){\bf 1}_{\sigma\geq ((C-b)/2)}(u') \Xi^H(h)\Xi^G(u'm)$$
 $$\Xi^M(h^{-1}m)\kappa_{N}(m)\sigma(u')^{D_{5}}\sigma(h)^{D_{5}}\sigma(m)^{D_{5}}\delta_{P}(m)^{-1/2}du'\,dh\,dm,$$
 pour un r\'eel $D_{5}$ convenable. D'apr\`es 3.3(5), il existe $c_{1}>0$ tel que $\Xi^M(h^{-1}m)<<exp(c_{1}\sigma(h))\Xi^M(m)$. Puisque $\sigma(h)<b$, on obtient
 $$I'_{<}(c',c',N,C,D)<<c^{_{'}D_{5}} b^{D_{5}}exp(c_{1}b)\int_{H(F)}{\bf 1}_{\sigma< b}dh\int_{M(F)}\int_{U(F)_{c'}}{\bf 1}_{\sigma\geq ((C-b)/2)}(u') \Xi^G(u'm)$$
 $$\Xi^M(m)\kappa_{N}(m)\sigma(u')^{D_{5}}\sigma(m)^{D_{5}}\delta_{P}(m)^{-1/2}du'\,dm.$$
 L'int\'egrale sur $H(F)$ est major\'ee par 4.3(1). Le produit des quatre premiers termes ci-dessus est born\'e par $c^{_{'}D_{5}}exp(c_{2}b)$ pour un r\'eel $c_{2}>0$ convenable. La seconde int\'egrale est major\'ee par le lemme 4.12. Pr\'ecis\'ement, fixons un r\'eel $R'$, soit $\alpha$ la constante que ce lemme lui associe. Alors on obtient une majoration
 $$(2) \qquad I'_{<}(c',c',N,C,D)<<c^{_{'}D_{5}}exp(c_{2}b)exp(-R'c')N^{-R'},$$
 pourvu que $(C-b)/2\geq \alpha(log(N)+c')$. Soit $R>0$ un r\'eel, choisissons $b=2\epsilon^{-1}(R+D_{3})log(N)+Rc'(2c_{2})^{-1}$. La majoration (1) devient
 $$I_{\geq  b}(c',c',N,C,D)<<c^{_{'}D_{3}}exp(-\epsilon Rc'(4c_{2})^{-1})N^{-R}<<N^{-R}.$$
 Choisissons $R'=R+2c_{2}\epsilon^{-1}(R+D_{3})$. La majoration (2) devient
 $$I_{<b}'(c',c',N,C,D)<<c^{_{'}D_{5}}exp(-Rc'/2)N^{-R}<<N^{-R}.$$
 Cela vaut pour $(C-b)/2\geq \alpha(log(N)+c')$. Mais il existe $\alpha'>0$ tel que cette condition soit v\'erifi\'ee pour $C\geq \alpha'(log(N)+c')$. C'est ce qu'on voulait d\'emontrer. $\square$
 
 \bigskip
 
 \subsection{Preuve de la relation 4.3(6)}
 
 L'expression $I(c,N,D)$ est \'egale \`a la somme de $I(c,c,N,D)$ et de $I(c,c,N,1,D)$.  Il suffit de reprendre les preuves des deux paragraphes pr\'ec\'edents pour montrer que ces deux termes sont born\'es par une puissance de $N$. $\square$ 
 
 \bigskip

\section{Entrelacements temp\'er\'es}

\bigskip

\subsection{D\'efinition des entrelacements temp\'er\'es}

Soient $(V,q_{V})$ et $(W,q_{W})$ deux espaces quadratiques compatibles, avec $d_{W}<d_{V}$. On utilise les notations de la section 4. Soit $\pi$, resp. $\rho$, une repr\'esentation temp\'er\'ee de $G(F)$, resp. $H(F)$.   On note $Hom_{H,\xi}(\pi,\rho)$ l'espace des applications lin\'eaires $l:E_{\pi}\to E_{\rho}$ telles que $l(\pi(hu)e)=\xi(u)\rho(h)l(e)$ pour tous $h\in H(F)$, $u\in U(F)$, $e\in E_{\pi}$.  Dans le cas o\`u $\pi$ et $\rho$ sont irr\'eductibles, on sait que $Hom_{H,\xi}(\pi,\rho)$ est de dimension au plus $1$. En tout cas, cet espace est de dimension finie. On note $m(\rho,\pi)$ cette dimension.  Ce nombre ne d\'epend pas des diff\'erents choix que l'on a effectu\'es. Pour simplifier la r\'edaction, on note aussi $m(\pi,\rho)=m(\rho,\pi)$.On munit $E_{\pi}$ et $E_{\rho}$ de produits scalaires invariants.    On d\'efinit une forme sesquilin\'eaire ${\cal L}_{\pi,\rho,c}$ sur $E_{\rho}\otimes_{{\mathbb C}}E_{\pi}$ par
$${\cal L}_{\pi,\rho,c}(\epsilon'\otimes e',\epsilon\otimes e)=\int_{H(F)U(F)_{c}}(\rho(h)\epsilon',\epsilon)(e',\pi(hu)e)\bar{\xi}(u)dudh$$
pour tous $\epsilon,\epsilon'\in E_{\rho}$ et $e,e'\in E_{\pi}$. Cette expression est absolument convergente. En effet, elle est major\'ee par
$$\int_{H(F)U(F)_{c}}\Xi^H(h)\Xi^G(hu)dudh$$
qui est convergente d'apr\`es 4.3(4).

Pour tout entier $c'\geq1$, notons $\omega_{A}(c')$ le sous-groupe des $a\in A(F)$ tels que $val_{F}(a_{i}-1)\geq c'$ pour tout $i=1,...,r$.  

\ass{Lemme }{Pour tous $\epsilon,\epsilon',e,e'$, il existe un entier $c_{0}$ tel que ${\cal L}_{\pi,\rho,c}(\epsilon'\otimes e',\epsilon\otimes e)$ soit ind\'ependant de $c$ pour $c\geq c_{0}$. Plus pr\'ecis\'ement, soit $c'$ un entier $\geq1$. Alors il existe un entier $c_{0}$ tel que cette conclusion soit v\'erifi\'ee pour tous $\epsilon,\epsilon'$ et tous $e,e'\in E_{\pi}^{\omega_{A}(c')}$.}

 La preuve est similaire \`a celle du lemme 3.5. $\square$
 
 On d\'efinit une forme sesquilin\'eaire ${\cal L}_{\pi,\rho}$ sur  $E_{\rho}\otimes_{{\mathbb C}}E_{\pi}$ par
 $${\cal L}_{\pi,\rho}(\epsilon'\otimes e',\epsilon\otimes e)=lim_{c\to \infty}{\cal L}_{\pi,\rho,c}(\epsilon'\otimes e',\epsilon\otimes e).$$
 Cette forme v\'erifie les relations
 $${\cal L}_{\pi,\rho}(\rho(h)\epsilon'\otimes e',\epsilon\otimes\pi(hu) e)=\xi(u){\cal L}_{\pi,\rho}(\epsilon'\otimes e',\epsilon\otimes e)$$
 $${\cal L}_{\pi,\rho}(\epsilon'\otimes \pi(hu)e',\rho(h)\epsilon\otimes e)=\bar{\xi}(u){\cal L}_{\pi,\rho}(\epsilon'\otimes e',\epsilon\otimes e)$$
 pour tous $h\in H(F)$, $u\in U(F)$. Elle est donc combinaison lin\'eaire de fonctions
 $$(\epsilon'\otimes e',\epsilon\otimes e)\mapsto (\epsilon',l(e))(l'(e'),\epsilon)$$
 o\`u $l,l'\in Hom_{H,\xi}(\pi,\rho)$. En particulier, elle est nulle si $Hom_{H,\xi}(\pi,\rho)=\{0\}$.

Supposons $\pi$ et $\rho$ irr\'eductibles et ${\cal L}_{\pi,\rho}$ non nulle. L'espace $Hom_{H,\xi}(\pi,\rho)$ est de dimension $1$. On peut fixer un \'el\'ement non nul $l$ de cet espace et un nombre complexe $c\not=0$ tels que
$${\cal L}_{\pi,\rho}(\epsilon'\otimes e',\epsilon\otimes e)=c(\epsilon',l(e))(l(e'),\epsilon)$$
pour tous $\epsilon,\epsilon',e,e'$. Il en r\'esulte l'\'egalit\'e
$$(1) \qquad \vert {\cal L}_{\pi,\rho}(\epsilon'\otimes e',\epsilon\otimes e)\vert ^2=\vert {\cal L}_{\pi,\rho}(\epsilon'\otimes e,\epsilon'\otimes e)\vert \vert {\cal L}_{\pi,\rho}(\epsilon\otimes e',\epsilon\otimes e')\vert .$$
 En cons\'equence
 
 (2) il existe $\epsilon$ et $e$ tels que
$${\cal L}_{\pi,\rho}(\epsilon\otimes e,\epsilon\otimes e)\not=0.$$

On a suppos\'e $d_{W}<d_{V}$. Pour simplifier la r\'edaction, dans le cas  $d_{W}>d_{V}$, on pose
${\cal L}_{\pi,\rho}={\cal L}_{\rho,\pi}$. 
\bigskip

\subsection{Induction d'entrelacements temp\'er\'es}

Soient $(V,q_{V})$ et $(W,q_{W})$ deux espaces quadratiques compatibles, avec $d_{W}<d_{V}$. Soit $V=\tilde{V}\oplus Y$ une d\'ecomposition orthogonale, avec $Y$ hyperbolique. On fixe une base hyperbolique $(y_{\pm j})_{j=1,...,k}$ de $Y$, on note $Y^+$, resp. $Y^-$, le lagrangien engendr\'e par $(y_{j})_{j=1,...,k}$, resp. $(y_{-j})_{j=1,...,k}$. On suppose $k\geq1$. On note $Q$ le sous-groupe parabolique de $G$ form\'e des \'el\'ements qui conservent$Y^+$, $L$ son sous-groupe de L\'evi form\'e des \'el\'ements qui conservent de plus $Y^-$, $U_{Q}$ le radical unipotent de $Q$, $GL_{k}$ le groupe des automorphismes lin\'eaires de $Y^+$ et $\tilde{G}$ le groupe sp\'ecial orthogonal de $\tilde{V}$. On a $L=GL_{k}\times \tilde{G}$. On fixe un sous-groupe compact sp\'ecial $K$ de $G(F)$ en bonne position relativement \`a $L$.

Soit $\tilde{\pi}$, resp. $\mu$, $\rho$, une repr\'esentation admissible irr\'eductible et temp\'er\'ee de $\tilde{G}(F)$, resp. $GL_{k}(F)$, $H(F)$. On fixe des produits scalaires invariants sur les espaces de ces repr\'esentations. Introduisons la repr\'esentation induite $\pi=Ind_{Q}^G(\mu\otimes \tilde{\pi})$, que l'on r\'ealise dans l'espace $E_{\pi}=E_{Q,\mu\otimes \tilde{\pi}}^G$.   On munit $E_{\pi}$ du produit scalaire invariant
$$(e',e)=\int_{Q(F)\backslash G(F)}(e'(g),e(g))dg.$$
On d\'efinit les formes sesquilin\'eaires ${\cal L}_{\pi,\rho}$ sur $E_{\rho}\otimes_{{\mathbb C}}E_{\pi}$ et ${\cal L}_{\tilde{\pi},\rho}$ sur $E_{\rho}\otimes_{{\mathbb C}}E_{\tilde{\pi}}$.

\ass{Proposition}{La forme ${\cal L}_{\pi,\rho}$ est non nulle si et seulement si la forme ${\cal L}_{\tilde{\pi},\rho}$ est non nulle.}

Preuve.  Soient $e,e'\in E_{\pi}$ et $\epsilon,\epsilon'\in E_{\rho}$, choisissons un entier $c$ assez grand. On a l'\'egalit\'e
$${\cal L}_{\pi,\rho}(\epsilon'\otimes e',\epsilon\otimes e)={\cal L}_{\pi,\rho,c}(\epsilon'\otimes e',\epsilon\otimes e)=\int_{H(F)U(F)_{c}}(\rho(h)\epsilon',\epsilon)$$
$$\qquad\int_{Q(F)\backslash G(F)}(e'(g),e(ghu))dg\,\bar{\xi}(u)\,du\,dh.$$
Montrons que

(1) cette expression est absolument convergente.

On a
$$\int_{Q(F)\backslash G(F)}\vert (e'(g),e(ghu))\vert dg=\int_{K}\vert (e'(k),e(khu))\vert dk$$
$$\qquad =\int_{K}\vert (e'(k),(\mu\otimes\tilde{\pi})(l_{Q}(khu))e(k_{Q}(khu)))\vert \delta_{Q}(l_{Q}(khu))^{1/2}dk.$$
Parce que $\mu$ et $\tilde{\pi}$ sont temp\'er\'ees, il s'ensuit que
$$\int_{Q(F)\backslash G(F)}\vert (e'(g),e(ghu))\vert dg<<\int_{K}\delta_{Q}(l_{Q}(khu))^{1/2}\Xi^L(l_{Q}(khu))dk=\Xi^G(hu),$$
la derni\`ere \'egalit\'e r\'esultant de [W2] lemme II.1.6. L'expression \`a \'etudier est donc born\'ee par un multiple de
$$\int_{H(F)U(F)_{c}}\Xi^H(h)\Xi^G(hu)dh\,du$$
qui est convergente d'apr\`es 4.3(4). Cela prouve (1).

Traitons d'abord le cas o\`u $k\leq r$. On peut supposer $y_{j}=v_{k-r-j}$ pour $j=1,...,k$. Alors $P\subset \bar{Q}$. On a
$$\int_{Q(F)\backslash G(F)}(e'(g),e(ghu))dg=\int_{U_{\bar{Q}}(F)}(e'(\bar{v}'),e(\bar{v}'uh))d\bar{v}',$$
d'o\`u, par le changement de variables $u\mapsto \bar{v}^{_{'}-1}u$,
$$(2) \qquad {\cal L}_{\pi,\rho}(\epsilon'\otimes e',\epsilon\otimes e)=\int_{H(F)}(\rho(h)\epsilon',\epsilon)$$
$$\qquad \int_{(\bar{v}',u)\in U_{\bar{Q}}(F)\times U(F); \bar{v}^{_{'}-1}u\in U(F)_{c}}(e'(\bar{v}'),e(uh))\bar{\xi}(\bar{v}')^{-1}\bar{\xi}(u)du\,d\bar{v}'\,dh.$$
Posons $U_{k}=U\cap GL_{k}$ et $\tilde{U}=U\cap \tilde{G}$.  Remarquons que la restriction de $\xi$ \`a $U_{k}(F)$ est le caract\`ere de ce groupe d\'efini en 3.5 tandis que la restriction de $\xi$ \`a $\tilde{U}(F)$ est le caract\`ere analogue \`a $\xi$ quand on remplace $G$ par $\tilde{G}$ dans les d\'efinitions.  D\'ecomposons $u$ en $u=n \tilde{u}\bar{v}$, avec $n\in U_{k}(F)$, $\tilde{u}\in \tilde{U}(F)$ et $\bar{v}\in U_{\bar{Q}}(F)$.  La condition $\bar{v}^{_{'}-1}u\in U(F)_{c}$ \'equivaut \`a $n\in U_{k}(F)_{c}$, $\tilde{u}\in \tilde{U}(F)_{c}$ et $\bar{v}^{_{'}-1}\bar{v}\in U(F)_{c}\cap U_{\bar{Q}}(F)$. Cette derni\`ere condition s'\'ecrit concr\`etement 
$$val_{F}(q_{V}(\bar{v}v_{r-k},v_{k-r-1})-q_{V}(\bar{v}'v_{r-k},v_{k-r-1}))\geq-c.$$
De m\^eme, on a l'\'egalit\'e
$$\bar{\xi}(\bar{v}^{_{'}-1}\bar{v})=\psi(q_{V}(\bar{v}v_{r-k},v_{k-r-1})-q_{V}(\bar{v}'v_{r-k},v_{k-r-1})).$$
Dans l'expression (2), effectuons le changement de variables $\bar{v}\mapsto h\bar{v}h^{-1}$. cela ne modifie pas les conditions ci-dessus et le jacobien de cette transformation vaut $1$. On obtient
$${\cal L}_{\pi,\rho}(\epsilon'\otimes e',\epsilon\otimes e)=\int_{(\bar{v}',\bar{v})\in U_{\bar{Q}}(F)^2; \bar{v}^{_{'}-1}\bar{v}\in U(F)_{c}}\int_{H(F)\tilde{U}(F)_{c}}\int_{U_{k}(F)_{c}}(\rho(h)\epsilon',\epsilon)$$
$$\qquad(e'(\bar{v}'),\mu(n)\tilde{\pi}(h\tilde{u})e(\bar{v}))\bar{\xi}(n)\bar{\xi}(\tilde{u})\bar{\xi}(\bar{v}^{_{'}-1}\bar{v})dn\,dh\,d\tilde{u}\,d\bar{v}\,d\bar{v}'.$$
Cette expression est encore absolument convergente. Pour un entier $c'\geq 1$, introduisons le sous-groupe $\omega_{A}(c')\subset A(F)$ de 5.1. Choisissons $c'$ tel que $e$ et $e'$ soient invariants par $\omega_{A}(c')$ et que le caract\`ere central de $\mu$ soit trivial sur $1+\mathfrak{p}_{F}^{c'}$. Soient $a,a'\in \omega_{A}(c')$. Consid\'erons l'expression
$$(3)\qquad \int_{(\bar{v}',\bar{v})\in U_{\bar{Q}}(F)^2; \bar{v}^{_{'}-1}\bar{v}\in U(F)_{c}}\int_{H(F)\tilde{U}(F)_{c}}\int_{U_{k}(F)_{c}}(\rho(h)\epsilon',\epsilon)$$
$$\qquad((\mu\otimes\tilde{\pi})(a')e'(\bar{v}'),\mu(n)\tilde{\pi}(h\tilde{u})(\mu\otimes\tilde{\pi})(a)e(\bar{v}))\bar{\xi}(n)\bar{\xi}(\tilde{u})\bar{\xi}(\bar{v}^{_{'}-1}\bar{v})dn\,dh\,d\tilde{u}\,d\bar{v}\,d\bar{v}'.$$
Introduisons le sous-groupe $\omega'_{A}(c')$ form\'e des $\alpha\in \omega_{A}(c')$ tels que $\alpha_{r-k}=\alpha_{r-k+1}$.On a l'\'egalit\'e $\omega_{A}(c')=(Z_{k}(F)\cap \omega_{A}(c'))\times \omega'_{A}(c')$. Ecrivons $a=z\alpha$, $a'=z'\alpha'$ conform\'ement \`a cette d\'ecomposition. On a $\mu(z)=1=\mu(z')$ et on a l'\'egalit\'e
$$((\mu\otimes\tilde{\pi})(a')e'(\bar{v}'),\mu(n)\tilde{\pi}(h\tilde{u})(\mu\otimes\tilde{\pi})(a)e(\bar{v}))=(e'(\alpha'\bar{v}'),\mu(n)\tilde{\pi}(h\tilde{u})e(\alpha\bar{v})).$$
On effectue les changements de variables $\bar{v}\mapsto \alpha^{-1}\bar{v}\alpha$, $\bar{v}'\mapsto \alpha^{_{'}-1}\bar{v}'\alpha'$. Cela ne modifie pas le domaine d'int\'egration ni le terme $\bar{\xi}(\bar{v}^{_{'}-1}\bar{v})$. Puisque $e$, resp. $e'$, est invariant par $\alpha$, resp. $\alpha'$, les termes $\alpha$ et $\alpha'$ disparaissent. Donc (3) est ind\'ependant de $a$ et $a'$. On peut int\'egrer cette expression sur $\omega_{A}(c')^2$. L'expression obtenue reste absolument convergente et on obtient
$${\cal L}_{\pi,\rho}(\epsilon'\otimes e',\epsilon\otimes e)=\int_{(\bar{v}',\bar{v})\in U_{\bar{Q}}(F)^2; \bar{v}^{_{'}-1}\bar{v}\in U(F)_{c}}\bar{\xi}(\bar{v}^{_{'}-1}\bar{v})\int_{H(F)\tilde{U}(F)_{c}} (\rho(h)\epsilon',\epsilon) \bar{\xi}(\tilde{u})$$
$$\qquad\int_{U_{k}(F)_{c}} (e_{c'}'(\bar{v}'),\mu(n)\tilde{\pi}(h\tilde{u}) e_{c'}(\bar{v}))\bar{\xi}(n)dn\,dh\,d\tilde{u}\,d\bar{v}\,d\bar{v}',$$
o\`u
$$e'_{c'}(\bar{v'})=mes(\omega_{A}(c'))^{-1}\int_{\omega_{A}(c')}(\mu\otimes\tilde{\pi})(a)e'(\bar{v}')da,$$
$$e_{c'}(\bar{v})=mes(\omega_{A}(c'))^{-1}\int_{\omega_{A}(c')}(\mu\otimes\tilde{\pi})(a)e(\bar{v})da.$$
Introduisons une fonctionnelle de Whittaker non nulle $\phi$ sur l'espace $E_{\mu}$. Notons $\Phi:E_{\mu}\otimes E_{\tilde{\pi}}\to E_{\tilde{\pi}}$ l'application $\phi\otimes id$. D'apr\`es le lemme 3.7(ii), si $c$ est assez grand,  on a l'\'egalit\'e
$$\int_{U_{k}(F)_{c}} (e_{c'}'(\bar{v}'),\mu(n)\tilde{\pi}(h\tilde{u})e_{c'}(\bar{v}))\bar{\xi}(n)dn=C  (\Phi e'_{c'}(\bar{v}'),\tilde{\pi}(h\tilde{u})\Phi e_{c'}(\bar{v})),$$
o\`u $C$ est une constante non nulle. D'apr\`es le lemme 5.1 appliqu\'e au groupe $\tilde{G}$, si $c$ est assez grand, on a l'\'egalit\'e
$$\int_{H(F)\tilde{U}(F)_{c}} (\rho(h)\epsilon',\epsilon) \bar{\xi}(\tilde{u}) (\Phi e_{c'}'(\bar{v}'),\tilde{\pi}(h\tilde{u})\Phi e_{c'}(\bar{v}))dh\,d\tilde{u}={\cal L}_{\tilde{\pi},\rho}(\epsilon'\otimes \Phi e_{c'}'(\bar{v}'),\epsilon\otimes \Phi e_{c'}(\bar{v})).$$
On obtient
$$(4) \qquad {\cal L}_{\pi,\rho}(\epsilon'\otimes e',\epsilon\otimes e)=$$
$$\qquad C \int_{(\bar{v}',\bar{v})\in U_{\bar{Q}}(F)^2; \bar{v}^{_{'}-1}\bar{v}\in U(F)_{c}}{\cal L}_{\tilde{\pi},\rho}(\epsilon'\otimes \Phi e_{c'}'(\bar{v}'),\epsilon\otimes \Phi e_{c'}(\bar{v}))\bar{\xi}(\bar{v}^{_{'}-1}\bar{v})d\bar{v}\,d\bar{v}'.$$
Cette \'egalit\'e entra\^{\i}ne que, si ${\cal L}_{\tilde{\pi},\rho}$ est nulle, ${\cal L}_{\pi,\rho}$ l'est aussi. Inversement, supposons ${\cal L}_{\tilde{\pi},\rho}$ non nulle. On peut fixer des \'el\'ements $\epsilon,\epsilon'\in E_{\rho}$, $\tilde{e},\tilde{e}'\in E_{\tilde{\pi}}$, $\eta,\eta'\in E_{\mu}$ de sorte que
$${\cal L}_{\tilde{\pi},\rho}(\epsilon'\otimes \tilde{e}',\epsilon\otimes \tilde{e})\not=0,$$
$$\phi(\eta)=\phi(\eta')=1.$$
Choisissons un sous-groupe ouvert compact $\Omega$ de $U_{\bar{Q}}(F)$ sur lequel $\bar{\xi}$ soit trivial. Introduisons l'unique fonction $e\in E_{\pi}$ telle que, pour $\bar{v}\in U_{\bar{Q}}(F)$, on ait l'\'egalit\'e
$$e(\bar{v})=\left\lbrace\begin{array}{cc}\eta\otimes \tilde{e},&\text{  si  }\bar{v}\in \Omega,\\ 0,&\text{  sinon.}\\ \end{array}\right.$$
Si $c'$ est assez grand, on a:
$$\Phi e_{c'}(\bar{v})=\left\lbrace\begin{array}{cc} \tilde{e},&\text{  si  }\bar{v}\in \Omega,\\ 0,&\text{  sinon.}\\ \end{array}\right.$$
Introduisons la fonction similaire $e'\in E_{\pi}$. L'\'egalit\'e (4) entra\^{\i}ne
$${\cal L}_{\pi,\rho}(\epsilon'\otimes e',\epsilon\otimes e)=Cmes(\Omega)^2{\cal L}_{\tilde{\pi},\rho}(\epsilon'\otimes \tilde{e}',\epsilon\otimes \tilde{e}).$$
Ce terme n'est pas nul, donc ${\cal L}_{\pi,\rho}$ n'est pas nulle.

Traitons maintenant le cas o\`u $r<k$. On peut supposer $v_{i}=y_{-k+r-i}$ pour $i=1,...,r$ et $v_{0}=y_{r-k}+\nu_{0}y_{k-r}$, o\`u $\nu_{0}=q_{V}(v_{0})$. Alors $U_{\bar{Q}}\subset P$ et $QP$est un ouvert de Zariski de $G$.  On a l'\'egalit\'e
$$(5) \qquad \int_{Q(F)\backslash G(F)}(e'(g),e(ghu))dg=\int_{(M(F)\cap Q(F))\backslash M(F)}$$
$$\qquad \int_{(U(F)\cap Q(F))\backslash U(F)}(e'(x'm),e(x'mhu))dx' \,dm.$$
Posons $Y^+_{0}=Y^+\cap V_{0}$,  $Q_{0}=Q\cap G_{0}$. L'espace $Y^+_{0}$ a pour base $(y_{j})_{j=1,...,k-r}$ et le groupe $Q_{0}$ est le sous-groupe parabolique des \'el\'ements de $G_{0}$ qui conservent $Y^+_{0}$. Posons $Y^+_{W}=Y^+\cap W$, $Y^-_{W}=Y^-\cap W$, notons $Q^H$ le sous-groupe parabolique de $H$ form\'e des \'el\'ements qui conservent $Y^+_{W}$, $L^H$ son sous-groupe de L\'evi form\'e des \'el\'ements qui conservent de plus $Y^-_{W}$ et $U_{Q^H}$ le radical unipotent de $Q^H$.  L'espace $Y^+_{W}$ a pour base $(y_{j})_{j=1,...,k-r-1}$. Notons $W_{0}$ l'orthogonal de $Y^+_{W}\oplus Y^-_{W}$ dans $W$, $H_{0}$ son groupe sp\'ecial orthogonal et $GL_{k-r-1}$ le groupe des automorphismes lin\'eaires de $Y^+_{W}$. On a l'\'egalit\'e $L^H=GL_{k-r-1}\times H_{0}$. Posons $w_{0}=-\frac{1}{2\nu_{0}}y_{r-k}+\frac{1}{2}y_{k-r}$ et $D^H=Fw_{0}$. Alors $W_{0}$ est la somme directe orthogonale de $\tilde{V}$ et de $D^H$. Le groupe $Q_{0}\cap H$ est le sous-groupe des \'el\'ements de $H$ qui conservent $Y^+_{0}$ (disons que, pour quelques instants, on \'etend les scalaires \`a $\bar{F}$). Un tel \'el\'ement conserve $Y^+_{0}\cap W=Y^+_{W}$, donc appartient \`a $Q^H$. Soit $h\in Q^H$. Pour qu'il conserve $Y^+_{0}$, il doit envoyer $y_{k-r}$ dans $Y^+_{0}$. Mais il fixe $v_{0}$ et on a $w_{0}=-\frac{1}{2\nu_{0}}v_{0}+y_{k-r}$.  On a donc $hw_{0}\in w_{0}+Y^+_{0}$. Puisque $hw_{0}$ et $w_{0}$ appartiennent \`a $W$, cela force $hw_{0}\in w_{0}+Y^+_{W}$. La r\'eciproque est similaire. Donc $Q_{0}\cap H$ est le sous-groupe des  $h\in Q^H$ tels que $hw_{0}\in w_{0 }+Y^+_{W}$. Autrement dit
$$Q_{0}\cap H=GL_{k-r-1}\tilde{G}U_{Q^H}.$$
On v\'erifie que la restriction  du module $\delta_{Q_{0}}$  au groupe $Q_{0}(F)\cap H(F)$ est \'egale au module  de ce groupe: si on \'ecrit un \'el\'ement $h\in Q_{0}(F)\cap H(F)$ sous la forme $\delta\tilde{g}n$, avec $\delta\in GL_{k-r-1}(F)$, $\tilde{g}\in \tilde{G}(F)$, $n\in U_{Q^H}(F)$, ces modules co\"{\i}ncident avec $\vert det(\delta)\vert _{F}^{d_{\tilde{V}}+k-r-1}$.  L'application naturelle $Q_{0}(F)\backslash G_{0}(F)\to (M(F)\cap Q(F))\backslash M(F)$ est un isomorphisme. Le groupe $H(F)$ agit sur l'ensemble des sous-espaces isotropes de $V_{0}$ de dimension $k-r$. Il y a deux orbites: l'orbite ouverte des sous-espaces dont l'intersection avec $W$ est de dimension $k-r-1$ et l'orbite ferm\'ee des sous-espaces contenus dans $W$. L'espace $Y^+_{0}$ appartient \`a l'orbite ouverte. Il en r\'esulte que l'application naturelle $(Q_{0}(F)\cap H(F))\backslash H(F) \to Q_{0}(F)\backslash G_{0}(F)$ est injective et a pour image un ouvert  de l'espace d'arriv\'ee dont le compl\'ementaire est de mesure nulle. On en d\'eduit ais\'ement l'assertion suivante. Soit $\varphi:G_{0}(F)\to {\mathbb C}$ une fonction telle que $\varphi(qg)=\delta_{Q_{0}}(q)\varphi(g)$ pour tous $q\in Q_{0}(F)$, $g\in G_{0}(F)$. Supposons $\varphi$ absolument int\'egrable sur $Q_{0}(F)\backslash G_{0}(F)$. Alors la restriction de $\varphi$ \`a $H(F)$ est absolument int\'egrable sur $(Q_{0}(F)\cap H(F))\backslash H(F)$ et on a l'\'egalit\'e
$$\int_{Q_{0}(F)\backslash G_{0}(F)}\varphi(g)dg=\int_{(Q_{0}(F)\cap H(F))\backslash H(F)}\varphi(h) dh.$$
Dans l'\'egalit\'e (5), on peut donc remplacer l'int\'egration sur $(M(F)\cap Q(F))\backslash M(F)$ par une int\'egration sur $(Q_{0}(F)\cap H(F))\backslash H(F)$. On obtient
$$ \int_{Q(F)\backslash G(F)}(e'(g),e(ghu))dg=\int_{(Q_{0}(F)\cap H(F)\backslash H(F)}\int_{(U(F)\cap Q(F))\backslash U(F)} (e'(x'h'),e(x'h'hu))dx'\,dh',$$
puis
$${\cal L}_{\pi,\rho}(\epsilon'\otimes e',\epsilon\otimes e)=\int_{H(F)U(F)_{c}}\int_{(Q_{0}(F)\cap H(F))\backslash H(F)}$$
$$\qquad \int_{(U(F)\cap Q(F))\backslash U(F)} (\rho(h)\epsilon',\epsilon)(e'(x'h'),e(x'h'hu))\bar{\xi}(u)dx'\,dh'\,du\,dh.$$
 On effectue les changements de variables $u\mapsto (h'h)^{-1}uh'h$, puis $h\mapsto h^{_{'}-1}h$, on d\'ecompose ensuite l'int\'egrale sur $H(F)$ en une int\'egrale compos\'ee d'une int\'egrale sur $Q_{0}(F)\cap H(F)$ et d'une int\'egrale sur $(Q_{0}\cap H(F))\backslash H(F)$ (la mesure sur $Q_{0}(F)\cap H(F)$ doit \^etre une mesure de Haar \`a gauche). On obtient
$${\cal L}_{\pi,\rho}(\epsilon'\otimes e',\epsilon\otimes e)=\int_{((Q_{0}(F)\cap H(F))\backslash H(F))^2}\int_{Q_{0}(F)\cap H(F)} \int_{ (U(F)\cap Q(F))\backslash U(F)}$$
$$\qquad \int_{U(F)_{c}}(\rho(qh)\epsilon',\rho(h')\epsilon)(e'(x'h'),e(x'uqh))\bar{\xi}(u)du\,dx'\,dq\,dh\,dh',$$
 On effectue le changement de variable $u\mapsto x^{_{'}-1}u$ puis on d\'ecompose l'int\'egrale en $u\in U(F)$ en compos\'ee d'une int\'egrale sur $u\in U(F)\cap Q(F)$ et d'une int\'egrale sur $x\in (U(F)\cap Q(F))\backslash U(F)$. Remarquons que $U(F)\cap Q(F)=U(F)\cap L(F)=U(F)\cap GL_{k}(F)$.
  La condition initiale $u\in U(F)_{c}$ est remplac\'ee par $x^{_{'}-1}ux\in U(F)_{c}$. Notons $\varphi_{c}(u,x,x')$ la fonction caract\'eristique de l'ensemble des $(u,x,x')$  v\'erifiant cette condition. On obtient
$${\cal L}_{\pi,\rho}(\epsilon'\otimes e',\epsilon\otimes e)=\int_{((Q_{0}(F)\cap H(F))\backslash H(F))^2}\int_{Q_{0}(F)\cap H(F)}\int_{((U(F)\cap Q(F))\backslash U(F))^2}$$
$$\qquad \int_{U(F)\cap GL_{k}(F)}(\rho(qh)\epsilon',\rho(h')\epsilon)(e'(x'h'),e(uxqh))\varphi_{c}(u,x,x')\bar{\xi}(x^{_{'}-1}ux)du \,dx\,dx'\,dq\,dh\,dh'.$$
On effectue le changement de variable $x\mapsto qxq^{-1}$: cette conjugaison pr\'eserve \`a la fois $U(F)$ et $U(F)\cap Q(F)$. Les termes $\bar{\xi}(x)$ et $\varphi_{c}(u,x,x')$ ne d\'ependent de $x$ que par l'int\'erm\'ediaire des coefficients $q_{V}(xv_{i},v_{-i-1})$pour $i=1,...,r-1$. Puisque $q$ fixe  les vecteurs $v_{i}$, la conjugaison par $q$ ne change pas ces termes. Par contre, elle introduit un module. Pour l'exprimer commod\'ement et pour poursuivre notre calcul, on d\'ecompose $q$ en $\delta n \tilde{g}$, o\`u $\delta\in GL_{k-r-1}(F)$, $n\in U_{Q^H}(F)$ et $\tilde{g}\in \tilde{G}(F)$. Le module est alors $\vert det(\delta)\vert _{F}^{-r}$.   On obtient
$$ {\cal L}_{\pi,\rho}(\epsilon'\otimes e',\epsilon\otimes e)=\int_{((Q_{0}(F)\cap H(F))\backslash H(F))^2}\int_{((U(F)\cap Q(F))\backslash U(F))^2}  \int_{\tilde{G}(F)}\int_{U_{Q^H}(F)}\int_{GL_{k-r-1}(F)}$$
$$\int_{U(F)\cap GL_{k}(F)}(\rho(\delta n \tilde{g}h)\epsilon',\rho(h')\epsilon)(e'(x'h'),e(u\delta n\tilde{g}xh))\varphi_{c}(u,x,x')$$
$$\bar{\xi}(x^{_{'}-1}ux)\vert det(\delta)\vert _{F}^{-r}du\,d\delta\,dn\,d\tilde{g}\,dx\,dx'\,dh\,dh'.$$
  On a l'\'egalit\'e $U(F)=(U(F)\cap Q(F))\times (U(F)\cap U_{\bar{Q}}(F)$ qui permet de remplacer l' int\'egration sur $(U(F)\cap Q(F))\backslash U(F)$ par une int\'egration sur $U(F)\cap U_{\bar{Q}}(F)$. On va l\'eg\`erement modifier cet ensemble de repr\'esentants. Soit $x\in U(F)\cap U_{\bar{Q}}(F)$. Pour $i=1,...,r-1$, on a 
 $$q_{V}(xv_{i},v_{-i-1})=q_{V}(xy_{-k+r-i},y_{k-r+i+1})=q_{V}(y_{-k+r-i},y_{k-r+i+1})=0$$
puisque $x $ fixe $y_{-k+r-i}$. Par contre, $q_{V}(xv_{0},v_{-1})$ n'est en g\'en\'eral pas nul. Exprimons matriciellement les \'el\'ements de $GL_{k}$ dans la base $(y_{j})_{j=1,...,k}$ de $Y^+$. Remarquons que le groupe $U\cap GL_{k}$ est le radical unipotent du sous-groupe parabolique de $GL_{k}$ triangulaire sup\'erieur par blocs, de blocs $k-r$, $1$,...,$1$. Pour $u\in U(F)\cap GL_{k}(F)$, on calcule ais\'ement
$$q_{V}(uv_{i},v_{-i-1})=-u_{k-r+i,k-r+i+1}$$
pour $i=0,...,r-1$. Pour $x\in U(F)\cap U_{\bar{Q}}(F)$, notons $u(x)$ l'\'el\'ement de $U(F)\cap GL_{k}(F)$ dont toutes les coordonn\'ees non diagonales sont nulles, sauf $u(x)_{k-r,k-r+1}$ qui vaut $q_{V}(xv_{0},v_{-1})$. Posons $x_{*}=u(x)x$. Alors $q_{V}(x_{*}v_{i},v_{-i-1})=0$ pour tout $i=0,...,r-1$ et $\{x_{*}; x\in U(F)\cap U_{\bar{Q}}(F)\}$ est encore un ensemble de repr\'esentants de $(U(F)\cap Q(F))\backslash U(F)$.   Pour $x,x'\in U(F)\cap U_{\bar{Q}}(F)$ et $u\in U(F)\cap GL_{k}(F)$, on a $\bar{\xi}(x_{*}^{_{'}-1}ux_{*})=\bar{\xi}(u)$ et $\varphi_{c}(u,x_{*},x'_{*})=1$ si et seulement si $u\in U_{k}(F)_{c}$. On obtient
$$(6) \qquad {\cal L}_{\pi,\rho}(\epsilon'\otimes e',\epsilon\otimes e)=\int_{((Q_{0}(F)\cap H(F))\backslash H(F))^2}\int_{(U(F)\cap U_{\bar{Q}}(F))^2} $$
$$I(x',h',x,h)dx\,dx'\,dh\,dh',$$
o\`u 
$$I(x',h',x,h)=\int_{\tilde{G}(F)}\int_{U_{Q^H}(F)}\int_{GL_{k-r-1}(F)}\int_{U(F)\cap U_{k}(F)_{c}}(\rho(\delta n \tilde{g}h)\epsilon',\rho(h')\epsilon)$$
$$(e'(x_{*}'h'),e(u\delta n\tilde{g}x_{*}h))\bar{\xi}(u)\vert det(\delta)\vert _{F}^{-r}du\,d\delta\,dn\,d\tilde{g}.$$
Toutes ces expressions sont  absolument convergentes d'apr\`es (1): on n'a jusqu'ici effectu\'e que des changements de variables et des permutations d'int\'egrales.

On va calculer $I(x',h',x,h)$. Fixons un sous-groupe compact sp\'ecial $K^H$de $H(F)$ en bonne position relativement au sous-groupe parabolique $Q^H$. L'application naturelle de $K^H$ dans $Q^H(F)\backslash H(F)$ est surjective. D'apr\`es la description que l'on a donn\'ee ci-dessus du groupe $Q_{0}\cap H$,  tout \'el\'ement de$(Q_{0}(F)\cap H(F))\backslash H(F)$ a un repr\'esentant qui appartient \`a $H_{0}(F)K^H$. On peut donc se limiter \`a calculer $I(x',h',x,h)$ pour des \'el\'ements $x,x'\in U(F)\cap U_{\bar{Q}}(F)$ et $h,h'\in H_{0}(F)K^H$. Dans l'expression de $I(x',h',x,h)$, on peut remplacer $e(u\delta n\tilde{g}x_{*}h)$ par $\mu(u)e(\delta n\tilde{g}x_{*}h)$. D'apr\`es les formules \'ecrites ci-dessus, $\bar{\xi}(u)=\psi(-\sum_{j=k-r,...,k-1}u_{j,j+1})$. Fixons une fonctionnelle de Whittaker $\phi$ non nulle sur $E_{\mu}$ et notons comme plus haut $\Phi:E_{\mu}\otimes_{{\mathbb C}}E_{\tilde{\pi}}\to E_{\tilde{\pi}}$ l'application $\phi\otimes id$. On peut appliquer le lemme 3.7(ii):  il existe une constante $C\not=0$ telle que 
 $$\int_{U(F)\cap U_{k}(F)_{c}}(e'(x_{*}'h'),\mu(u)e(\delta n\tilde{g}x_{*}h))\bar{\xi}(u)du= $$
 $$C\int_{U_{k-r-1}(F)\backslash GL_{k-r-1}(F)\times \omega_{[k-r]}(c+c_{\psi})}(\Phi \mu(\gamma a)e'(x_{*}'h'),\Phi\mu(\gamma a)e(\delta n\tilde{g}x_{*}h))\vert det(\gamma)\vert _{F}^{-r}da\,d\gamma,$$
 pourvu que $c+c_{\psi}\geq1$. On peut remplacer $ \mu(\gamma a)e(\delta n \tilde{g}x_{*}h)$ par 
 $$\mu(\gamma \delta a)e(n \tilde{g}x_{*}h)\delta_{Q}(\delta)^{1/2}=\mu(\gamma \delta a)e(n \tilde{g}x_{*}h)\vert  det(\delta)\vert _{F}^{(d_{\tilde{V}}+k-1)/2}.$$
 On obtient
 $$I(x',h',x,h)=C \int_{\tilde{G}(F)}\int_{U_{Q^H}(F)}\int_{GL_{k-r-1}(F)}$$
 $$\int_{U_{k-r-1}(F)\backslash GL_{k-r-1}(F)\times \omega_{[k-r]}(c+c_{\psi})}(\rho(\delta n \tilde{g}h)\epsilon',\rho(h')\epsilon)(\Phi \mu(\gamma a)e'(x_{*}'h'),\Phi\mu(\gamma \delta a)e( n\tilde{g}x_{*}h))$$
 $$\vert det(\gamma)\vert _{F}^{-r}\vert det(\delta)\vert _{F}^{-r+(d_{\tilde{V}}+k-1)/2}da\,d\gamma\,d\delta\,dn\,d\tilde{g}.$$
 Montrons que
 
 (7) pour $\tilde{g}$ et $n$ fix\'es, l'int\'egrale int\'erieure sur $GL_{k-r-1}(F)\times( U_{k-r-1}(F)\backslash GL_{k-r-1}(F))\times \omega_{[k-r]}(c+c_{\psi})$ est absolument convergente.
 
 La variable $a\in \omega_{[k-r]}(c+c_{\psi})$ dispara\^{\i}t tout de suite: la fonction que l'on int\`egre est localement constante en cette variable et le domaine d'int\'egration est compact. On a une majoration
 $$\vert (\rho(\delta n \tilde{g}h)\epsilon',\rho(h')\epsilon)\vert <<\Xi^H(\delta).$$
 On peut d\'ecomposer $e'(x_{*}'h')$ et $e(n\tilde{g}x_{*}h)$ en combinaisons lin\'eaires de produits $\eta\otimes \tilde{e}$, o\`u $\eta\in E_{\mu}$ et $\tilde{e}\in E_{\tilde{\pi}}$. Cela nous ram\`ene \`a montrer que, pour $\eta,\eta'\in E_{\mu}$, l'int\'egrale
 $$\int_{GL_{k-r-1}(F)}\int_{U_{k-r-1}(F)\backslash GL_{k-r-1}(F)}\Xi^H(\delta)\vert \phi\mu(\gamma)\eta'\vert \vert \phi\mu(\gamma \delta)\eta\vert$$
 $$ \vert det(\gamma)\vert _{F}^{-r}\vert det(\delta)\vert _{F}^{-r+(d_{\tilde{V}}+k-1)/2}d\gamma\,d\delta$$
 est convergente.  On effectue le changement de variable $\delta\mapsto \gamma^{-1}\delta$. On remplace ensuite la variable $\gamma$ par $t'k'$, avec $t'\in A_{k-r-1}(F)$ et $k'\in K_{k-r-1}$ et $\delta$ par $tuk$, avec $t\in A_{k-r-1}(F)$, $u\in U_{k-r-1}(F)$, $k\in K_{k-r-1}$.  Cela remplace $d\gamma\, d\delta$ par $\delta_{B_{k-r-1}}(t')^{-1}dt'\,dk'\,dt\,du\,dk$.   De nouveau, les int\'egrales en $k$ et $k'$ sont inessentielles et on est ramen\'e \`a l'int\'egrale
 $$(8) \qquad\int_{U_{k-r-1}(F)}\int_{A_{k-r-1}(F)^2}\Xi^H(t^{_{'}-1}tu)\vert \phi\mu(t')\eta'\vert \vert \phi\mu(t)\eta\vert \delta_{B_{k-r-1}}(t')^{-1}$$
 $$\vert det(t')\vert _{F}^{(1-d_{\tilde{V}}-k)/2}\vert det(t)\vert _{F}^{-r+(d_{\tilde{V}}+k-1)/2}dt\,dt'\,du.$$
   Montrons que
 
 (9) pour tous r\'eels $R>0$ et $\epsilon$ avec $0<\epsilon<1/2$, on a une majoration
 $$\Xi^H(g)<<\Xi^{GL_{k-r-1}}(g)\sigma(g)^{-R}\vert det(g)\vert _{F}^{\epsilon+(r+1-d_{\tilde{V}}-k)/2}$$
 pour tout $g\in GL_{k-r-1}(F)$.
 
 On peut supposer $g=a\in A_{k-r-1}(F)$. Notons $a_{j}$,  pour $j=1,...,k-r-1$, les coefficients diagonaux de $a$. Choisissons un sous-groupe de L\'evi minimal $M_{min}$ de $H$ contenant $A_{k-r-1}$ et un sous-groupe parabolique minimal $P_{min}\in {\cal P}(M_{min})$  tel que $a$ soit "n\'egatif " pour $P_{min}$, c'est-\`a-dire que $\vert \alpha(a)\vert \leq 1$ pour toute racine $\alpha$ de $A_{M_{min}}$ dans $\mathfrak{u}_{P_{min}}$.  D'apr\`es [W2] lemme II.1.1, on a des in\'egalit\'es
 $$(10) \qquad \delta_{P_{min}}(a)<<\Xi^H(a)^2<<\delta_{P_{min}}(a)\sigma(a)^D$$
 o\`u $D$ est un certain entier. On \'enum\`ere les valeurs de $\alpha(a)$ pour toutes les racines $\alpha$ de $A_{M_{min}}$ dans $\mathfrak{h}$: ce sont $a_{j}a_{j'}^{-1}$ pour $j\not=j'$, $a_{j}a_{j'}$ et $(a_{j}a_{j'})^{-1}$ pour $j<j'$, qui interviennent avec multiplicit\'e $1$, et $a_{j}$ et $a_{j}^{-1}$ qui interviennent avec multiplicit\'e $d_{W_{0}}=d_{\tilde{V}}+1$ (\'evidemment, les $j$, $j'$ parcourent $\{1,...,k-r-1\}$). Le module $\delta_{P_{min}}(a)$ est le produit  de celles des valeurs absolues de ces termes qui sont inf\'erieures ou \'egales \`a $1$. Donc
 $$\delta_{P_{min}}(a)=I_{1}I_{2}I_{3},$$
 o\`u
 $$I_{1}=\prod_{j\not=j'; val_{F}(a_{j})\geq val_{F}(a_{j'})}\vert a_{j}a_{j'}^{-1}\vert _{F};$$
 $$I_{2}=(\prod_{j<j'; val_{F}(a_{j}a_{j'})\geq0}\vert a_{j}a_{j'}\vert _{F})(\prod_{j<j'; val_{F}(a_{j}a_{j'})<0}\vert a_{j}a_{j'}\vert _{F}^{-1});$$
 $$I_{3}=(\prod_{j; val_{F}(a_{j})\geq 0}\vert a_{j}\vert ^{d_{\tilde{V}}+1})(\prod_{j; val_{F}(a_{j})\leq 0}\vert a_{j}\vert _{F}^{-d_{\tilde{V}}-1}).$$
 On peut majorer le premier produit de $I_{2}$ par son inverse. On obtient
 $$I_{2}\leq \prod_{j<j'}\vert a_{j}a_{j'}\vert_{F} ^{-1}=\vert det(a)\vert _{F}^{r+2-k}.$$
 On a $I_{3}=I_{4}I_{5}$, o\`u
 $$I_{4}=(\prod_{j; val_{F}(a_{j})\geq 0}\vert a_{j}\vert_{F} ^{2\epsilon})(\prod_{j; val_{F}(a_{j})\leq 0}\vert a_{j}\vert _{F}^{-2\epsilon}),$$
 $$I_{5}=(\prod_{j; val_{F}(a_{j})\geq 0}\vert a_{j}\vert_{F} ^{d_{\tilde{V}}+1-2\epsilon})(\prod_{j; val_{F}(a_{j})\leq 0}\vert a_{j}\vert _{F}^{-d_{\tilde{V}}-1+2\epsilon}).$$
 Comme ci-dessus, le premier produit de $I_{5}$ est major\'e par son inverse, d'o\`u
 $$I_{5}\leq\prod_{j}\vert a_{j}\vert _{F}^{-d_{\tilde{V}}-1+2\epsilon}=\vert det(a)\vert _{F}^{-d_{\tilde{V}}-1+2\epsilon}.$$
 On a $I_{4}=q^{-2\epsilon b}$ o\`u
 $$b=\sum_{j}\vert val_{F}(a_{j})\vert >> \sigma(a),$$
 donc $I_{4}<<q^{-2\epsilon \sigma(a)}$. D'o\`u, en utilisant la seconde majoration de (9),
 $$\Xi^H(a)^2<<\sigma(a)^Dq^{-2\epsilon\sigma(a)}\vert det(a)\vert _{F}^{r+1-d_{\tilde{V}}-k+2\epsilon}I_{1}.$$
 Un calcul similaire vaut en rempla\c{c}ant le groupe $H$ par $GL_{k-r-1}$. En utilisant cette fois la premi\`ere majoration de (10), on obtient simplement
 $$I_{1}<< \Xi^{GL_{k-r-1}}(a)^2.$$
 De ces deux majorations se d\'eduit l'assertion (9). 
 
 De (9) et de la proposition II.4.5 de [W2] se d\'eduit que l'int\'egrale
 $$\int_{U_{k-r-1}(F)}\Xi^H(t^{_{'}-1}tu)du$$
 est convergente et que, pour tout  $\epsilon$ tel que $0<\epsilon<1/2$,elle est born\'ee par 
 $$ \delta_{B_{k-r-1}}(t^{_{'}-1}t)^{-1/2}\vert det(t^{_{'}-1}t)\vert _{F}^{\epsilon+(r+1-d_{\tilde{V}}-k)/2}.$$
Pour tout entier $c'\in {\mathbb N}$,  on a introduit en 3.5 la fonction $\iota_{c'}$  sur $A_{k}(F)$.  D'apr\`es le lemme 3.7(i), il existe  un entier $c'\in {\mathbb N}$ et un r\'eel $R\geq 0$ tels  que l'on ait la majoration
$$\vert \phi\mu(a)\eta\vert<<\iota_{c'}(a)\delta_{B_{k}}(a)^{1/2}\sigma(a)^{R}$$
pour tout $a\in A_{k}(F)$. Pour $t\in A_{k-r-1}(F)$, on a $\delta_{B_{k}}(t)= \delta_{B_{k-r-1}}(t)\vert det(t)\vert _{F}^{r+1}$, d'o\`u
$$(11) \qquad \vert \phi\mu(t)\eta\vert<<\iota_{c'}(t)\delta_{B_{k-r-1}}(t)^{1/2}\vert det(t)\vert _{F}^{(r+1)/2}\sigma(t)^{R}.$$
 La fonction $g\mapsto \phi\mu(g)\eta'$ v\'erifie une majoration analogue. L'expression (8) est donc major\'ee par le produit des int\'egrales
$$\int_{A_{k-r-1}(F)}\iota_{c'}(t') \vert det(t')\vert _{F}^{-\epsilon+1/2} \sigma(t')^{R} dt'$$
et
$$\int_{A_{k-r-1}(F)}\iota_{c'}(t)\vert det(t)\vert _{F}^{\epsilon+1/2} \sigma(t)^R dt.$$ 
Consid\'erons par exemple la premi\`ere. On remplace les variables $t'_{1},...,t'_{k-r-1}$ par $a_{1}....a_{k-r-1}$,
 $a_{2}...a_{k-r-1},...,a_{k-r-1}$. L'int\'egrale est alors essentiellement born\'ee par le produit sur $j=1,...,k-r-1$ des int\'egrales
$$\int_{a_{j}\in F^{\times}; val_{F}(a_{j})\geq -c'}\vert a_{j}\vert _{F}^{j(-\epsilon+1/2)}(1+\vert val_{F}(a_{j})\vert)^R da_{j}.$$ 
Chacune d'elles est convergente, ce qui ach\`eve la preuve de (7).

On utilise (7) pour permuter les deux int\'egrales int\'erieures dans l'expression qui pr\'ec\`ede cette assertion. On effectue ensuite le changement de variables $\delta\mapsto \gamma^{-1}\delta$, puis on d\'ecompose $\gamma$ en $t'k'$ et $\delta$ en $utk$, avec $u\in U_{k-r-1}(F)$, $t,t'\in A_{k-r-1}(F)$ et $k,k'\in K_{k-r-1}$. Cela change $d\gamma\,d\delta$ en $\delta_{B_{k-r-1}}(tt')^{-1}du\,dt\,dt'\,dk\,dk'$. On obtient
$$I(x',h',x,h)=C \int_{\tilde{G}(F)}I(x',h',x,h,\tilde{g})d\tilde{g}$$
o\`u
$$I(x',h',x,h,\tilde{g})=\int_{U_{Q^H}(F)}\int_{\omega_{[k-r]}(c+c_{\psi})}\int_{K_{k-r-1}^2}\int_{A_{k-r-1}(F)^2}
\int_{U_{k-r-1}(F)}$$
$$(\rho(utkn\tilde{g}h)\epsilon',\rho(t'k'h')\epsilon)(\Phi\mu(t'k'a)e'(x'_{*}h'),\Phi\mu(utka)e(n\tilde{g}x_{*}h))$$
$$\vert det(t')\vert _{F}^{(1-d_{\tilde{V}}-k)/2}\vert det(t)\vert _{F}^{-r+(d_{\tilde{V}}+k-1)/2}\delta_{B_{k-r-1}}(tt')^{-1}du \,dt\,dt'\,dk\,dk'\,da\,dn.$$
Fixons une suite $(\Omega_{l})_{l\in {\mathbb N}}$ de sous-groupes ouverts compacts de $U_{Q^H}(F)$ telle que $\Omega_{l}\subset \Omega_{l+1}$ pour tout $l$, 
$$\bigcup_{l\in {\mathbb N}}\Omega_{l}=U_{Q^H}(F)$$
et tout $\Omega_{l}$ soit invariant par conjugaison par $K_{k-r-1}$. Notons $I_{l}(x',h',x,h,\tilde{g})$ l'expression obtenue en rempla\c{c}ant dans $I(x',h',x,h,\tilde{g})$ la premi\`ere int\'egrale sur $U_{Q^H}(F)$ par l'int\'egrale sur $\Omega_{l}$. On a
$$I(x',h',x,h,\tilde{g})=lim_{l\to \infty}I_{l}(x',h',x,h,\tilde{g}).$$
Fixons $l\in {\mathbb N}$. En (7), on avait fix\'e $n$. Mais, les fonctions consid\'er\'ees \'etant localement constantes en cette variable, on aurait aussi bien pu l'autoriser \`a varier dans un sous-ensemble compact de $U_{Q^H}(F)$ et on aurait obtenu une majoration uniforme en $n$ de l'int\'egrale consid\'er\'ee. Donc l'expression $I_{l}(x',h',x,h,\tilde{g})$ est absolument convergente. On peut changer l'ordre des int\'egrales puis effectuer le changement de variable $n\mapsto (tk)^{-1}ntk$. Cela introduit le module $\delta_{Q^H}(t)^{-1}=\vert det(t)\vert _{F}^{r+1-k-d_{\tilde{V}}}$.  On obtient
$$I_{l}(x',h',x,h,\tilde{g})=\int_{\omega_{[k-r]}(c+c_{\psi})}\int_{K_{k-r-1}^2}\int_{A_{k-r-1}(F)^2}\int_{t\Omega_{l}t^{-1}}
\int_{U_{k-r-1}(F)}$$
$$(\rho(untk\tilde{g}h)\epsilon',\rho(t'k'h')\epsilon)(\Phi\mu(t'k'a)e'(x'_{*}h'),\Phi\mu(utka)e((tk)^{-1}ntk\tilde{g}x_{*}h))$$
$$\vert det(tt')\vert _{F}^{(1-d_{\tilde{V}}-k)/2}\delta_{B_{k-r-1}}(tt')^{-1}du\,dn \,dt\,dt'\,dk\,dk'\,da.$$
 Soit $n\in U_{Q^H}(F)$. On a $nw_{0}\in w_{0}+Y^+_{W}$. Introduisons l'\'el\'ement $u(n)\in U_{k}(F)$ qui fixe tous les vecteurs $y_{j}$ pour $j=1,...,k$, $j\not=k-r$, et envoie $y_{k-r}$ sur $y_{k-r}+ nw_{0}-w_{0}$. On a 
 
 (12) $u(n)^{-1}n\in U_{Q_{0}}(F)$. 
 
 En effet, $u(n)^{-1}n$ fixe tout $y_{j}$ pour $j=1,...,k-r-1$ et envoie $\tilde{V}$ dans $Y^+_{W}$, a fortiori dans $Y^+_{0}$. Pour qu'il appartienne \`a $U_{Q_{0}}(F)$, il suffit qu'il fixe de plus $y_{k-r}$. On a $y_{k-r}=w_{0}+\frac{1}{2\nu_{0}}v_{0}$, donc $ny_{k-r}=nw_{0}+\frac{1}{2\nu_{0}}v_{0}=y_{k-r}+nw_{0}-w_{0}$, puis $u(n)^{-1}ny_{k-r}=u(n)^{-1}y_{k-r}+nw_{0}-w_{0}=y_{k-r}$. Cela prouve (12).
 
 Puisque $e$ est invariante \`a gauche par $U_{Q}(F)$, a fortiori par $U_{Q_{0}}(F)$, on a
 $$\mu(utka)e((tk)^{-1}ntk\tilde{g}x_{*}h)=\delta_{Q}(tk)^{-1/2}\mu(uau(n))e(u(n)^{-1}ntk\tilde{g}x_{*}h)$$
 $$\qquad =\delta_{Q}(tk)^{-1/2}\mu(uau(n))e(tk\tilde{g}x_{*}h)=\mu(uau(n)tk)e(\tilde{g}x_{*}h).$$
 Introduisons le sous-groupe parabolique $P^H$ de $H$, contenu dans $Q^H$, dont l'intersection avec $L^H$ est $B_{k-r-1}\times H_{0}$. Son radical unipotent est $U_{P^H}=U_{k-r-1}\times U_{Q^H}$. On d\'efinit un caract\`ere $\bar{\xi}^H$ de $U_{P^H}(F)$ par la formule
 $$\bar{\xi}^H(\nu)=\psi(q_{W}(\nu w_{0},y_{r-k+1})+\sum_{j=2,...,k-r-1}q_{W}(\nu y_{j},y_{-j+1}))$$
 pour $\nu\in U_{P^H}(F)$. Remarquons que les donn\'ees $H$, $P^H$, $w_{0}$, $\bar{\xi}^H$, $\tilde{G}$ sont d'exactes similaires des donn\'ees $G$, $P$, $v_{0}$, $\bar{\xi}$, $H$. Plus g\'en\'eralement, pour $a\in \omega_{[k-r]}(c+c_{\psi})$, d\'efinissons un caract\`ere $\bar{\xi}^H_{a}$ de $U_{P^H}(F)$ par
 $$\bar{\xi}_{a}^H(\nu)=\psi(a_{k-r}^{-1}q_{W}(\nu w_{0},y_{r-k+1})+\sum_{j=2,...,k-r-1}q_{W}(\nu y_{j},y_{-j+1})).$$
Si $\nu=un$, avec $u\in U_{k-r-1}(F)$ et $n\in U_{Q^H}(F)$, on a l'\'egalit\'e
$$\bar{\xi}^H_{a}(un)=\psi(a_{k-r}^{-1}u(n)_{k-r-1,k-r}+\sum_{j=1,...k-r-2}u_{j,j+1}).$$ 
 D'apr\`es la d\'efinition de $\Phi$, on en d\'eduit
 $$\Phi\mu(uau(n)tk)e(\tilde{g}x_{*}h)=\bar{\xi}_{a}^H(un)\Phi\mu(tka)e(\tilde{g}x_{*}h)=\bar{\xi}^H_{a}(un)\tilde{\pi}(\tilde{g})\Phi\mu(tka)e(x_{*}h).$$
 puis
 $$(13) \qquad I_{l}(x',h',x,h,\tilde{g})=\int_{\omega_{[k-r]}(c+c_{\psi})}\int_{K_{k-r-1}^2}\int_{A_{k-r-1}(F)^2}\int_{t\Omega_{l}t^{-1}}
\int_{U_{k-r-1}(F)}$$
$$(\rho(untk\tilde{g}h)\epsilon',\rho(t'k'h')\epsilon)(\Phi\mu(t'k'a)e'(x'_{*}h'),\tilde{\pi}(\tilde{g})\Phi\mu(tka)e(x_{*}h))\bar{\xi}^H_{a}(un)$$
$$\vert det(tt')\vert _{F}^{(1-d_{\tilde{V}}-k)/2}\delta_{B_{k-r-1}}(tt')^{-1}du\,dn \,dt\,dt'\,dk\,dk'\,da.$$
Pour un entier $c'\geq1$, on introduit le sous-groupe $\omega_{A_{k-r-1}}(c')\subset A_{k-r-1}(F)$. Montrons que

(14) il existe un entier $c'\geq1$ tel que, pour tous $t,t'\in A_{k-r-1}(F)$, $k,k'\in K_{k-r-1}$, $h,h'\in H_{0}(F)K^H$ et $\tilde{g}\in \tilde{G}(F)$, les \'el\'ements $\rho(tk\tilde{g}h)\epsilon'$ et $\rho(t'k'h')\epsilon$
de $E_{\rho}$ soient invariants par $\rho(\alpha)$ pour tout $\alpha\in\omega_{A_{k-r-1}}(c')$.

Ecrivons $h=h_{0}\kappa$, avec $h_{0}\in H_{0}(F)$ et $\kappa\in K^H$. Soit $\alpha\in A_{k-r-1}(F)$. En utilisant le fait que $G_{k-r-1}$ commute \`a $H_{0}$, a fortiori \`a $\tilde{G}$, on a
$$\alpha tk\tilde{g}h=tk\tilde{g}h\gamma$$
o\`u $\gamma=(k\kappa)^{-1}\alpha k\kappa$. Puisque $k\kappa$ reste dans un ensemble compact, il existe $c'$ tel que la condition $\alpha\in \omega_{A_{k-r-1}}(c')$ entra\^{\i}ne que $\rho(\gamma)$  fixe $\epsilon'$, donc que $\rho(\alpha)$ fixe $\rho(tk\tilde{g}h)\epsilon'$. La preuve est la m\^eme pour l'\'el\'ement $\rho(t'k'h')\epsilon$. Cela prouve (14).

On fixe $c'$ comme en (14). Les deux int\'egrales int\'erieures de la formule (13) se regroupent en une int\'egrale sur le sous-ensemble $U_{k-r-1}(F)\times t\Omega_{l}t^{-1}$ de $U_{P^H}(F)$. Introduisons le sous-groupe $U_{P^H}(F)_{c'}$ de ce dernier groupe. Posons $\Delta_{l}(t)=U_{P^H}(F)_{c'}\cap (U_{k-r-1}(F)\times t\Omega_{l}t^{-1})$.  La preuve du lemme 3.5 montre que l'int\'egrale sur le compl\'ementaire de $ \Delta_{l}(t)$ dans $U_{k-r-1}(F)\times t\Omega_{l}t^{-1}$ est nulle. On peut donc remplacer les deux int\'egrales int\'erieures de (13) par l'int\'egrale sur $ \Delta_{l}(t)$. Remarquons que $c'$ est ind\'ependant de $c$, on peut donc supposer que $c$ est assez grand relativement \`a $c'$. Alors, pour $a\in \omega_{[k-r]}(c+c_{\psi})$, le caract\`ere $\bar{\xi}^H_{a}$ co\"{\i}ncide avec $\bar{\xi}^H$ sur $U_{P^H}(F)_{c'}$. On obtient
$$I_{l}(x',h',x,h)=\int_{\tilde{G}(F)}\int_{\omega_{[k-r]}(c+c_{\psi})}\int_{K_{k-r-1}^2}\int_{A_{k-r-1}(F)^2}\int_{\Delta_{l}(t)}$$
$$(\rho(utk\tilde{g}h)\epsilon',\rho(t'k'h')\epsilon)(\Phi\mu(t'k'a)e'(x'_{*}h'),\tilde{\pi}(\tilde{g})\Phi\mu(tka)e(x_{*}h))\bar{\xi}^H(u)$$
$$\vert det(tt')\vert _{F}^{(1-d_{\tilde{V}}-k)/2}\delta_{B_{k-r-1}}(tt')^{-1}du\,dt\,dt'\,dk\,dk'\,da\,d\tilde{g}.$$
Montrons que

(15) l'expression obtenue en rempla\c{c}ant ci-dessus l'int\'egrale sur $\Delta_{l}(t)$ par l'int\'egrale sur $U_{P^H}(F)_{c'}$ est absolument convergente.

Comme toujours, les int\'egrales sur $\omega_{[k-r]}(c+c_{\psi})\times K_{k-r-1}^2$ sont inessentielles. Oublions-les. On peut supposer $e(x_{*}h)=\eta\otimes \tilde{e}$ et $e'(x'_{*}h')=\eta'\otimes \tilde{e}'$, avec $\eta,\eta'\in E_{\mu}$ et $\tilde{e}, \tilde{e}'\in E_{\tilde{\pi}}$. On a les majorations
$$\vert (\rho(ut\tilde{g}h)\epsilon',\rho(t'h')\epsilon)\vert<<\Xi^H(t^{_{'}-1}ut\tilde{g}),$$
$$\vert  (\Phi\mu(t'k'a)e'(x'_{*}h'),\tilde{\pi}(\tilde{g})\Phi\mu(tka)e(x_{*}h))\vert <<\Xi^{\tilde{G}}(\tilde{g})\vert \phi\mu(t')\eta'\vert \vert \phi\mu(t)\eta\vert .$$
On peut encore majorer cette derni\`ere expression gr\^ace \`a (11): il existe un entier $c''$ et un r\'eel $R\geq0$ tels que
$$\vert  (\Phi\mu(t'k'a)e'(x'_{*}h'),\tilde{\pi}(\tilde{g})\Phi\mu(tka)e(x_{*}h))\vert <<\Xi^{\tilde{G}}(\tilde{g})\iota_{c''}(t)\iota_{c''}(t')$$
$$\delta_{B_{k-r-1}}(tt')^{1/2}\vert det(tt')\vert _{F}^{(r+1)/2}\sigma(t)^{R}\sigma(t')^R.$$

On effectue le changement de variable $u\mapsto t'ut^{_{'}-1}$, qui introduit le module $\delta_{P^H}(t') $, on est ramen\'e \`a l'expression
$$\int_{\tilde{G}(F)}\int_{A_{k-r-1}(F)^2}\int_{t^{_{'}-1}U_{P^H}(F)_{c'}t'}\Xi^H(ut^{_{'}-1}t\tilde{g})\Xi^{\tilde{G}}(\tilde{g})\iota_{c''}(t)\iota_{c''}(t')$$
$$\delta_{B_{k-r-1}}(tt')^{-1/2}\delta_{P^H}(t')\vert det(tt')\vert _{F}^{(1+(r-d_{\tilde{V}}-k)/2} \sigma(t)^R\sigma(t')^Rdu\,dt\,dt'\,d\tilde{g}.$$
 Pour $t'\in A_{k-r-1}(F)$, dont on \'ecrit les coefficients diagonaux $t'_{1},...,t'_{k-r-1}$, posons
 $$c'(t')=sup\{c', c'+val_{F}(t_{1})-val_{F}(t_{2}),...,c'+val_{F}(t_{k-r-2})-val_{F}(t_{k-r-1}), c'+val_{F}(t_{k-r-1})\}.$$
 Alors $t^{_{'}-1}U_{P^H}(F)_{c'}t'\subset U_{P^H}(F)_{c'(t')}$ et on peut remplacer dans l'expression ci-dessus l'int\'egrale sur $t^{_{'}-1}U_{P^H}(F)_{c'}t'$ par l'int\'egrale sur $U_{P^H}(F)_{c'(t')}$. D'apr\`es 4.3(3), il existe un r\'eel $R'\geq0$ tel que
$$ \int_{U_{P^H}(F)_{c'(t')}}\Xi^H(ut^{_{'}-1}t\tilde{g})du<<\delta_{P^H}(tt^{_{'}-1})^{1/2}\Xi^{H_{0}}(\tilde{g})\sigma(\tilde{g})^{R'}\sigma(tt^{_{'}-1})^{R'}c'(t')^{R'}.$$
On peut aussi bien remplacer $\sigma(tt^{_{'}-1})^{R'}c'(t')^{R'}$ par $\sigma(t)^{R'}\sigma(t')^{R'}$. On a $$\delta_{P^H}(tt')=\delta_{B_{k-r-1}}(tt')\vert det(tt')\vert _{F}^{d_{\tilde{V}}+k-r-1}.$$
 Alors  notre expression est major\'ee par le produit de trois int\'egrales. La premi\`ere est
$$\int_{\tilde{G}(F)}\Xi^{H_{0}}(\tilde{g})\Xi^{\tilde{G}}(\tilde{g})\sigma(\tilde{g})^{R'}d\tilde{g},$$
qui est convergente d'apr\`es 4.3(4). Les deux autres sont toutes deux \'egales \`a
$$\int_{A_{k-r-1}(F)}\iota_{c''}(t)\vert det(t)\vert _{F}^{1/2}\sigma(t)dt.$$
On a vu dans la preuve de (7) que cette int\'egrale \'etait convergente. Cela prouve (15).

Gr\^ace \`a (15), le th\'eor\`eme de convergence domin\'ee nous permet de calculer la limite de $I_{l}(x',h',x,h)$ quand $l$ tend vers l'infini, c'est-\`a-dire $I(x',h',x,h)$, comme l'int\'egrale obtenue en rempla\c{c}ant, dans la formule qui pr\'ec\`ede (15), l'int\'egrale sur $\Delta_{l}(t)$ par l'int\'egrale sur la limite de cet ensemble, c'est-\`a-dire sur $U_{P^H}(F)_{c'}$ tout entier. On reconna\^{\i}t la double int\'egrale sur $\tilde{G}(F)\times U_{P^H}(F)_{c'}$: elle donne 
$${\cal L}_{\rho,\tilde{\pi},c'}(\rho(tkh)\epsilon'\otimes \Phi\mu(t'k'a)e'(x'_{*}h'),\rho(tk'h')\epsilon\otimes \Phi\mu(tka)e(x_{*}h)),$$
 ce terme \'etant calcul\'e \`a l'aide du caract\`ere $\xi^H$. D'apr\`es le choix de $c'$ et le lemme 5.1, c'est aussi
 $${\cal L}_{\rho,\tilde{\pi}}(\rho(tkh)\epsilon'\otimes \Phi\mu(t'k'a)e'(x'_{*}h'),\rho(tk'h')\epsilon\otimes \Phi\mu(tka)e(x_{*}h)).$$
 D'o\`u
 $$I(x',h',x,h)=\int_{\omega_{[k-r]}(c+c_{\psi})}\int_{K_{k-r-1}^2}\int_{A_{k-r-1}(F)^2}\vert det(tt')\vert _{F}^{(1-d_{\tilde{V}}-k)/2}\delta_{B_{k-r-1}}(tt')^{-1}$$
 $${\cal L}_{\rho,\tilde{\pi}}(\rho(tkh)\epsilon'\otimes \Phi\mu(t'k'a)e'(x'_{*}h'),\rho(t'k'h')\epsilon\otimes \Phi\mu(tka)e(x_{*}h))dt\,dt'\,dk\,dk'\,da.$$
 Remarquons que les int\'egrales sur $K_{k-r-1}$ et $A_{k-r-1}(F)$ se regroupent en int\'egrales sur $U_{k-r-1}(F)\backslash GL_{k-r-1}(F)$, d'o\`u
 $$(16) \qquad I(x',h',x,h)=\int_{\omega_{[k-r]}(c+c_{\psi})}\int_{(U_{k-r-1}(F)\backslash GL_{k-r-1}(F))^2}\vert det(\gamma\gamma')\vert _{F}^{(1-d_{\tilde{V}}-k)/2}$$
 $${\cal L}_{\rho,\tilde{\pi}}(\rho(\gamma h)\epsilon'\otimes \Phi\mu(\gamma' a)e'(x'_{*}h'),\rho(\gamma'h')\epsilon\otimes \Phi\mu(\gamma a)e(x_{*}h))d\gamma\,d\gamma'\,da.$$
 Cette \'egalit\'e et (6) entra\^{\i}nent que, si ${\cal L}_{\rho,\tilde{\pi}}$ est nulle, ${\cal L}_{\pi,\rho}$ l'est aussi. Inversement, supposons ${\cal L}_{\rho,\tilde{\pi}}$ non nulle. Fixons $\tilde{e},\tilde{e}'\in E_{\tilde{\pi}}$, $\epsilon,\epsilon'\in E_{\rho}$ tels que ${\cal L}_{\rho,\tilde{\pi}}(\epsilon'\otimes \tilde{e}',\epsilon\otimes\tilde{e})\not=0$. Fixons un sous-groupe ouvert compact $K'\subset GL_{k-r-1}(F)$ tel que $\epsilon$ et $\epsilon'$ soient invariant par $K'$ et $\bar{\xi}$ soit trivial sur $U_{k-r-1}(F)\cap K'$. L'espace des fonctions $\gamma\mapsto \phi\mu(\gamma)\eta$ sur $GL_{k-1}(F)$, quand $\eta$ parcourt $E_{\mu}$, est le mod\`ele de Kirillov de $\mu$. Ce mod\`ele contient toutes les fonctions localement constantes, se transformant \`a gauche par $U_{k-1}(F)$ selon le caract\`ere $\xi$ et de support d'image compacte dans $U_{k-1}(F)\backslash GL_{k-1}(F)$. On peut donc choisir $\eta$ tel que  la fonction  $\gamma\mapsto \phi\mu(\gamma)\eta$  sur $GL_{k-r-1}(F)$ soit \`a support dans $U_{k-r-1}(F)K'$ et vaille $1$ sur $K'$. On choisit $\eta'=\eta$.  Fixons une sous-vari\'et\'e analytique $\Lambda$ de $K^H$, contenant $1$, telle que l'application produit $(Q_{0}(F)\cap H(F))\times \Lambda\to H(F)$ soit un hom\'eomorphisme  au voisinage des \'el\'ements neutres (c'est-\`a-dire un hom\'eomorphisme d'un voisinage de $(1,1)$ dans l'espace de d\'epart sur un voisinage de $1$ dans l'espace d'arriv\'ee). Les arguments que l'on a utilis\'es ci-dessus pour d\'ecomposer les int\'egrales 
 montrent que l'application
 $$\begin{array}{ccc}(U(F)\cap U_{\bar{Q}}(F))\times \Lambda&\to&Q(F)\backslash G(F)\\ (x,h)&\mapsto& x_{*}h\\ \end{array}$$
 est aussi un hom\'eomorphisme au voisinage des \'el\'ements neutres. On peut donc choisir un voisinage ouvert compact $\omega_{U}$ de $1$ dans $U(F)\cap U_{\bar{Q}}(F)$ et un  voisinage ouvert compact $\omega_{H}\subset \Lambda$ de sorte qu'il existe deux fonctions $e,e'\in E_{\pi}$ v\'erifiant les propri\'et\'es suivantes: l'image dans $Q(F)\backslash G(F)$ du support de $e$, resp. $e'$, est \'egale \`a l'image de $\omega_{U}\times \omega_{H}$ par l'application pr\'ec\'edente; pour $(x,h)\in \omega_{U}\times \omega_{H}$, $e(x_{*}h)=\eta\otimes \tilde{e}$, resp. $e'(x_{*}h)=\eta\otimes \tilde{e}$.   
 Quitte \`a restreindre $\omega_{H}$, on peut supposer que $\omega_{H}$ fixe $\epsilon$ et $\epsilon'$. Calculons ${\cal L}_{\pi,\rho}(\epsilon'\otimes e',\epsilon\otimes e)$. On a
 
(17)  $I(x',h',x,h)=0$ si $h$ ou $h'$ n'appartient pas \`a $(Q_{0}(F)\cap H(F))\omega_{H}$; pour $h,h'\in \omega_{H}$, $I(x',h',x,h)=0$ si $x$ ou $x'$ n'appartient pas \`a $\omega_{U}$. 

En effet, modulo multiplication \`a gauche par des \'el\'ements de $Q_{0}(F)\cap H(F)$, on peut supposer que $h,h'\in H_{0}(F)K^H$. Alors $I(x',h',x,h)$ est calcul\'e par la formule (16). Il r\'esulte de cette formule et de la d\'efinition de $e$ que, si $I(x',h',x,h)\not=0$, il existe $x''\in \omega_{U}$ et $h''\in \omega_{H}$ de sorte que $xh\in Q(F)x''h''$. Mais cette relation entra\^{\i}ne que les images dans $(Q_{0}(F)\cap H(F))\backslash H(F)$ de $h$ et $h''$ sont \'egales. Donc $h\in (Q_{0}(F)\cap H(F))\omega_{H}$.  Supposons maintenant $h\in \omega_{H}$. Si $I(x',h',x,h)\not=0$, on a encore des $x''$ et $h''$ comme ci-dessus, et les images de $h$ et $h''$ dans $(Q_{0}(F)\cap H(F))\backslash H(F)$ sont \'egales.  Mais l'application naturelle de $\omega_{H}$ dans ce quotient est injective. Donc $h=h''$ puis $x\in Q(F)x''$. Puisque $x,x''\in U(F)$, cela entra\^{\i}ne $x\in (U(F)\cap Q(F))x''$. L'application naturelle de $U(F)\cap U_{\bar{Q}}(F)$ dans $(U(F)\cap Q(F))\backslash U(F)$ est injective, donc $x=x''\in \omega_{U}$. Les m\^emes arguments s'appliquent \`a $x'$ et $h'$.
 
 Cette propri\'et\'e nous permet de remplacer dans la formule (6) les int\'egrales sur
 $(Q_{0}(F)\cap H(F))\backslash H(F)$ et $U(F)\cap U_{\bar{Q}}(F)$ par des int\'egrales sur $\omega_{H}$ et $\omega_{U}$, pour une mesure convenable sur $\omega_{H}$. Sur ces ensembles d'int\'egration, $I(x',h',x,h)$ est constante, \'egale \`a $I(1,1,1,1)$. Donc ${\cal L}_{\pi,\rho}(\epsilon'\otimes e',\epsilon\otimes e)$ est un multiple non nul de ce terme. On peut choisir $c$ assez grand pour que $\omega_{[k-r]}(c+c_{\psi})$ fixe $\eta$ et $\eta'$. Alors la formule (16) devient
 $$I(1,1,1,1)=mes(\omega_{[k-r]}(c+c_{\psi})\int_{(U_{k-r-1}(F)\backslash GL_{k-r-1}(F))^2}\vert det(\gamma\gamma')\vert _{F}^{(1-d_{\tilde{V}}-k)/2}$$
 $${\cal L}_{\rho,\tilde{\pi}}(\rho(\gamma )\epsilon'\otimes\tilde{e}',\rho(\gamma')\epsilon\otimes \tilde{e})(\phi\mu(\gamma )\eta)( \phi\mu(\gamma' )\eta') d\gamma\,d\gamma'.$$
 D'apr\`es les choix de $\eta$ et $\eta'$, cette expression est un multiple non nul de ${\cal L}_{\rho,\tilde{\pi}}(\epsilon'\otimes \tilde{e}',\epsilon\otimes \tilde{e})$. Donc ${\cal L}_{\pi,\rho}(\epsilon'\otimes e',\epsilon\otimes e)\not=0$, ce qui ach\`eve la d\'emonstration. $\square$
   
 \bigskip
 
 \subsection{Variable d'induction et entrelacements temp\'er\'es}
 
 Soient $(V,q_{V})$ et $(W,q_{W})$ deux espaces quadratiques compatibles, avec $d_{W}<d_{V}$. Soit $Q=LU_{Q}$ un sous-groupe parabolique de $G$.  Le groupe $L$ s'identifie \`a un produit
 $$L=GL_{k_{1}}\times...\times GL_{k_{s}}\times \tilde{G},$$
 o\`u $\tilde{G}$ est le groupe sp\'ecial orthogonal d'un sous-espace $\tilde{V}$ de $V$. Le groupe $A_{L}(F)$ s'identifie \`a $F^{\times s}$.  On suppose que $K$ est en bonne position relativement \`a $L$. Soient $\mu_{1}$, resp. $\mu_{2},...,\mu_{s},\tilde{\pi}$ des repr\'esentations admissibles irr\'eductibles et temp\'er\'ees de $GL_{k_{1}}(F)$, resp. $GL_{k_{2}}(F)$,...,$GL_{k_{s}}(F)$, $\tilde{G}(F)$. Posons $\tau=\mu_{1}\otimes...\otimes \mu_{s}\otimes \tilde{\pi}$. C'est une repr\'esentation de $L(F)$. Pour $\lambda\in i{\cal A}_{L,F}^*$, on d\'efinit la repr\'esentation   induite $\pi_{\lambda}=Ind_{Q}^G(\tau_{\lambda})$, que l'on r\'ealise dans  l'espace ${\cal K}_{Q,\tau}^G$. Fixons des produits scalaires invariants sur les espaces de $\mu_{1}$,...,$\mu_{s}$, $\tilde{\pi}$. On construit la forme hermitienne
 $$(e',e)=\int_{K}(e'(k),e(k))dk$$
 sur ${\cal K}_{Q,\tau}^G$.   Soit $\rho\in Temp(H)$. Pour $\lambda\in i{\cal A}_{L,F}^*$, $\epsilon,\epsilon'\in E_{\rho}$ et $e,e'\in {\cal K}_{Q,\tau}^G$, on d\'efinit ${\cal L}_{\pi_{\lambda},\rho}(\epsilon'\otimes e',\epsilon\otimes e)$.
 
 \ass{Lemme}{(i) L'application $\lambda\mapsto {\cal L}_{\pi_{\lambda},\rho}(\epsilon'\otimes e',\epsilon\otimes e)$ est $C^{\infty}$ sur $i{\cal A}_{L,F}^*$.
 
 (ii) Les trois conditions suivantes sont \'equivalentes:
 
 - il existe $\lambda$ tel que ${\cal L}_{\pi_{\lambda},\rho}$ ne soit pas nul;
 
 - pour tout $\lambda$, ${\cal L}_{\pi_{\lambda},\rho}$ n'est pas nul;
 
 - ${\cal L}_{\tilde{\pi},\rho}$ n'est pas nul.
 
 (iii) Si ces conditions sont v\'erifi\'ees, on peut choisir des familles finies $(\epsilon'_{i})_{i=1,...,n}$, $(\epsilon_{i})_{i=1,...,n}$ d'\'el\'ements de $E_{\rho}$, des familles finies $(e'_{i})_{i=1,...,n}$, $(e_{i})_{i=1,...,n}$ d'\'el\'ements de ${\cal K}_{Q,\tau}^G$ et une famille finie $(\varphi_{i})_{i=1,...,n}$ de fonctions $C^{\infty}$ sur $i{\cal A}_{L,F}^*$ de sorte que
 $$\sum_{i=1,...,n}\varphi_{i}(\lambda){\cal L}_{\pi_{\lambda},\rho}(\epsilon'_{i}\otimes e'_{i},\epsilon_{i}\otimes e_{i})=1$$
 pour tout $\lambda$.} 
 
 Preuve. Pour $\epsilon,\epsilon',e,e'$ fix\'es, la preuve du lemme 3.5 montre   qu'il existe $c_{0}$ tel que, pour tout $c\geq c_{0}$ et tout $\lambda\in i{\cal A}_{L,F}^*$, on ait l'\'egalit\'e
 $${\cal L}_{\pi_{\lambda},\rho}(\epsilon'\otimes e',\epsilon\otimes e)={\cal L}_{\pi_{\lambda},\rho,c}(\epsilon'\otimes e',\epsilon\otimes e).$$
 Soit $D$ un op\'erateur diff\'erentiel sur $i{\cal A}_{L,F}^*$, \`a coefficients $C^{\infty}$. Il existe un r\'eel $R\geq0$ tel que l'on ait une majoration
 $$\vert D(e',\pi_{\lambda}(g)e)\vert<<\sigma(g)^R\Xi^G(g)$$ 
 pour tous $g\in G(F)$ et $\lambda\in i{\cal A}_{L,F}^*$. Alors $D{\cal L}_{\pi_{\lambda},\rho,c}(\epsilon'\otimes e',\epsilon\otimes e)$ est uniform\'ement major\'ee par l'int\'egrale
 $$\int_{H(F)U(F)_{c}}\Xi^H(h)\Xi^G(hu)\sigma(hu)^Rdu\,dh$$
 qui est convergente d'apr\`es 4.3(4). Le th\'eor\`eme usuel de d\'erivation d'une int\'egrale d\'ependant d'un param\`etre entra\^{\i}ne que $\lambda\mapsto {\cal L}_{\pi_{\lambda},\rho,c}(\epsilon'\otimes e',\epsilon\otimes e)$ est $C^{\infty}$. D'o\`u (i).
 
 Introduisons le sous-groupe parabolique $Q'$ de $G$ contenant $Q$, dont la composante de L\'evi est $L'=GL_{k}\times \tilde{G}$, o\`u $k=k_{1}+...+k_{s}$. Soit $\lambda\in i{\cal A}_{L,F}^*$. Posons $\tau^{L'}_{\lambda}=Ind_{L'\cap Q}^{L'}(\tau_{\lambda})$. C'est une repr\'esentation de $L'(F)$ de la forme $\mu_{\lambda}\otimes \tilde{\pi}$. La repr\'esentation $\mu_{\lambda}$ est une repr\'esentation  temp\'er\'ee et irr\'eductible de $GL_{k}(F)$. On sait en effet que, dans les groupes lin\'eaires,  l'induite d'une repr\'esentation temp\'er\'ee irr\'eductible est irr\'eductible. On a $\pi_{\lambda}=Ind_{Q'}^G(\tau^{L'}_{\lambda})$. Appliquant la proposition 5.2, on voit que la non-nullit\'e de ${\cal L}_{\pi_{\lambda},\rho}$est \'equivalente \`a la non-nullit\'e de ${\cal L}_{\tilde{\pi},\rho}$. L'assertion (ii) s'en d\'eduit.
 
 Supposons satisfaites les conditions de (ii). Pour tout $\lambda$, on choisit des \'el\'ements $\epsilon_{\lambda}$,$\epsilon'_{\lambda}$, $e_{\lambda}$ et $e'_{\lambda}$ tels que 
 $${\cal L}_{\pi_{\lambda},\rho}(\epsilon'_{\lambda}\otimes e'_{\lambda},\epsilon_{\lambda}\otimes e_{\lambda})\not=0.$$
 Gr\^ace \`a (i), il existe un voisinage $\omega_{\lambda}$ de $\lambda$ dans $i{\cal A}_{L,F}^*$tel que 
 $${\cal L}_{\pi_{\lambda'},\rho}(\epsilon'_{\lambda}\otimes e'_{\lambda},\epsilon_{\lambda}\otimes e_{\lambda})\not=0$$
 pour tout $\lambda'\in \omega_{\lambda}$. Puisque $i{\cal A}_{L,F}^*$ est compact, on peut choisir une famille finie $(\lambda_{i})_{i=1,...,n}$ de sorte que 
 $$\bigcup_{i=1,...,n}\omega_{\lambda_{i}}=i{\cal A}_{L,F}^*.$$
 On d\'efinit $\epsilon_{i},\epsilon'_{i},e_{i},e'_{i}$ par $\epsilon_{i}=\epsilon_{\lambda_{i}}$ etc... Notons $\varphi'_{i}$ la conjugu\'ee complexe de la fonction
 $$\lambda\mapsto {\cal L}_{\pi_{\lambda},\rho}(\epsilon'_{i}\otimes e'_{i},\epsilon_{i}\otimes e_{i}).$$
 et posons
$$\varphi(\lambda)=\sum_{i=1,...,n}\varphi'_{i}(\lambda) {\cal L}_{\pi_{\lambda},\rho}(\epsilon'_{i}\otimes e'_{i},\epsilon_{i}\otimes e_{i}).$$
Alors $\varphi$ est une fonction $C^{\infty}$ \`a valeurs strictement positives. On satisfait la condition (iii) en d\'efinissant $\varphi_{i}=\varphi'_{i}/\varphi$. $\square$

 \bigskip
 
 \subsection{Le cas des induites r\'eductibles}
 
 On conserve les donn\'ees de 5.3. Posons $\pi=\pi_{0}$. Cette repr\'esentation est  semi-simple et toute sous-repr\'esentation irr\'eductible y intervient avec multiplicit\'e $1$.  
  
  \ass{Lemme}{Supposons ${\cal L}_{\pi,\rho}$ non nulle. Alors il existe une unique sous-repr\'esentation irr\'eductible $\pi'$ de $\pi$ telle que ${\cal L}_{\pi',\rho}$ soit non nulle.}
  
  Preuve. Par construction, pour toute sous-repr\'esentation irr\'eductible $\pi'$ de $\pi$, la forme ${\cal L}_{\pi',\rho}$ est la restriction de ${\cal L}_{\pi,\rho}$ \`a l'espace $E_{\rho}\otimes E_{\pi'}$, o\`u $E_{\pi'}\subset {\cal K}_{Q,\tau}^G$ est l'espace de $\pi'$.  D'autre part, soient $\epsilon',\epsilon''\in E_{\rho}$ et $e',e''\in {\cal K}_{Q,\tau}^G$. Pour $\lambda\in i{\cal A}_{L,F}^*$ en position g\'en\'erale, la repr\'esentation $\pi_{\lambda}$ est irr\'eductible. D'apr\`es 5.1(1), on a l'\'egalit\'e
 $$\vert {\cal L}_{\pi_{\lambda},\rho}(\epsilon'\otimes e'',\epsilon''\otimes e')\vert ^2= \vert {\cal L}_{\pi_{\lambda},\rho}(\epsilon'\otimes e',\epsilon'\otimes e')\vert\vert {\cal L}_{\pi_{\lambda},\rho}(\epsilon''\otimes e'',\epsilon''\otimes e'')\vert.$$
 D'apr\`es le lemme 5.3(i), tous les termes sont continus en $\lambda$. L'\'egalit\'e est donc vraie pour tout $\lambda$. Pour $\lambda=0$, elle donne
  $$(1) \qquad \vert {\cal L}_{\pi,\rho}(\epsilon'\otimes e'',\epsilon''\otimes e')\vert ^2= \vert {\cal L}_{\pi,\rho}(\epsilon'\otimes e',\epsilon'\otimes e')\vert\vert {\cal L}_{\pi,\rho}(\epsilon''\otimes e'',\epsilon''\otimes e'')\vert.$$
  Puisque ${\cal L}_{\pi,\rho}$ n'est pas nulle, on peut choisir $\epsilon',\epsilon'',e',e''$ tels que ${\cal L}_{\pi,\rho}(\epsilon'\otimes e'',\epsilon''\otimes e')\not=0$. Par lin\'earit\'e, on peut supposer qu'il existe des 
sous-repr\'esentations irr\'eductibles $\pi'$ et $\pi''$ de $\pi$ telles que $e'\in E_{\pi'}$ et $e''\in E_{\pi''}$. 
La relation (1) entra\^{\i}ne que ${\cal L}_{\pi,\rho}(\epsilon'\otimes e',\epsilon'\otimes e')\not=0$, donc ${\cal L}_{\pi',\rho}$ n'est pas nulle. D'o\`u l'existence affirm\'ee dans l'\'enonc\'e. Soient maintenant  $\pi'$ et $\pi''$ des sous-repr\'esentations irr\'eductibles de $\pi$ telles que ${\cal L}_{\pi',\rho}$ et ${\cal L}_{\pi'',\rho}$ ne soient pas nulles. En appliquant  5.1(2), on peut fixer $\epsilon',\epsilon''\in E_{\rho}$, $e'\in E_{\pi'}$ et $e''\in E_{\pi''}$ de sorte que ${\cal L}_{\pi,\rho}(\epsilon'\otimes e',\epsilon'\otimes e')\not=0$ et ${\cal L}_{\pi,\rho}(\epsilon''\otimes e'',\epsilon''\otimes e'')\not=0$. D'apr\`es (1), on a aussi ${\cal L}_{\pi,\rho}(\epsilon'\otimes e'',\epsilon''\otimes e')\not=0$. Par construction de ${\cal L}_{\pi,\rho}$, cela entra\^{\i}ne que le coefficient $g\mapsto (e'',\pi(g)e')$ n'est pas nul. Les sous-espaces $E_{\pi'}$ et $E_{\pi''}$ ne sont donc pas orthogonaux, ce qui implique $\pi'=\pi''$.  $\square$

  \bigskip
  
  \subsection{Le cas $r=0$: majoration des entrelacements}
  
    Soient $(V,q_{V})$ et $(W,q_{W})$ deux espaces quadratiques compatibles. On suppose $d_{V}=d_{W}+1$. Soit $\pi$, resp. $\rho$, une repr\'esentation  temp\'er\'ee de $G(F)$, resp. $H(F)$. Soient $l\in Hom_{H,\xi}(\pi,\rho)$ et $\underline{K}$ un sous-groupe ouvert compact de $G(F)$.
  
  \ass{Proposition}{Pour tous $\epsilon\in E_{\rho}$, $e\in E_{\pi}$, il existe un r\'eel  $D\geq0$ tels que
  $$\int_{\underline{K}}\vert (\epsilon,l(\pi(kg)e))\vert  dk <<\sigma(g)^D\Xi^G(g)$$
  pour tout $g\in G(F)$. }
  
  Preuve.   La proposition d\'epend du groupe $\underline{K}$. Soit $\underline{K}'$ un autre sous-groupe ouvert compact de $G(F)$. Montrons d'abord que la proposition relative \`a $\underline{K}$ est \'equivalente \`a celle relative \`a $\underline{K}'$. On peut d\'emontrer que la proposition relative \`a $\underline{K}$ est \'equivalente \`a celle relative \`a $\underline{K}\cap \underline{K}'$, puis que celle-ci est \'equivalente \`a celle relative \`a $\underline{K}'$. Cela nous ram\`ene au cas o\`u $\underline{K}'\subset \underline{K}$. Notons $I_{\underline{K}}(\epsilon,e,g)$ l'int\'egrale de l'\'enonc\'e et $I_{\underline{K}'}(\epsilon,e,g)$ son analogue pour le groupe $\underline{K}'$. On a les in\'egalit\'es
  $$I_{\underline{K}'}(\epsilon,e,g)\leq I_{\underline{K}}(\epsilon,e,g)\leq \sum_{k\in \underline{K}'\backslash \underline{K}}I_{\underline{K}'}(\epsilon,e,kg),$$
  dont on d\'eduit la propri\'et\'e voulue.
  
  On a fix\'e un \'el\'ement $v_{0}\in V$, non nul et orthogonal \`a $W$. On peut multiplier $q_{V}$ et $q_{W}$ par une constante, cela ne change rien aux  groupes $G$ et $H$ ni aux donn\'ees de notre probl\`eme. On peut donc supposer $q_{V}(v_{0})=1$, ce qui simplifiera la r\'edaction. La proposition est triviale si $q_{V}$ est anisotrope: dans ce cas, $G(F)$ est compact et la fonction $g\mapsto (\epsilon,l(\pi(g)e)$ ne prend qu'un nombre fini de valeurs. Supposons que $q_{V}$ n'est pas anisotrope. Fixons un syst\`eme hyperbolique maximal $(v_{\pm i})_{i=1,...,n}$ de $V$. On a $n\geq 1$. Notons $V_{an}$ l'orthogonal du sous-espace engendr\'e par ces \'el\'ements et  posons:
  $$\underline{R}_{an}=\{v\in V_{an}; q_{V}(v)\in \mathfrak{o}_{F}\}.$$
  Notons $\underline{R}$ la somme de $\underline{R}_{an}$ et du $\mathfrak{o}_{F}$-module engendr\'e par les $v_{\pm i}$, $i=1,...,n$. Ce r\'eseau poss\`ede un \'el\'ement $v$ tel que $q_{V}(v)=1$, par exemple $v=v_{1}+v_{-1}$. Quitte \`a transformer le syst\`eme hyperbolique et le r\'eseau $\underline{R}$ par l'action d'un \'el\'ement de $G(F)$, on peut donc supposer $v_{0}\in \underline{R}$. Notons $\underline{K}$ le stabilisateur de $\underline{R}$ dans $G(F)$. C'est un sous-groupe compact maximal de $G(F)$, pas forc\'ement sp\'ecial (il l'est si $V_{an}=\{0\}$ ou si $V_{an}$ poss\`ede un \'el\'ement $v$ tel que $q_{V}(v)\in \mathfrak{o}_{F}^{\times}$). Le groupe $\underline{K}$ contient en tout cas comme sous-groupe d'indice fini un sous-groupe parahorique de $G(F)$. Soit $A_{min}$ le sous-tore d\'eploy\'e maximal de $G$ form\'e des \'el\'ements qui conservent chaque droite $Fv_{\pm i}$ et agissent trivialement sur $V_{an}$. Notons $M_{min}$ le commutant de $A_{min}$ dans $G$, qui est un L\'evi minimal de $G$. Le groupe $\underline{K}$ est en bonne position relativement \`a $M_{min}$. Pour obtenir la majoration  de l'\'enonc\'e pour le groupe $\underline{K}$, il suffit de majorer $I_{\underline{K}}(\epsilon,e,a)$ pour $a\in A_{min}(F)$. En effet, d'apr\`es Bruhat-Tits, il existe un sous-ensemble compact $\Gamma$ de $G(F)$ tel que tout \'el\'ement $g\in G(F)$ s'\'ecrive $g=ka\gamma$, avec $k\in \underline{K}$, $a\in A_{min}(F)$ et $\gamma\in \Gamma$. On a alors
  $$I_{\underline{K}}(\epsilon,e,g)=I_{\underline{K}}(\epsilon,\pi(\gamma)e,a)$$
  et la majoration de ce dernier terme implique la majoration cherch\'ee de $I_{\underline{K}}(\epsilon,e,g)$. Comme dans la preuve de 4.9, introduisons le groupe de permutations $\mathfrak{S}$ et, pour $s\in \mathfrak{S}$, le sous-ensemble $A_{min}(F)_{s}^-$ de $A_{min}(F)$.
   Le groupe $A_{min}(F)$ est r\'eunion des sous-ensembles $A_{min}(F)_{s}^-$, quand $s$ parcourt $\mathfrak{S}$.   On peut donc fixer $s$ et se limiter aux $a\in A_{min}(F)_{s}^-$. Quitte \`a r\'eindexer notre syst\`eme hyperbolique, on peut supposer que  $s$ est l'identit\'e. Notons simplement $A_{min}(F)^-$ l'ensemble des $a\in A_{min}(F)$ tels que $val_{F}(a_{n})\geq...\geq val_{F}(a_{1})\geq 0$. On peut se limiter \`a majorer $I_{\underline{K}}(\epsilon,e,a)$ pour $a\in A_{min}(F)^-$.   
  
  Fixons un entier $l\geq0$ tel que
  
 (1) $l\geq val_{F}(2)$;
  
  (2) pour tout $x\in F$, la condition $val_{F}(x)\geq l+1$ entra\^{\i}ne que $1+x$ est un carr\'e et poss\`ede une racine $\sqrt{1+x}$ telle que $val_{F}(\sqrt{1+x}-1)=val_{F}(x/2)$.
  
  Si $p\not=2$, on peut prendre $l=0$. Pour $i=1,...,n$, posons $v'_{i}=\varpi_{F}^{il}v_{i}$ et $v'_{-i}=\varpi_{F}^{-il}v_{-i}$. Notons $\underline{R}'$ le r\'eseau somme de $\underline{R}_{an}$ et du $\mathfrak{o}_{F}$-module engendr\'e par les $ v'_{\pm i}$ pour $i=1,...,l$. Notons $\underline{K}'$ le stabilisateur de $\underline{R}'$ dans $G(F)$ ($\underline{R}'=\underline{R}$ et $\underline{K}'=\underline{K}$ si $l=0$). S'il existe $v\in V_{an}$ tel que $q_{V}(v)=1$, on fixe un tel \'el\'ement que l'on note $v'_{0}$ et on pose $E=\{\pm v'_{0}\}$, $\iota=1$. Sinon, on pose $v'_{0}=v_{1}+v_{-1}$, $E=\{\varpi_{F}^{-c}v'_{1}+\varpi_{F}^cv'_{-1}; c=0,...,l\}$, $\iota=2$. Remarquons que $E\subset \underline{R}'$ dans le premier cas.  Fixons $a\in A_{min}(F)^-$ et posons $\alpha_{i}=val_{F}(a_{i})$ pour tout $i=1,...,n$.    Notons ${\cal B}'$ le sous-ensemble des suites ${\bf b}=(b_{1},...,b_{n})\in {\mathbb Z}^{n}$ telles que
  
  $b_{n}\geq b_{n-1}\geq...\geq b_{1}$;
  
  $\alpha_{n}-b_{n}\geq \alpha_{n-1}-b_{n-1}\geq...\geq \alpha_{1}-b_{1}\geq 0$.

   Notons $X'=\{x\in V; q_{V}(x)=1\}$.   Soit $x\in X'$. Pour $i=1,...,n$, introduisons l'ensemble
  $$B_{i}(x)=\{b\in {\mathbb Z}; b\leq\alpha_{i},\, x\in \varpi_{F}^{b-\alpha_{i}}a\underline{R}'+\varpi_{F}^b\underline{R}' \}.$$
  L'ensemble $B_{i}(x)$ n'est pas vide: il contient tout \'el\'ement suffisamment n\'egatif. Posons $b_{i}(x)=sup(B_{i}(x))$, puis ${\bf b}(x)=(b_{1}(x),...,b_{n}(x)\}$. Montrons que
  
  (3) ${\bf b}(x)$ appartient \`a ${\cal B}'$.
  
  Pour $i\leq n-1$ et $b\in B_{i}(x)$, on a $b\leq \alpha_{i}\leq \alpha_{i+1}$ et
  $$x\in \varpi_{F}^{b-\alpha_{i}}a\underline{R}'+\varpi_{F}^b\underline{R}' \subset \varpi_{F}^{b-\alpha_{i+1}}a\underline{R}' +\varpi_{F}^{b}\underline{R}' .$$
  Donc $b\in B_{i+1}(x)$. Ainsi $B_{i}(x)\subset B_{i+1}(x)$, d'o\`u $b_{i+1}(x)\geq b_{i}(x)$. Posons 
  $$C_{i}(x)=\{c\in {\mathbb Z}; c\leq0,\, x\in \varpi_{F}^ca\underline{R}'+\varpi_{F}^{c+\alpha_{i}}\underline{R}' \}.$$
  On a $b_{i}(x)-\alpha_{i}=sup(C_{i}(x))$. On voit comme ci-dessus que $C_{i+1}(x) \subset C_{i}(x)$, d'o\`u $b_{i+1}(x)-\alpha_{i+1}\leq b_{i}(x)-\alpha_{i}$. Enfin, $\alpha_{1}-b_{1}(x)\geq0$ par d\'efinition de $b_{1}(x)$. Cela prouve (3). 
  
  Le groupe $\underline{K}'\cap a\underline{K}'a^{-1}$ agit sur $X'$. Par construction, les ensembles $B_{i}(x)$ ne d\'ependent que de l'orbite de $x$ pour cette action. L'\'el\'ement ${\bf x}$ ne d\'epend donc lui-aussi que de cette orbite.  L'application
  $$\begin{array}{ccc}H(F)\backslash G(F)&\to&X' \\ g&\mapsto &g^{-1}v_{0}\\ \end{array}$$
  est un isomorphisme. Le quotient des mesures de Haar sur $H(F)$ et $G(F)$ est une mesure sur l'ensemble de d\'epart que l'on transporte \`a l'ensemble d'arriv\'ee.  Posons $X=\underline{R}\cap X'$.  Pour ${\bf b}\in {\cal B}'$, notons $X({\bf b})$ l'ensemble des $x\in X$ tels que ${\bf b}(x)={\bf b}$. Notons ${\cal B}$ le sous-ensemble des ${\bf b}\in {\cal B}'$ qui v\'erifient de plus $b_{1}\geq-nl$ et, dans le cas o\`u $\iota=2$, $b_{1}\leq0$. Montrons que
  
  (4) pour  ${\bf b}\in {\cal B}'$, $X({\bf b})$ est vide si ${\bf b}\not\in {\cal B}$; on a une majoration
  $$mes(X({\bf b}))<<q^{-\sum_{i=1,...,n}b_{i}}$$
  pour tout ${\bf b}=(b_{1},...,b_{n})\in {\cal B}$.

 Soit $x\in X$. On a l'inclusion $\underline{R}\subset \varpi_{F}^{-nl}\underline{R}'$, donc $-nl\in B_{1}(x)$ et $b_{1}(x)\geq-nl$. Supposons $b_{1}(x)>0$, donc aussi $\alpha_{1}>0$. Ecrivons $x=x_{an}+\sum_{i=\pm 1,...,\pm n}x_{ i}v'_{i}$, avec $x_{an}\in V_{an}$ et $x_{ i}\in F$ pour tout $i$. Puisque $x\in \underline{R}$, on a $x_{an}\in \underline{R}_{an}$ et $val_{F}(x_{-i})\geq l$ pour tout $i=1,...,n$. La condition
  $$x\in \varpi_{F}^{1-\alpha_{1}}a\underline{R}'+\varpi_{F}\underline{R}' $$
  entra\^{\i}ne   $val_{F}(x_{i})\geq 1$ pour tout $i$. La condition $q_{V}(x)=1$ entra\^{\i}ne alors $val_{F}(q_{V}(x_{an})-1)\geq l+1$. Donc $q_{V}(x_{an})$ est un carr\'e et l'\'el\'ement $v=x_{an}/\sqrt{q_{V}(x_{an})}$ de $V_{an}$ v\'erifie $q_{V}(v)=1$. Donc $\iota=1$. Cela prouve la premi\`ere assertion de (4).
 
 Soit ${\bf b}\in {\cal B}$. Comme ci-dessus, \'ecrivons tout \'el\'ement $x\in V$ sous la forme $x=x_{an}+\sum_{i=\pm 1,...\pm n}x_{i}v'_{i}$. Notons $\underline{R}({\bf b})$ l'ensemble des $x\in \underline{R}$ tels que $val_{F}(x_{i})\geq b_{i}$ pour tout $i=1,...,n$. On v\'erifie que $X({\bf b})\subset \underline{R}({\bf b})$.
  Consid\'erons l'application
  $$\begin{array}{ccc}V&\to&F\\ v&\mapsto &q_{V}(v).\\ \end{array}$$
  Soit $\nu\in F^{\times}$, notons $V_{\nu}$ la fibre de cette application au-dessus de $\nu$. Au-dessus de $F^{\times}$, l'application est une fibration localement triviale. Des mesures de Haar sur $V$ et $F$ se d\'eduit donc une mesure sur $V_{\nu}$, notons-la $d_{\nu}v$. Pour $f\in C_{c}^{\infty}(V)$, on a l'\'egalit\'e
  $$\int_{V}f(v)dv=\int_{F}\int_{V_{\nu}}f(v)d_{\nu}v\,d\nu.$$
  La mesure $d_{\nu}$ est invariante par l'action de $G(F)$ et c'est la seule mesure v\'erifiant cette condition, \`a une constante. Pour $\nu=1$, on a $V_{1}=X'$ et on peut  supposer que $d_{1}$ co\"{\i}ncide avec la mesure que l'on a pr\'ec\'edemment introduite sur cet ensemble. D'autre part, pour $\lambda\in F^{\times}$, on v\'erifie que 
  $$\int_{V_{\nu}}f^{\lambda}( v)d_{\nu}v=\vert \lambda\vert _{F}^{2-dim(V)}\int_{V_{\lambda^2\nu}}f(v)d_{\lambda^2\nu},$$
  o\`u, comme d'habitude, $f^{\lambda}$ est d\'efinie par $f^{\lambda}(v)=f(\lambda v)$. Supposons $f^{\lambda}=f$ pour $\lambda\in \mathfrak{o}_{F}^{\times}$. On d\'eduit des relations  ci-dessus l'\'egalit\'e  
  $$mes(\mathfrak{o}_{F}^{\times,2})\int_{V_{\nu}}f(v)d_{\nu}v=\int_{v\in V; q_{V}(v)\in \nu\mathfrak{o}_{F}^{\times,2}}f(v)dv.$$
  On applique cela \`a $\nu=1$ et \`a la fonction caract\'eristique $f$ de l'ensemble $\underline{R}({\bf b})$. On obtient une majoration
  $$mes(X({\bf b}))<<mes(\{x\in \underline{R}({\bf b}); q_{V}(x)\in \mathfrak{o}_{F}^{\times,2}\})<<mes( \underline{R}({\bf b}))<<q^{-\sum_{i=1,...,n}b_{i}}.$$
  Cela prouve  la seconde assertion de (4).
  
   Soit $x\in \underline{R}$. Montrons que
  
  (5) dans l'orbite de $x$ sous l'action de $\underline{K}'\cap a\underline{K}'a^{-1}$, il existe un \'el\'ement $y$ (non n\'ecessairement dans $\underline{R}$)  de la forme $y=e+\sum_{i=\iota,...,n}y_{i}v'_{i}$, o\`u $e\in E$.
 
   Pour cela, on \'ecrit $x=x_{an}+\sum_{\pm 1,...,\pm n}x_{i}v'_{i}$ et on montre par r\'ecurrence  sur $j=1,...,n$ que
  
  (6) la conclusion de (5) est v\'erifi\'ee si $x_{-i}=0$ pour $i=j+1,...,n$.
  
  Pour $v',v''\in V$, on note $c(v',v'')$ l'\'el\'ement de $\mathfrak{g}(F)$ d\'efini par
  $$c(v',v'')(v)=q_{V}(v,v')v''-q_{V}(v,v'')v'.$$
   Notons $V_{j}$, resp. $V_{j-1}$, le sous-espace de $V$ engendr\'e par $V_{an}$ et les $v'_{\pm i}$ pour $i=1,...,j$, resp. $i=1,...,j-1$. Pour $v\in V_{j-1}$, on pose $k_{v}=exp(c(v'_{-j},v))$. On a
  $$k_{v}v'_{-j}=v'_{-j};$$
  $$k_{v}v'_{j}=v'_{j}+v-q_{V}(v)v'_{-j};$$
  $$k_{v}(z)=z-q_{V}(z,v)v'_{-j},\text{ pour }z\in V_{j-1};$$
  $$k_{v}z=z, \text{ pour }z\text{ orthogonal \`a }V_{j}.$$
  En particulier, si $v\in \underline{R}'\cap V_{j-1}$, on a $k_{v}\in \underline{K}'$. Mais $a^{-1}k_{v}a=exp(c(v_{-j},a_{j}a^{-1}v))$ et $a_{j}a^{-1}(\underline{R}'\cap V_{j-1})\subset \underline{R}'\cap V_{j-1}$. Donc, si $v\in \underline{R}'\cap V_{j-1}$, on a $k_{v}\in \underline{K}'\cap a\underline{K}'a^{-1}$. Soit $x$ v\'erifiant l'hypoth\`ese de (6). Ecrivons $x=u+x_{j}v'_{j}+z+x_{-j}v'_{-j}$, avec $u=\sum_{i=j+1,...,n}x_{i}v'_{i}$ et $z=x_{an}+\sum_{i=\pm 1,...\pm (j-1)}x_{i}v'_{i}$. L'\'el\'ement $u$ est orthogonal \`a $V_{j}$ et appartient \`a $\underline{R}$. L'\'el\'ement $z$ appartient \`a $\underline{R}\cap V_{j-1 }\subset \varpi_{F}^{(1-j)l}\underline{R}'\cap V_{j-1}$. On a aussi $val_{F}(x_{j})\geq -jl$, $val_{F}(x_{-j})\ge jl$. Supposons d'abord $val_{F}(x_{j})\geq (1-j)l+1$, donc $val_{F}(x_{j}x_{-j})\geq l+1$. L'\'egalit\'e $q_{V}(x)=1$ \'equivaut \`a $q_{V}(z)=1-x_{j}x_{-j}$. D'apr\`es le choix de $l$, on peut introduire une racine carr\'ee $\sqrt{q_{V}(z)}$ telle que $val_{F}(1-\sqrt{q_{V}(z)})=val_{F}(x_{j}x_{-j}/2)\geq1$. Posons 
  $$\lambda=\frac{1-\sqrt{q_{V}(z)}}{x_{j}\sqrt{q_{V}(z)}}.$$
  Alors $val_{F}(\lambda)= val_{F}(x_{-j}/2)\geq (j-1)l$. Donc $\lambda z\in \underline{R}'\cap V_{j-1}$ et $k_{\lambda z}\in \underline{K}'\cap a\underline{K}'a^{-1}$. On a
  $$k_{\lambda z}x=u+x_{j}v'_{j}+(1-\lambda x_{j})z+(x_{-j}-2\lambda q_{V}(z)-\lambda^2x_{j}q_{V}(z))v'_{-j}.$$
  On v\'erifie que le coefficient de $v'_{-j}$ est nul: on a choisi $\lambda$ pour cela. D'autre part $val_{F}(\lambda x_{j})\geq1$, donc $(1-\lambda x_{j})z\in \underline{R}$. Supposons $j\geq 2$. Alors $k_{\lambda z}x$ v\'erifie les m\^emes conditions que $x$, l'entier $j$ \'etant remplac\'e par $j-1$. L'hypoth\`ese de r\'ecurrence permet de conclure. Si $j=1$, l'\'el\'ement $k_{\lambda z}x$ est de la forme $k_{\lambda z}x=y_{an}+\sum_{i=1,...,n}y_{i}v'_{i}$. On a $q_{V}(y_{an})=1$, donc $\iota=1$. Il existe un \'el\'ement $g_{an}\in G_{an}(F)$ tel que $g_{an}(y_{an})=\pm v'_{0}$. Rappelons que $G_{an}(F)$ est contenu dans $\underline{K}'\cap a\underline{K}'a^{-1}$.  En posant $y=g_{an}k_{\lambda z}x$, cet \'el\'ement v\'erifie les conditions requises. Cela conclut le cas o\`u $val_{F}(x_{j})\geq (1-j)l+1$. Supposons maintenant $val_{F}(x_{j})\leq (1-j)l$. Alors $-z/x_{j}\in \underline{R}'\cap V_{j-1}$. L'\'el\'ement $k_{-z/x_{j}}x$ est \'egal \`a $u+x_{j}v'_{j}+x_{j}^{-1}v'_{-j}$. Si $j\geq \iota$, consid\'erons l'\'el\'ement $v'_{0}\in \underline{R}\cap V_{j-1}$. On a  $v'_{0}/x_{j}\in \underline{R}'\cap V_{j-1}$, $k_{v'_{0}/x_{j}}k_{-z/x_{j}}x=u+x_{j}v'_{j}+v'_{0}$ et cet \'el\'ement a la forme voulue. Si $j=1<\iota$, on a   $-l\leq val_{F}(x_{1})\leq 0$. En faisant agir un \'el\'ement de $A_{min}(F)\cap \underline{K}'\cap a\underline{K}'a^{-1}$, on peut remplacer $k_{-z/x_{1}}x$ par un \'el\'ement $u+\varpi_{F}^{-c}v'_{1}+\varpi_{F}^{c}v'_{-1}$, avec $0\leq c\leq l$, qui est de la forme voulue. Cela prouve (6) et (5).
  
   Soient ${\bf b}=(b_{1},...,b_{n})\in {\cal B}$ et $x\in X({\bf b})$. Montrons que
  
  (7) dans l'orbite de $x$ sous l'action de $\underline{K}'\cap a\underline{K}'a^{-1}$, il existe un \'el\'ement $y$ de la forme $y=e+\sum_{i=\iota,...,n}\varpi_{F}^{\beta_{i}}v'_{i}$, avec $e\in E$, $b_{i}\leq \beta_{i}\leq b_{i}+l$ pour tout $i=\iota,...,n$.
  
  D'apr\`es (5), on peut remplacer $x$ par un \'el\'ement de la forme $x=e+\sum_{i=\iota,...,n}x_{i}v'_{i}$ avec $e\in E$. On va prouver que, pour tout $i=\iota,...,n$, il existe un \'el\'ement $y$ de l'orbite de $x$ sous l'action de $\underline{K}'\cap a\underline{K}'a^{-1}$ tel que $y=e+\varpi_{F}^{\beta_{i}}v'_{i}+\sum_{j=\iota,...,n, j\not=i}x_{i}v'_{i}$, avec $b_{i}\leq \beta_{i}\leq b_{i}+l$. L'assertion (7)  r\'esultera de cette propri\'et\'e appliqu\'ee successivement \`a tout $i=\iota,...,n$. Dans le cas o\`u $\iota=2$, on introduit la coordonn\'ee suppl\'ementaire $x_{1}$ en posant $e=x_{1}v'_{1}+x_{1}^{-1}v'_{-1}$. En tout cas $x\in a\underline{R}'+\sum_{i=1,...,n}x_{i}v'_{i}$. Pour un tel \'el\'ement, on calcule ais\'ement ${\bf b}(x)$: pour tout $i=1,...,n$, $b_{i}(x)$ est le plus petit des entiers
  $$\alpha_{i},\,\,val_{F}(x_{j})\text{ pour }j=i,...,n,\,\,val_{F}(x_{j})-\alpha_{j}+\alpha_{i}\text{ pour }j=1,...,i-1.$$
  Fixons $i\in \{\iota,...,n\}$. Puisque $x\in X({\bf b})$, on a $b_{i}(x)=b_{i}$. En particulier, $val_{F}(x_{i})\geq b_{i}$. Si  $val_{F}(x_{i})\leq b_{i}+l$,   il suffit de remplacer $x$ par $a'x$, o\`u $a'$ est un \'el\'ement convenable de $A_{min}(F)\cap \underline{K}'\cap a\underline{K}'a^{-1}$ pour obtenir l'\'el\'ement $y$ cherch\'e. Supposons $val_{F}(x_{i})> b_{i}+l$. Si le plus petit des  entiers ci-dessus est $val_{F}(x_{j})$ pour un $j$ tel que $i+1\leq j\leq n$, on introduit l'\'el\'ement $k=exp(c(v'_{-j},v'_{i}))$. Il appartient \`a $\underline{K}'$. On a $a^{-1}ka=exp(a_{j}a_{i}^{-1}c(v'_{-j},v'_{i}))$. Puisque $j>i$, $val_{F}(a_{j}a_{i}^{-1})\geq0$, donc $a^{-1}ka\in \underline{K}'$ et $k\in \underline{K}'\cap a\underline{K}'a^{-1}$. Posons $y=kx$. Les coordonn\'ees de $y$ sont les m\^emes que celles de $x$, sauf $y_{i}$ qui vaut $x_{i}+x_{j}$. Alors $val_{F}(y_{i})=b_{i}$ et on conclut en appliquant encore un \'el\'ement convenable de $A_{min}(F)\cap \underline{K}'\cap a\underline{K}'a^{-1}$. Si le plus petit des entiers ci-dessus est $val_{F}(x_{j})-\alpha_{j}+\alpha_{i}$, pour un $j$ tel que $1\leq j\leq i-1$, on introduit l'\'el\'ement $k=exp(a_{i}a_{j}^{-1}c(v'_{-j},v'_{i}))$ et on pose  $y=kx$. On v\'erifie de m\^eme que $k$ appartient \`a $\underline{K}'\cap a\underline{K}'a^{-1}$ et que les coordonn\'ees de $y$ sont les m\^emes que celles de $x$, sauf $y_{i}$ qui v\'erifie $val_{F}(y_{i})=b_{i}$. On conclut comme pr\'ec\'edemment. Reste le cas o\`u le plus petit des entiers ci-dessus est $\alpha_{i}$.  On peut m\^eme supposer que tous les autres sont strictement plus grands. En particulier $val_{F}(x_{n})>\alpha_{n}\geq 0$. D'apr\`es la d\'efinition de $E$, cela entra\^{\i}ne  $\iota=1$. On introduit alors l'\'el\'ement $k=exp(c(e,a_{i}v'_{i}))$ et on pose $y=kx$. On v\'erifie encore que $k\in \underline{K}'\cap a\underline{K}'a^{-1}$ et que $y$ a les m\^emes coordonn\'ees que $x$, sauf $y_{i}$ qui v\'erifie $val_{F}(y_{i})=b_{i}+val_{F}(2)$. On conclut comme pr\'ec\'edemment. Cela prouve (7).
  
  Notons $X''$ l'ensemble des $x\in X'$ de la forme $x=e+\sum_{i=1,...,n}x_{i}v_{i}$ si $\iota=1$, $x=x_{1}^{-1}v_{-1}+\sum_{i=1,...,n}x_{i}v_{i}$ si $\iota=2$, o\`u $e\in E$ dans le premier cas et o\`u, dans les deux cas et pour tout $i=1,...,n$, $x_{i}=\varpi_{F}^{c_{i}}$, avec $0\leq c_{i}\leq(n+1)l$. C'est un ensemble fini.
  Pour tout $x\in X''$, fixons $\gamma_{x}\in G(F)$ tel que $\gamma_{x}^{-1}v_{0}=x$. Posons $\Gamma=\{\gamma_{x}; x\in X''\}$.  Pour ${\bf b}=(b_{1},...,b_{n})\in {\cal B}$, introduisons l'\'el\'ement $a({\bf b})\in A_{min}(F)$ tel que $a({\bf b})_{i}=\varpi_{F}^{b_{i}}$ pour tout $i=1,...,n$.  Un \'el\'ement $y$ v\'erifiant la conclusion de l'assertion (7) est de la forme $a({\bf b})x''$ pour un $x''\in X''$, autrement dit de la forme $a({\bf b})\gamma^{-1}v_{0}$ pour un $\gamma\in \Gamma$. On a donc
  
  (8) pour tout $x\in X({\bf b})$, il existe $k'\in \underline{K}'\cap a\underline{K}'a^{-1}$ et $\gamma\in \Gamma$ tels que $k'x=a({\bf b})\gamma^{-1} v_{0}$.

   Pour tout $x\in X''$, notons $H^x$ le sous-groupe des \'el\'ements de $G$ qui fixent $x$.   Introduisons le sous-groupe parabolique $P_{min}=M_{min}U_{min}\in {\cal P}(M_{min})$ form\'e des \'el\'ements qui conservent le drapeau de sous-espaces
  $$Fv_{n}\subset Fv_{n}\oplus Fv_{n-1}\subset ...\subset Fv_{n}\oplus...\oplus Fv_{1}.$$
  Introduisons aussi le sous-groupe parabolique $P'=M'U'\in {\cal F}(M_{min})$ form\'e des \'el\'ements qui conservent le drapeau des $n+1-\iota$ premiers sous-espaces ci-dessus. On a $P'=P_{min}$ si $\iota=1$ ou si $V_{an}=\{0\}$.  Montrons que
  
  (9) pour tout $x\in X''$, l'application
  $$\begin{array}{ccc}(P'(F)\cap H^x(F))\times \bar{P}_{min}(F)&\to&G(F)\\ (p',\bar{p})&\mapsto&p'\bar{p}\\ \end{array}$$
  est submersive \`a l'origine.
  
  Il suffit de prouver l'\'egalit\'e
  $$(\mathfrak{p}'(F)\cap \mathfrak{h}^x(F))+\bar{\mathfrak{p}}_{min}(F)=\mathfrak{g}(F).$$
  Notons $\mathfrak{h}^{x,\perp}$ l'orthogonal de $\mathfrak{h}^x$ dans $\mathfrak{g}$ pour la forme $(X,Y)\mapsto trace(XY)/2$ et notons $\bar{U}_{min}$ le radical unipotent de $\bar{P}_{min}$. Il suffit encore de prouver l'\'egalit\'e
  $$(\mathfrak{u}'(F)+\mathfrak{h}^{x,\perp}(F))\cap \bar{\mathfrak{u}}_{min}(F)=\{0\}.$$
  Notons $W^x$ l'orthogonal de $x$ dans $V$. Tout \'el\'ement de $\mathfrak{h}^{x,\perp}(F)$ est de la forme $c(x,y)$, pour un $y\in W^x$. On doit donc prouver que, pour $y\in W^x$, $N'\in \mathfrak{u}'(F)$, $\bar{N}\in \bar{\mathfrak{u}}_{min}(F)$, l'\'egalit\'e $N'+c(x,y)=\bar{N}$ entra\^{\i}ne $\bar{N}=0$. Consid\'erons de tels \'el\'ements. Soit $i=\iota,...,n$. On a $c(x,y)v_{i}=q_{V}(v_{i},x)y-q_{V}(v_{i},y)x=-q_{V}(v_{i},y)x$. On a $q_{V}(N'v_{i},v_{-i})=q_{V}(\bar{N}v_{i},v_{-i})=0$. Donc $q_{V}(c(x,y)v_{i},v_{-i})=0$, c'est-\`a-dire $-x_{i}q_{V}(v_{i},y)=0$. Puisque $x_{i}\not=0$, on a $q_{V}(v_{i},y)=0$, donc $c(x,y)v_{i}=0$. Alors $N'v_{i}=\bar{N}v_{i}$. Ces \'el\'ements appartiennent \`a des sous-espaces de $V$ d'intersection nulle, donc $\bar{N}v_{i}=0$. Si $\iota=1$, ou si $V_{an}=\{0\}$, cela suffit pour conclure: un \'el\'ement de $\bar{\mathfrak{u}}(F)$ qui annule tous les $v_{i}$, pour $i=\iota,...,n$, est nul. Supposons $\iota=2$ et $V_{an}\not=\{0\}$. Il faut montrer de plus que $ \bar{N}$ annule $V_{an}$. Soit $v_{an}\in V_{an}$. On a $c(x,y)v_{an}=q_{V}(v_{an},x)y-q_{V}(v_{an},y)x=-q_{V}(v_{an},y)x$. On a $q_{V}(N'v_{an},v_{-1})=q_{V}(\bar{N}v_{an},v_{-1})=0$, donc $q_{V}(c(x,y)v_{an},v_{-1})=0$, c'est-\`a-dire $-x_{1}q_{V}(v_{an},y)=0$. Puisque $x_{1}\not=0$, on a $q_{V}(v_{an},y)=0$, donc $c(x,y)v_{an}=0$. Alors $N'v_{an}=\bar{N}v_{an}$. Ces \'el\'ements appartiennent encore \`a des sous-espaces de $V$ d'intersection nulle, donc $\bar{N}v_{an}=0$. Cela d\'emontre (9).
  
    Apr\`es ces pr\'eparatifs, passons \`a la majoration de $I_{\underline{K}}(\epsilon,e,a)$. Pour tout $k\in \underline{K}$, on a $k^{-1}v_{0}\in X$. Pour ${\bf b}\in {\cal B}$, posons $\underline{K}({\bf b})=\{k\in \underline{K}; k^{-1}v_{0}\in X({\bf b})\}$. On  a l'\'egalit\'e
  $$I_{\underline{K}}(\epsilon,e,a)=\sum_{{\bf b}\in {\cal B}}I_{\underline{K}({\bf b})}(\epsilon,e,a),$$
  o\`u
  $$I_{\underline{K}({\bf b})}(\epsilon,e,a)=\int_{\underline{K}({\bf b})}\vert (\epsilon,l(\pi(ka)e))\vert  dk.$$
  Fixons ${\bf b}=(b_{1},...,b_{n})\in {\cal B}$.  Soit $k\in \underline{K}({\bf b})$, appliquons (8) \`a $x=k^{-1}v_{0}$. Il y a $k'\in \underline{K}'\cap a\underline{K}'a^{-1}$ et $\gamma\in \Gamma$ tels que $k'k^{-1}v_{0}=a({\bf b})\gamma^{-1}v_{0}$. Fixons de tels \'el\'ements et posons $h=kk^{_{'}-1}a({\bf b})\gamma^{-1}$. Alors $hv_{0}=v_{0}$, c'est-\`a-dire $h\in H(F)$. On a $k=h\gamma a({\bf b})^{-1}k'$, donc $\pi(ka)e=\pi(h\gamma a({\bf b})^{-1}a)e'$, o\`u $e'=\pi(a^{-1}k'a)e$. Notons $x$ l'\'el\'ement de $X''$ tel que $\gamma=\gamma_{x}$. Montrons que
   
   (10) il existe des sous-groupes ouverts compacts $\underline{K}_{1}\subset P'(F)\cap H^x(F)$ et $\underline{K}_{2}\subset \bar{P}(F)$, ind\'ependants des variables $a$, ${\bf b}$ et $k$, tels que la fonction  
   $$(k_{1},k_{2})\mapsto  (\epsilon,l(\pi(h\gamma k_{1}k_{2}a({\bf b})^{-1}a)e'))$$
   soit constante sur $\underline{K}_{1}\times \underline{K}_{2}$.
   
   Puisque $a^{-1}k'a\in \underline{K}'$, le vecteur $e'$ appartient \`a un ensemble fini ind\'ependant des variables.
  D'autre part, par d\'efinition de ${\cal B}'$, on a $a({\bf b})^{-1}a\in A_{min}(F)^-$. La conjugaison par $a({\bf b})a^{-1}$ contracte $\bar{P}$. Il existe donc $\underline{K}_{2}$ comme ci-dessus tel que $\pi(k_{2}a({\bf b})^{-1}a)e'=\pi(a({\bf b})^{-1}a)e'$ pour tout $k_{2}\in \underline{K}_{2}$. Consid\'erons l'application $k_{1}\mapsto h\gamma k_{1}\gamma^{-1}h^{-1}$ sur $P'(F)\cap H^x(F)$. Elle prend ses valeurs dans $H(F)$ car $\gamma H^x \gamma^{-1}=H$. D'autre part, elle est "born\'ee" en un sens facile \`a pr\'eciser. En effet, on a $h\gamma=kk^{_{'}-1}a({\bf b})$; les racines de $a({\bf b})$ dans $\mathfrak{p}'$ sont born\'ees par d\'efinition de ${\cal B}$ et $k$ et $k'$ appartiennent \`a des compacts. Il existe donc $\underline{K}_{1}$ comme dans l'\'enonc\'e tel que $\rho(h\gamma k_{1}^{-1}\gamma^{-1}h^{-1})\epsilon=\epsilon$ pour tout $k_{1}\in \underline{K}_{1}$. D'o\`u (10).
  
  Gr\^ace \`a (9) et (10), et \`a la finitude de $\Gamma$, il existe un sous-groupe ouvert compact $\underline{K}''$ de $G(F)$, ind\'ependant des variables, tel que la fonction
  $$k''\mapsto (\epsilon,l(\pi(h\gamma k''a({\bf b})^{-1}a)e'))$$
  soit constante sur $\underline{K}''$. Fixons un tel $\underline{K}''$. On a alors
  $$(\epsilon,l(\pi(ka)e))=(\epsilon,l(\pi(h\gamma a({\bf b})^{-1}a)e'))=(\epsilon,l(\pi(h\gamma)e'')),$$
  o\`u
  $$e''=mes(\underline{K}'')^{-1}\int_{\underline{K}''}\pi(k''a({\bf b})^{-1}a)e'\,dk''.$$
   Fixons une base orthonorm\'ee $(e_{j})_{j=1,...,k}$ du sous-espace $E_{\pi}^{\underline{K}''}$. On a
  $$e''=\sum_{j=1,...,k}(e_{j},e'')e_{j},$$
d'o\`u
  $$(\epsilon,l(\pi( ka)e))=\sum_{j=1,...,k}(\epsilon,l(\pi(h\gamma)e_{j}))(e_{j},e'').$$
  Pour tout $j$, on a $(e_{j},e'')=(e_{j},\pi(a({\bf b})^{-1}a(a^{-1}k'a))e)$. Rappelons que $a^{-1}k'a\in \underline{K}'$. D'o\`u une majoration
  $$\vert (e_{j},e'')\vert <<\Xi^G(a({\bf b})^{-1}a).$$
  On a aussi
  $$\vert (\epsilon,l(\pi(h\gamma)e_{j}))\vert=\vert (\rho(h^{-1})\epsilon,l(\pi(\gamma)e_{j}))\vert <<\Xi^H(h).$$ 
  Rappelons que $h=kk^{_{'}-1}a({\bf b})\gamma^{-1}$ et que l'on a not\'e $x$ l'\'el\'ement de $X''$ tel que $\gamma=\gamma_{x}$. Ecrivons $x=y+\sum_{i=\iota,...,n}x_{i}v_{i}$ avec $y\in E$ si $\iota=1$, $y=x_{1}^{-1}v_{-1}+x_{1}v_{1}$ si $\iota=2$. On a $x=exp(c(y,(y-x)/2))y$. Posons $\gamma_{0}=\gamma exp(c(y,(x-y)/2))$. Alors $\gamma_{0}^{-1}v_{0}=y$. Notons $H^y$ le sous-groupe des \'el\'ements de $G$ qui fixent $y$. On a $\gamma_{0}^{-1}H\gamma_{0}=H^y$ et une majoration
  $$\Xi^H(h)<<\Xi^{H^y}(\gamma_{0}^{-1}h\gamma_{0}) .$$
  On a $\gamma_{0}^{-1}h\gamma_{0}=\gamma_{0}^{-1}kk^{_{'}-1}a({\bf b})exp(c(y,(x-y)/2))$. Introduisons l'\'el\'ement $a({\bf b})'\in A_{min}(F)$ tel que
  
  - si $\iota=1$, $a({\bf b})'_{i}=\varpi_{F}^{b_{i}+nl}$ pour $i=1,...,n$;
  
  - si $\iota=2$, $a({\bf b})'_{i}=\varpi_{F}^{b_{i}-b_{1}}$ pour $i=1,...,n$.
  
  On a
  $$\gamma_{0}^{-1}h\gamma_{0}= k_{1}a({\bf b})''u'a({\bf b})',$$
  o\`u $k_{1}=\gamma_{0}kk^{_{'}-1}$, $a({\bf b})''=a({\bf b})a({\bf b})^{_{'}-1}$, $u'=a({\bf b})'exp(c(y,(x-y)/2))a({\bf b})^{_{'}-1}$. Remarquons que $a({\bf b})'$ appartient \`a $H^y(F)$, donc aussi $k_{1}a({\bf b})''u'\in H^y(F)$. L'\'el\'ement $k_{1}$ reste dans un compact. D'apr\`es la d\'efinition de ${\cal B}$, l'\'el\'ement $a({\bf b})''$ reste lui-aussi dans un compact. Enfin $c(y,(x-y)/2)$ appartient \`a $\mathfrak{u}'(F)$ et la conjugaison par $a({\bf b})'$ contracte cet ensemble. Donc $u'$ reste dans un compact. On en d\'eduit une majoration
  $$\Xi^{H^y}(\gamma_{0}^{-1}h\gamma_{0})<<\Xi^{H^y}(a({\bf b})'),$$
  puis
   $$\vert (\epsilon,l(\pi( ka)e))\vert <<\Xi^{H^y}(a({\bf b})')\Xi^G(a({\bf b})^{-1}a).$$
   Les \'el\'ements $a({\bf b})^{-1}a$ et $a({\bf b})'$ sont "n\'egatifs" pour $P_{min}$, resp. $P'\cap H^y$. D'apr\`es [W2] lemme II.1.1, il existe un r\'eel $D$ tel que le membre de droite ci-dessus soit essentiellement born\'e par
   $$\delta_{P'\cap H^y}(a({\bf b})')^{1/2}\delta_{P_{min}}(a({\bf b})^{-1}a)^{1/2}\sigma(a({\bf b})')^D\sigma(a({\bf b})a)^D.$$
   Les coefficients de $a({\bf b})'$ sont essentiellement les m\^emes que ceux de $a({\bf b})$. En calculant explicitement l'expression ci-dessus, on obtient   la majoration
   $$\vert  (\epsilon,l(\pi(ka)e))\vert <<q^{\sum_{i=1,...,n}b_{i}/2}\sigma(a({\bf b}))^D\Xi^G(a)\sigma(a)^D.$$
   L'application $k\mapsto k^{-1}v_{0}$ de $(\underline{K}\cap H(F))\backslash \underline{K}({\bf b})$ dans $X({\bf b})$ est injective et pr\'eserve les mesures. D'apr\`es (4), on a donc
   $$mes(\underline{K}({\bf b}))<<q^{-\sum_{i=1,...,n}b_{i}},$$
   puis
   $$I_{\underline{K}({\bf b})}(\epsilon,e,a)<<q^{-\sum_{i=1,...,n}b_{i}/2}\sigma(a({\bf b}))^D\Xi^G(a)\sigma(a)^D,$$
   et enfin
   $$I_{\underline{K}}(\epsilon,e,a)<<\Xi^G(a)\sigma(a)^D\sum_{{\bf b}\in {\cal B}}q^{-\sum_{i=1,...,n}b_{i}/2}\sigma(a({\bf b}))^D.$$
   L'ensemble ${\cal B}$ d\'epend de $a$ mais est contenu dans l'ensemble des $(b_{1},...,b_{n})\in {\mathbb Z}^n$ tels que $b_{i}\geq-nl$ pour tout $i$. On peut remplacer ${\cal B}$ par cet ensemble, la s\'erie ci-dessus y est convergente et on obtient   la majoration
$$I_{\underline{K}}(\epsilon,e,a)<<\Xi^G(a)\sigma(a)^D$$
que l'on voulait d\'emontrer. Cela ach\`eve la preuve. $\square$

\bigskip

\subsection{Le cas $r=0$: tout entrelacement est temp\'er\'e}

Soient $(V,q_{V})$ et $(W,q_{W})$ deux espaces quadratiques compatibles.  Soient $\pi\in Temp(G)$ et   $\rho\in Temp(H)$.

\ass{Proposition}{Supposons $d_{V}=d_{W}+1$. Alors $m(\pi,\rho)=1$ si et seulement si ${\cal L}_{\pi,\rho}$ n'est pas nulle.}

Preuve. Pour une raison qui va appara\^{\i}tre, modifions la notation en notant $\pi'$ plut\^ot que $\pi$ la repr\'esentation de $G(F)$. Un sens de l'\'equivalence ("si") est clair d'apr\`es 5.1. On doit prouver que, si $Hom_{H,\xi}(\pi',\rho)$ n'est pas nul, ${\cal L}_{\pi',\rho}$ ne l'est pas non plus.  Puisque $\pi'$ est temp\'er\'ee, on peut fixer des donn\'ees comme en 5.3, avec de plus $\tilde{\pi}$, $\mu_{1}$,....,$\mu_{s}$ de la s\'erie discr\`ete, de sorte que $\pi'$ soit une sous-repr\'esentation de la repr\'esentation induite $\pi=\pi_{0}$. On peut supposer que $E_{\pi'}\subset {\cal K}_{Q,\tau}^G$. Soient $e,e'\in {\cal K}_{Q,\tau}^G$ et $\varphi$ une fonction $C^{\infty}$ sur $i{\cal A}_{L,F}^*$. Comme en 1.6, d\'efinissons une fonction $f=f_{e,e',\varphi}$ sur $G(F)$ par
$$f(g)=\int_{i{\cal A}_{L,F}^*}\varphi(\lambda)(\pi_{\lambda}(g)e',e)m(\tau_{\lambda})d\lambda.$$
 Elle appartient \`a l'espace de Schwartz-Harish-Chandra ${\cal S}(G(F))$ et agit donc dans $\pi'$. Par d\'efinition de cette action, on a
$$(e_{0}',\pi'(f)e_{0})=\int_{G(F)}(e_{0}',\pi'(g)e_{0})f(g)dg$$
pour tous $e_{0},e_{0}'\in E_{\pi'}$.  Soit $l\in Hom_{H,\xi}(\pi',\rho)$, supposons $l\not=0$. Soient $e_{0}\in E_{\pi'}$ et $\epsilon\in E_{\rho}$. Posons
$$I(\epsilon,e_{0},f)=\int_{G(F)}(\epsilon,l(\pi(g)e_{0}))f(g)dg.$$
Gr\^ace \`a la proposition 5.5, cette int\'egrale est absolument convergente. On peut la calculer de deux fa\c{c}ons. La premi\`ere consiste \`a fixer un sous-groupe ouvert compact $K_{f}$ de $G(F)$ tel que $f$ soit biinvariante par $K_{f}$ et une suite $(\Omega_{n})_{n\geq1}$ de sous-ensembles ouverts compacts de $G(F)$ biinvariants par $K_{f}$, telle que
$$\Omega_{n}\subset \Omega_{n+1},\text{ et  }\bigcup_{n\geq1}\Omega_{n}=G(F).$$
Alors
$$I(\epsilon,e_{0},f)=lim_{n\to \infty}\int_{\Omega_{n}}(\epsilon,l(\pi'(g)e_{0}))f(g)dg.$$
Consid\'erons cette derni\`ere int\'egrale. Puisqu'elle est limit\'ee \`a un compact, on a
$$\int_{\Omega_{n}}(\epsilon,l(\pi'(g)e_{0}))f(g)dg=(\epsilon,l(e_{n})),$$
o\`u
$$e_{n}=\int_{\Omega_{n}}\pi'(g)f(g)dg.$$
Ces vecteurs restent dans le sous-espace de dimension finie $E_{\pi'}^{K_{f}}$. De plus, dans ce sous-espace, $lim_{n\to \infty}e_{n}=\pi'(f)e_{0}$. On en d\'eduit
$$lim_{n\to \infty}(\epsilon,l(e_{n}))=(\epsilon,l(\pi'(f)e_{0})),$$
puis
$$(1) \qquad I(\epsilon,e_{0},f)=(\epsilon,l(\pi'(f)e_{0})).$$
D'autre part, on a
$$I(\epsilon,e_{0},f)=\int_{H(F)\backslash G(F)}I(\epsilon,e_{0},f,g)dg,$$
o\`u
$$I(\epsilon,e_{0},f,g)=\int_{H(F)}(\epsilon,l(\pi'(hg)e_{0}))f(hg)dh$$
$$=\int_{H(F)}(\epsilon,l(\pi'(hg)e_{0}))\int_{i{\cal A}_{L,F}^*}\varphi(\lambda)(\pi_{\lambda}(hg)e',e)m(\tau_{\lambda})d\lambda \,dh.$$
Fixons $g$. La derni\`ere int\'egrale est absolument convergente. En effet, quand on remplace tous les termes par leurs valeurs absolues, il existe un entier $D$ tel que l'int\'egrale int\'erieure soit   $<<\Xi^G(h)\sigma(h)^D$. Le premier terme est $<<\Xi^H(h)$ et l'assertion r\'esulte de 4.3(4). On permute les deux int\'egrales, on utilise l'\'egalit\'e $(\epsilon,l(\pi'(hg)e_{0}))=(\rho(h)^{-1}\epsilon,l(\pi(g)e_{0}))$ et on change $h$ en $h^{-1}$. On obtient
$$I(\epsilon,e_{0},f,g)=\int_{i{\cal A}_{L,F}^*}\varphi(\lambda)m(\tau_{\lambda})\int_{H(F)}(\rho(h)\epsilon,l(\pi'(g)e_{0}))(\pi_{\lambda}(g)e',\pi_{\lambda}(h)e)dh\,d\lambda.$$
On reconna\^{\i}t  l'int\'egrale int\'erieure: c'est ${\cal L}_{\pi_{\lambda},\rho}(\epsilon\otimes \pi_{\lambda}(g)e',l(\pi'(g)e_{0})\otimes e)$. On obtient
$$(2)\qquad I(\epsilon,e_{0},f)=\int_{H(F)\backslash G(F)}\int_{i{\cal A}_{L,F}^*}\varphi(\lambda)m(\tau_{\lambda}){\cal L}_{\pi_{\lambda},\rho}(\epsilon\otimes \pi_{\lambda}(g)e',l(\pi'(g)e_{0})\otimes e)d\lambda\,dg.$$

Fixons un voisinage $\omega$ de $0$ dans $i{\cal A}_{L,F}^*$ tel que la relation 1.6(1) soit v\'erifi\'ee. Fixons $\epsilon\in E_{\rho}$ et $e_{0}\in E_{\pi'}$ tels que $(\epsilon,l(e_{0}))\not=0$. Appliquons les constructions ci-dessus \`a $e=e'=e_{0}$ et \`a une fonction $\varphi$ \`a support dans $\omega$ telle que $\varphi(0)\not=0$. La relation 1.6(1) nous dit que $\pi'(f)(e_{0})$ est un multiple non nul de $e_{0}$. D'apr\`es (1), on a donc $I(\epsilon,e_{0},f)\not=0$.  Alors (2) implique qu'il existe $\lambda$ tel que ${\cal L}_{\pi_{\lambda},\rho}$ ne soit pas nul. D'apr\`es le lemme 5.3(ii), ${\cal L}_{\pi,\rho}$ n'est pas nul. Si $\pi$ est irr\'eductible, on a $\pi'=\pi$ et c'est termin\'e. Sinon, d'apr\`es le lemme 5.4, il y a une sous-repr\'esentation irr\'eductible $\pi''$ de $\pi$ telle que ${\cal L}_{\pi'',\rho}$ ne soit pas nul. D'apr\`es  5.1(2) on peut fixer $e_{1}\in E_{\pi''}\subset {\cal K}_{Q,\tau}^G$ et $\epsilon_{1}\in E_{\rho}$ de sorte que ${\cal L}_{\pi,\rho}(\epsilon_{1}\otimes e_{1},\epsilon_{1}\otimes e_{1})\not=0$. Notons $c$ la valeur non nulle de ce terme.  D'apr\`es le lemme 5.3(i),  quitte \`a restreindre $\omega$, on peut supposer que $\vert {\cal L}_{\pi_{\lambda},\rho}(\epsilon_{1}\otimes e_{1},\epsilon_{1}\otimes e_{1})\vert \geq\vert c\vert /2$ pour $\lambda\in \omega$.  Soit $\varphi'$ la fonction \`a support dans $\omega$ telle que, pour $\lambda\in \omega$,
$$\varphi'(\lambda)=\varphi(\lambda){\cal L}_{\pi_{\lambda},\rho}(\epsilon\otimes e_{1},\epsilon_{1}\otimes e_{0}){\cal L}_{\pi_{\lambda},\rho}(\epsilon_{1}\otimes e_{1},\epsilon_{1}\otimes e_{1})^{-1}.$$
C'est une fonction $C^{\infty}$. Posons $f'=f_{e_{1},e_{0},\varphi'}$. Montrons que l'on a l'\'egalit\'e
$$(3) \qquad I(\epsilon,e_{0},f) =I(\epsilon_{1},e_{0},f').$$
D'apr\`es (2), il suffit de prouver que, pour tous $\lambda$ et $g$, on a l'\'egalit\'e
$$\varphi(\lambda){\cal L}_{\pi_{\lambda},\rho}(\epsilon\otimes \pi_{\lambda}(g)e_{0},l(\pi'(g)e_{0})\otimes e_{0})=\varphi'(\lambda){\cal L}_{\pi_{\lambda},\rho}(\epsilon_{1}\otimes\pi_{\lambda}(g)e_{0},l(\pi'(g)e_{0})\otimes e_{1}).$$
Tous les termes \'etant $C^{\infty}$, on peut supposer $\lambda$ en position g\'en\'erale, donc $\pi_{\lambda}$ irr\'eductible. On peut aussi supposer $\lambda\in \omega$. Fixons un tel $\lambda$ et un \'el\'ement non nul $\underline{l}\in Hom_{H}(\pi_{\lambda},\rho)$. D'apr\`es 5.1, il existe $c\in{\mathbb  C}^{\times}$ tel que
$${\cal L}_{\pi_{\lambda},\rho}(\underline{\epsilon}'\otimes \underline{e}',\underline{\epsilon}\otimes \underline{e})=c(\underline{\epsilon}',\underline{l}(\underline{e}))(\underline{l}(\underline{e}'),\underline{\epsilon})$$
pour tous $\underline{e},\underline{e}'\in {\cal K}_{Q,\tau}^G$ et $\underline{\epsilon},\underline{\epsilon}'\in E_{\rho}$. Les deux membres de l'\'egalit\'e \`a prouver valent
$$c\varphi(\lambda)(\epsilon,\underline{l}(e_{0}))(\underline{l}(\pi_{\lambda}(g)e_{0}),l(\pi'(g)e_{0})).$$
Cela prouve cette \'egalit\'e et (3). D'apr\`es (3), $I(\epsilon_{1},e_{0},f')\not=0$. D'apr\`es (1), $\pi'(f')e_{0}\not=0$. Alors, d'apr\`es 1.6(1), le produit scalaire $(e_{0},e_{1})$ n'est pas nul. Puisque $e_{0}\in E_{\pi'}$ et $e_{1}\in E_{\pi''}$, on a donc $\pi'=\pi''$ et  ${\cal L}_{\pi',\rho}$ est non nul par d\'efinition de $\pi''$. $\square$

\subsection{Tout entrelacement est temp\'er\'e}

Soient $(V,q_{V})$ et $(W,q_{W})$ deux espaces quadratiques compatibles. Soit $\pi\in Temp(G)$ et $\rho\in Temp(H)$.

\ass{Proposition}{On a $m(\pi,\rho)=1$ si et seulement si ${\cal L}_{\pi,\rho}$ est non nul.}

Preuve. On peut supposer $d_{V}>d_{W}$. Le cas $d_{V}=d_{W}+1$ est trait\'e par la proposition pr\'ec\'edente. Supposons $d_{V}\geq d_{W}+3$. Comme dans la preuve pr\'ec\'edente, on peut supposer $m(\pi,\rho)=1$ et on doit montrer que ${\cal L}_{\pi,\rho}$ n'est pas nul. Posons $k=(d_{V}-d_{W}+1)/2$. Soit $(Z',q_{Z'})$ un espace hyperbolique de dimension $2k$, notons $(V',q_{V'})$ la somme directe orthogonale de $W$ et $Z'$. Alors $(V',q_{V'})$ et $(V,q_{V})$ sont compatibles et $d_{V'}=d_{V}+1$. Le groupe sp\'ecial orthogonal $G'$ de $V'$ contient un groupe de L\'evi $L'$ isomorphe \`a $GL_{k}\times H$. Soient $P'\in {\cal P}(L')$ et $\gamma$ une repr\'esentation admissible irr\'eductible et cuspidale de $GL_{k}(F)$. Posons $\rho'=Ind_{P'}^{G'}(\gamma\otimes \rho)$. D'apr\`es le th\'eor\`eme 20.1 de [GGP], on peut choisir $\gamma$ de sorte que d'une part $\rho'$ soit irr\'eductible, d'autre part l'hypoth\`ese $m(\pi,\rho)=1$ entra\^{\i}ne $m(\rho',\pi)=1$. On fixe un tel $\gamma$. Gr\^ace \`a la proposition 5.6, la forme ${\cal L}_{\rho',\pi}$ n'est pas nulle. Gr\^ace au lemme 5.3(ii), la forme ${\cal L}_{\pi,\rho}$ est elle-aussi non nulle. C'est ce qu'il fallait prouver. $\square$

 \bigskip
 
\section{Expression spectrale de la limite d'une int\'egrale}

\bigskip

\subsection{Le th\'eor\`eme}

 Soient $(V,q_{V})$ et $(W,q_{W})$ deux espaces quadratiques compatibles. On suppose $d_{V}> d_{W}$. On utilise les constructions et notations de 4.2.  On fixe un L\'evi minimal $M_{min}$ de $G$ contenu dans $M$.  On suppose, ainsi qu'il est loisible, que le groupe $K$ est en bonne position relativement \`a $M_{min}$. On fixe des mesures de Haar sur $G(F)$ et $H(F)$. Les autres mesures que l'on utilisera sont normalis\'ees comme en 1.2. On fixe une repr\'esentation   $\rho\in Temp(H)$  et on note $\theta_{\rho}$ son caract\`ere.  Soit $f\in C_{c}^{\infty}(G(F))$. Pour $g\in G(F)$, on d\'efinit une fonction $^gf^{\xi}$ sur $H(F)$ par
 $$^gf^{\xi}(h)=\int_{U(F)}f(g^{-1}hug)\xi(u)du$$
 et on pose
$$I(\theta_{\rho},f,g)=\int_{H(F)}\theta_{\rho}(h){^gf}^{\xi}(h)dh.$$
Pour un entier $N\geq1$, on pose
$$I_{N}(\theta_{\rho},f)=\int_{U(F)H(F)\backslash G(F)}I(\theta_{\rho},f,g)\kappa_{N}(g)dg.$$

Pour $L\in {\cal L}(M_{min})$, notons $\Pi_{ell}(L)$ l'ensemble des classes d'isomorphie de repr\'esentations admissibles irr\'eductibles temp\'er\'ees et elliptiques de $L(F)$. Cet ensemble se d\'ecompose en orbites pour l'action $\pi\mapsto \pi_{\lambda}$ de $i{\cal A}_{L}^*$. On note $\{\Pi_{ell}(L)\}$ l'ensemble des orbites. Soient ${\cal O}$ une telle orbite et $\pi\in {\cal O}$.   Ecrivons
$$L=GL_{k_{1}}\times...\times GL_{k_{s}}\times \tilde{G},$$
$$\pi=\mu_{1}\otimes...\otimes \mu_{s}\otimes \tilde{\pi}.$$
La repr\'esentation $\tilde{\pi}$ ne d\'epend pas du choix de $\pi$ dans ${\cal O}$. D'autre part, $\tilde{G}$ est le groupe sp\'ecial orthogonal d'un sous-espace $\tilde{V}$ de $V$ qui est compatible \`a $W$.  On d\'efinit comme en 4.1 et 5.1 les nombres $t(\tilde{\pi})$ et $m(\tilde{\pi},\rho)$. On  pose $t(\pi)=t(\tilde{\pi})$ et $m({\cal O},\rho)=m(\pi,\rho)=m(\tilde{\pi},\rho)$. Notons aussi $i{\cal A}_{{\cal O}}^{\vee}$ le groupe des $\lambda\in i{\cal A}_{L}^*$ tels que, pour tout $\pi\in {\cal O}$, $\pi_{\lambda}$ soit \'equivalente \`a $\pi$.    On pose
$$I_{spec}(\theta_{\rho},f)=\sum_{L\in {\cal L}(M_{min})}\vert W^L\vert \vert W^G\vert ^{-1}(-1)^{a_{L}}\sum_{{\cal O}\in \{\Pi_{ell}(L)\}; m({\cal O},\rho)=1} $$
$$[i{\cal A}_{{\cal O}}^{\vee}:i{\cal A}_{L,F}^{\vee}]^{-1}t(\pi)^{-1}\int_{i{\cal A}_{L,F}^*}J_{L}^G(\pi_{\lambda},f)d\lambda,$$
o\`u, pour toute orbite ${\cal O}$, on a fix\'e un \'el\'ement $\pi\in {\cal O}$.  La fonction $\lambda\mapsto J_{L}^G(\pi_{\lambda},f)$ est $C^{\infty}$. Si l'on fixe un sous-groupe ouvert compact $K_{f}$ de $G(F)$ tel que $f$ soit biinvariante par $K_{f}$, elle n'est non nulle que si $\pi$ admet des invariants non nuls par $K_{f}$. Il n'y a qu'un nombre fini d'orbites ${\cal O}$ v\'erifiant cette condition.  L'expression ci-dessus est donc absolument convergente.

\ass{Th\'eor\`eme}{Soit $f\in C_{c}^{\infty}(G(F))$. Si $f$ est tr\`es cuspidale, on a l'\'egalit\'e
$$lim_{N\to \infty}I_{N}(\theta_{\rho},f)=I_{spec}(\theta_{\rho},f).$$}

Toute la section est consacr\'ee \`a la preuve de ce th\'eor\`eme. Pour toute cette section, on fixe une fonction $f\in C_{c}^{\infty}(G(F))$, que l'on suppose tr\`es cuspidale. On fixe un sous-groupe ouvert compact $K_{f}$ de $K$, distingu\'e dans $K$, tel que $f$ soit biinvariante par $K_{f}$.

\bigskip

\subsection{Utilisation de la formule de Plancherel}

Exprimons $f$ \`a l'aide de la formule de Plancherel. Comme on l'a dit en 1.6, pour tout $L\in {\cal L}(M_{min})$, on peut fixer un sous-ensemble fini ${\Pi_2(L)}_{f}\subset {\Pi_2(L)}$ de sorte que pour tout $g\in G(F)$, 
$$f(g)=\sum_{L\in {\cal L}(M_{min})}\vert W^L\vert \vert W^G\vert ^{-1}\sum_{{\cal O}\in {\Pi_2(L)}_{f}} f_{{\cal O}}(g),$$
o\`u
$$f_{{\cal O}}(g)=[i{\cal A}_{{\cal O}}^{\vee}:i{\cal A}_{L,F}^{\vee}]^{-1}\int_{i{\cal A}_{L,F}^*}m(\tau_{\lambda})trace(Ind_{Q}^G(\tau_{\lambda},g^{-1})Ind_{Q}^G(\tau_{\lambda},f))\,d\lambda.$$
On a remplac\'e $M$ et $P$ par $L$ et $Q$ dans la formule de 1.6. Pour $g\in G(F)$ et $h\in H(F)$, on a donc
$${^gf}^{\xi}(h)=\int_{U(F)}\sum_{L\in {\cal L}(M_{min})}\vert W^L\vert \vert W^G\vert ^{-1}\sum_{{\cal O}\in {\Pi_2(L)}_{f}}  f_{{\cal O}}(g^{-1}hug)\xi(u)\,du.$$

On fixe un produit scalaire invariant sur $E_{\rho}$. Soit $N\geq1$. Le support de la fonction $\kappa_{N}$ restreinte \`a $M(F)$ est d'image compacte dans $H(F)\backslash M(F)$. On fixe un sous-ensemble ouvert compact $\Gamma_{N}$ de $M(F)$ tel que ce support soit contenu dans $H(F)\Gamma_{N}$. On fixe un sous-groupe ouvert compact $K_{N}$ de $H(F)$ tel que, pour tout $g\in \Gamma_{N}$, $K_{N}
$ soit contenu dans $gK_{f}g^{-1}$. Pour $g\in \Gamma_{N}K$, la fonction ${^gf}^{\xi}$ sur $H(F)$ est biinvariante par $K_{N}$. Fixons une base orthonorm\'ee ${\cal B}_{\rho}^{K_{N}}$ du sous-espace $E_{\rho}^{K_{N}}$ des \'el\'ements de $E_{\rho}$ invariants par $K_{N}$. Alors, pour $g\in \Gamma_{N} K$, on a l'\'egalit\'e
$$I(\theta_{\rho},f,g)=\sum_{\epsilon\in {\cal B}_{\rho}^{K_{N}}}\int_{H(F)}(\epsilon,\rho(h)\epsilon){^gf}^{\xi}(h)dh,$$
d'o\`u
$$I(\theta_{\rho},f,g)=\sum_{\epsilon\in {\cal B}_{\rho}^{K_{N}}}\int_{H(F)}(\epsilon,\rho(h)\epsilon)\int_{U(F)}\sum_{L\in {\cal L}(M_{min})}\vert W^L\vert \vert W^G\vert ^{-1}$$
$$\sum_{{\cal O}\in {\Pi_2(L)}_{f}}  f_{{\cal O}}(g^{-1}hug)\xi(u)\,du\,dh.$$
Pour $\epsilon\in E_{\rho}$, $L\in {\cal L}(M_{min})$, ${\cal O}\in  {\Pi_2(L)}_{f}$, et $g\in G(F)$, posons
$$(1) \qquad I_{L,{\cal O}}(\epsilon,f,g)=\int_{H(F)\times U(F)}(\epsilon,\rho(h)\epsilon)f_{{\cal O}}(g^{-1}hug)\xi(u)\,du\,dh.$$
On a

(2) cette expression est absolument convergente.

  D'apr\`es Harish-Chandra ([W2] proposition VI.3.1), la fonction $f_{{\cal O}}$ appartient \`a l'espace de Schwartz-Harish-Chandra ${\cal S}(G(F))$. D'apr\`es [W2] proposition II.4.5, pour tout entier $D$, on a une majoration
$$\int_{U(F)}\vert f_{{\cal O}}(g^{-1}hug)\vert du<<\delta_{P}(h)^{-1/2}\Xi^M(h)\sigma(h)^{-D}$$
pour tout $h\in H(F)$. Sur $H(F)$, le module $\delta_{P}$ est trivial et $\Xi^M$ co\"{\i}ncide avec $\Xi^{G_{0}}$  . D'autre part, on a une majoration
$$\vert (\epsilon,\rho(h)\epsilon)\vert <<\Xi^H(h)$$
pour tout $h\in H(F)$. Enfin, l'int\'egrale
$$\int_{H(F)}\Xi^H(h)\Xi^{G_{0}}(h)dh$$
est convergente d'apr\`es 4.3(4). Cela prouve (2).

Pour $g\in \Gamma_{N}K$, on a donc l'\'egalit\'e
$$(3) \qquad I(\theta_{\rho},f,g)=\sum_{\epsilon\in {\cal B}_{\rho}^{K_{N}}}\sum_{L\in {\cal L}(M_{min})}\vert W^L\vert \vert W^G\vert ^{-1}\sum_{{\cal O}\in {\Pi_2(L)}_{f}}I_{L,{\cal O}}(\epsilon,f,g).$$

\bigskip

\subsection{Apparition des entrelacements temp\'er\'es}

On poursuit le calcul pr\'ec\'edent. Fixons $L\in {\cal L}(M_{min})$ et ${\cal O}\in {\Pi_2(L)}_{f}$.   Pour $c\in {\mathbb N}$, introduisons le sous-groupe $U(F)_{c}$ de $U(F)$, cf. 4.3. 
On a

(1) il existe $c_{0}\in {\mathbb N}$ tel que, pour tout $c\geq c_{0}$, tout $g\in M(F)K$ et tout $h\in H(F)$, on ait l'\'egalit\'e
$$\int_{U(F)}f_{{\cal O}}(g^{-1}hug)\xi(u)\,du\,=\int_{U(F)_{c}}f_{{\cal O}}(g^{-1}hug)\xi(u)\,du.$$

La preuve est similaire \`a celle du lemme 3.5. Pour $c\geq c_{0}$ notons $U(F)_{c}-U(F)_{c_{0}}$ le compl\'ementaire de $U(F)_{c_{0}}$ dans $U(F)_{c}$. Il existe $c_{0}$ tel que pour tout $c\geq c_{0}$ et tout $u\in U(F)_{c}-U(F)_{c_{0}}$, l'int\'egrale
$$ \int_{A(F)\cap K_{f}}\xi(aua^{-1})da$$
soit nulle. Choisissons un tel $c_{0}$, soit $c\geq c_{0}$. Parce que $A$ commute \`a $M$ et \`a $H\subset M$, parce que $K_{f}$ est distingu\'e dans $K$ et parce que $f_{{\cal O}}$ est, comme $f$, biinvariante par $K_{f}$, on a l'\'egalit\'e
$$f_{{\cal O}}(g^{-1}hug)=f_{{\cal O}}(g^{-1}ha^{-1}uag)$$
pour tous $g\in M(F)K$, $h\in H(F)$, $a\in A(F)\cap K_{f}$. Alors
$$(2)\qquad \int_{U(F)_{c}-U(F)_{c_{0}}}f_{{\cal O}}(g^{-1}hug)\xi(u)\,du=(mes(A(F)\cap K_{f}))^{-1}$$
$$\int_{U(F)_{c}-U(F)_{c_{0}}}\int_{A(F)\cap K_{f}}f_{{\cal O}}(g^{-1}ha^{-1}uag)\xi(u)\,da\,du.$$
Cette expression est absolument convergente. On peut effectuer le changement de variable $u\mapsto aua^{-1}$ et l'expression ci-dessus devient
$$(mes(A(F)\cap K_{f}))^{-1}\int_{U(F)_{c}-U(F)_{c_{0}}}f_{{\cal O}}(g^{-1}hug)\int_{A(F)\cap K_{f}}\xi(aua^{-1})da\,du.$$
L'int\'egrale int\'erieure est nulle donc aussi le membre de gauche de (2). Cela suffit \`a prouver (1). 

Comme en 1.6, on fixe $\tau\in {\cal O}$ et un \'el\'ement $Q\in {\cal P}(L)$. Pour simplifier les notations, on pose
$$\pi_{\lambda}=Ind_{Q}^G(\tau_{\lambda})$$
pour tout $\lambda\in i{\cal A}_{L,F}^*$. On r\'ealise toutes les repr\'esentations $\pi_{\lambda}$ dans l'espace commun ${\cal K}^G_{Q,\tau}$.  Fixons une base orthonorm\'ee ${\cal B}_{{\cal O}}^{K_{f}}$ du sous-espace $({\cal K}^G_{Q,\tau})^{K_{f}}$. Pour tout $g\in G(F)$, on a l'\'egalit\'e
$$f_{{\cal O}}(g)=[i{\cal A}_{{\cal O}}^{\vee}:i{\cal A}_{L,F}^{\vee}]^{-1}\sum_{e\in {\cal B}_{{\cal O}}^{K_{f}}}\int_{i{\cal A}_{L,F}^*}m(\tau_{\lambda})(e,\pi_{\lambda}(g^{-1})\pi_{\lambda}(f)e)\,d\lambda.$$
Soient $c_{0}$ v\'erifiant (1), $c\geq c_{0}$, $g\in M(F)K$ et $\epsilon\in E_{\rho}$. D'apr\`es (1) et la d\'efinition 6.2(1), on a l'\'egalit\'e
$$I_{L,{\cal O}}(\epsilon,f,g)=[i{\cal A}_{{\cal O}}^{\vee}:i{\cal A}_{L,F}^{\vee}]^{-1}\int_{H(F)\times U(F)_{c}}(\epsilon,\rho(h)\epsilon)$$
$$ \sum_{e\in {\cal B}_{{\cal O}}^{K_{f}}}\int_{i{\cal A}_{L,F}^*}m(\tau_{\lambda})(\pi_{\lambda}(g)e,\pi_{\lambda}((hu)^{-1}g)\pi_{\lambda}(f)e)\xi(u)\,d\lambda\,du\,dh.$$
En changeant $h$ et $u$ en leurs inverses, on obtient 
$$I_{L,{\cal O}}(\epsilon,f,g)=[i{\cal A}_{{\cal O}}^{\vee}:i{\cal A}_{L,F}^{\vee}]^{-1}\int_{H(F)\times U(F)_{c}}(\rho(h)\epsilon,\epsilon)$$
$$ \sum_{e\in {\cal B}_{{\cal O}}^{K_{f}}}\int_{i{\cal A}_{L,F}^*}m(\tau_{\lambda})(\pi_{\lambda}(g)e,\pi_{\lambda}(hug)\pi_{\lambda}(f)e)\bar{\xi}(u)\,d\lambda\,du\,dh.$$
Remarquons que, pour $g$ fix\'e, on a une majoration
$$\vert (\pi_{\lambda}(g)e,\pi_{\lambda}(hug)\pi_{\lambda}(f)e)\vert << \Xi^G(hu)$$
pour tous $\lambda$, $h$, $u$. Gr\^ace \`a 4.3(4), on en d\'eduit que l'expression ci-dessus est absolument convergente. On peut donc permuter les int\'egrales:
$$I_{L,{\cal O}}(\epsilon,f,g)=[i{\cal A}_{{\cal O}}^{\vee}:i{\cal A}_{L,F}^{\vee}]^{-1} \sum_{e\in {\cal B}_{{\cal O}}^{K_{f}}}\int_{i{\cal A}_{L,F}^*}m(\tau_{\lambda})$$
$$\int_{H(F)\times U(F)_{c}}(\rho(h)\epsilon,\epsilon))(\pi_{\lambda}(g)e,\pi_{\lambda}(hug)\pi_{\lambda}(f)e)\bar{\xi}(u)du\,dh\,d\lambda.$$
On reconna\^{\i}t l'int\'egrale int\'erieure: c'est ${\cal L}_{\pi_{\lambda},\rho,c}(\epsilon\otimes \pi_{\lambda}(g)e,\epsilon\otimes \pi_{\lambda}(g)\pi_{\lambda}(f)e)$. D'apr\`es le lemme 5.1, quitte \`a accro\^{\i}tre $c_{0}$ (en fait, la preuve de (1) montre que ce n'est pas n\'ecessaire), c'est aussi ${\cal L}_{\pi_{\lambda},\rho}(\epsilon\otimes \pi_{\lambda}(g)e,\epsilon\otimes \pi_{\lambda}(g)\pi_{\lambda}(f)e)$. On obtient
$$(3) \qquad I_{L,{\cal O}}(\epsilon,f,g)=[i{\cal A}_{{\cal O}}^{\vee}:i{\cal A}_{L,F}^{\vee}]^{-1} \sum_{e\in {\cal B}_{{\cal O}}^{K_{f}}}\int_{i{\cal A}_{L,F}^*}m(\tau_{\lambda})$$
$${\cal L}_{\pi_{\lambda},\rho}(\epsilon\otimes \pi_{\lambda}(g)e,\epsilon\otimes \pi_{\lambda}(g)\pi_{\lambda}(f)e)d\lambda$$
pour tout $g\in M(F)K$.

 On peut \'ecrire $L$ et $\tau$ comme en 5.3. C'est-\`a-dire que 
$$L=GL_{k_{1}}\times...\times GL_{k_{s}}\times \tilde{G}$$
et $\tau=\mu_{1}\otimes...\otimes \mu_{s}\otimes \tilde{\pi}$. Si $m(\tilde{\pi},\rho)=0$, on a aussi ${\cal L}_{\tilde{\pi},\rho}=0$ d'apr\`es la proposition 5.7 et ${\cal L}_{\pi_{\lambda},\rho}=0$ pour tout $\lambda$ d'apr\`es le lemme 5.3(ii). Posons $m({\cal O},\rho)=m(\tilde{\pi},\rho)$. Alors

(4) si $m({\cal O},\rho)=0$, $I_{L,{\cal O}}(\epsilon,f,g)=0$ pour tout $g\in M(F)K$ et tout $\epsilon\in E_{\rho}$.

Supposons d\'esormais $m({\cal O},\rho)=1$. On fixe des familles $(\epsilon'_{j})_{j=1,...,n}$, $(\epsilon_{j})_{j=1,...,n}$, $(e'_{j})_{j=1,...,n}$, $(e_{j})_{j=1,...,n}$, $(\varphi_{j})_{j=1,...,n}$ v\'erifiant le lemme 5.3(iii). Soit $N\geq1$. Pour $\lambda\in i{\cal A}_{L,F}^*$, $g\in M(F)K$ et $e\in {\cal K}^G_{Q,\tau}$, consid\'erons la somme
$$X_{\lambda}(e,g)=\sum_{\epsilon\in {\cal B}_{\rho}^{K_{N}}}{\cal L}_{\pi_{\lambda},\rho}(\epsilon\otimes \pi_{\lambda}(g)e,\epsilon\otimes \pi_{\lambda}(g)\pi_{\lambda}(f)e).$$
Supposons $\lambda$ en position g\'en\'erale. Alors $\pi_{\lambda}$ est irr\'eductible. Fixons un \'el\'ement non nul $l_{\lambda}\in Hom_{H,\xi}(\pi_{\lambda},\rho)$. Comme on l'a dit en 5.1, il existe un nombre complexe non nul $c_{\lambda}$ tel que 
$${\cal L}_{\pi_{\lambda},\rho}(\epsilon'\otimes e',\epsilon\otimes e)=c_{\lambda}(\epsilon',l_{\lambda}(e))(l_{\lambda}(e'),\epsilon)$$
pour tous $\epsilon,\epsilon'\in E_{\rho}$, $e,e'\in {\cal K}^G_{Q,\tau}$. La propri\'et\'e (iii) du lemme 5.3 s'\'ecrit
$$(5) \qquad \sum_{j=1,...,n}c_{\lambda}\varphi_{j}(\lambda)(\epsilon_{j}',l_{\lambda}(e_{j}))(l_{\lambda}(e_{j}'),\epsilon_{j})=1.$$
 On a
$$X_{\lambda}(e,g)=c_{\lambda}\sum_{\epsilon\in {\cal B}_{\rho}^{K_{N}}}(\epsilon,l_{\lambda}(\pi_{\lambda}(g)\pi_{\lambda}(f)e))(l_{\lambda}(\pi_{\lambda}(g)e),\epsilon).$$
 L'\'el\'ement $l_{\lambda}(\pi_{\lambda}(g)\pi_{\lambda}(f)e)$ est invariant par $K_{N}$ pour tout $g\in \Gamma_{N} K$ et tout $e\in {\cal K}^G_{Q,\tau}$. Supposons $g\in \Gamma_{N}K$. Alors
$$l_{\lambda}(\pi_{\lambda}(g)\pi_{\lambda}(f)e)=\sum_{\epsilon\in {\cal B}_{\rho}^{K_{N}}}(\epsilon,l_{\lambda}(\pi_{\lambda}(g)\pi_{\lambda}(f)e))\epsilon,$$
et
$$X_{\lambda}(e,g)=c_{\lambda}(l_{\lambda}(\pi_{\lambda}(g)e),l_{\lambda}(\pi_{\lambda}(g)\pi_{\lambda}(f)e)).$$
On peut multiplier $X_{\lambda}(e,g)$ par le membre de gauche de (5) et on obtient
$$(6) \qquad X_{\lambda}(e,g)=\sum_{j=1,...,n}\varphi_{j}(\lambda)X_{\lambda,j}(e,g),$$
o\`u, pour tout $g\in M(F)K$, on a pos\'e
$$X_{\lambda,j}(e,g)=c_{\lambda}^2(l_{\lambda}(\pi_{\lambda}(g)e),l_{\lambda}(\pi_{\lambda}(g)\pi_{\lambda}(f)e))(\epsilon_{j}',l_{\lambda}(e_{j}))(l_{\lambda}(e_{j}'),\epsilon_{j}).$$
Fixons $j$. Le produit de l'un des facteurs $c_{\lambda}$ et des deux premiers produits scalaires est \'egal \`a
$${\cal L}_{\pi_{\lambda},\rho}(\epsilon_{j}'\otimes \pi_{\lambda}(g)e,l_{\lambda}(\pi_{\lambda}(g)\pi_{\lambda}(f)e)\otimes e_{j}).$$
En choisissant un entier $c_{0}$ assez grand, on peut ici remplacer ${\cal L}_{\pi_{\lambda},\rho}$ par ${\cal L}_{\pi_{\lambda},\rho,c}$ pour tout $c\geq c_{0}$. Remarquons que $c_{0}$ est ind\'ependant de  de $\lambda$: cela r\'esulte de la preuve du lemme 3.5. Il est aussi ind\'ependant de $N$ qui n'intervient pas ici. On a donc
$$X_{\lambda,j}(e,g)=\int_{H(F)U(F)_{c}}(\rho(h)\epsilon'_{j},l_{\lambda}(\pi_{\lambda}(g)\pi_{\lambda}(f)e))(\pi_{\lambda}(g)e,\pi_{\lambda}(hu)e_{j})c_{\lambda}(l_{\lambda}(e_{j}'),\epsilon_{j})\bar{\xi}(u)du\,dh$$
 pourvu que $c\geq c_{0}$. Dans l'int\'egrale, le produit de $c_{\lambda}$ et des produits scalaires extr\^emes est \'egal \`a ${\cal L}_{\pi_{\lambda},\rho}(\rho(h)\epsilon'_{j}\otimes e'_{j},\epsilon_{j}\otimes \pi_{\lambda}(g)\pi_{\lambda}(f)e)$. On peut encore remplacer ${\cal L}_{\pi_{\lambda},\rho}$ par ${\cal L}_{\pi_{\lambda},\rho,c}$ pourvu que $c\geq c_{0}$. On obtient
 $$(7) \qquad X_{\lambda,j}(e,g)=\int_{H(F)U(F)_{c}}\int_{H(F)U(F)_{c}}(\pi_{\lambda}(g)e,\pi_{\lambda}(hu)e_{j})(\rho(h'h)\epsilon'_{j},\epsilon_{j})$$
 $$(e'_{j},\pi_{\lambda}(h'u'g)\pi_{\lambda}(f)e)\bar{\xi}(u')\bar{\xi}(u)du'\,dh'\,du\,dh.$$
 On a
 
 (8) pour $g$ fix\'e, cette expression est absolument convergente, uniform\'ement en $\lambda$.
 
 En effet, elle est major\'ee en valeur absolue par
 $$\int_{H(F)U(F)_{c}}\int_{H(F)U(F)_{c}}\Xi^G(hu)\Xi^H(h'h)\Xi^G(h'u')du'\,dh'\,du\,dh$$
 qui est convergente d'apr\`es 4.3(5).
 
 On peut maintenant lever l'hypoth\`ese que $\lambda$ est en position g\'en\'erale. Gr\^ace \`a (8), la formule (7) d\'efinit une fonction $C^{\infty}$ de $\lambda$ et l'\'egalit\'e (6) se prolonge par continuit\'e \`a tout $\lambda$. Pour deux entiers $c,c'\in {\mathbb N}$, et pour $g\in M(F)K$, posons
 $$X_{\lambda,j,c,c'}(e,g)=\int_{H(F)U(F)_{c}}\int_{H(F)U(F)_{c'}}(\rho(h)\epsilon'_{j},\epsilon_{j})(\pi_{\lambda}(h'u'g)e,\pi_{\lambda}(hu)e_{i})$$
 $$(e'_{j},\pi_{\lambda}(h'u'g)\pi_{\lambda}(f)e)\bar{\xi}(u)du'\,dh'\,du\,dh.$$
 Comme (7), cette expression est absolument convergente. On a
 
 (9) il existe $c_{0}$ ind\'ependant de $N$ et $\lambda$ tel que, si $ c\geq c_{0}$ et $c'\geq c_{0}$, alors $X_{\lambda,j}(e,g)=X_{\lambda,j,c,c'}(e,g)$ pour tout $g\in M(F)K$.
 
 Pour $a$ appartenant \`a un sous-groupe ouvert compact assez petit de $A(F)$, le changement de variables $u\mapsto aua^{-1}$, $u'\mapsto au'a^{-1}$ dans la d\'efinition ci-dessus de $X_{\lambda,i,c,c'}(e,g)$ revient \`a y remplacer $\bar{\xi}(u)$ par $\bar{\xi}(aua^{-1})$.  Comme dans la preuve de (1), on en d\'eduit que, si $c_{0}$ est assez grand, $X_{\lambda,j,c,c'}(e,g)$ ne d\'epend pas de $c$, pourvu que $c\geq c_{0}$. Pour $c,c'\geq c_{0}$, on a donc $X_{\lambda,j,c,c'}(e,g)=X_{\lambda,j,c',c'}$. On peut remplacer $c$ par $c'$ dans la formule (7). Dans cette formule, effectuons le changement de variables $h\mapsto h^{_{'}-1}h$, $u\mapsto h^{-1}h'u^{_{'}-1}h^{_{'}-1}hu$. Alors le membre de droite de (7) devient $X_{\lambda,j,c',c'}$. Cela prouve (9).
 
 En rassemblant l'\'egalit\'e (2), la d\'efinition de $X_{\lambda}(e,g)$, l'\'egalit\'e (6) et la propri\'et\'e (9), on obtient le r\'esultat suivant. Rappelons que l'on a suppos\'e $m({\cal O},\rho)=1$. Il existe $c_{0}$ tel que, pour tout $N$, tout $g\in \Gamma_{N}K$, tous $c,c'$ tels que $c\geq c_{0}$, $c'\geq c_{0}$, on a l'\'egalit\'e
 $$(10) \qquad  \sum_{\epsilon\in {\cal B}_{\rho}^{K_{N}}}I_{L,{\cal O}}(\epsilon,f,g)=[i{\cal A}_{{\cal O}}^{\vee}:i{\cal A}_{L,F}^{\vee}]^{-1} \sum_{e\in {\cal B}_{{\cal O}}^{K_{f}}}\sum_{j=1,...,n}$$
 $$\int_{i{\cal A}_{L,F}^*}m(\tau_{\lambda})\varphi_{j}(\lambda)X_{\lambda,j,c,c'}(e,g) \,d\lambda.$$
  On note $I_{L,{\cal O},N}(\theta_{\rho},f,g) $ le membre de droite de cette \'egalit\'e, qui est d\'efini pour tout $g\in M(F)K$.
 \bigskip
 
 \subsection{Une premi\`ere approximation}
 
 D'apr\`es 6.2(3), 6.3(4) et 6.3(10), on a l'\'egalit\'e
 $$I(\theta_{\rho},f,g)= \sum_{L\in {\cal L}(M_{min})}\vert W^L\vert \vert W^G\vert ^{-1}\sum_{{\cal O}\in {\Pi_2(L)}_{f}, m({\cal O},\rho)=1}I_{L,{\cal O},N}(\theta_{\rho},f,g),$$
 pour tout $g\in \Gamma_{N}K$. Comme fonctions de $g$, les deux membres sont d\'efinis  sur $ M(F)K$ et sont invariants \`a gauche par $H(F)$. Ils sont donc \'egaux pour $g\in H(F)\Gamma_{N}K$. L'intersection de $M(F)K$ avec le support de $\kappa_{N}$ \'etant contenu dans cet ensemble,  si l'on multiplie les deux membres de l'\'egalit\'e par $\kappa_{N}(g)$, l'\'egalit\'e obtenue devient vraie pour tout $g\in M(F)K$. Par d\'efinition, $I_{N}(\theta_{\rho},f)$ est l'int\'egrale de $I(\theta_{\rho},f,g)\kappa_{N}(g)$ sur $g\in H(F)U(F)\backslash G(F)$ ou, ce qui revient au m\^eme, l'int\'egrale de $I(\theta_{\rho},f,mk)\kappa_{N}(mk)\delta_{P}(m)^{-1}$ sur $m\in H(F)\backslash M(F)$ et $k\in K$. Donc
 $$(1) \qquad I_{N}(\theta_{\rho},f)=\int_{H(F)\backslash M(F)}\int_{K}\sum_{L\in {\cal L}(M_{min})}\vert W^L\vert \vert W^G\vert ^{-1}$$
 $$\sum_{{\cal O}\in {\Pi_2(L)}_{f}, m({\cal O},\rho)=1}I_{L,{\cal O},N}(f,mk)\kappa_{N}(mk)\delta_{P}(m)^{-1}dk\, dm.$$
  Soient $L\in {\cal L}(M_{min})$ et ${\cal O}\in {\Pi_2(L)}_{f}$  tel que $m({\cal O},\rho)=1$. On reprend les notations du paragraphe pr\'ec\'edent. On fixe $c_{0}$ v\'erifiant 6.3(9) et $c\geq c_{0}$.   Pour tout entier $C\in  {\mathbb N}$, posons
  $$I_{L,{\cal O},N,C}(\theta_{\rho},f)=[i{\cal A}_{{\cal O}}^{\vee}:i{\cal A}_{L,F}^{\vee}]^{-1} \sum_{e\in {\cal B}_{{\cal O}}^{K_{f}}}\sum_{j=1,...,n}\int_{i{\cal A}_{L,F}^*}m(\tau_{\lambda})\varphi_{j}(\lambda)$$
  $$\int_{H(F)U(F)_{c}}{\bf 1}_{\sigma< Clog(N)}(hu)(\rho(h)\epsilon'_{j},\epsilon_{j})\bar{\xi}(u)\int_{G(F)}(\pi_{\lambda}(g)e,\pi_{\lambda}(hu)e_{j})$$
 $$(e'_{j},\pi_{\lambda}(g)\pi_{\lambda}(f)e)\kappa_{N}(g)dg\,du\,dh\,d\lambda.$$
 
 \ass{Lemme}{(i) Cette expression est absolument convergente.
 
 (ii) Il existe $C$ tel que l'on ait la majoration
 $$\vert I_{N}(\theta_{\rho},f)-\sum_{L\in {\cal L}(M_{min})}\vert W^L\vert \vert W^G\vert ^{-1}\sum_{{\cal O}\in {\Pi_2(L)}_{f}, m({\cal O},\rho)=1}I_{L,{\cal O},N,C}(\theta_{\rho},f)\vert <<N^{-1}$$
 pour tout entier $N\geq 2$.}
 
 Preuve. Soient $L$ et ${\cal O}$ comme avant l'\'enonc\'e. Notons $ I _{L,{\cal O},N,\infty}(\theta_{\rho},f)$ l'expression obtenue en supprimant la fonction ${\bf 1}_{\sigma< Clog(N)}$ dans la d\'efinition de $I_{L,{\cal O},N,C}(\theta_{\rho},f)$. Pour $m\in M(F)$, $k\in K$, $u,u'\in U(F)$, $h\in H(F)$, $\lambda\in i{\cal A}_{L,F}^*$, posons
$$\Phi(m,k,u,u',h,\lambda)=[i{\cal A}_{{\cal O}}^{\vee}:i{\cal A}_{L,F}^{\vee}]^{-1} \sum_{e\in {\cal B}_{{\cal O}}^{K_{f}}}\sum_{j=1,...,n}m(\tau_{\lambda})\varphi_{j}(\lambda)(\rho(h)\epsilon'_{j},\epsilon_{j})\bar{\xi}(u)$$
$$(\pi_{\lambda}(u'mk)e,\pi_{\lambda}(hu)e_{i})(e'_{j},\pi_{\lambda}(u'mk)\pi_{\lambda}(f)e)\kappa_{N}(mk)\delta_{P}(m)^{-1}.$$
L'expression $I_{L,{\cal O},N,\infty}(\theta_{\rho},f)$ contient une int\'egrale sur $ G(F)$. On y d\'ecompose la variable $g\in G(F)$ en $g=u'mk$, avec $u'\in U(F)$, $m\in M(F)$, $k\in K$. On obtient formellement l'\'egalit\'e
$$I_{L,{\cal O},N,\infty}(\theta_{\rho},f)=\int_{i{\cal A}_{L,F}^*}\int_{M(F)}\int_{K}\int_{H(F)U(F)_{c}}\int_{U(F)}\Phi(m,k,u,u',h,\lambda)du'\,du\,dh\,dk\,dm\,d\lambda.$$
Cette \'egalit\'e est justifi\'ee par 

(2) cette expression est absolument convergente.

Les int\'egrales sur les ensembles compacts $K$ et $i{\cal A}_{L,F}^*$ sont insignifiantes.
On a
$$(3) \qquad \vert \Phi(m,k,u,u',h,\lambda)\vert<< \Xi^H(h)\Xi^G(u^{-1}h^{-1}u'm)\Xi^G(u'm)\kappa_{N}(m)\delta_{P}(m)^{-1}.$$
L'int\'egrale de cette fonction en $m$, $u$, $u'$, $h$ est convergente d'apr\`es 4.3(6). Cela d\'emontre (2).

En appliquant la m\^eme proc\'edure, on obtient
$$I_{L,{\cal O},N,C}(\theta_{\rho},f)=\int_{i{\cal A}_{L,F}^*}\int_{M(F)}\int_{K}\int_{H(F)U(F)_{c}}\int_{U(F)}{\bf 1}_{\sigma< Clog(N)}(hu)$$
$$\Phi(m,k,u,u',h,\lambda)du'\,du\,dh\,dk\,dm\,d\lambda,$$ 
$$\int_{H(F)\backslash M(F)}\int_{K}I_{L,{\cal O},N}(f,mk)\kappa_{N}(mk)\delta_{P}(m)^{-1}dk\, dm=\int_{i{\cal A}_{L,F}^*}\int_{M(F)}\int_{K}$$
$$\int_{H(F)U(F)_{c}}\int_{U(F)_{c'}}\Phi(m,k,u,u',h,\lambda)du'\,du\,dh\,dk\,dm\,d\lambda,$$
pourvu que $c'\geq c_{0}$. Ces expressions sont absolument convergentes d'apr\`es ce que l'on vient de prouver. La convergence de la premi\`ere est l'assertion (i) de l'\'enonc\'e. Pour prouver l'assertion (ii), il suffit d'apr\`es (1) de montrer que l'on peut choisir $C$ tel que la diff\'erence entre les deux expressions ci-dessus soit essentiellement major\'ee par $N^{-1}$.  D'apr\`es (3), cette diff\'erence est essentiellement major\'ee par la somme de l'int\'egrale
$$\int_{M(F)}\int_{H(F)U(F)_{c}}\int_{U(F)-U(F)_{c'}}{\bf 1}_{\sigma< Clog(N)}(hu) \Xi^H(h)\Xi^G(u^{-1}h^{-1}u'm)$$
$$\Xi^G(u'm)\kappa_{N}(m)\delta_{P}(m)^{-1}du'\,du\,dh\,dm,$$
o\`u $U(F)-U(F)_{c'}$ est le compl\'ementaire de $U(F)_{c'}$ dans $U(F)$, et de l'int\'egrale
$$\int_{M(F)}\int_{H(F)U(F)_{c}}\int_{U(F)_{c'}}(1-{\bf 1}_{\sigma< Clog(N)}(hu)) \Xi^H(h)\Xi^G(u^{-1}h^{-1}u'm)$$
$$\Xi^G(u'm)\kappa_{N}(m)\delta_{P}(m)^{-1}du'\,du\,dh\,dm.$$
Dans la premi\`ere, on peut oublier le terme ${\bf 1}_{\sigma< Clog(N)}(hu)$. La relation 4.3(7) nous dit qu'il existe $c_{1}$ tel que, si l'on prend $c'=c_{1}log(N)+c_{0}$, l'int\'egrale est essentiellement born\'ee par $N^{-1}$. L'entier $c'$ \'etant ainsi choisi, la relation 4.3(8) nous dit qu'il existe $c_{2}$ tel que la seconde int\'egrale soit essentiellement major\'ee par $N^{-1}$ pourvu que $Clog(N)\geq c_{2}(log(N)+c')$. On peut choisir $C$ tel qu'il en soit ainsi pour tout $N\geq 2$. Cela ach\`eve la preuve. $\square$

On fixe d\'esormais $C$ v\'erifiant l'assertion (ii) du lemme pr\'ec\'edent. Dans les paragraphes suivants  et jusqu'en 6.9, on fixe $L\in {\cal L}(M_{min})$ et ${\cal O}\in {\Pi_2(L)}_{f}$ tel que $m({\cal O},\rho)=1$.

\bigskip

\subsection{Rappels sur les termes constants faibles des coefficients des repr\'esentations temp\'er\'es}

Fixons un \'el\'ement $P_{min}=M_{min}U_{min}\in {\cal P}(M_{min})$. On note $\Delta$ l'ensemble de racines simples associ\'e.  On note $M_{min}(F)^+$ le sous-ensemble des $m\in M_{min}(F)$ tels que $H_{M_{min}}(m)\in {\cal A}_{P_{min}}^+$.Soit $Q_{1}=L_{1}U_{1}\in {\cal F}(M_{min})$ tel que $P_{min}\subset Q_{1}$. On note $\Delta^{L_{1}}$ le sous-ensemble de $\Delta$ associ\'e au sous-groupe parabolique minimal $P_{min}\cap L_{1}$ de $L_{1}$. Posons
$$W(L_{1}\vert G\vert L)=\{s\in W^G; sLs^{-1}\subset L_{1}, P_{min}\cap L_{1}\subset sQs^{-1}\}.$$
C'est un ensemble de repr\'esentants du quotient
$$ W^{L_{1}}\backslash\{s\in W^G; sLs^{-1}\subset L_{1}\}.$$
On identifie cet ensemble \`a un ensemble de repr\'esentants dans $K$. Pour $s\in W(L_{1}\vert G\vert L)$, notons $\gamma(s):{\cal K}^G_{Q,\tau}\to {\cal K}^G_{sQs^{-1},\tau}$ l'op\'erateur d\'efini par $(\gamma(s)e)(k)=e(s^{-1}k)$. Posons $Q_{1,s}=(L_{1}\cap sQs^{-1})U_{1}$, $\tilde{Q}_{1,s}=(L_{1}\cap sQs^{-1})\bar{U}_{1}$, o\`u $\bar{U}_{1}$ est le radical unipotent de $\bar{Q}_{1}$. Pour $\lambda\in i{\cal A}_{L,F}^*$, on d\'efinit les op\'erateurs:
$$J_{Q_{1,s}\vert sQs^{-1}}((s\tau)_{s\lambda})\circ \gamma(s):{\cal K}^G_{Q,\tau}\to {\cal K}^G_{Q_{1,s},s\tau}$$
et
$$J_{\tilde{Q}_{1,s}\vert sQs^{-1}}((s\tau)_{s\lambda}\circ\gamma(s):{\cal K}^G_{Q,\tau}\to {\cal K}^G_{\tilde{Q}_{1,s},s\tau}.$$
Posons $K_{1}=K\cap L_{1}(F)$. La restriction de $K$ \`a $K_{1}$ d\'efinit des homomorphismes
$$\begin{array}{ccc}{\cal K}^G_{Q_{1,s},s\tau}&&\\ &\searrow&\\ &&{\cal K}^{L_{1}}_{L_{1}\cap sQs^{-1},s\tau}\\ &\nearrow&\\ {\cal K}^G_{\tilde{Q}_{1,s},s\tau}&&\\ \end{array}$$
Pour $e'\in {\cal K}^G_{Q_{1,s},s\tau}$ et $e\in {\cal K}^G_{\tilde{Q}_{1,s},s\tau}$, on pose
$$(e',e)^{L_{1}}=\int_{K_{1} }(e'(k_{1}),e(k_{1}))dk_{1},$$
la mesure \'etant de masse totale $1$. Pour $e,e'\in {\cal K}^G_{Q,\tau}$, posons
$$(e',e)_{Q_{1},\lambda}=\sum_{s\in W(L_{1}\vert G\vert L)}(J_{Q_{1,s}\vert sQs^{-1}}((s\tau)_{s\lambda})\circ \gamma(s)e',J_{\tilde{Q}_{1,s}\vert sQs^{-1}}((s\tau)_{s\lambda})\circ\gamma(s)e)^{L_{1}}.$$

\ass{Lemme}{Soient $e,e'\in {\cal K}^G_{Q,\tau}$. Il existe $c>0$, $R\geq 0$ et $\epsilon>0$ tels que
$$\delta_{Q_{1}}(m)^{1/2}\vert  (e',\pi_{\lambda}(m)e)- (e',\pi_{\lambda}(m)e)_{Q_{1},\lambda}\vert <<$$
$$ \Xi^{L_{1}}(m)\sigma(m)^Rsup\{exp(-\epsilon \alpha( H_{M_{1}}(m)) ); \alpha\in \Delta-\Delta^{L_{1}}\}$$
pour tout $\lambda\in i{\cal A}_{L,F}^*$ et tout $m\in M_{min}(F)^+$ v\'erifiant les conditions
$\alpha(H_{M_{min}}(m))>c$ pour tout $\alpha\in \Delta- \Delta^{L_{1}}$.}

Preuve. D'apr\`es un r\'esultat de Casselman, on peut trouver $c$ tel que, pour $\lambda$ et $m$ comme dans l'\'enonc\'e, $\delta_{Q_{1}}(m)^{1/2}(e',\pi_{\lambda}(m)e)$ soit \'egal au "terme constant" de $(e',\pi_{\lambda}(m)e)$ au sens d'Harish-Chandra  (cf. [W2] corollaire I.4.4; l'uniformit\'e en $\lambda$ est expliqu\'ee en [W2] I.5).  L'expression $ \delta_{Q_{1}}(m)^{1/2}(e',\pi_{\lambda}(m)e)_{Q_{1},\lambda}$ est le "terme constant faible" de $(e',\pi_{\lambda}(m)e)$, cf. [W2] proposition V.1.1). On doit donc prouver que la diff\'erence entre le terme constant et le terme constant faible v\'erifie la majoration indiqu\'ee. Cette diff\'erence est exactement une fonction not\'ee $m\mapsto (E_{Q}^G\psi_{\omega})_{\bar{Q}_{1}}^+(m)$ dans le lemme VI.2.3 de [W2]. Ce lemme donne une majoration de cette fonction qui n'est pas tout-\`a-fait celle dont nous avons besoin. Il suffit de modifier sa d\'emonstration. Notons simplement $m\mapsto E_{\lambda}(m)$ notre fonction. D'apr\`es [W2] VI.2(2), on peut trouver des r\'eels $c>0$ et $R\geq0$ et  un \'el\'ement $a\in A_{L_{1}}(F)\cap M_{min}(F)^+$ de sorte que l'on ait la majoration
  $$\vert E_{\lambda}(a^nm)\vert <<2^{-n}\Xi^{L_{1}}(m)\sigma(m)^R$$
pour $\lambda$ et $m$ comme dans l'\'enonc\'e et $n\in {\mathbb N}$. Cela \'etant, pour $m$ comme dans l'\'enonc\'e, prenons pour $n$ la partie enti\`ere du plus petit des nombres
$$\frac{\alpha(H_{M_{min}}(m))-c}{\alpha(H_{M_{min}}(a))+1}$$
pour $\alpha\in \Delta-\Delta^{L_{1}}$. Posons $m'=a^{-n}m$. Cet \'el\'ement v\'erifie encore les conditions de l'\'enonc\'e. En appliquant la majoration ci-dessus \`a $n$ et $m'$, on obtient
$$(1) \qquad \vert E_{\lambda}(m)\vert <<2^{-n}\Xi^{L_{1}}(m')\sigma(m')^R.$$
Parce que $a\in A_{L_{1}}(F)$, $\Xi^{L_{1}}(m')=\Xi^{L_{1}}(m)$. On a 
$$\sigma(m')<<\sigma(m)+\sigma(a^{-n})<<\sigma(m)+n\sigma(a),$$
d'o\`u
$$2^{-n/2}\sigma(m')^R<<\sigma(m)^R.$$
D'autre part il existe $\epsilon>0$ tel que
$$2^{-n/2}<<sup\{exp(-\epsilon \alpha( H_{M_{1}}(m)) ); \alpha\in \Delta-\Delta^{L_{1}}\}.$$
Alors (1) entra\^{\i}ne la majoration cherch\'ee. $\square$

\bigskip

\subsection{Changement de fonction de troncature}

  Soit $Y\in {\cal A}_{P_{min}}^+$. Comme en [A3] paragraphe 3, on en d\'eduit une famille ${\cal Y}=(Y_{P'})_{P'\in {\cal P}(M_{min})}$ de points de ${\cal A}_{M_{min}}$, qui est $(G,M_{min})$-orthogonale et positive. On note $\zeta\mapsto\sigma_{M_{min}}^G(\zeta,{\cal Y})$ la fonction caract\'eristique de son enveloppe convexe. On note $g\mapsto u(g,Y)$ la fonction caract\'eristique de l'ensemble des $g\in G(F)$ qui s'\'ecrivent $g=kmk'$, o\`u $k,k'\in K$, $m\in M_{min}(F)$, avec $\sigma_{M_{min}}^G(H_{M_{min}},{\cal Y})=1$. Fixons $e',e''\in {\cal K}^G_{Q,\tau}$ et une fonction $\varphi$ sur $i{\cal A}_{L,F}^*$, que l'on suppose $C^{\infty}$. Pour $e\in {\cal K}^G_{Q,\tau}$, $g,g'\in G(F)$ et $\lambda\in i{\cal A}_{L,F}^*$, posons
$$\Phi(e,g,g',\lambda)=(\pi_{\lambda}(g)e,\pi_{\lambda}(g')e')(e'',\pi_{\lambda}(g)\pi_{\lambda}(f)e).$$
Posons
$$\Phi_{N}(g')=[i{\cal A}_{{\cal O}}^{\vee}:i{\cal A}_{L,F}^{\vee}]^{-1} \sum_{e\in {\cal B}_{{\cal O}}^{K_{f}}}\int_{i{\cal A}_{L,F}^*}m(\tau_{\lambda})\varphi(\lambda)\int_{G(F)}\Phi(e,g,g',\lambda)\kappa_{N}(g)dg\,d\lambda,$$
$$\Phi_{Y}(g')=[i{\cal A}_{{\cal O}}^{\vee}:i{\cal A}_{L,F}^{\vee}]^{-1} \sum_{e\in {\cal B}_{{\cal O}}^{K_{f}}}\int_{i{\cal A}_{L,F}^*}m(\tau_{\lambda})\varphi(\lambda)\int_{G(F)}\Phi(e,g,g',\lambda)u(g,Y)dg\,d\lambda.$$
Ces expressions sont absolument convergentes. Pour la seconde, c'est \'evident: les int\'egrales sont \`a support compact. Pour la premi\`ere, pour $e$ et $g'$ fix\'es, on a une majoration
$$\vert \Phi(e,g,g',\lambda)\vert <<\Xi^G(g)^2$$
pour tous $\lambda$, $g$. Or l'int\'egrale
$$(1)\qquad \int_{G(F)}\kappa_{N}(g)\Xi^G(g)^2dg$$
est convergente d'apr\`es 4.3(2).

\ass{Proposition}{ Soit $R$ un r\'eel. Il existe des r\'eels $c_{1},c_{2}>0$ tels que l'on ait la majoration
$$\vert \Phi_{N}(g')-\Phi_{Y}(g')\vert <<N^{-R}$$
pour tout $N\geq2$,   tout $g'\in G(F)$ tel que $\sigma(g')\leq Clog(N)$ et tout $Y\in {\cal A}_{P_{min}}^+$ v\'erifiant les in\'egalit\'es $c_{1}log(N)\leq \alpha(Y)\leq c_{2}N$ pour tout $\alpha\in \Delta_{min}$.}

Preuve.  Commen\c{c}ons par pr\'eciser le calcul de convergence que l'on a fait avant l'\'enonc\'e. On a

(2) il existe $R_{1}\geq0$ tel que
$$\vert \Phi_{N}(g')\vert<< N^{R_{1}}$$
pour tout $N\geq2 $ et tout $g'\in G(F)$ tel que $\sigma(g')\leq Clog(N)$.

En effet, pour $e\in {\cal K}^G_{Q,\tau}$ on a plus pr\'ecis\'ement la majoration
$$ \vert \Phi(e,g,g',\lambda)\vert <<\Xi^G(g^{_{'}-1}g)\Xi^G(g)$$
pour tous $\lambda$, $g$, $g'$. Gr\^ace \`a 3.3(5), il existe $R_{2}>0$ tel que 
$$\Xi^G(g^{_{'}-1}g)<<exp(R_{2}\sigma(g'))\Xi^G(g).$$
Pour $g'$ tel que $\sigma(g')\leq Clog(N)$,on obtient
$$ \vert \Phi(e,g,g',\lambda)\vert <<N^{CR_{2}}\Xi^G(g)^2.$$
D'apr\`es 4.3(2), il existe $R_{3}>0$ tel que l'int\'egrale (1) soit essentiellement major\'ee par $N^{R_{3}}$. On en d\'eduit (2) avec $R_{1}=CR_{2}+R_{3}$.

D\'efinissons une fonction $D$ sur $M_{min}(F)^+$ par
$$D(m)=mes(KmK)mes(K\cap M_{min}(F))^{-1}.$$
D'apr\`es [W2], on a une majoration $D(m)<<\delta_{P_{min}}(m)$. Pour toute fonction $f'\in C_{c}^{\infty}(G(F))$, on a l'\'egalit\'e
$$\int_{G(F)}f'(g)dg=\int_{K}\int_{K}\int_{M_{min}(F)^+}D(m)f'(k_{1}mk_{2})dm\,dk_{1}\,dk_{2},$$
cf. [A3] (1.3). Pour tout $Q_{1}=L_{1}U_{1}\in {\cal F}(M_{min})$, la famille ${\cal Y}$ permet de d\'efinir des fonctions $\zeta\mapsto \sigma_{M_{min}}^{Q_{1}}(\zeta,{\cal Y})$ et $\zeta\mapsto \tau_{Q_{1}}(\zeta-Y_{Q_{1}})$ sur ${\cal A}_{M_{min}}$ (cf.[A3]; on a repris les d\'efinitions en [W1] 10.3). On v\'erifie que l'\'egalit\'e
$$\sigma_{M_{min}}^{Q_{1}}(\zeta,{\cal Y})\tau_{Q_{1}}(\zeta-Y_{Q_{1}})=1$$
entraine

(3) $\beta(\zeta)>inf\{\alpha(Y); \alpha\in \Delta\}\geq0$ pour toute racine $\beta$ de $A_{M_{min}}$ dans $\mathfrak{u}_{1}$. 

On a l'\'egalit\'e
$$\sum_{Q_{1}\in {\cal F}(M_{min})}\sigma_{M_{min}}^{Q_{1}}(\zeta,{\cal Y})\tau_{Q_{1}}(\zeta-Y_{Q_{1}})=1$$
pour tout $\zeta$. L'in\'egalit\'e pr\'ec\'edente entra\^{\i}ne que, pour $\zeta\in {\cal A}_{P_{min}}^+$, seuls interviennent de fa\c{c}on non nulle les $Q_{1}$ contenant $P_{min}$. C'est-\`a-dire que, pour $\zeta\in {\cal A}_{P_{min}}^+$, on a l'\'egalit\'e
$$\sum_{Q_{1}\in {\cal F}(M_{min}), P_{min}\subset Q_{1}}\sigma_{M_{min}}^{Q_{1}}(\zeta,{\cal Y})\tau_{Q_{1}}(\zeta-Y_{Q_{1}})=1.$$
Pour toute fonction $f'\in C_{c}^{\infty}(G(F))$, on a donc l'\'egalit\'e
$$(4)\qquad \int_{G(F)}f'(g)dg=\sum_{Q_{1}\in {\cal F}(M_{min}), P_{min}\subset Q_{1}}\int_{K}\int_{K}\int_{M_{min}(F)^+}$$
$$f'(k_{1}mk_{2})D(m)\sigma_{M_{min}}^{Q_{1}}(H_{M_{min}}(m),{\cal Y})\tau_{Q_{1}}(H_{M_{min}}(m)-Y_{Q_{1}})dm\,dk_{1}\,dk_{2}.$$
Pour $Q_{1}\in {\cal F}(M_{min})$ tel que $P_{min}\subset Q_{1}$,  posons
 
$$\Phi_{N,Y,Q_{1}}(g')=[i{\cal A}_{{\cal O}}^{\vee}:i{\cal A}_{L,F}^{\vee}]^{-1} \sum_{e\in {\cal B}_{{\cal O}}^{K_{f}}}\int_{i{\cal A}_{L,F}^*}m(\tau_{\lambda})\varphi(\lambda)\int_{K}\int_{K}\int_{M_{min}(F)^+}$$
$$\Phi(e,k_{1}mk_{2},g',\lambda)\kappa_{N}(k_{1}mk_{2})$$
$$D(m)\sigma_{M_{min}}^{Q_{1}}(H_{M_{min}}(m),{\cal Y})\tau_{Q_{1}}(H_{M_{min}}(m)-Y_{Q_{1}})dm\,dk_{1}\,dk_{2}\,d\lambda.$$
En l'appliquant (4), on a l'\'egalit\'e
$$\Phi_{N}(g')=\sum_{Q_{1}\in {\cal F}(M_{min}), P_{min}\subset Q_{1}}\Phi_{N,Y,Q_{1}}(g').$$
 
Consid\'erons d'abord le sous-groupe parabolique $Q_{1}=G$. Dans ce cas, pour $g=k_{1}mk_{2}$, avec $k_{1},k_{2}\in K$ et $m\in M_{min}(F)^+$, on a $\sigma_{M_{min}}^G(H_{M_{min}}(m),{\cal Y})\tau_{G}(H_{M_{min}}(m)-Y_{G})=1$ si et seulement si $u(g,Y)=1$.  On a donc
$$\Phi_{N,Y,G}(g')=[i{\cal A}_{{\cal O}}^{\vee}:i{\cal A}_{L,F}^{\vee}]^{-1} \sum_{e\in {\cal B}_{{\cal O}}^{K_{f}}}\int_{i{\cal A}_{L,F}^*}m(\tau_{\lambda})\varphi(\lambda)$$
$$\int_{G(F)}\Phi(e,g,g',\lambda)\kappa_{N}(g)u(g,Y)dg\,d\lambda.$$
Montrons que

(5) il existe $c_{2}>0$ tel que, si $\alpha(Y)\leq c_{2}N$ pour tout $\alpha\in \Delta_{min}$, on a l'\'egalit\'e $u(g,Y)\kappa_{N}(g)=u(g,Y)$ pour tout $g\in G(F)$.

Ecrivons $g=muk$, avec $m\in M(F)$, $u\in U(F)$ et $k\in K$. On a $\kappa_{N}(g)=\kappa_{N}(m)$ et une majoration $\sigma(m)<<\sigma(g)$.  Par d\'efinition de la fonction $\kappa_{N}$, il existe $c_{3}$ tel que la majoration $\sigma(m)\leq c_{3}N$ entra\^{\i}ne $\kappa_{N}(m)=1$. On a d'autre part une majoration $\sigma(g)<<sup\{\alpha(Y); \alpha\in \Delta\}$ pour tout $g\in G(F)$ tel que $u(g,Y)=1$. La combinaison de ces propri\'et\'es entra\^{\i}ne (5).

On d\'eduit  de (5) que, si $Y$ v\'erifie les conditions de cette assertion, on a l'\'egalit\'e $\Phi_{N,Y,G}(g')=\Phi_{Y}(g')$ pour tout $g'$. Pour d\'emontrer la proposition, il suffit donc de trouver $c_{1}$ tel que, si $Y$ v\'erifie les minorations de l'\'enonc\'e, on a la majoration
$$(6)\qquad \vert  \Phi_{N,Y,Q_{1}}(g')\vert <<N^{-R}$$
pour $g'$ comme dans l'\'enonc\'e et tout $Q_{1}\not=G$.

On fixe d\'esormais $Q_{1}=L_{1}U_{1}\in {\cal F}(M_{min})$ tel que $P_{min}\subset Q_{1}$ et $Q_{1}\not=G$. Pour simplifier la r\'edaction, on va consid\'erer un r\'eel $c_{1}$ et supposer $\alpha(Y)\geq c_{1}log(N)$ pour tout $\alpha\in \Delta$.  On montrera que toutes les propri\'et\'es dont on a besoin sont v\'erifi\'ees si $c_{1}$ est assez grand. Soient $g'\in G(F)$ tel que $\sigma(g')\leq Clog(N)$, $k_{1},k_{2}\in K$ et $m\in M_{min}(F)^+$. On pose $\zeta=H_{M_{min}}(m)$ et on suppose  $\sigma_{M_{min}}^{Q_{1}}(\zeta,{\cal Y})\tau_{Q_{1}}(\zeta-Y_{Q_{1}})=1$. Ecrivons $g^{_{'}-1}k_{1}=k^{_{'}-1}l'u'$, avec $u'\in U_{1}(F)$, $l'\in L_{1}(F)$, $k'\in K$. On a

(7) si $c_{1}$ est assez grand, $k_{2}^{-1}m^{-1}u'mk_{2}$ appartient \`a $K_{f}$.

Posons $u'=exp(X')$, avec $X'\in \mathfrak{u}_{1}(F)$ et fixons une norme $\vert .\vert $ sur $\mathfrak{g}(F)$. On a $log(\vert X'\vert )<<\sigma(u')<<\sigma(g')\leq Clog(N)$. Il existe donc $c_{4},c_{5}>0$ tels que
$$log(\vert m^{-1}X'm\vert )\leq c_{4}log(N)-c_{5 }inf\{\beta(\zeta); \beta\text{ racine de }A_{M_{min}}\text{ dans }\mathfrak{u}_{1}\}.$$
Gr\^ace \`a (3), $m^{-1}X'm$ est aussi petit que l'on veut pourvu que $c_{1}$ soit assez grand. On peut en particulier imposer que $m^{-1}u'm=exp(m^{-1}X'm)$ appartienne \`a $K_{f}$. Puisque $K_{f}$ est distingu\'e dans $K$, (7) en r\'esulte.

Notons $M_{min}(F)^{L_{1},+}$ l'ensemble des $m'\in M_{min}(F)$ tels que $\alpha(H_{M_{min}}(m'))\geq0$ pour tout $\alpha\in \Delta^{L_{1}}$. Posons $K_{1}=K\cap L_{1}(F)$. L'\'el\'ement $l'm$ de $L_{1}(F)$ s'\'ecrit $l'm=k_{3}m'k_{4}$, avec $k_{3},k_{4}\in K _{1}$ et $m'\in M_{min}(F)^{L_{1},+}$. Posons $\zeta'=H_{M_{min}}(m')$. On a

(8) pour tout $c>0$,   $m'$ appartient \`a $M_{min}(F)^+$ et v\'erifie $\alpha(\zeta')>clog(N)$ pour tout $\alpha\in \Delta-\Delta^{L_{1}}$, pourvu que $c_{1}$ soit assez grand.

Soit $\alpha$. Notons $\mathfrak{u}_{\alpha}$ le sous-espace radiciel de $\mathfrak{u}_{min}$ associ\'e \`a $\alpha$. Il existe $c_{6},c_{7}>0$ tels que, pour tout \'el\'ement non nul $X\in \mathfrak{u}_{\alpha}(F)$ et tout $x\in M_{min}(F)$, on ait les in\'egalit\'es
$$(9) \qquad \alpha(H_{M_{min}}(x))-c_{6}\leq log(\frac{\vert xXx^{-1}\vert }{\vert X\vert })\leq \alpha(H_{M_{min}}(x))+c_{7}.$$
Supposons $\alpha\in \Delta-\Delta^{L_{1}}$, soit $X$ un \'el\'ement non nul de $\mathfrak{u}_{\alpha}(F)$. On a $m'Xm^{_{'}-1}=k_{3}^{-1}l'mk_{4}^{-1}Xk_{4}m^{-1}l^{_{'}-1}k_{3}$.  Puisque $k_{4}$ appartient \`a $L_{1}(F)$, $k_{4}^{-1}Xk_{4}$ appartient \`a $\mathfrak{u}_{1}(F)$. En utilisant (9) pour $m$, on voit qu'il existe des constantes $c_{8},c_{9},c_{10}>0$ telles que l'on ait l'in\'egalit\'e
$$log(\vert m'Xm^{_{'}-1}\vert )\geq log(\vert X\vert )-c_{8}\sigma(l')+c_{9}inf\{\beta(\zeta); \beta\text{ racine de }A_{M_{min}}\text{ dans }\mathfrak{u}_{1}\}-c_{10}.$$
On a $\sigma(l')<<\sigma(g')\leq C log(N)$. Gr\^ace \`a (3), on a donc
$$log(\vert m'Xm^{_{'}-1}\vert )> log(\vert X\vert )+clog(N)+c_{7}$$
pourvu que $c_{1}$ soit assez grand.  Alors (9) entra\^{\i}ne la minoration cherch\'ee de $\alpha(\zeta')$. En particulier $\alpha(\zeta')>0$. Puisque l'on a aussi $\alpha(\zeta')\geq0$ pour $\alpha\in \Delta^{L_{1}}$ par d\'efinition de $m'$, on a bien $m'\in M_{min}(F)^+$. Cela prouve (8).

Pour tout $\lambda$ et tout $e\in {\cal B}_{{\cal O}}^{K_{f}}$, on a l'\'egalit\'e
$$\Phi(e,k_{1}mk_{2},g',\lambda)=(\pi_{\lambda}(l'u'mk_{2})e,\pi_{\lambda}(k')e')(\pi_{\lambda}(k_{1}^{-1})e'',\pi_{\lambda}(mk_{2})\pi_{\lambda}(f)e).$$
Gr\^ace \`a (7), on peut supprimer $u'$ de cette expression pourvu que $c_{1}$ soit assez grand. On obtient
$$\Phi(e,k_{1}mk_{2},g',\lambda)=(\pi_{\lambda}(m'k_{4}k_{2})e,\pi_{\lambda}(k_{3}^{-1}k')e')(\pi_{\lambda}(k_{1}^{-1})e'',\pi_{\lambda}(mk_{2})\pi_{\lambda}(f)e).$$
Posons
$$\Phi_{w}(e,k_{1}mk_{2},g',\lambda)= (\pi_{\lambda}(m'k_{4}k_{2})e,\pi_{\lambda}(k_{3}^{-1}k')e')_{Q_{1},\lambda}(\pi_{\lambda}(k_{1}^{-1})e'',\pi_{\lambda}(mk_{2})\pi_{\lambda}(f)e)_{Q_{1},\lambda}.$$
 Le lemme 6.5 affirme l'existence de r\'eels $R_{4}\geq0$ et $\epsilon>0$ tels que la valeur absolue de la diff\'erence
$$\Phi(e,k_{1}mk_{2},g',\lambda)- \Phi_{w}(e,k_{1}mk_{2},g',\lambda)$$
soit  born\'ee par la somme de
$$\delta_{Q_{1}}(m')^{-1/2}\Xi^{L_{1}}(m')\sigma(m')^{R_{4}}sup\{exp(-\epsilon\alpha(\zeta'); \alpha\in \Delta-\Delta^{L_{1}}\}\vert (\pi_{\lambda}(k_{1}^{-1})e'',\pi_{\lambda}(mk_{2})\pi_{\lambda}(f)e)\vert $$
et de
$$\delta_{Q_{1}}(m')^{-1/2}\vert (\pi_{\lambda}(m'k_{4}k_{2})e,\pi_{\lambda}(k_{3}^{-1}k')e')_{Q_{1},\lambda}\vert \delta_{Q_{1}}(m)^{-1/2}\Xi^{L_{1}}(m)$$
$$\sigma(m)^{R_{4}}sup\{exp(-\epsilon\alpha(\zeta); \alpha\in \Delta-\Delta^{L_{1}}\}.$$
On sait qu'il existe $R_{5}\geq0$ tel que $\delta_{Q_{1}}(x)^{-1/2}\Xi^{L_{1}}(x)<<\Xi^G(x)\sigma(x)^{R_{5}}$ pour tout $x\in M_{min}(F)^+$. D'autre part, $\Xi^G(m')=\Xi^G(l'm)$. En utilisant (8), on voit qu'il existe $R_{6}\geq0$ tel que, pour tout $c>0$, les expressions ci-dessus soient essentiellement major\'ees par
$$N^{-c}\Xi^G(l'm)\Xi^G(m)\sigma(l')^{R_{6}}\sigma(m)^{R_{6}}$$
pourvu que $c_{1}$ soit assez grand. D\'efinissons des termes $\Phi_{N,Y,Q_{1},w}(e,g',\lambda)$ et $\Phi_{N,Y,Q_{1},w}(g')$ en rempla\c{c}ant $\Phi(e,k_{1}mk_{2},g',\lambda)$ par $ \Phi_{w}(e,k_{1}mk_{2},g',\lambda)$ dans les d\'efinitions de $\Phi_{N,Y,Q_{1}}(e,g',\lambda)$ et $\Phi_{N,Y,Q_{1}}(g')$. Si on oublie le facteur $N^{-c}$ de la majoration ci-dessus, le m\^eme calcul qu'en (1) montre l'existence de $R_{7}\geq0$ tel que
$$\vert \Phi_{N,Y,Q_{1}}(g')-\Phi_{N,Y,Q_{1},w}(g')\vert<<N^{R_{7}}.$$
En r\'eintroduisant le facteur $N^{-c}$, on voit que, pour tout $c>0$, la diff\'erence ci-dessus est essentiellement major\'ee par $N^{-c}$ pourvu que $c_{1}$ soit assez grand.  Pour prouver (6), il suffit donc de prouver la majoration
$$ \qquad \vert \Phi_{N,Y,Q_{1},w}(g')\vert<<N^{-R}.$$ 

On poursuit le calcul pr\'ec\'edent, avec les m\^emes notations. D'apr\`es la d\'efinition de 6.5, on peut d\'ecomposer $\Phi_{w}(e,k_{1}mk_{2},g',\lambda)$ en
$$\sum_{s_{1},s_{2}\in W(L_{1}\vert G\vert L)}\Phi_{s_{1},s_{2}}(e,k_{1}mk_{2},g',\lambda)$$
o\`u
$$\Phi_{s_{1},s_{2}}(e,k_{1}mk_{2},g',\lambda)= $$
$$(J_{\tilde{Q}_{1,s_{1}}\vert s_{1}Qs_{1}^{-1}}((s_{1}\tau)_{s_{1}\lambda})\circ\gamma(s_{1})\circ\pi_{\lambda}(m'k_{4}k_{2})e,J_{Q_{1,s_{1}}\vert s_{1}Qs_{1}^{-1}}((s_{1}\tau)_{s_{1}\lambda})\circ\gamma(s_{1})\circ\pi_{\lambda}(k_{3}^{-1}k')e')^{L_{1}}$$
$$(J_{Q_{1,s_{2}}\vert s_{2}Qs_{2}^{-1}}((s_{2}\tau)_{s_{2}\lambda})\circ\gamma(s_{2})\circ\pi_{\lambda}(k_{1}^{-1})e'',J_{\tilde{Q}_{1,s_{2}}\vert s_{2}Qs_{2}^{-1}}((s_{2}\tau)_{s_{2}\lambda})\circ\gamma(s_{2})\circ\pi_{\lambda}(mk_{2})\pi_{\lambda}(f)e)^{L_{1}}.$$
Cette d\'efinition n'a bien s\^ur de sens que "presque partout" en $\lambda$, les op\'erateurs d'entrelacement pouvant avoir des p\^oles. Au moins formellement, on peut d\'ecomposer de m\^eme   $\Phi_{N,Y,Q_{1},w}(g')$ en somme de termes  $\Phi_{N,Y,Q_{1},s_{1},s_{2}}(g')$.  Il y a un probl\`eme de convergence  \`a cause des p\^oles des op\'erateurs d'entrelacement.  L'assertion suivante va r\'esoudre ce probl\`eme.  Soient $s_{1},s_{2}\in W(L_{1}\vert G\vert L)$. Montrons que

(10) pour $e$ fix\'e, la fonction $m(\tau_{\lambda})\Phi_{s_{1},s_{2}}(e,k_{1}mk_{2},g',\lambda)$ est une combinaison lin\'eaire de fonctions qui sont elles-m\^emes des produits  $f_{1}(m',\lambda)f_{2}(m,\lambda)f_{3}(k_{1},k_{2},k_{3},k_{4},k')f_{4}(\lambda)$, o\`u:

$$f_{1}(m',\lambda)=\delta_{Q_{1}}(m')^{-1/2}(Ind_{L_{1}\cap s_{1}Qs_{1}^{-1}}^{L_{1}}((s_{1}\tau)_{s_{1}\lambda},m')e'_{1},e_{1})$$
pour des \'el\'ements $e_{1}$ et $e'_{1}$ de ${\cal K}^{L_{1}}_{L_{1}\cap s_{1}Qs_{1}^{-1},s_{1}\tau}$;

$$f_{2}(m,\lambda)=\delta_{Q_{1}}(m)^{-1/2}(e'_{2},Ind_{L_{1}\cap s_{2}Qs_{2}^{-1}}^{L_{1}}((s_{2}\tau)_{s_{2}\lambda},m)e_{2})$$
 pour des \'el\'ements $e_{2}$ et $e'_{2}$ de ${\cal K}^{L_{1}}_{L_{1}\cap s_{2}Qs_{2}^{-1},s_{2}\tau}$;
 
 $f_{3}$ est une fonction localement constante des variables $k_{1},k_{2},k_{3},k_{4}$ et $k'$;

$f_{4}$ est une fonction $C^{\infty}$ de $\lambda$.
\bigskip

Fixons un sous-groupe ouvert compact $K_{0}$ de $ K$, distingu\'e dans $K$, inclus dans $K_{f}$ et tel que $e$, $e'$ et $e''$ soient invariants par $K_{0}$. Fixons des bases ${\cal B}$, resp. ${\cal B}_{1}$, $\tilde{{\cal B}}_{1}$, ${\cal B}_{2}$, $\tilde{{\cal B}}_{2}$, des sous-espaces des \'el\'ements invariants par $K_{0}$ dans ${\cal K}^G_{Q,\tau}$, resp. ${\cal K}^G_{Q_{1,s_{1}},s_{1}\tau}$, ${\cal K}^G_{\tilde{Q}_{1,s_{1}},s_{1}\tau}$, ${\cal K}^G_{Q_{1,s_{2}},s_{2}\tau}$, ${\cal K}^G_{\tilde{Q}_{1,s_{2}},s_{2}\tau}$. Consid\'erons par exemple le terme
$$J_{\tilde{Q}_{1,s_{2}}\vert s_{2}Qs_{2}^{-1}}((s_{2}\tau)_{s_{2}\lambda})\circ\gamma(s_{2})\circ\pi_{\lambda}(mk_{2})\pi_{\lambda}(f)e$$
qui intervient dans la d\'efinition de $\Phi_{s_{1},s_{2}}(e,k_{1}mk_{2},g',\lambda)$. D'apr\`es les propri\'et\'es d'entrelacement de nos op\'erateurs, il est \'egal \`a
$$Ind_{\tilde{Q}_{1,s_{2}}}^G((s_{2}\tau)_{s_{2}\lambda},m)\circ J_{\tilde{Q}_{1,s_{2}}\vert s_{2}Qs_{2}^{-1}}((s_{2}\tau)_{s_{2}\lambda})\circ\gamma(s_{2})\circ\pi_{\lambda}(k_{2})\pi_{\lambda}(f)e.$$
  On peut remplacer $\pi_{\lambda}(k_{2})\pi_{\lambda}(f)e$ par son expression dans la base ${\cal B}$. Les coefficients sont combinaisons lin\'eaires de produits d'une fonction localement constante en $k_{2}$ et d'une fonction $C^{\infty}$ en $\lambda$. Pour tout $b\in {\cal B}$, on peut ensuite remplacer $J_{\tilde{Q}_{1,s_{2}}\vert s_{2}Qs_{2}^{-1}}((s_{2}\tau)_{s_{2}\lambda})\circ\gamma(s_{2})b$ par son expression dans la base $\tilde{{\cal B}}_{2}$. Les coefficients sont  des fonctions de $\lambda$ de la forme
$$(\tilde{b}_{2},J_{\tilde{Q}_{1,s_{2}}\vert s_{2}Qs_{2}^{-1}}((s_{2}\tau)_{s_{2}\lambda})\circ\gamma(s_{2})b)$$
o\`u $b\in {\cal B}$, $\tilde{b}_{2}\in \tilde{{\cal B}}_{2}$. Notons $j_{\tilde{Q}_{1,s_{2}}}(\lambda)$ une telle fonction. En appliquant le m\^eme calcul aux autres termes, l'expression  $m(\tau_{\lambda})\Phi_{s_{1},s_{2}}(e,k_{1}mk_{2},g',\lambda)$  appara\^{\i}t comme une combinaison lin\'eaire de fonctions qui sont produits de fonctions des deux derniers  types de l'assertion et de fonctions des types

- une fonction $(Ind_{\tilde{Q}_{1,s_{1}}}^G((s_{1}\tau)_{s_{1}\lambda},m')\tilde{b}_{1},b_{1})^{L_{1}}$ o\`u $\tilde{b}_{1}\in \tilde{{\cal B}}_{1}$, $b_{1}\in {\cal B}_{1}$;

- une fonction $(b_{2},(Ind_{\tilde{Q}_{1,s_{2}}}^G((s_{2}\tau)_{s_{2}\lambda},m)\tilde{b}_{2})^{L_{1}}$ o\`u $\tilde{b}_{2}\in \tilde{{\cal B}}_{2}$, $b_{2}\in {\cal B}_{2}$;

- une fonction $m(\tau_{\lambda})\overline{{j}_{\tilde{Q}_{1,s_{1}}}(\lambda)j_{Q_{1,s_{2}}}(\lambda)}j_{\tilde{Q}_{1,s_{2}}}(\lambda)j_{Q_{1,s_{1}}}(\lambda)$;

Consid\'erons le premier type ci-dessus. Notons $e'_{1}$, resp. $e_{1}$, la restriction de $\tilde{b}_{1}$, resp. $b_{1}$, \`a $K_{1}$. Ce sont des \'el\'ements de ${\cal K}^{L_{1}}_{L_{1}\cap s_{1}Qs_{1}^{-1},s_{1}\tau}$. En explicitant les d\'efinitions, on voit que, pour tout $x\in L_{1}(F)$, a fortiori pour $x\in M_{min}(F)$, on a l'\'egalit\'e
$$(Ind_{\tilde{Q}_{1,s_{1}}}^G((s_{1}\tau)_{s_{1}\lambda},x)\tilde{b}_{1},b_{1})^{L_{1}}=\delta_{Q_{1}}(x)^{-1/2}(Ind_{L_{1}\cap s_{1}Ls_{1}^{-1}}^{L_{1}}((s_{1}\tau)_{s_{1}\lambda},x)e'_{1},e_{1}).$$
La fonction du premier type ci-dessus est donc aussi du premier type de l'assertion (10). De m\^eme pour les fonctions du deuxi\`eme type. Reste celles du troisi\`eme type ci-dessus. Or la d\'emonstration du corollaire V.2.3 de [W2] s'applique \`a ces fonctions et montre qu'elles sont des restrictions \`a $i{\cal A}_{L,F}^*$ de fonctions rationnelles sur $({\cal A}_{L}^*\otimes_{{\mathbb R}}{\mathbb C})/i{\cal A}_{L,F}^{\vee}$, holomorphes au voisinage de $i{\cal A}_{L,F}^*$. Ces fonctions sont donc du quatri\`eme type de l'assertion (10). Cela d\'emontre cette assertion.

A l'aide de (10), on voit que l'expression qui d\'efinit  $\Phi_{N,Y,Q_{1},s_{1},s_{2}}(g')$ est absolument convergente: puisqu'il n'y a pas de singularit\'e en $\lambda$, il suffit de reprendre la preuve d\'ej\`a faite pour $\Phi_{N}(g')$. On a alors
$$ \Phi_{N,Y,Q_{1},w}(g')=\sum_{s_{1},s_{2}\in W(L_{1}\vert G\vert L)}\Phi_{N,Y,Q_{1},s_{1},s_{2}}(g').$$
Il  suffit de prouver que, pour tous $s_{1},s_{2}$, on a une majoration
$$\vert \Phi_{N,Y,Q_{1},s_{1},s_{2}}(g')\vert <<N^{-R}.$$

Fixons $s_{1},s_{2}$. Posons $s=s_{1}s_{2}^{-1}$. Supposons d'abord v\'erifi\'ee l'hypoth\`ese

{\it {\bf (Hyp)} il n'existe pas de sous-groupe parabolique $Q_{2}=L_{2}U_{2}\in {\cal F}_{M_{min}}$ tel que $Q_{1}\subset Q_{2}\not=G$ et que $s$ fixe tout point de ${\cal A}_{L_{2}}$.}

Dans ce cas, introduisons une fonction
$$\Psi_{N,Y,Q_{1},s_{1},s_{2}}(g')=\int_{i{\cal A}_{L,F}^*}f_{4}(\lambda)\int_{K}\int_{K}\int_{M_{min}(F)^+}f_{1}(m',\lambda)f_{2}(m,\lambda)f_{3}(k_{1},k_{2},k_{3},k_{4},k')$$
$$D(m)\kappa_{N}(k_{1}mk_{2})\sigma_{M_{min}}^{Q_{1}}(\zeta,{\cal Y})\tau_{Q_{1}}(\zeta-Y_{Q_{1}})dm\,dk_{1}\,dk_{2}\,d\lambda$$
o\`u $f_{1}$, $f_{2}$, $f_{3}$ et $f_{4}$ v\'erifient les conditions de l'assertion (10). Cette assertion nous dit que $\Phi_{N,Y,Q_{1},s_{1},s_{2}}(g')$ est combinaison lin\'eaire de fonctions de ce type. On va majorer $\vert \Psi_{N,Y,Q_{1},s_{1},s_{2}}(g')\vert $. Pour tout $x\in L_{1}(F)$, on choisit des \'el\'ements $l_{s_{1}}(x)\in s_{1}L(F)s_{1}^{-1}$, $u_{s_{1}}(x)\in L_{1}(F)\cap s_{1}U_{Q}(F)s_{1}^{-1}$ et $k_{s_{1}}(x)\in K_{1}$ de sorte que $x=l_{s_{1}}(x)u_{s_{1}}(x)k_{s_{1}}(x)$. On a l'\'egalit\'e
$$f_{1}(m',\lambda)=\int_{K_{1}}f'_{1}(m',x)exp(-(s_{1}\lambda)(H_{L_{1}\cap s_{1}Qs_{1}^{-1}}(xm')))dx,$$
o\`u
$$f'_{1}(m',x)=\delta_{Q_{1}}(m')^{-1/2}((s_{1}\tau)(l_{s_{1}}(xm'))(e'_{1}(k_{s_{1}}(xm'))),e_{1}(x))\delta_{L_{1}\cap s_{1}Qs_{1}}^{L_{1}}(l_{s_{1}}(xm'))^{1/2 }.$$
Le produit scalaire figurant dans cette expression est celui de deux \'el\'ements de $E_{s_{1}\tau}=E_{\tau}$. On \'ecrit $f_{2}(m,\lambda)$ de la m\^eme fa\c{c}on et on obtient
$$\Psi_{N,Y,Q_{1},s_{1},s_{2}}(g')=\int_{K}\int_{K}\int_{M_{min}(F)^+}\int_{K_{1}}\int_{K_{1}}f_{3}(k_{1},k_{2},k_{3},k_{4},k')D(m)\kappa_{N}(k_{1}mk_{2})$$
$$\sigma_{M_{min}}^{Q_{1}}(\zeta,{\cal Y})\tau_{Q_{1}}(\zeta-Y_{Q_{1}})f'_{1}(m',x)f'_{2}(m,y)f_{5}(xm',ym)dy\,dx\,dm\,dk_{1}\,dk_{2}$$
o\`u
$$f_{5}(xm',ym)=\int_{i{\cal A}_{L,F}^*}f_{4}(\lambda)exp(-(s_{1}\lambda)(H_{L_{1}\cap s_{1}Qs_{1}^{-1}}(xm'))+(s_{2}\lambda)(H_{L_{1}\cap s_{2}Qs_{2}^{-1}}(ym)))d\lambda.$$
Par le changement de variable $\lambda\mapsto s_{1}^{-1}\lambda$, on a
$$f_{5}(xm',ym)=\int_{i{\cal A}_{s_{1}Ls_{1},F}^*}f_{4}(s_{1}^{-1}\lambda)exp(-\lambda(\zeta(xm',ym)))d\lambda,$$
o\`u
$$\zeta(xm',ym)=H_{L_{1}\cap s_{1}Qs_{1}^{-1}}(xm'))-s(H_{L_{1}\cap s_{2}Qs_{2}^{-1}}(ym)).$$
On suppose toujours, comme il est loisible, que $m$ v\'erifie la relation $\sigma_{M_{min}}^{Q_{1}}(\zeta,{\cal Y})\tau_{Q_{1}}(\zeta-Y_{Q_{1}})=1$. Munissons ${\cal A}_{M_{min}}$ d'une norme euclidienne $\vert .\vert $ invariante par l'action de $W^G$. Montrons que

(11) on a la minoration
$$\vert \zeta(xm',ym)\vert >> log(N)$$
pour tous $x,y\in K_{1}$ pourvu que $c_{1}$ soit assez grand.

Posons $\zeta(xm')=H_{P_{min}}(xm')$, $\zeta(ym)=H_{P_{min}}(ym)$. Rappelons que, par d\'efinition de $W(L_{1}\vert G\vert L)$, on a $L_{1}\cap P_{min}\subset L_{1}\cap s_{1}Qs_{1}^{-1}$, $L_{1}\cap P_{min}\subset L_{1}\cap s_{2}Qs_{2}^{-1}$. On en d\'eduit que
$$\zeta(xm',ym)=(\zeta(xm')-s\zeta(ym))_{ s_{1}Ls_{1}^{-1}},$$
o\`u, comme toujours, $\zeta''\mapsto \zeta''_{ s_{1}Ls_{1}^{-1}}$ d\'esigne la projection orthogonale de ${\cal A}_{M_{min}}$ sur ${\cal A}_{ s_{1}Ls_{1}^{-1}}$. 
  Il suffit de minorer $\vert (\zeta(xm')-s\zeta(ym))_{L_{1}}\vert $. Parce que $x,y\in K_{1}$ et $k_{3}m'k_{4}=l'm$, avec $k_{3},k_{4}\in K_{1}$, on a les \'egalit\'es
$$\zeta(xm')_{L_{1}}=\zeta'_{L_{1}}=\zeta_{L_{1}}+H_{L_{1}}(l')=\zeta(ym)_{L_{1}}+H_{L_{1}}(l').$$
Puisque $\sigma(l')<<\sigma(g')\leq Clog(N)$,  on a la majoration $\vert H_{L_{1}}(l')\vert <<log(N)$. Il nous suffit  de montrer que, pour tout $c>0$, on a la minoration
$$\vert (\zeta(ym)-s\zeta(ym))_{L_{1}}\vert> c log(N) $$
pourvu que $c_{1}$ soit assez grand. Soit $c>0$. 
Le m\^eme calcul qu'en (8) montre que l'on a une majoration
$$\beta(\zeta(ym))> c log(N)$$
pour toute racine $\beta$ de $A_{M_{min}}$ dans $\mathfrak{u}_{1}$, pourvu que $c_{1}$ soit assez grand. A fortiori $\beta(\zeta(ym))>0$. Introduisons un \'el\'ement $P'_{min}\in {\cal P}(M_{min})$ tel que $\zeta(ym)\in {\cal A}_{P'_{min}}^+$.  La relation pr\'ec\'edente entra\^{\i}ne que $P'_{min}\subset Q_{1}$. Notons $\Delta'$ la base de racines simples associ\'ee \`a $P'_{min}$, $(\Delta')^{L_{1}}$ le sous-ensemble associ\'e \`a $L_{1}\cap P'_{min}$, $\{\check{\alpha}; \alpha\in \Delta'\}$ l'ensemble de coracines associ\'e \`a $\Delta'$, et $\{\varpi_{\alpha}; \alpha\in \Delta'\}$ la base de ${\cal A}_{M_{min}}$ duale de $\Delta'$. Ecrivons
$$\zeta(ym)=\sum_{\alpha\in \Delta'}\alpha(\zeta(ym))\varpi_{\alpha}.$$
Tous les coefficients sont positifs ou nuls. On a
$$(\zeta(ym)-s\zeta(ym))_{L_{1}}=\sum_{\alpha\in \Delta'}\alpha(\zeta(ym))(\varpi_{\alpha}-s\varpi_{\alpha})_{L_{1}}.$$
On sait que, pour tout $\alpha\in \Delta'$, $\varpi_{\alpha}-s\varpi_{\alpha}$ est combinaison lin\'eaire \`a coefficients positifs ou nuls de coracines $\check{\beta}$ pour $\beta\in \Delta'$. Donc $(\varpi_{\alpha}-s\varpi_{\alpha})_{L_{1}}$ appartient au c\^one ferm\'e engendr\'e par les $\check{\beta}_{L_{1}}$ pour $\beta\in \Delta'-(\Delta')^{L_{1}}$. Si $\alpha\in \Delta'-(\Delta')^{L_{1}}$, l'\'el\'ement  $(\varpi_{\alpha}-s\varpi_{\alpha})_{L_{1}}$ n'est pas nul. En effet, s'il l'\'etait, on aurait
$$(s\varpi_{\alpha})_{L_{1}}=(\varpi_{\alpha})_{L_{1}}=\varpi_{\alpha}.$$
Puisque $s\varpi_{\alpha}$ est de m\^eme norme que $\varpi_{\alpha}$, cela entra\^{\i}nerait
$$s\varpi_{\alpha}=(s\varpi_{\alpha})_{L_{1}}=\varpi_{\alpha}.$$
Mais cette \'egalit\'e est interdite par l'hypoth\`ese (Hyp), d'o\`u la conclusion. Des propri\'et\'es ci-dessus r\'esulte une minoration
$$\vert (\zeta(ym)-s\zeta(ym))_{L_{1}}\vert >>\sum_{\alpha\in \Delta'-(\Delta')^{L_{1}}}\alpha(\zeta(ym)).$$
Mais on a dit ci-dessus que $\alpha(\zeta(ym))>clog(N)$ pour tous les $\alpha$ qui interviennent ici. La minoration cherch\'ee en r\'esulte, d'o\`u (11).

La fonction $f_{5}$ est la transform\'ee de Fourier d'une fonction $C^{\infty}$ \'evalu\'ee en $\zeta(xm',ym)$. Elle est donc \`a d\'ecroissance rapide en cette variable. Gr\^ace \`a (11), pour tout entier $D\geq0$, on a une majoration
$$\vert f_{5}(xm',ym)\vert<< N^{-D} $$
pourvu que $c_{1}$ soit assez grand. On a aussi
$$\vert \Psi_{N,Y,Q_{1},s_{1},s_{2}}(g')\vert<<N^{-D}\int_{K}\int_{K}\int_{M_{min}(F)^+}\int_{K_{1}}\int_{K_{1}} D(m)\kappa_{N}(k_{1}mk_{2})$$
$$\sigma_{M_{min}}^{Q_{1}}(\zeta,{\cal Y})\tau_{Q_{1}}(\zeta-Y_{Q_{1}})\vert f'_{1}(m',x)f'_{2}(m,y)\vert dy\,dx\,dm\,dk_{1}\,dk_{2}.$$
On a
$$\vert f'_{2}(m,y)\vert <<\delta_{Q_{1}}(m)^{-1/2}\delta_{L_{1}\cap s_{2}Qs_{2}^{-1}}^{L_{1}}(l_{s_{2}}(ym))^{1/2}\Xi^{s_{2}Ls_{2}^{-1}}(l_{s_{2}}(ym)).$$ 
D'apr\`es [W2] lemmes II.1.6 et II.1.1, on en d\'eduit l'existence d'un r\'eel $R_{8}\geq0$ tel que
$$\int_{K_{1}}\vert f'_{2}(m,y)\vert <<\delta_{Q_{1}}(m)^{-1/2}\Xi^{L_{1}}(m)<<\Xi^G(m)\sigma(m)^{R_{8}}.$$
On a une majoration analogue pour la fonction $f'_{1}(m',x)$ et on obtient une majoration
$$\vert \Psi_{N,Y,Q_{1},s_{1},s_{2}}(g')\vert<<N^{-D}\int_{K}\int_{K}\int_{M_{min}(F)^+}D(m)\kappa_{N}(k_{1}mk_{2})$$
$$\Xi^G(l'm)\Xi^G(m)\sigma(l'm)^{R_{8}}\sigma(m)^{R_{8}}\,dm\,dk_{1}\,dk_{2}.$$
On a d\'ej\`a major\'e une int\'egrale de ce type: il existe $R_{9}\geq0$ tel qu'elle soit essentiellement major\'ee par $N^{R_{9}}$. En tenant compte du facteur $N^{-D}$ et en se rappelant que $D$ est quelconque, on obtient
$$\vert \Psi_{N,Y,Q_{1},s_{1},s_{2}}(g')\vert<<N^{-R}$$
pourvu que $c_{1}$ soit assez grand. C'est ce que l'on voulait d\'emontrer.

Supposons maintenant que l'hypoth\`ese (Hyp) n'est pas v\'erifi\'ee. Dans ce cas, on va montrer que $\Phi_{N,Y,Q_{1},s_{1},s_{2}}(g')=0$ pourvu que $c_{1}$ soit assez grand. En se rappelant la d\'efinition de ce terme, on voit qu'il suffit de prouver que, si $c_{1}$ est assez grand, l'assertion suivante est v\'erifi\'ee:

(12) soient $\lambda\in i{\cal A}_{L,F}^*$, $m\in M_{min}(F)^+$  et $k_{1},k_{2}\in K$; alors
$$\sum_{e\in {\cal B}_{{\cal O}}^{K_{f}}}\Phi_{s_{1},s_{2}}(e,k_{1}mk_{2},g',\lambda)=0.$$ 

Notons $X(\lambda)$ la somme ci-dessus. C'est une fonction m\'eromorphe en $\lambda$ (plus exactement, c'est la restriction \`a $i{\cal A}_{L,F}^*$ d'une fonction m\'eromorphe sur $({\cal A}_{L}^*\otimes_{{\mathbb R}}{\mathbb C})/i{\cal A}_{L,F}^{\vee}$). Il suffit de montrer qu'elle est nulle pour presque tout $\lambda$. On peut donc supposer que tous les op\'erateurs d'entrelacement qui vont intervenir n'ont pas de p\^oles en $\lambda$ et sont inversibles. Revenons \`a la d\'efinition de $\Phi_{s_{1},s_{2}}(e,k_{1}mk_{2},g',\lambda)$.  On a une \'egalit\'e 
$$\Phi_{s_{1},s_{2}}(e,k_{1}mk_{2},g',\lambda)= $$
$$(J_{\tilde{Q}_{1,s_{1}}\vert s_{1}Qs_{1}^{-1}}((s_{1}\tau)_{s_{1}\lambda})\circ\gamma(s_{1})\circ\pi_{\lambda}(m'k_{4}k_{2})e, e_{1})^{L_{1}}$$
$$( e_{2},J_{\tilde{Q}_{1,s_{2}}\vert s_{2}Qs_{2}^{-1}}((s_{2}\tau)_{s_{2}\lambda})\circ\gamma(s_{2})\circ\pi_{\lambda}(mk_{2})\pi_{\lambda}(f)e)^{L_{1}},$$
o\`u $e_{1}\in {\cal K}^G_{Q_{1,s_{1}},s_{1}\tau}$ et $e_{2}\in {\cal K}^G_{Q_{1,s_{2}},s_{2}\tau}$.  On $\pi_{\lambda}(k_{2})\pi_{\lambda}(f) e=\pi_{\lambda}(f')\pi_{\lambda}(k_{2})e$, o\`u $f'={^{k_{2}}f}$. Posons ${\cal B}_{\natural}^{K_{f}}=\{\pi_{\lambda}(k_{2})e;e\in {\cal B}_{{\cal O}}^{K_{f}}\}$. C'est encore une base orthonorm\'ee de $(E_{Q,\tau}^{G})^{K_{f}}$ et on a
$$X(\lambda)=\sum_{e\in {\cal B}_{\natural}^{K_{f}}}(J_{\tilde{Q}_{1,s_{1}}\vert s_{1}Qs_{1}^{-1}}((s_{1}\tau)_{s_{1}\lambda})\circ\gamma(s_{1})\circ\pi_{\lambda}(m'k_{4})e, e_{1})^{L_{1}}$$
$$( e_{2},J_{\tilde{Q}_{1,s_{2}}\vert s_{2}Qs_{2}^{-1}}((s_{2}\tau)_{s_{2}\lambda})\circ\gamma(s_{2})\circ\pi_{\lambda}(m)\pi_{\lambda}(f')e)^{L_{1}}.$$
 Il existe une fonction $j_{1}(\lambda)$ qui est m\'eromorphe, au m\^eme sens que ci-dessus, telle que
$$J_{\tilde{Q}_{1,s_{1}}\vert s_{1}Qs_{1}^{-1}}((s_{1}\tau)_{s_{1}\lambda})\circ\gamma(s_{1})=j_{1}(\lambda)J_{\tilde{Q}_{1,s_{1}}\vert s\tilde{Q}_{1,s_{2}}s^{-1}}((s_{1}\tau)_{s_{1}\lambda})\circ\gamma(s)\circ
J_{\tilde{Q}_{1,s_{2}}\vert s_{2}Qs_{2}^{-1}}((s_{2}\tau)_{s_{2}\lambda})\circ\gamma(s_{2}).$$
 L'ensemble 
 $$\{J_{\tilde{Q}_{1,s_{2}}\vert s_{2}Qs_{2}^{-1}}((s_{2}\tau)_{s_{2}\lambda})\circ\gamma(s_{2})e; e\in {\cal B}_{\natural}^{K_{f}}\}$$
 est une base de $({\cal K}_{\tilde{Q}_{1,s_{2}},\tau}^G)^{K_{f}}$. Les propri\'et\'es d'adjonction et de composition des op\'erateurs d'entrelacement entra\^{\i}ne qu'elle est orthogonale et que tous ses \'el\'ements ont la m\^eme norme. Notons $j_{2}(\lambda)$ cette norme et divisons tout \'el\'ement de cette base par $\sqrt{j_{2}(\lambda)}$. On obtient  une base orthonorm\'ee  de $({\cal K}_{\tilde{Q}_{1,s_{2}},\tau}^G)^{K_{f}}$ que l'on note ${\cal B}_{\sharp}^{K_{f}}$. On a l'\'egalit\'e
 $$X(\lambda)=j_{1}(\lambda)j_{2}(\lambda)\sum_{e\in {\cal B}_{\sharp}^{K_{f}}}(J_{\tilde{Q}_{1,s_{1}}\vert s\tilde{Q}_{1,s_{2}}s^{-1}}((s_{1}\tau)_{s_{1}\lambda})\circ\gamma(s)\circ Ind_{\tilde{Q}_{1,s_{2}}}^G((s_{2}\tau)_{s_{2}\lambda},m'k_{4})e,e_{1})^{L_{1}}$$
$$(e_{2}, Ind_{\tilde{Q}_{1,s_{2}}}^G((s_{2}\tau)_{s_{2}\lambda},m)Ind_{\tilde{Q}_{1,s_{2}}}^G((s_{2}\tau)_{s_{2}\lambda},f')e)^{L_{1}}.$$
A ce point, le sous-groupe parabolique $Q=LU_{Q}$ n'intervient plus (sauf via les fonctions $j_{1}$ et $j_{2}$). Pour simplifier les notations, on peut supposer $s_{2}=1$ et $Q=\tilde{Q}_{1,s_{2}}$. Auquel cas $s_{1}=s$ et  l'expression pr\'ecedente se simplifie en
 $$X(\lambda)=j_{1}(\lambda)j_{2}(\lambda)\sum_{e\in {\cal B}_{\sharp}^{K_{f}}}(J_{\tilde{Q}_{1,s}\vert sQs^{-1}}((s\tau)_{s\lambda})\circ\gamma(s)\circ\pi_{\lambda}(m'k_{4})e,e_{1})^{L_{1}}$$
$$(e_{2}, \pi_{\lambda}(m)\pi_{\lambda}(f')e)^{L_{1}}.$$
Puisque l'hypoth\`ese (Hyp) n'est pas v\'erifi\'ee, on peut fixer un sous-groupe parabolique $Q'=L'U'\in {\cal F}(M_{min})$ tel que $Q_{1}\subset Q'\not=G$ et $s$ fixe tout \'el\'ement de ${\cal A}_{L'}$. Cela entra\^{\i}ne $s\in W^{L'}$. On a $Q\subset \bar{Q}_{1}\subset \bar{Q}'$, donc aussi $\tilde{Q}_{1,s}\subset \bar{Q}'$. Introduisons les espaces ${\cal K}_{L'\cap Q,\tau}^{L'}$ et ${\cal K}_{L'\cap \tilde{Q}_{1,s},\tau}^{L'}$. On dispose de l'op\'erateur
$$  J_{L'\cap\tilde{Q}_{1,s}\vert L'\cap sQs^{-1}}^{L'}((s\tau)_{s\lambda})\circ\gamma(s):{\cal K}_{L'\cap Q,\tau}^{L'}\to{\cal K}_{L'\cap \tilde{Q}_{1,s},\tau}^{L'}.$$
Posons $\pi'=Ind_{L'\cap Q}^{L'}(\tau_{\lambda})$ et r\'ealisons cette repr\'esentation dans ${\cal K}_{L'\cap Q,\tau}^{L'}$. Pour $e\in {\cal K}_{\bar{Q}',\pi'}^G$, posons
$$b_{1}(e)= \delta_{Q'}(m'k_{4})^{-1/2}\int_{K_{1}}((J_{L'\cap\tilde{Q}_{1,s}\vert L'\cap sQs^{-1}}^{L'}((s\tau)_{s\lambda})\circ\gamma(s)\circ\pi'(m'k_{4})(e(1)))(x),e_{1}(x))dx,$$
$$b_{2}(e)=\delta_{Q'}(m)^{-1/2}\int_{K_{1}}(e_{2}(x),(\pi'(m)(e(1)))(x))dx.$$
Expliquons par exemple la signification de
$$(J_{L'\cap\tilde{Q}_{1,s}\vert L'\cap sQs^{-1}}^{L'}((s\tau)_{s\lambda})\circ\gamma(s)\circ\pi'(m'k_{4})(e(1)))(x).$$
On \'evalue $e$ au point $1$. On obtient un \'el\'ement $e(1)$ de ${\cal K}_{L'\cap Q,\tau}^{L'}$. On applique successivement \`a cet \'el\'ement les op\'erateurs $\pi'(m'k_{4})$ (notons que $m'k_{4}\in L_{1}(F)\subset L'(F)$) puis $J_{L'\cap\tilde{Q}_{1,s}\vert L'\cap sQs^{-1}}^{L'}((s\tau)_{s\lambda})\circ\gamma(s)$. On obtient un \'el\'ement de ${\cal K}_{L'\cap \tilde{Q}_{1,s},s\tau}^{L'}$, que l'on \'evalue au point $x\in K_{1}\subset K\cap L'(F)$.   D\'efinissons une forme sesquilin\'eaire $B$ sur ${\cal K}_{\bar{Q}',\pi'}^G$ par
$$B(e',e)=b_{1}(e')b_{2}(e).$$
 Identifions ${\cal K}_{Q,\tau}^G$ \`a  ${\cal K}_{\bar{Q}',\pi'}^G$. Modulo cette identification, on a les \'egalit\'es
$$(e_{2},\pi_{\lambda}(m)e)^{L_{1}}= b_{2}(e),$$
$$(J_{\tilde{Q}_{1,s}\vert sQs^{-1}}((s\tau)_{s\lambda})\circ\gamma(s)\circ\pi_{\lambda}(m'k_{4})e,e_{1})^{L_{1}}= b_{1}(e),$$
pour tout $e\in {\cal K}_{\bar{Q}',\pi'}^G$.  Alors
$$X(\lambda)=j_{1}(\lambda)j_{2}(\lambda)trace_{B}(Ind_{\bar{Q}'}^G(\pi',f')).$$
Comme $f$, la fonction $f'$ est tr\`es cuspidale. Puisque $b_{1}(e)$ et $b_{2}(e)$ ne d\'ependent que de $e(1)$, la forme $B$ v\'erifie la condition (3) de 2.1. Le lemme de ce paragraphe entra\^{\i}ne $X(\lambda)=0$. Cela prouve (12) et ach\`eve la preuve de la proposition. $\square$
    
  \bigskip
  
  \subsection{Utilisation des calculs spectraux d'Arthur}
  
  Les donn\'ees sont les m\^emes que dans le paragraphe pr\'ec\'edent. Pour tout $\epsilon>0$, on note ${\cal D}(\epsilon)$ l'ensemble des \'el\'ements $Y\in {\cal A}_{P_{min}}^+$ tels que
  $$\inf\{\alpha(Y); \alpha\in \Delta\}>\epsilon sup\{\alpha(Y); \alpha\in \Delta\}.$$
 Pour $L'\in {\cal L}(L)$ et $t\in W^{L'}(L)_{reg}$, notons $\Lambda_{{\cal O}}^{L'}(t)$ l'ensemble des $\lambda\in i{\cal A}_{L}^*$ tels que $t(\tau_{\lambda})\simeq \tau_{\lambda}$. Cet ensemble est stable par translation par $i{\cal A}_{L,F}^{\vee}+i{\cal A}_{L'}^*$. L'ensemble des orbites est fini. Soit $\lambda$ un \'el\'ement de cet ensemble. Arthur d\'efinit  en [A3] p.87 un signe, qu'il note $\epsilon_{\bar{\sigma}}(t)$ et que nous notons $\epsilon_{\tau_{\lambda}}(t)$. Sa d\'efinition ne nous importe pas. Il ne d\'epend que de l'orbite de $\lambda$. On dispose de l'op\'erateur $R_{Q}(t,\tau_{\lambda})$ de  ${\cal K}_{Q,\tau}^G$.  Soit $Q'=L'U'\in {\cal P}(L')$. Posons $Q(Q')=(L'\cap Q)U'$. C'est un \'el\'ement de ${\cal P}(L)$. D\'efinissons sur $i{\cal A}_{L'}^*$  une fonction $\nu\mapsto j_{Q'}(t,\lambda,\nu)$ et, pour $g'\in G(F)$, une fonction $\nu\mapsto d_{Q'}(t,\lambda,g',\nu)$  par
 $$j_{Q'}(t,\lambda,\nu)=\sum_{e\in {\cal B}_{{\cal O}}^{K_{f}}}(e,J_{Q(\bar{Q}')\vert Q}(\tau_{\lambda})^{-1}J_{Q(\bar{Q}')\vert Q}(\tau_{\lambda+\nu})R_{Q}(t,\tau_{\lambda})\pi_{\lambda}(f)e),$$
 $$d_{Q'}(t,\lambda,g',\nu)=(J_{Q(Q')\vert Q}(\tau_{\lambda})^{-1}J_{Q(Q')\vert Q}(\tau_{\lambda+\nu})R_{Q}(t,\tau_{\lambda})e'',\pi_{\lambda}(g')e').$$
 Les op\'erateurs d'entrelacement peuvent avoir des p\^oles, mais, au moins pour $\lambda$ dans un ouvert dense de $\Lambda_{{\cal O}}^{L'}(t)$, les fonctions ci-dessus sont bien d\'efinies pour $\nu$ proche de $0$. Nous utiliserons une troisi\`eme fonction $\nu\mapsto c_{Q'}(\nu)$: c'est celle qui est not\'ee $c_{\bar{Q}'}(\nu)$ en [A3] (12.12). Il nous est inutile de rappeler sa d\'efinition. Les trois familles de fonctions $(j_{Q'}(t,\lambda))_{Q'\in {\cal P}(L')}$, $(d_{Q'}(t,\lambda,g'))_{Q'\in {\cal P}(L')}$ et $(c_{Q'})_{Q'\in {\cal P}(L')}$ sont des $(G,L')$-familles. De plus $c_{Q'}(0)=1$ pour tout $Q'$.  Posons $(jdc)_{Q'}(t,\lambda,g')=j_{Q'}(t,\lambda)d_{Q'}(t,\lambda,g')c_{Q'}$. La famille $((jdc)_{Q'}(t,\lambda,g'))_{Q'\in {\cal P}(L')}$ est encore une $(G,L')$-famille, \`a laquelle on associe un nombre $(jdc)_{L'}(t,\lambda,g')$. Admettons un instant que ce nombre soit fonction $C^{\infty}$ de $\lambda$. Rappelons par ailleurs que, puisque $t\in W^{L'}(L)_{reg}$, l'op\'erateur $t-1$ sur ${\cal A}_{L}/{\cal A}_{L'}$ est inversible et a donc un d\'eterminant $det(t-1)_{{\cal A}_{L}/{\cal A}_{L'}}$ non nul. Posons
 $$\Phi(g')=[i{\cal A}_{{\cal O}}^{\vee};i{\cal A}_{L,F}^{\vee}]^{-1}\sum_{L'\in {\cal L}(L)}\sum_{t\in W^{L'}(L)_{reg}}\sum_{\lambda\in \Lambda_{{\cal O}}^{L'}(t)/(i{\cal A}_{L,F}^{\vee}+i{\cal A}_{L'}^*)}\epsilon_{\tau_{\lambda}}(t)\vert det(t-1)_{{\cal A}_{L}/{\cal A}_{L'}}\vert ^{-1}$$
 $$\int_{i{\cal A}_{L',F}^*}(jdc)_{L'}(t,\lambda+\mu,g')\varphi(\lambda+\mu)d\mu.$$
 
 \ass{Proposition}{(i) Pour tout $L'\in {\cal L}(L)$, la fonction $\lambda\mapsto (jdc)_{L'}(t,\lambda,g')$ est $C^{\infty}$.
 
 (ii) Soit $\epsilon>0$ et $R\geq1$ un entier. On a une majoration
 $$\vert \Phi_{Y}(g')-\Phi(g')\vert<<\sigma(g')^R\Xi^G(g')\vert Y\vert ^{-R}$$
 pour tout $g'\in G(F)$ et tout $Y\in {\cal D}(\epsilon)\cap {\cal A}_{M_{min},F}$.}
 
 Preuve. Fixons $\epsilon$ et $R$. Oublions d'abord $g'$, c'est-\`a-dire supposons $g'=1$. Supposons aussi que la fonction $\varphi$ est constante de valeur $1$. Dans [A3] p.69, Arthur \'etudie une expression qu'il note $K^T(f)$. C'est une somme sur $M\in {\cal L}(M_{min})$, $\sigma\in \{\Pi_{2}(M(F))\}$ d'expressions qui sont des int\'egrales sur $i{\cal A}_{M,F}^*\times G(F)$.   Dans le cas o\`u le $T$ d'Arthur est notre $Y$ et o\`u le couple $(M,\sigma)$ d'Arthur est \'egal \`a notre couple $(L,\tau)$, cette expression est presque notre terme $\Phi_{Y}(1)$. Plus pr\'ecis\'ement, notre terme d\'epend d'\'el\'ements fix\'es $e',e''$ et le terme d'Arthur est \'egal \`a la somme de nos termes $\Phi_{Y}(1)$ associ\'es \`a un nombre fini de couples $(e',e'')$.  Dans [A3], proposition 11.3, p.88, figure une expression $\tilde{J}(f)$. C'est une somme sur les couples $(M,\sigma)$ comme ci-dessus d'expressions qui, pour le couple $(M,\sigma)=(L,\tau)$, sont presqu'\'egales \`a notre terme $\Phi(1)$. Plus exactement, ce terme est la somme de termes $\Phi(1)$ associ\'es aux m\^emes couples $(e',e'')$ que ci-dessus.   Arthur d\'emontre  en [A3] corollaire 10.4 que les fonctions $\lambda\mapsto (jcd)_{L'}(t,\lambda,1)$ sont $C^{\infty}$. Entre les pages 69 et 88 de [A3], il d\'emontre le r\'esultat suivant. Il existe une fonction $T\mapsto J^T(f)$ qui v\'erifie les trois conditions
 
 (1) il existe un entier $D\geq1$ et, pour tout $\xi\in (\frac{1}{D}{\cal A}_{M_{min},F})/{\cal A}_{M_{min},F}$, un polyn\^ome $T\mapsto q_{\xi}(T)$, de sorte que $J^T(f)=\sum_{\xi \in (\frac{1}{D}{\cal A}_{M_{min},F})/{\cal A}_{M_{min},F}}exp(\xi(T))q_{\xi}(T)$ pour tout $T\in {\cal A}_{M_{min},F}$;
 
 (2) $\tilde{J}(f)=q_{0}(0)$;
 
 (3)   on a une majoration $\vert K^T(f)-J^T(f)\vert <<\vert T\vert ^{-R}$ pour tout $T\in {\cal D}(\epsilon)\cap {\cal A}_{M_{min},F}$.
  
  Cf. la preuve du lemme 11.1 de [A3]. En inspectant la preuve, on voit d'une part que la somme sur les couples $(M,\sigma)$ ne sert \`a rien dans ce passage: la preuve se fait terme par terme. D'autre part sommer sur un ensemble fini de couples ne sert \`a rien non plus, la m\^eme preuve s'applique pour chaque couple. En revenant \`a nos notations, cette preuve montre donc qu'il existe une fonction $Y\mapsto \Phi^Y$ sur ${\cal A}_{M_{min},F}$, de la forme
  $$(4) \qquad \Phi^Y= \sum_{\xi \in (\frac{1}{D}{\cal A}_{M_{min},F})/{\cal A}_{M_{min},F}}exp(\xi(Y))q_{\xi}(Y),$$
  telle que $\Phi(1)=p_{0}(0)$ et telle que l'on ait une majoration
  $$(5) \qquad \vert \Phi_{Y}(1)-\Phi^Y\vert <<\vert Y\vert ^{-R}$$
  pour tout $Y\in {\cal D}(\epsilon)\cap {\cal A}_{M_{min},F}$. Montrons qu'en fait 
 
 (6) $p_{0}$ est constant et  $p_{\xi}=0$ si $\xi\not=0$.
 
  Remarquons que, pour $Y\in {\cal D}(\epsilon)$, on a des majorations
  $$\vert Y\vert << sup\{\alpha(Y); \alpha\in \Delta\}<<\vert Y\vert .$$
  On peut donc remplacer $sup\{\alpha(Y); \alpha\in \Delta\}$ par $\vert Y\vert $ dans l'\'enonc\'e de la proposition 6.6. Fixons $c_{1}$ et $c_{2}$ v\'erifiant cet \'enonc\'e modifi\'e. Soit $Y\in {\cal D}(\epsilon)\cap {\cal A}_{M_{min},F}$. Notons $N_{Y}$ la partie enti\`ere de $2c_{2}\vert Y\vert +1$. Si $\vert Y\vert $ est assez grand, le couple $(N_{Y},Y)$ v\'erifie les conditions de la proposition 6.6. Plus g\'en\'eralement, il en est de m\^eme du couple $(N_{Y},Y')$ pour tout $Y'\in {\cal D}(\epsilon)\cap {\cal A}_{M_{min}}$ tel que $\vert Y-Y'\vert \leq\vert Y\vert /2$. Pour un tel $Y'$, on a donc
  $$\vert \Phi_{Y'}(1)-\Phi_{Y}(1)\vert<<\vert \Phi_{Y'}(1)-\Phi_{N_{Y}}(1)\vert +\vert \Phi_{N_{Y}}(1)-\Phi_{Y}(1)\vert <<N_{Y}^{-R}<<\vert Y\vert ^{-R}.$$
  En appliquant (5), on a aussi
  $$\vert \Phi^{Y'} -\Phi^Y\vert <<\vert Y\vert ^{-R}.$$
 Si  $Y\mapsto \Phi^Y$ ne v\'erifie pas (6), c'est-\`a-dire n'est pas constante, on peut fixer  $Y_{0}\in {\cal A}_{M_{min},F}$ tel que la fonction $Y\mapsto \Phi^{Y+Y_{0}}-\Phi^Y$ soit non nulle. Cette fonction est encore de la forme (4). Pour $Y$ assez grand, le point $Y'=Y+Y_{0}$ v\'erifie $\vert Y-Y'\vert \leq\vert Y\vert /2$. La fonction est donc essentiellement major\'ee par $\vert Y\vert ^{-R}$. Mais une fonction de la forme (4) ne peut v\'erifier cette majoration que si elle est nulle. Cela contredit le choix de $Y_{0}$, d'o\`u (6).
 
 Gr\^ace \`a (6), on a $\Phi(1)=\Phi^Y$ et la majoration (5) est celle que l'on voulait prouver.
 
 Si $\varphi$  n'est plus la fonction constante de valeur $1$, on inspecte la preuve d'Arthur et on voit que l'on peut glisser la fonction $\varphi$ tout le long de cette preuve. Le r\'esultat est celui annonc\'e. Bien s\^ur, la constante qui se trouve implicitement dans la majoration de l'\'enonc\'e d\'epend de $\varphi$. Cela ne nous g\^ene pas pourvu que $\varphi$ soit fix\'ee mais, pour la suite du raisonnement, pr\'ecisons tout-de-m\^eme cette d\'ependance. Pour tout entier $k\geq0$, fixons une base ${\cal X}_{k}$ de l'espace des op\'erateurs diff\'erentiels sur $i{\cal A}_{L}^*$, \`a coefficients constants et de degr\'e $\leq k$. Posons
 $$\vert \varphi\vert _{k}=sup\{\vert (X\varphi)(\lambda); \lambda\in i{\cal A}_{L,F}^*, X\in {\cal X}_{k}\}.$$
 Dans la preuve d'Arthur, les approximations interviennent \`a deux endroits. D'abord dans l'utilisation du th\'eor\`eme 8.1. Cette approximation porte uniquement sur l'int\'egrale int\'erieure sur $G(F)$. Quand on int\`egre ensuite sur $i{\cal A}_{L,F}^*$, la constante implicite est simplement multipli\'ee par $\vert \varphi\vert _{0}$. Il y a ensuite les approximations de la page 80 qui se r\'ef\`erent elles-m\^emes \`a [A6], p. 1306,1307. D'apr\`es cette derni\`ere r\'ef\'erence, la constante implicite est de la forme
 $$sup\{(\Psi\varphi)\hat{}(Z)\vert Z\vert ^{R};Z\in {\cal A}_{L,F}\},$$
 o\`u $\Psi$ est une certaine fonction $C^{\infty}$ sur $i{\cal A}_{L,F}^*$ ind\'ependante de $\varphi$ et $(\Psi\varphi)\hat{}$ est la transform\'ee de Fourier de $\Psi\varphi$. Comme on le sait bien, le terme ci-dessus est essentiellement born\'e par $\vert \varphi\vert _{R}$. Finalement, la majoration de l'\'enonc\'e se pr\'ecise en
  $$(7)\qquad \vert \Phi_{Y}(1)-\Phi(1)\vert\leq c\vert \varphi\vert _{R}\vert Y\vert ^{-R}$$
  o\`u $c$ est ind\'ependant de $\varphi$.

 Passons au cas g\'en\'eral o\`u $g'$ est quelconque. Fixons un sous-groupe ouvert compact $K_{0}$ de $K$ tel que $e''$ soit invariant par $K_{0}$. On utilise le fait que la fonction $g\mapsto u(g,Y)$ est invariante \`a gauche par $K$ (ce que n'\'etait pas la fonction $g\mapsto \kappa_{N}(g)$ des paragraphes pr\'ec\'edents). Dans la d\'efinition de $\Phi_{Y}(g')$ donn\'ee au d\'ebut du paragraphe 6.6, on remplace la variable $g$ par $kg$ et on int\`egre sur $k\in K_{0}$, en divisant le tout par $mes(K_{0})$. On obtient une expression analogue \`a $\Phi_{Y}(g')$, o\`u le terme $\pi_{\lambda}(g')e'$ est remplac\'e par 
 $$mes(K_{0})^{-1}\int_{K_{0}}\pi_{\lambda}(kg')e'.$$
 Fixons une base orthonorm\'ee ${\cal B}_{{\cal O}}^{K_{0}}$ de ${\cal K}_{Q,\tau}^{K_{0}}$. On peut encore remplacer le terme ci-dessus par
 $$\sum_{e_{0}\in {\cal B}_{{\cal O}}^{K_{0}}}(e_{0},\pi_{\lambda}(g')e')e_{0}.$$
 Alors $\Phi_{Y}(g')$ est une somme de termes analogues o\`u le triplet $(g',e',\varphi)$ est remplac\'e par $(1,e_{0},\varphi')$, o\`u
 $$\varphi'(\lambda)=\varphi(\lambda)(e_{0},\pi_{\lambda}(g')e').$$
 Une d\'ecomposition analogue vaut pour le terme $\Phi(g')$. En utilisant la majoration (7), on voit que, pour obtenir la majoration de l'\'enonc\'e, il nous reste \`a prouver que l'on a une majoration
 $$\vert \varphi'\vert _{R}<<\sigma(g')^R\Xi^G(g')$$
 pour toute fonction $\varphi'$ comme ci-dessus. Il suffit de prouver une majoration analogue pour la fonction $\varphi''(\lambda)=(e_{0},\pi_{\lambda}(g')e')$. On a
 $$\varphi''(\lambda)=\int_{K}(e_{0}(k),\tau(l_{Q}(kg'))e'(k_{Q}(kg'))\delta_{Q}(l_{Q}(kg'))^{1/2}exp(\lambda(H_{Q}(kg')))dk.$$
 Appliquer un op\'erateur diff\'erentiel $X\in {\cal X}_{R}$ revient \`a glisser dans l'int\'egrale un terme $P(H_{Q}(kg'))$, o\`u $P$ est un polyn\^ome de degr\'e $\leq R$. On a une majoration
 $$\vert H_{Q}(kg')\vert <<\sigma(kg')=\sigma(g').$$
 D'o\`u les majorations
 $$\vert \varphi''\vert <<\sigma(g')^R\int_{K}\vert (e_{0}(k),\tau(l_{Q}(kg'))e'(k_{Q}(kg'))\vert \delta_{Q}(l_{Q}(kg'))^{1/2}dk$$
$$<<\sigma(g')^R\int_{K} \Xi^L(l_{Q}(kg'))\delta_{Q}(l_{Q}(kg'))^{1/2}dk$$   
  $$<<\sigma(g')^R\Xi^G(g')$$
  d'apr\`es [W2] lemme II.1.6. C'est ce que l'on voulait d\'emontrer. $\square$
  
  \bigskip
  
  \subsection{Simplification de $\Phi(g')$}
  
  Pour $L'\in {\cal L}(L)$ et $\lambda\in i{\cal A}_{L}^*$, on d\'efinit les groupes $W^{L'}(\tau_{\lambda})$, $(W^{L'})'(\tau_{\lambda})$ et $R^{L'}(\tau_{\lambda})$: ce sont les analogues de $W(\tau_{\lambda})$, $W'(\tau_{\lambda})$ et $R(\tau_{\lambda})$ quand on remplace $G$ par $L'$. Notons $\Lambda_{{\cal O},ell}^{L'}$ l'ensemble des $\lambda\in i{\cal A}_{L}^*$ tels que $R^{L'}(\tau_{\lambda})\cap W^{L'}(L)_{reg}\not=\emptyset$. Cet ensemble est stable par translations par $i{\cal A}_{L,F}^{\vee}+i{\cal A}_{L'}^*$. Soit $\lambda$ un \'el\'ement de cet ensemble. On a d\'ecrit le groupe  $R^{L'}(\tau_{\lambda})$ en 4.1. Notons $R^{L'}(\tau_{\lambda})^{\vee}$ son groupe dual.  La repr\'esentation $Ind_{L'\cap Q}^{L'}(\tau_{\lambda})$ se d\'ecompose en somme de sous-repr\'esentations irr\'eductibles $Ind_{L'\cap Q}^{L'}(\tau_{\lambda},\zeta)$ pour $\zeta\in R^{L'}(\tau_{\lambda})^{\vee}$.  On a conform\'ement la d\'ecomposition en somme orthogonale
  $${\cal K}_{L'\cap Q,\tau}^{L'}=\oplus_{\zeta\in R^{L'}(\tau_{\lambda})^{\vee}}{\cal K}_{L'\cap Q,\tau_{\lambda},\zeta}^{L'}.$$
  Fixons $S'=L'U'\in {\cal P}(L')$ et, comme dans le paragraphe pr\'ec\'edent, notons $Q(S')=(L'\cap Q)U'\in {\cal P}(L)$. L'op\'erateur $R_{Q(S')\vert Q}(\tau_{\lambda})$ est une isom\'etrie de ${\cal K}_{Q,\tau}^G$ sur ${\cal K}_{Q(S'),\tau}^G$. D'autre part, ce dernier espace s'identifie \`a un espace de fonctions de $K$ dans ${\cal K}_{L'\cap Q,\tau}^{L'}$. La d\'ecomposition ci-dessus de cet espace induit une d\'ecomposition orthogonale
  $$(1) \qquad {\cal K}_{Q(S'),\tau}^G=\oplus_{\zeta\in R^{L'}(L)}{\cal K}_{Q(S'),\tau_{\lambda},\zeta}^G.$$
 Notons $proj_{\lambda,\zeta}$ la projection de ${\cal K}_{Q(S'),\tau}^G$ sur ${\cal K}_{Q(S'),\tau_{\lambda},\zeta}^G$. Remarquons que ces sous-espaces et ces projections ne d\'ependent que de l'orbite de $\lambda$. Rappelons d'autre part que l'on note $a_{L'}$ la dimension de ${\cal A}_{L'}$.
 
 \ass{Lemme}{Pour tout $g'\in G(F)$, on a l'\'egalit\'e
 $$\Phi(g')=[i{\cal A}_{{\cal O}}^{\vee}:i{\cal A}_{L,F}^{\vee}]^{-1}\sum_{L'\in {\cal L}(L)}(-1)^{a_{L'}}\sum_{\lambda\in \Lambda^{L'}_{{\cal O},ell}/(i{\cal A}_{L,F}^{\vee}+i{\cal A}_{L'}^*)}\vert R^{L'}(\tau_{\lambda})\vert 2^{a_{L'}-a_{L}}\sum_{\zeta\in R^{L'}(\tau_{\lambda})^{\vee}}$$
 $$\int_{i{\cal A}^*_{L',F}}(proj_{\lambda,\zeta}\circ R_{Q(S')\vert Q}(\tau_{\lambda+\mu})e'',proj_{\lambda,\zeta}\circ R_{Q(S')\vert Q}(\tau_{\lambda+\mu})\circ Ind_{Q}^G(\tau_{\lambda+\mu},g')e')$$
 $$J_{L'}^G(Ind_{L'\cap Q}^{L'}(\tau_{\lambda+\mu},\zeta),f)\varphi(\lambda+\mu)d\mu.$$}
 
 Preuve. Fixons $L'\in {\cal L}(L)$, $t\in W^{L'}(L)_{reg}$ et $\lambda\in \Lambda^{L'}_{{\cal O}}(t)$. Consid\'erons le terme $(jdc)_{L'}(t,\lambda,g')$. On commence par remplacer, dans les d\'efinitions des $(G,L')$-familles $(j_{Q'}(t,\lambda))_{Q'\in {\cal P}(L')}$ et $(d_{Q'}(t,\lambda,g'))_{Q'\in {\cal P}(L')}$, les op\'erateurs d'entrelacement par des op\'erateurs normalis\'es. Cela multiplie ces familles par des $(G,L')$-familles \`a valeurs scalaires. Quitte \`a multiplier la  $(G,L')$-famille $(c_{Q'})_{Q'\in {\cal P}(L')}$ 
par ces familles, on retrouve une expression similaire \`a celle de d\'epart (la famille qui remplace $(c_{Q'})_{Q'\in {\cal P}(L')}$ d\'epend de $\lambda$).  On peut donc consid\'erer que l'on a
$$j_{Q'}(t,\lambda,\nu)=\sum_{e\in {\cal B}_{{\cal O}}^{K_{f}}}(e,R_{Q(\bar{Q}')\vert Q }(\tau_{\lambda})^{-1}R_{Q(\bar{Q}')\vert Q}(\tau_{\lambda+\nu})R^{L'}_{Q}(t,\tau_{\lambda})\pi_{\lambda}(f)e).$$
L'ensemble $\{R_{Q(S')\vert Q}(\tau_{\lambda})(e); e\in {\cal B}_{{\cal O}}^{K_{f}}\}$ est une base orthonorm\'ee de ${\cal K}_{Q(S'),\tau}^G$. Notons-la ${\cal B}_{{\cal O},Q(S')}^{K_{f}}$. Posons
$$j'_{Q'}(t,\lambda,\nu)=\sum_{e\in {\cal B}_{{\cal O},Q(S')}^{K_{f}}}(e,R_{Q(\bar{Q}')\vert Q(S')}(\tau_{\lambda})^{-1}R_{Q(\bar{Q}')\vert Q(S')}(\tau_{\lambda+\nu})R^{L'}_{Q(S')}(t,\tau_{\lambda})Ind_{Q(S')}^G(\tau_{\lambda},f)e).$$
On a l'\'egalit\'e
$$j'_{Q'}(t,\lambda,\nu)=\sum_{e\in {\cal B}_{{\cal O}}^{K_{f}}}(e,R_{Q(\bar{Q}')\vert Q }(\tau_{\lambda})^{-1}R_{Q(\bar{Q}')\vert Q}(\tau_{\lambda+\nu})\pi_{\lambda}(f)r(\lambda,\nu)e),$$
o\`u 
$$r(\lambda,\nu)=R_{Q(S')\vert Q}(\lambda+\nu)^{-1}R_{Q(S')\vert Q}(\lambda).$$
Le point est que cet op\'erateur ne d\'epend pas de $Q'$ et que $r(\lambda,0)$ est l'identit\'e. Une propri\'et\'e famili\`ere des $(G,L')$-familles entra\^{\i}ne que  l'on a l'\'egalit\'e
$$(j')_{L'}^{Q''}(t,\lambda)=j_{L'}^{Q''}(t,\lambda)$$
pour tout $Q''\in {\cal F}(L')$. D'apr\`es les formules de descente, on peut donc remplacer la famille $(j_{Q'}(t,\lambda))_{Q'\in {\cal P}(L')}$ par  $(j'_{Q'}(t,\lambda))_{Q'\in {\cal P}(L')}$. De la m\^eme fa\c{c}on, on peut remplacer la famille $(d_{Q'}(t,\lambda,g'))_{Q'\in {\cal P}(L')}$ par $(d'_{Q'}(t,\lambda,g'))_{Q'\in {\cal P}(L')}$, o\`u
$$d'_{Q'}(t,\lambda,g',\nu)=(R_{Q(Q')\vert Q(S')}(\tau_{\lambda})^{-1}R_{Q(Q')\vert Q(S')}(\tau_{\lambda+\nu})R_{Q(S')}(t,\tau_{\lambda})R_{Q(S')\vert Q}(\tau_{\lambda})e'', $$
$$R_{Q(S')\vert Q}(\tau_{\lambda})Ind_{Q}^G(\tau_{\lambda},g')e').$$
Autrement dit, quitte \`a remplacer les \'el\'ements $Ind_{Q}^G(\tau_{\lambda})e'$ et $e''$ par $R_{Q(S')\vert Q}(\tau_{\lambda})Ind_{Q}^G(\tau_{\lambda})e'$ et $R_{Q(S')\vert Q}(\tau_{\lambda})e''$, on peut remplacer le parabolique $Q$ par $Q(S')$.   Les consid\'erations qui pr\'ec\`edent l'\'enonc\'e sont valables ind\'ependamment de l'hypoth\`ese $R^{L'}(\tau_{\lambda})\cap W^{L'}(L)_{reg}\not=\emptyset$. On a la d\'ecomposition (1). On peut supposer que la base ${\cal B}_{{\cal O},S'}^{K_{f}}$ est r\'eunion de bases des diff\'erents sous-espaces. L'op\'erateur $R_{Q(S')}(t,\tau_{\lambda})$ agit par $\zeta(t)$ sur le sous-espace ${\cal K}_{Q(S'),\tau,\zeta}^G$ (o\`u le caract\`ere $\zeta$ de $R^{L'}(\tau_{\lambda})$ est \'etendu en un caract\`ere de $W^{L'}(\tau_{\lambda})$ trivial sur $(W^{L'})'(\tau_{\lambda})$). On en d\'eduit une \'egalit\'e
$$j'_{Q'}(t,\lambda,\nu)=\sum_{\zeta\in R^{L'}(\tau_{\lambda})^{\vee}}\zeta(t)j'_{Q'}(\lambda,\zeta,\nu),$$
avec une d\'efinition plus ou moins \'evidente de ce dernier terme. On a alors
$$(j')_{L'}^{Q''}(t,\lambda)=(-1)^{a_{L'}-a_{L''}}\sum_{\zeta\in R^{L'}(\tau_{\lambda})^{\vee}}\zeta(t)J_{L'}^{\bar{Q}''}(Ind_{L'\cap Q}^{L'}(\tau_{\lambda},\zeta),f)$$
pour tout $Q''=L''U''\in {\cal F}(L')$. Le signe $(-1)^{a_{L'}-a_{L''}}$ ainsi que le $\bar{Q}''$ en exposant viennent de ce que c'est le parabolique $Q(\bar{Q}')$ qui intervient dans la d\'efinition de $j'_{Q'}(t,\lambda,\nu)$, cf. [A3] p.92. Supposons $R^{L'}(\tau_{\lambda})\cap W^{L'}(L)_{reg}=\emptyset$, c'est-\`a-dire $\lambda\not\in \Lambda^{L'}_{{\cal O},ell}$. Dans ce cas, toutes les repr\'esentations $Ind_{L'\cap Q}^{L'}(\tau_{\lambda},\zeta)$ sont induites et les termes ci-dessus sont nuls (lemme 2.2(ii)). Gr\^ace \`a la formule de descente, on obtient
 
 (2) $(jdc)(t,\lambda,g')=0$ si $\lambda\not\in \Lambda^{L'}_{{\cal O},ell}$.
 
  Supposons maintenant $R^{L'}(\tau_{\lambda})\cap W^{L'}(L)_{reg}\not=\emptyset$, c'est-\`a-dire $\lambda\in \Lambda^{L'}_{{\cal O},ell}$. Le  lemme 2.2(i) et la formule de descente entra\^{\i}ne l'\'egalit\'e
$$ (jdc)(t,\lambda,g')=(j')_{L'}^G(t,\lambda)d'_{Q'}(t,\lambda,g',0)c_{Q'}(\lambda,0)$$
o\`u $Q'$ est un \'el\'ement quelconque de ${\cal P}(L')$. On a $c_{Q'}(\lambda,0)=1$, 
$$d'_{Q'}(t,\lambda,g',0)=(R_{Q(S')}(t,\tau_{\lambda})R_{Q(S')\vert Q}(\tau_{\lambda})e'', R_{Q(S')\vert Q}(\tau_{\lambda})Ind_{Q}^G(\tau_{\lambda})e')$$
et
$$(j')_{L'}^G(t,\lambda)=(-1)^{a_{L'}}\sum_{\zeta\in R^{L'}(\tau_{\lambda})^{\vee}}\zeta(t)J_{L'}^{G}(Ind_{L'\cap Q}^{L'}(\tau_{\lambda},\zeta),f).$$
Rappelons que l'hypoth\`ese $R^{L'}(\tau_{\lambda})\cap W^{L'}(L)_{reg}\not=\emptyset$ entra\^{\i}ne que $(W^{L'})'(\tau_{\lambda})=\{1\}$, donc $R^{L'}(\tau_{\lambda})=W^{L'}(\tau_{\lambda})$. De plus, $R^{L'}(\tau_{\lambda})\cap W^{L'}(L)_{reg}$ poss\`ede un unique \'el\'ement. Puisque $t$ appartient \`a cette intersection, cet unique \'el\'ement est $t$. Soit $x\in R^{L'}(\tau_{\lambda})$, $x\not=t$. Consid\'erons la repr\'esentation virtuelle
$$\sum_{\zeta\in R^{L'}(\tau_{\lambda})^{\vee}}\zeta(x)Ind_{L'\cap Q}^{L'}(\tau_{\lambda}).$$
D'apr\`es [A5] proposition 2.1(b), c'est une somme, \`a coefficients dans ${\mathbb Z}$, de repr\'esentations induites. D'apr\`es le lemme 2.2(ii), on a donc
$$\sum_{\zeta\in R^{L'}(\tau_{\lambda})^{\vee}}\zeta(x)J_{L'}^{G}(Ind_{L'\cap Q}^{L'}(\tau_{\lambda},\zeta),f)=0.$$
Il en r\'esulte que $\zeta(t)J_{L'}^{G}(Ind_{L'\cap Q}^{L'}(\tau_{\lambda},\zeta),f)$ est ind\'ependant de $\zeta$.  On en d\'eduit l'\'egalit\'e
$$(j')_{L'}^G(t,\lambda)=(-1)^{a_{L'}}\vert R^{L'}(\tau_{\lambda})\vert \zeta(t)J_{L'}^{G}(Ind_{L'\cap Q}^{L'}(\tau_{\lambda},\zeta),f)$$
pour tout $\zeta\in R^{L'}(\tau_{\lambda})^{\vee}$. D'autre part, on peut d\'ecomposer les \'el\'ements $R_{Q(S')\vert Q}(\tau_{\lambda})Ind_{Q}^G(\tau_{\lambda})e'$ et $R_{Q(S')\vert Q}(\tau_{\lambda})e''$ selon la d\'ecomposition (1). Il en r\'esulte l'\'egalit\'e
$$d'_{Q'}(t,\lambda,g',0)=\sum_{\zeta\in R^{L'}(\tau_{\lambda})^{\vee}}\zeta(t)(proj_{\lambda,\zeta}\circ R_{Q(S')\vert Q}(\tau_{\lambda})e'',proj_{\lambda,\zeta}\circ R_{Q(S')\vert Q}(\tau_{\lambda})\circ Ind_{Q}^G(\tau_{\lambda},g')e').$$
Des deux \'egalit\'es pr\'ec\'edentes r\'esulte la relation

(3) si $\lambda\in \Lambda^{L'}_{{\cal O},ell}$,
$$(jdc)_{L'}(t,\lambda,g')=(-1)^{a_{L'}}\vert R^{L'}(\tau_{\lambda})\vert \sum_{\zeta\in R^{L'}(\tau_{\lambda})^{\vee}}$$
$$(proj_{\lambda,\zeta}\circ R_{Q(S')\vert Q}(\tau_{\lambda})e'',proj_{\lambda,\zeta}\circ R_{Q(S')\vert Q}(\tau_{\lambda})\circ Ind_{Q}^G(\tau_{\lambda},g')e')J_{L'}^G(Ind_{L'\cap Q}^{L'}(\tau_{\lambda},\zeta),f).$$

D'autre part, le signe $\epsilon_{\tau_{\lambda}}(t)$ vaut $1$ parce que $(W^{L'})'(\tau_{\lambda})=\{1\}$, cf. [A3] p.95. Enfin, on a d\'ecrit $t$ en 4.2 et on voit que $\vert det(t-1)_{{\cal A}_{L}/{\cal A}_{L'}}\vert ^{-1}=2^{a_{L'}-a_{L}}$. Il suffit de reporter ces \'egalit\'es et celles des relations (2) et (3) dans la d\'efinition de $\Phi(g')$ pour obtenir l'\'egalit\'e de l'\'enonc\'e. $\square$

 \bigskip 
  
  \subsection{Evaluation d'une limite}

  \ass{Lemme}{On a l'\'egalit\'e
  $$lim_{N\to \infty}I_{L,{\cal O},N,C}(\theta_{\rho},f)=[i{\cal A}_{{\cal O}}^{\vee}:i{\cal A}_{L,F}^{\vee}]^{-1}\sum_{L'\in {\cal L}(L)}(-1)^{a_{L'}}\sum_{ \lambda\in \Lambda^{L'}_{{\cal O},ell} /(i{\cal A}_{L,F}^{\vee}+i{\cal A}_{L'}^*)}$$
  $$\vert R^{L'}(\tau_{\lambda})\vert 2^{a_{L'}-a_{L}}\sum_{\zeta\in R^{L'}(\tau_{\lambda})^{\vee}; m(Ind_{L'\cap Q}^{L'}(\tau_{\lambda},\zeta),\rho)=1} \int_{i{\cal A}_{L',F}^*}J_{L'}^G(Ind_{L'\cap Q}^{L'}(\tau_{\lambda+\mu},\zeta),f)d\mu.$$}
  
  Remarque. Le nombre $m(Ind_{L'\cap Q}^{L'}(\tau_{\lambda},\zeta),\rho)$ a \'et\'e d\'efini en 6.1.
  
  Preuve. Consid\'erons la d\'efinition de $I_{L,{\cal O},N,C}(\theta_{\rho},f)$ donn\'ee avant le lemme 6.4. Il y intervient des objets $e_{j}$, $e_{j}$ et $\varphi_{j}$ pour $j=1,...,n$. Dans les paragraphes pr\'ec\'edents, on a introduit des fonctions $\Phi_{N}(g')$, $\Phi_{Y}(g')$ et $\Phi(g')$ qui d\'ependaient de choix d'\'el\'ements $e'$, $e''$ et d'une fonction $\varphi$. On note $\Phi_{N,j}(g')$, $\Phi_{Y,j}(g')$, $\Phi_{j}(g')$ ces fonctions relatives \`a $e'=e_{j}$, $e''=e'_{j}$, $\varphi=\varphi_{j}$. On a alors
    $$I_{L,{\cal O},N,C}(\theta_{\rho},f)=\sum_{j=1,...,n}\int_{H(F)U(F)_{c}}{\bf 1}_{\sigma< Clog(N)}(hu)(\rho(h)\epsilon'_{j},\epsilon_{j})\bar{\xi}(u)\Phi_{N,j}(hu)du\,dh.$$
  Fixons un r\'eel $\epsilon$ tel que $0<\epsilon<1$ et consid\'erons un entier $R>0$ que nous pr\'eciserons par la suite. Introduisons des constantes $c_{1},c_{2}$  qui v\'erifient les conditions de la proposition 6.6 pour chaque couple de fonctions $(\Phi_{N,j}(g'),\Phi_{Y,j}(g'))$.  Il y a une constante $c_{3}>0$ telle que, pour tout $N$, il existe $Y\in {\cal A}_{M_{min},F}$ tel que $c_{3}N<\alpha(Y)<c_{2}N$ pour tout $\alpha\in \Delta$. Fixons un tel $c_{3} $ et, pour tout $N$, un \'el\'ement $Y_{N}$ v\'erifiant ces in\'egalit\'es. Si $N$ est assez grand, $Y_{N}$ v\'erifie les hypoth\`eses de la proposition 6.6 et celles de la proposition 6.7 (pour chacune de nos fonctions $\Phi_{N,j}(g')$ etc...). Ces propositions entra\^{\i}nent que l'on a une majoration
  $$\vert \Phi_{N,j}(g')-\Phi_{j}(g')\vert<< (1+\sigma(g')^R\Xi^G(g'))N^{-R}$$
  pour tout $j$, tout $N$ assez grand et  tout $g'\in G(F)$ tel que ${\bf 1}_{\sigma< Clog(N)}(g')=1$. On peut oublier le terme $\Xi^G(g')$ qui est born\'e. Posons
  $$X_{N}=\sum_{j=1,...,n}\int_{H(F)U(F)_{c}}{\bf 1}_{\sigma< Clog(N)}(hu)(\rho(h)\epsilon'_{j},\epsilon_{j})\bar{\xi}(u)\Phi_{j}(hu)du\,dh.$$
  Alors
  $$\vert I_{L,{\cal O},N,C}(\theta_{\rho},f)-X_{N}\vert <<N^{-R}\int_{H(F)U(F)_{c}}{\bf 1}_{\sigma< Clog(N)}(hu)\Xi^H(h)\sigma(hu)^Rdu\,dh$$
  $$<<log(N)^RN^{-R}\int_{H(F)U(F)_{c}}{\bf 1}_{\sigma< Clog(N)}(hu)du\,dh.$$
  On peut choisir $C'>0$ tel que la condition $\sigma(hu)< Clog(N)$ entra\^{\i}ne $\sigma(h)< C'log(N)$ et $\sigma(u)< C'log(N)$. L'expression ci-dessus est essentiellement major\'ee par
  $$log(N)^RN^{-R}\int_{H(F)}{\bf 1}_{\sigma< C'log(N)}dh\int_{U(F)}{\bf 1}_{\sigma< C'log(N)}du.$$
  D'apr\`es 4.3(1) et sa preuve, il existe $R'>0$ tel que chacune de ces int\'egrales soit essentiellement major\'ee par $N^{R'}$.    On obtient
  $$\vert I_{L,{\cal O},N,C}(\theta_{\rho},f)-X_{N}\vert <<log(N)^RN^{-R+R'}.$$
  Le r\'eel $R'$ est ind\'ependant de $R$. On choisit $R=R'+2$  et on obtient que $\vert I_{L,{\cal O},N,C}(\theta_{\rho},f)-X_{N}\vert $  tend vers $0$ quand $N$ tend vers l'infini. Cela nous ram\`ene \`a calculer $lim_{N\to \infty}X_{N}$.
  
  On utilise le lemme 6.8 qui calcule les fonctions $\Phi_{j}$. On obtient une expression de $X_{N}$  que l'on peut au moins formellement \'ecrire
  $$X_{N}=[i{\cal A}_{{\cal O}}^{\vee}:i{\cal A}_{L,F}^{\vee}]^{-1}\sum_{L'\in {\cal L}(L)}(-1)^{a_{L'}}\sum_{\lambda\in \Lambda^{L'}_{{\cal O},ell}/(i{\cal A}_{L,F}^{\vee}+i{\cal A}_{L'}^*)}\vert R^{L'}(\tau_{\lambda})\vert 2^{a_{L'}-a_{L}}\sum_{\zeta\in R^{L'}(\tau_{\lambda})^{\vee}}$$
 $$\int_{i{\cal A}^*_{L',F}} X_{N}(L',\lambda,\mu,\zeta)J_{L'}^G(Ind_{L'\cap Q}^{L'}(\tau_{\lambda+\mu},\zeta),f)d\mu,$$
 o\`u
 $$X_{N}(L',\lambda,\mu,\zeta)=\sum_{j=1,...,n}\varphi_{j}(\lambda+\mu)\int_{H(F)U(F)_{c}}{\bf 1}_{\sigma< Clog(N)}(\rho(h)\epsilon'_{j},\epsilon_{j})$$
 $$(proj_{\lambda,\zeta}\circ R_{Q(S')\vert Q}(\tau_{\lambda+\mu})e'_{j},proj_{\lambda,\zeta}\circ R_{Q(S')\vert Q}(\tau_{\lambda+\mu})\circ Ind_{Q}^G(\tau_{\lambda+\mu}(hu))e_{j})\bar{\xi}(u)du\,dh.$$
 Cette expression est essentiellement major\'ee ind\'ependamment de $N$, $\lambda$ et $\mu$ par
 $$\int_{H(F)U(F)_{c}}\Xi^H(h)\Xi^G(hu)du\,dh$$
 et on sait que cette int\'egrale est convergente d'apr\`es 4.3(4). Cela justifie le calcul formel que l'on a fait ci-dessus et nous permet en m\^eme temps de calculer la limite quand $N$ tend vers l'infini. On obtient
 $$(1) \qquad lim_{N\to \infty}X_{N}=[i{\cal A}_{{\cal O}}^{\vee}:i{\cal A}_{L,F}^{\vee}]^{-1}\sum_{L'\in {\cal L}(L)}(-1)^{a_{L'}}\sum_{\lambda\in \Lambda^{L'}_{{\cal O},ell}/(i{\cal A}_{L,F}^{\vee}+i{\cal A}_{L'}^*)}\vert R^{L'}(\tau_{\lambda})\vert 2^{a_{L'}-a_{L}}$$
 $$\sum_{\zeta\in R^{L'}(\tau_{\lambda})^{\vee}}\int_{i{\cal A}^*_{L',F}} X(L',\lambda,\mu,\zeta)J_{L'}^G(Ind_{L'\cap Q}^{L'}(\tau_{\lambda+\mu},\zeta),f)d\mu,$$
  o\`u
  $$X(L',\lambda,\mu,\zeta)=\sum_{j=1,...,n}\varphi_{j}(\lambda+\mu)\int_{H(F)U(F)_{c}} (\rho(h)\epsilon'_{j},\epsilon_{j})$$
  $$(proj_{\lambda,\zeta}\circ R_{Q(S')\vert Q}(\tau_{\lambda+\mu})e'_{j},proj_{\lambda,\zeta}\circ R_{Q(S')\vert Q}(\tau_{\lambda+\mu})\circ Ind_{Q}^G(\tau_{\lambda+\mu}(hu))e_{j})\bar{\xi}(u)du\,dh.$$
Fixons $L'$, $\lambda$, $\mu$ et $\zeta$. Posons $\pi'=Ind_{L'\cap Q}^{L'}(\tau_{\lambda})$, $\pi'(\zeta)=Ind_{L'\cap Q}^{L'}(\tau_{\lambda},\zeta)$, $\tilde{e}'_{j}=R_{Q(S')\vert Q}(\tau_{\lambda+\mu})e'_{j}$, $\tilde{e}_{j}=R_{Q(S')\vert Q}(\tau_{\lambda+\mu})e_{j}$. On r\'ecrit
$$(2) \qquad X(L',\lambda,\mu,\zeta)=\sum_{j=1,...,n}\varphi_{j}(\lambda+\mu)\int_{H(F)U(F)_{c}} (\rho(h)\epsilon'_{j},\epsilon_{j})$$
$$(proj_{\lambda,\zeta}\tilde{e}'_{j},Ind_{S'}^G(\pi'(\zeta)_{\mu},hu)proj_{\lambda,\zeta}\tilde{e}_{j})\bar{\xi}(u)du\,dh.$$
On reconna\^{\i}t
$$X(L',\lambda,\mu,\zeta)=\sum_{j=1,...,n}\varphi_{j}(\lambda+\mu){\cal L}_{Ind_{S'}^G(\pi'(\zeta)_{\mu}),\rho,c}(\epsilon'_{j}\otimes proj_{\lambda,\zeta}\tilde{e}'_{j},\epsilon_{j}\otimes proj_{\lambda,\zeta}\tilde{e}_{j}).$$
On a fix\'e $c$, mais on peut le supposer assez grand. Puisque $L'$, $\lambda$ et $\zeta$ ne parcourent que des ensembles finis, la preuve du lemme 3.5 nous permet de remplacer ci-dessus ${\cal L}_{Ind_{S'}^G(\pi'(\zeta)_{\mu}),\rho,c}$ par ${\cal L}_{Ind_{S'}^G(\pi'(\zeta)_{\mu}),\rho}$. D'apr\`es le lemme 5.3(ii), ${\cal L}_{Ind_{S'}^G(\pi'(\zeta)_{\mu}),\rho}$ est non nul si et seulement si $m(\pi'(\zeta),\rho)=1$. Quand $\zeta$ parcourt $R^{L'}(\tau_{\lambda})^{\vee}$, les repr\'esentations $\pi'(\zeta)$ parcourent les diff\'erentes composantes irr\'eductibles de $\pi'$. D'apr\`es le lemme 5.4, il y a un unique $\zeta$ tel que $m(\pi'(\zeta),\rho)=1$. Notons $\zeta_{\lambda,\rho}$ cet \'el\'ement. On obtient

(3) si $\zeta\not=\zeta_{\lambda,\rho}$, $X(L',\lambda,\mu,\zeta)=0$.

Ce r\'esultat entra\^{\i}ne
$$X(L',\lambda,\mu,\zeta_{\lambda,\rho})=\sum_{\zeta\in R^{L'}(\tau_{\lambda})^{\vee}}X(L',\lambda,\mu,\zeta).$$
Gr\^ace \`a (2), cette derni\`ere somme est \'egale \`a
$$\sum_{j=1,...,n}\varphi_{j}(\lambda+\mu)\int_{H(F)U(F)_{c}} (\rho(h)\epsilon'_{j},\epsilon_{j})( \tilde{e}'_{j},Ind_{S'}^G(\pi'_{\mu},hu)\tilde{e}_{j})\bar{\xi}(u)du\,dh.$$
Puisque les op\'erateurs $R_{Q(S')\vert Q}(\tau_{\lambda+\mu})$ sont unitaires, on a
$$( \tilde{e}'_{j},Ind_{S'}^G(\pi'_{\mu},hu)\tilde{e}_{j})=(e'_{j},Ind_{Q}^G(\tau_{\lambda+\mu},hu)e_{j}).$$
De nouveau, on reconna\^{\i}t
$$X(L',\lambda,\mu,\zeta_{\lambda,\rho})=\sum_{j=1,...,n}\varphi_{j}(\lambda+\mu){\cal L}_{Ind_{Q}^G(\tau_{\lambda+\mu}),\rho,c}(\epsilon'_{j}\otimes e'_{j},\epsilon_{j}\otimes e_{j}).$$
On a ${\cal L}_{Ind_{Q}^G(\tau_{\lambda+\mu}),\rho,c}={\cal L}_{Ind_{Q}^G(\tau_{\lambda+\mu}),\rho}$ et on a justement choisi les donn\'ees $\varphi_{j}$, $\epsilon'_{j}$,$\epsilon_{j}$, $e'_{j}$ et $e_{j}$ pour que l'expression ci-dessus soit \'egale \`a $1$. Donc

(4) $X(L',\lambda,\mu,\zeta_{\lambda,\rho})=1$.

Gr\^ace \`a (3) et (4), la formule (1) devient
 $$lim_{N\to \infty}X_{N}=[i{\cal A}_{{\cal O}}^{\vee}:i{\cal A}_{L,F}^{\vee}]^{-1}\sum_{L'\in {\cal L}(L)}(-1)^{a_{L'}}\sum_{\lambda\in \Lambda^{L'}_{{\cal O},ell}/(i{\cal A}_{L,F}^{\vee}+i{\cal A}_{L'}^*)} \vert R^{L'}(\tau_{\lambda})\vert 2^{a_{L'}-a_{L}}$$
 $$\int_{i{\cal A}^*_{L',F}}  J_{L'}^G(Ind_{L'\cap Q}^{L'}(\tau_{\lambda+\mu},\zeta_{\lambda,\rho}),f)d\mu,$$
 ce qui est l'\'egalit\'e de l'\'enonc\'e. $\square$
 
 \bigskip
 
 \subsection{Preuve du th\'eor\`eme 6.1}
 
 Les lemmes 6.4 et 6.9 prouvent que $I_{N}(\theta_{\rho},f)$ a une limite quand $N$ tend vers l'infini et ils calculent cette limite. En intervertissant les sommations sur $L$ et $L'$, on a
 $$(1) \qquad lim_{N\to \infty}I_{N}(\theta_{\rho},f)=\sum_{L'\in {\cal L}(M_{min})}(-1)^{a_{L'}}\vert W^G\vert ^{-1}X(L'),$$
 o\`u
 $$X(L')=\sum_{L\in {\cal L}^{L'}(M_{min})}\vert W^L\vert \sum_{{\cal O}\in \{\Pi_2(L)\}_{f}; m({\cal O},\rho)=1}[i{\cal A}_{{\cal O}}^{\vee}:i{\cal A}_{L,F}^{\vee}]^{-1} \sum_{ \lambda\in \Lambda^{L'}_{{\cal O},ell} /(i{\cal A}_{L,F}^{\vee}+i{\cal A}_{L'}^*)}$$
  $$\vert R^{L'}(\tau_{\lambda})\vert 2^{a_{L'}-a_{L}}\sum_{\zeta\in R^{L'}(\tau_{\lambda})^{\vee}; m(Ind_{L'\cap Q}^{L'}(\tau_{\lambda},\zeta),\rho)=1} \int_{i{\cal A}_{L',F}^*}J_{L'}^G(Ind_{L'\cap Q}^{L'}(\tau_{\lambda+\mu},\zeta),f)d\mu.$$
 Fixons $L'$. Dans la formule ci-dessus, on peut remplacer la somme sur ${\cal O}\in \{\Pi_2(L)\}_{f}$ tel que $m({\cal O},\rho)=1$ par une somme sur ${\cal O}\in \{\Pi_2(L)\}$ tout entier. En effet, si ${\cal O}\in \{\Pi_2(L)\}$ v\'erifie $m({\cal O},\rho)=0$, la somme en $\zeta$ est vide d'apr\`es les propositions 5.2 et 5.7. Si ${\cal O}\not\in \{\Pi_2(L)\}_{f}$, les fonctions $J_{L'}^G(Ind_{L'\cap Q}^{L'}(\tau_{\lambda+\mu},\zeta),f)$ sont nulles.
Consid\'erons l'ensemble ${\cal Z}$ des quadruplets $z=(L,{\cal O},\lambda,\zeta)$ tels que $L\in {\cal L}^{L'}(M_{min})$, ${\cal O}\in \{\Pi_2(L)\}$, $\lambda\in \Lambda^{L'}_{{\cal O},ell} /(i{\cal A}_{L,F}^{\vee}+i{\cal A}_{L'}^*)$ et $\zeta\in R^{L'}(\tau_{\lambda})^{\vee}$. L'application
$$\iota:z=(L,{\cal O},\lambda,\zeta)\mapsto \{(Ind_{L'\cap Q}^{L'}(\tau_{\lambda},\zeta))_{\mu}; \mu\in i{\cal A}_{L'}^*\}$$
est une surjection de ${\cal Z}$ sur $\{\Pi_{ell}(L')\}$. Notons ${\cal Z}_{\rho}$ le sous-ensemble des $(L,{\cal O},\lambda,\zeta)$ tels que $m(Ind_{L'\cap Q}^{L'}(\tau_{\lambda},\zeta),\rho)=1$. Alors ${\cal Z}_{\rho}$ est l'image r\'eciproque par $\iota$ du sous-ensemble des ${\cal O}'\in \{\Pi_{ell}(L')\}$ tels que $m({\cal O}',\rho)=1$. On a donc
$$(2) \qquad X(L')=\sum_{{\cal O}'\in \{\Pi_{ell}(L')\}; m({\cal O}',\rho)=1} c({\cal O}')\int_{i{\cal A}_{L',F}^*}J_{L'}^G(\pi'_{\mu},f)d\mu,$$
o\`u, pour tout ${\cal O}'$, on a fix\'e un \'el\'ement $\pi'\in {\cal O}'$, et
$$c({\cal O}')=\sum_{z=(L,{\cal O},\lambda,\zeta)\in {\cal Z}; \iota(z)={\cal O}'} \vert W^L\vert [i{\cal A}_{{\cal O}}^{\vee}:i{\cal A}_{L,F}^{\vee}]^{-1}\vert R^{L'}(\tau_{\lambda})\vert 2^{a_{L'}-a_{L}}.$$
Fixons ${\cal O}'$. Notons ${\cal Z}'$ l'ensemble des quadruplets $z'=(L,{\cal O},\lambda,\zeta)$ tels que  $L\in {\cal L}^{L'}(M_{min})$, ${\cal O}\in \{\Pi_2(L)\}$, $\lambda\in \Lambda^{L'}_{{\cal O},ell} /i{\cal A}_{L,F}^{\vee}$  et $\zeta\in R^{L'}(\tau_{\lambda})^{\vee}$, qui v\'erifient l'\'egalit\'e $\pi'=Ind_{L'\cap Q}^{L'}(\tau_{\lambda},\zeta)$. Les  projections de $\Lambda^{L'}_{{\cal O},ell} /i{\cal A}_{L,F}^{\vee}$  sur $\Lambda^{L'}_{{\cal O},ell} /(i{\cal A}_{L,F}^{\vee}+i{\cal A}_{L'}^*)$ induisent une application de ${\cal Z}'$ dans ${\cal Z}$. Son image est pr\'ecis\'ement l'ensemble des $z\in {\cal Z}$ tels que $\iota(z)={\cal O}'$. Deux \'el\'ements de ${\cal Z}'$ ont m\^eme image par cette application si et seulement s'ils sont de la forme $(L,{\cal O},\lambda,\zeta)$, $(L,{\cal O},\lambda+\mu,\zeta)$, avec $ \mu\in i{\cal A}_{L'}^*$. Pour que les deux \'el\'ements appartiennent \`a ${\cal Z}'$, il faut et il suffit que $\pi'_{\mu}=\pi'$, autrement dit $\mu\in i{\cal A}_{{\cal O}'}^{\vee}$. La fibre de l'application au-dessus de l'image de nos deux \'el\'ements a donc le m\^eme nombre d'\'el\'ements qu'une orbite par translation par $i{\cal A}_{{\cal O}'}^{\vee}$ dans $i{\cal A}_{L,F}^*$. Puisque $i{\cal A}_{L,F}^{\vee}\cap i{\cal A}_{L'}^*=i{\cal A}_{L',F}^{\vee}$, ce nombre d'\'el\'ements est constant, \'egal \`a $[i{\cal A}_{{\cal O}'}^{\vee}/i{\cal A}_{L',F}^{\vee}]$. On obtient
$$c({\cal O}')=[i{\cal A}_{{\cal O}'}^{\vee}/i{\cal A}_{L',F}^{\vee}]^{-1}\sum_{z'=(L,{\cal O},\lambda,\zeta)\in {\cal Z}'}\vert W^L\vert [i{\cal A}_{{\cal O}}^{\vee}:i{\cal A}_{L,F}^{\vee}]^{-1}\vert R^{L'}(\tau_{\lambda})\vert 2^{a_{L'}-a_{L}}.$$
Fixons un \'el\'ement $z'=(L,{\cal O},\lambda,\zeta)\in {\cal Z}'$. Consid\'erons un autre \'el\'ement $\tilde{z}=(\tilde{L},\tilde{{\cal O}},\tilde{\lambda},\tilde{\zeta})\in {\cal Z}'$. Alors 
$$(3) \qquad Ind_{L'\cap \tilde{Q}}^{L'}(\tilde{\tau}_{\tilde{\lambda}},\tilde{\zeta}) =\pi'=Ind_{L'\cap Q}^{L'}(\tau_{\lambda},\zeta).$$
 Donc les induites $Ind_{L'\cap \tilde{Q}}^{L'}(\tilde{\tau}_{\tilde{\lambda}})$ et $Ind_{L'\cap Q}^{L'}(\tau_{\lambda})$ ont une composante irr\'eductible commune. D'apr\`es un r\'esultat d'Harish-Chandra, il existe $w\in W^{L'}$ tel que $wLw^{-1}=\tilde{L}$, $w{\cal O}=\tilde{{\cal O}}$ et $w(\tau_{\lambda})=\tilde{\tau}_{\tilde{\lambda}}$. D'o\`u les \'egalit\'es
$$\vert W^{\tilde{L}}\vert [i{\cal A}_{\tilde{{\cal O}}}^{\vee}:i{\cal A}_{\tilde{L},F}^{\vee}]^{-1}=\vert W^L\vert [i{\cal A}_{{\cal O}}^{\vee}:i{\cal A}_{L,F}^{\vee}]^{-1},$$
$$\vert R^{L'}(\tau_{\lambda})\vert =\vert R^{L'}(\tilde{\tau}_{\tilde{\lambda}})\vert,\,\,2^{a_{L'}-a_{L}}=2^{a_{L'}-a_{\tilde{L}}}.$$
Ces deux derniers nombres ne sont autres que $r(\pi')$ et $t(\pi')^{-1}$. D'o\`u 
$$c({\cal O}')=[i{\cal A}_{{\cal O}'}^{\vee}:i{\cal A}_{L',F}^{\vee}]^{-1} \vert W^L\vert [i{\cal A}_{{\cal O}}^{\vee}:i{\cal A}_{L,F}^{\vee}]^{-1}r(\pi')t(\pi')^{-1}\vert {\cal Z}'\vert .$$
Inversement, pour $w\in W^{L'}$, d\'efinissons $\tilde{L}$, $\tilde{{\cal O}}$  par les deux premi\`eres \'egalit\'es pr\'ec\'edentes et notons $\Lambda(w)$ l'ensemble des $\tilde{\lambda}\in i{\cal A}_{\tilde{L},F}^*$ tels que la troisi\`eme soit v\'erifi\'ee. Pour $\tilde{\lambda}\in \Lambda(w)$, les induites ci-dessus ont les m\^emes composantes irr\'eductibles et il existe un unique $\tilde{\zeta}$ tel que (3) soit v\'erifi\'e. L'ensemble ${\cal Z}'$ appara\^{\i}t comme l'image d'une application d'ensemble de d\'epart
$$ \{(w,\tilde{\lambda}); w\in W^{L'}, \tilde{\lambda}\in \Lambda(w)\}.$$
Tous les ensembles $\Lambda(w)$ ont m\^eme nombre d'\'el\'ements, qui est \'egal \`a $[i{\cal A}_{\tilde{{\cal O}}}^{\vee}:i{\cal A}_{\tilde{L},F}^{\vee}]$. D'autre part, deux \'el\'ements $(w_{1},\tilde{\lambda}_{1})$ et $(w_{2},\tilde{\lambda}_{2})$ ont m\^eme image dans ${\cal Z}'$ si et seulement si $\tilde{\lambda}_{1}=\tilde{\lambda}_{2}$, $w_{1}Lw_{1}^{-1}=w_{2}Lw_{2}^{-1}$, $w_{1}{\cal O}=w_{2}{\cal O}$ et $w_{1}(\tau_{\lambda})=w_{2}(\tau_{\lambda})$. Ces derni\`eres conditions sont \'equivalentes \`a ce que l'\'el\'ement$w=w_{2}^{-1}w_{1}$ conserve $L$ et ait pour image dans $W^{L'}(L)$ un \'el\'ement de $W^{L'}(\tau_{\lambda})$. Puisque $Ind_{L'\cap Q}^{L'}(\tau_{\lambda})$ a une composante elliptique, ce dernier groupe n'est autre que $R^{L'}(\tau_{\lambda})$. Autrement dit, les fibres de l'application pr\'ec\'edente ont pour nombre d'\'el\'ements 
$$\vert W^L\vert \vert R^{L'}(\tau_{\lambda})\vert=r(\pi')\vert W^L\vert  .$$
 Donc
  $$\vert {\cal Z}'\vert =\vert W^{L'}\vert \vert W^L\vert ^{-1}[i{\cal A}_{\tilde{{\cal O}}}^{\vee}:i{\cal A}_{\tilde{L},F}^{\vee}]r(\pi')^{-1},$$
  puis
  $$c({\cal O}')=[i{\cal A}_{{\cal O}'}^{\vee}:i{\cal A}_{L',F}^{\vee}]^{-1} \vert W^{L'}\vert t(\pi')^{-1}.$$
  En reportant  cette valeur dans (2) puis (1), on obtient l'\'egalit\'e du th\'eor\`eme. $\square$

 \bigskip

\section{Une formule int\'egrale calculant la multiplicit\'e; application}

\bigskip

\subsection{Le th\'eor\`eme principal}

Soient $(V,q_{V})$ et $(W,q_{W})$ deux espaces quadratiques compatibles. On plonge le plus petit dans le plus grand comme en 4.2 et on utilise les notations de ce paragraphe. Soient $\pi\in Temp(G)$ et $\rho\in Temp(H)$. On a d\'efini le nombre $m(\rho,\pi)$. Supposons $d_{V}>d_{W}$. On d\'efinit un autre nombre $m_{geom}(\rho,\pi)$  comme en [W1] 13.1, c'est-\`a-dire par l'\'egalit\'e
$$m_{geom}(\rho,\pi)=\sum_{T\in {\cal T}}\vert W(H,T)\vert ^{-1}\int_{T(F)}c_{\check{\rho}}(t)c_{\pi}(t)D^H(t)\Delta(t)^rdt.$$
Les ingr\'edients de cette formule ont \'et\'e d\'efinis en [W1] 7.3. Le terme $\nu(T)$ de [W1] 13.1 dispara\^{\i}t car nous utilisons ici la mesure sur $T$ d\'efinie en 1.2. On a mis un $\check{\rho}$ dans la formule parce que c'est le terme qui intervient naturellement mais  c'est inessentiel car

(1) on a les \'egalit\'es $m(\check{\rho},\check{\pi})=m(\rho,\pi)=m(\check{\rho},\pi)=m(\rho,\check{\pi})$ et $m_{geom}(\check{\rho},\check{\pi})=m_{geom}(\rho,\pi)=m_{geom}(\check{\rho},\pi)=m_{geom}(\rho,\check{\pi})$.

En effet, choisissons un \'el\'ement $\gamma$ du groupe orthogonal de $W$ tel que $det(\gamma)=-1$. On d\'efinit la repr\'esentation $\rho^{\gamma}$ par $\rho^{\gamma}(h)=\rho(\gamma h\gamma^{-1})$. Il est bien connu que $\rho^{\gamma}$ est \'equivalente \`a $\check{\rho}$. En particulier, si $d_{W}$ est impair, on peut choisir pour $\gamma$ la multiplication par $-1$. Cet \'el\'ement commute \`a $H$ et on obtient $\check{\rho}=\rho$. Les m\^emes propri\'et\'es valent pour $\pi$. Puisque l'un des deux nombres $d_{V}$ ou $d_{W}$ est impair, on voit qu'il suffit de prouver les deux premi\`eres \'egalit\'es de chaque s\'erie. Soit $\gamma$ comme ci-dessus. On peut consid\'erer $\gamma$ comme un \'el\'ement du groupe orthogonal de $V$ qui agit par l'identit\'e sur l'orthogonal de $W$ dans $V$. On a aussi $\check{\pi}\simeq \pi^{\gamma}$. On v\'erifie que $Hom_{H,\xi}(\pi^{\gamma},\rho^{\gamma})=Hom_{H,\xi}(\pi,\rho)$, d'o\`u la premi\`ere \'egalit\'e. L'ensemble ${\cal T}$ est un ensemble de repr\'esentants des classes  de conjugaison par $H(F)$ dans un ensemble de tores $\underline{{\cal T}}$. Ce dernier est stable par conjugaison par $\gamma$. Pour d\'emontrer la seconde \'egalit\'e, il suffit de prouver que, pour $T\in \underline{{\cal T}}$ et pour un \'el\'ement $t\in T(F)$ en position g\'en\'erale, on a l'\'egalit\'e $c_{\rho^{\gamma}}(t)=c_{\rho}(\gamma t\gamma^{-1})$ et une \'egalit\'e similaire pour la repr\'esentation $\pi$. Le terme $c_{\rho}(t)$ est le coefficient associ\'e \`a une certaine orbite, notons-la ici ${\cal O}^t\in Nil( \mathfrak{h}_{t})$ dans le d\'eveloppement du caract\`ere $\theta_{\rho}$ au voisinage de $t$. La relation \`a prouver est imm\'ediate, pourvu que l'image par la conjugaison par $\gamma$ de ${\cal O}^t$ soit \'egale \`a ${\cal O}^{\gamma t\gamma^{-1}}$. On le v\'erifie sur la d\'efinition de ces orbites, cf. [W1] 7.3. En fait, l'argument implicite est que ces orbites sont r\'eguli\`eres et que toute orbite nilpotente r\'eguli\`ere de l'alg\`ebre de Lie d'un groupe sp\'ecial orthogonal est conserv\'ee par le groupe orthogonal tout entier.

Si maintenant $d_{W}>d_{V}$, on pose $m_{geom}(\rho,\pi)=m_{geom}(\pi,\rho)$.  

\ass{Th\'eor\`eme}{Pour tout $\pi\in Temp(G)$ et tout $\rho\in Temp(H)$, on a l'\'egalit\'e $m(\rho,\pi)=m_{geom}(\rho,\pi)$.}

La preuve sera donn\'ee dans les paragraphes 7.7 \`a 7.9.

\bigskip

\subsection{Multiplicit\'es g\'eom\'etriques pour les quasi-caract\`eres et induction}

Soient $(V,q_{V})$ et $(W,q_{W})$ deux espaces quadratiques compatibles.   Soient $\rho\in Temp(H)$ et $\theta$ un quasi-caract\`ere de $G(F)$. Supposons d'abord $d_{V}>d_{W}$. On pose
$$m_{geom}(\rho,\theta)=m_{geom}(\theta,\rho)=\sum_{T\in {\cal T}}\vert W(H,T)\vert ^{-1}\int_{T(F)}c_{\check{\rho}}(t)c_{\theta}(t)D^H(t)\Delta(t)^rdt,$$
o\`u $c_{\theta}$ a \'et\'e d\'efini en [W1] 7.3. Supposons maintenant  $d_{V}<d_{W}$. On inverse les r\^oles de $V$ et $W$ en posant $V'=W$, $W'=V$ et en introduisant les objets relatifs au couple $(V',W')$, que l'on affecte d'un $'$. On pose
$$m_{geom}(\rho,\theta)=m_{geom}(\theta,\rho)=\sum_{T\in {\cal T}'}\vert W(G,T)\vert ^{-1}\int_{T(F)}c'_{\theta}(t)c'_{\check{\rho}}(t)D^G(t)\Delta(t)^{r'}dt.$$
Remarquons que, pour $\pi\in Temp(G)$, on a $m_{geom}(\rho,\pi)=m_{geom}(\rho,\theta_{\pi})$.

Soit $k\geq1$ un entier, posons $G=GL_{k}$ et soit $\theta$ un quasi-caract\`ere sur $G(F)$. Dans $\mathfrak{g}(F)$, il y a une unique orbite nilpotente r\'eguli\`ere, notons-la ${\cal O}_{GL_{k}}$. On pose
$$m_{geom}(\theta)=c_{\theta,{\cal O}_{GL_{k}}}(1).$$

Soient $(V,q_{V})$ et $(W,q_{W})$ comme ci-dessus et $L$ un L\'evi de $G$. Ecrivons 
$$L=GL_{k_{1}}\times ...\times GL_{k_{s}}\times \tilde{G},$$
o\`u $\tilde{G}$ est le groupe sp\'ecial orthogonal d'un sous-espace quadratique $\tilde{V}$ de $V$. Soient $\rho\in Temp(H)$, $\tilde{\theta}$ un quasi-caract\`ere sur $\tilde{G}(F)$ et, pour $j=1,...,s$, $\theta_{j}$ un quasi-caract\`ere sur $GL_{k_{j}}(F)$. Posons $\theta^L=\theta_{1}\otimes...\otimes \theta_{s}\otimes \tilde{\theta}$. Rappelons que les espaces quadratiques $\tilde{V}$ et $W$ sont compatibles. On peut donc poser
$$m_{geom}(\rho,\theta^L)=m_{geom}(\theta_{1})...m_{geom}(\theta_{s})m_{geom}(\rho,\tilde{\theta}).$$

{\bf Remarque.} La d\'efinition se g\'en\'eralise \'evidemment \`a tout quasi-caract\`ere sur $L(F)$ qui est combinaison lin\'eaire de quasi-caract\`eres comme ci-dessus. En fait, elle se g\'en\'eralise \`a tout quasi-caract\`ere sur $L(F)$, mais nous n'en aurons pas besoin.

On d\'efinit le quasi-caract\`ere induit $\theta=Ind_{L}^G(\theta^L)$.

\ass{Lemme}{Sous ces hypoth\`eses, on a l'\'egalit\'e $m_{geom}(\rho,\theta)=m_{geom}(\rho,\theta^L)$.}

 Preuve. Traitons le cas o\`u $d_{V}>d_{W}$, le cas oppos\'e \'etant similaire et en fait plus simple.
On suppose d'abord $d_{W}<d_{\tilde{V}}$. On peut alors supposer $W\subset \tilde{V}$. Consid\'erons les formules qui d\'efinissent $m_{geom}(\rho,\theta)$ et $m_{geom}(\rho,\tilde{\theta})$. En se reportant aux d\'efinitions de [W1] 7.3, on voit que les ensembles ${\cal T}$ qui y interviennent sont les m\^emes. Fixons $T\in {\cal T}$. Les fonctions $c_{\check{\rho}}$ sont aussi les m\^emes, ainsi que les fonctions $D^H$ et $\Delta$. L'entier $r$ qui intervient dans la premi\`ere formule est chang\'e en $r-k$ dans la seconde. Pour d\'emontrer l'\'egalit\'e cherch\'ee, il suffit donc de prouver que, pour $t\in T(F)$ en position g\'en\'erale, on a l'\'egalit\'e
 $$(1) \qquad c_{\theta}(t)\Delta(t)^k=c_{\tilde{\theta}}(t)\prod_{j=1,...,s} m_{geom}(\theta_{j}).$$
 Rappelons qu'\`a $T$ est attach\'e une d\'ecomposition orthogonale $W=W'\oplus W''$. L'espace $W'$ est de dimension paire et $T$ est un sous-tore maximal de $H'(F)$, o\`u $H'$ est le groupe sp\'ecial orthogonal de $V'$. De plus, $A_{T}=\{1\}$, c'est-\`a-dire que $T$ ne contient aucun sous-tore d\'eploy\'e non trivial. On note $V''$ l'orthogonal de $W'$ dans $V$ et $G''$ son groupe sp\'ecial orthogonal. Soit $t\in T(F)$ tel que toutes ses valeurs propres dans $V$ soient distinctes (donc aussi diff\'erentes de $1$). On a alors $Z_{G}(t)=T\times G''$.  Par d\'efinition, on a $c_{\theta}(t)=c_{\theta,{\cal O}}(t)$ pour une certaine orbite nilpotente r\'eguli\`ere ${\cal O}$ de $\mathfrak{g}_{t}(F)$. Utilisons le lemme 2.3 pour calculer ce terme.

Montrons que l'ensemble ${\cal X}^L(t)$ qui y intervient peut \^etre suppos\'e r\'eduit \`a $\{t\}$, autrement dit que tout \'el\'ement de $L(F)$ qui est conjugu\'e \`a $t$ par un \'el\'ement de $G(F)$ l'est par un \'el\'ement de $L(F)$. Soit $g\in G(F)$ tel que $gtg^{-1}\in L(F)$. Alors $g^{-1}A_{L}g$ est inclus dans le commutant de $t$, c'est-\`a-dire dans $T\times G''$. Sa projection dans $T$ ne peut \^etre que triviale, donc $g^{-1}A_{L}g\subset G''$. L'intersection des noyaux des op\'erateurs $a-1$ de $V$, pour $a\in A_{L}(F)$, est \'egale \`a $\tilde{V}$. L'inclusion pr\'ec\'edente entra\^{\i}ne $W'\subset g^{-1}\tilde{V}$ et $gW'\subset \tilde{V}$. Mais on a aussi $W'\subset W\subset \tilde{V}$. Les sous-espaces quadratiques $W'$ et $gW'$ \'etant isomorphes, le th\'eor\`eme de Witt entra\^{\i}ne que l'on peut choisir $\tilde{g}\in \tilde{G}(F)$ tel que $W'=\tilde{g}gW'$. Le groupe $\tilde{G}$ \'etant inclus dans $L$, cela montre que, quitte \`a changer $g$ par un \'el\'ement de sa classe $L(F)g$, on peut supposer $gW'=W'$. Alors $g$ induit des \'el\'ements $g'$ et $g''$ des groupes orthogonaux de $W'$ et $W''$. Si $g'\in H'(F)$, on peut consid\'erer $g'$ comme un \'el\'ement de $G(F)$, qui appartient en fait \`a $L(F)$. Alors $gtg^{-1}=g'tg^{_{'}-1}$ et cet \'el\'ement est conjugu\'e \`a $t$ par un \'el\'ement de $L(F)$. Si $det(g')=-1$, on remarque que l'orthogonal de $W'$ dans $\tilde{V}$ n'est pas nul. Fixons un \'el\'ement $\epsilon$ du groupe
 orthogonal de cet espace tel que $det(\epsilon)=-1$. Soit $g_{1}$ l'\'el\'ement de $G(F)$ qui agit par $g'$ sur $W'$, par $\epsilon$ sur l'orthogonal ci-dessus et par l'identit\'e sur l'orthogonal de $\tilde{V}$ dans $V$. On a encore $g_{1}\in L(F)$ et $gtg^{-1}=g_{1}tg_{1}^{-1}$, ce qui d\'emontre l'assertion.

Pour l'unique \'el\'ement $t$ de ${\cal X}^L(t)$, l'ensemble $\Gamma_{t}/G_{t}(F)$ du lemme 2.3 est r\'eduit \`a $\{1\}$ puisque $\Gamma_{t}=Z_{G}(t)(F)$ et $Z_{G}(t)$ est connexe. Le groupe $Z_{L}(t)$ est lui-aussi connexe et le lemme 2.3 se r\'eduit donc \`a l'\'egalit\'e
 $$(2) \qquad c_{\theta}(t)=c_{\theta,{\cal O}}(t)=D^G(t)^{-1/2}D^L(t)^{1/2}c_{\theta^L,{\cal O}^L}(t),$$
 o\`u ${\cal O}^L$ est l'unique \'el\'ement de $Nil(\mathfrak{l})$ tel que $[{\cal O}:{\cal O}^L]=1$. Pour calculer les deux premiers facteurs, on peut passer \`a la cl\^oture alg\'ebrique et supposer $T$ inclus dans un tore maximal $A_{0}$ de $G$. Le rapport $D^G(t)D^L(t)^{-1}$ est la valeur absolue du produit des $\alpha(t)-1$ sur les racines $\alpha$ de $A_{0}$ dans $G$ qui ne sont pas dans $L$ et qui sont telles que $\alpha(t)\not=1$. Notons $t_{1},...,t_{n},t_{n}^{-1},...,t_{1}^{-1}$ les valeurs propres de $t$ dans $W'$. La description habituelle des racines montre que le produit ci-dessus est \'egal \`a
 $$(\prod_{i=1,...,n}(t_{i}-1)(t_{i}^{-1}-1))^{2k}.$$
 En comparant avec la d\'efinition de la fonction $\Delta$, on obtient
 $$(3) \qquad D^G(t)^{-1/2}D^L(t)^{1/2}=\Delta(t)^{-k}.$$
 Supposons $d_{V}$ impair ou $d_{V}=2$. Alors ${\cal O}$ est l'unique orbite nilpotente r\'eguli\`ere. L'orbite ${\cal O}^L$ est aussi l'unique orbite nilpotente r\'eguli\`ere. Elle se d\'ecompose en la somme des uniques orbites nilpotentes r\'eguli\`eres ${\cal O}_{GL_{k_{j}}}$ des $\mathfrak{gl}_{k_{j}}(F)$ et de l'unique orbite r\'eguli\`ere $\tilde{{\cal O}}$ de $\tilde{\mathfrak{g}}(F)$. On a l'\'egalit\'e 
 $$(4) \qquad c_{\theta^L,{\cal O}^L}(t)= c_{\tilde{\theta},\tilde{{\cal O}}}(t)\prod_{j=1,...,s}c_{\theta_{j},{\cal O}_{GL_{k_{j}}}}(1).$$
  Le terme $c_{\tilde{\theta},\tilde{{\cal O}}}(t)$ n'est autre que $c_{\tilde{\theta}}(t)$ tandis que les termes $c_{\theta_{j},{\cal O}_{GL_{k_{j}}}}(1)$ ne sont autres que $m_{geom}(\theta_{j})$. Alors l'\'egalit\'e (2) devient l'\'egalit\'e (1) cherch\'ee. Supposons $d_{V}$ pair et $d_{V}\geq4$. Alors ${\cal O}$ est l'orbite param\'etr\'ee en [W1] 7.1 par $\nu_{0}$, o\`u $x\mapsto \nu_{0}x^2$ est le noyau anisotrope de l'orthogonal de $W$ dans $V$. Quand on remplace $V$ par $\tilde{V}$, ce $\nu_{0}$ ne change pas. On v\'erifie sur la d\'efinition des param\'etrages que ${\cal O}^L$ est somme des ${\cal O}_{GL_{k_{j}}}$ et, ou bien de l'orbite $\tilde{{\cal O}}$ param\'etr\'ee par $\nu_{0}$ si $d_{\tilde{V}}\geq4$, ou bien de l'unique orbite r\'eguli\`ere $\tilde{{\cal O}}$ si $d_{\tilde{V}}\leq 2$. On conclut comme pr\'ec\'edemment.

 On suppose maintenant que $d_{\tilde{V}}<d_{W}$. On peut supposer $\tilde{V}\subset W$. La formule de d\'efinition de $m_{geom}(\rho,\theta)$ s'\'ecrit encore
 $$m_{geom}(\rho,\theta)=\sum_{T\in {\cal T}}m_{T}(\rho,\theta),$$
 o\`u
 $$m_{T}(\rho,\theta)=\vert W(H,T)\vert ^{-1}\int_{T(F)}c_{\check{\rho}}(t)c_{\theta}(t)D^H(t)\Delta(t)^r dt.$$
 Dans celle qui d\'efinit $m_{geom}(\rho,\tilde{\theta})$, on se rappelle qu'il faut inverser les r\^oles de $H$ et $\tilde{G}$ puisque $d_{\tilde{V}}<d_{W}$.  D'o\`u
  $$m_{geom}(\rho,\tilde{\theta})=\sum_{T\in \tilde{{\cal T}}}m_{T}(\rho,\tilde{\theta}),$$
 o\`u
 $$m_{T}(\rho,\tilde{\theta})=\vert W(\tilde{G},T)\vert ^{-1}\int_{T(F)}\tilde{c}_{\check{\rho}}(t)c_{\tilde{\theta}}(t)D^{\tilde{G}}(t)\Delta(t)^{\tilde{r}}dt.$$
 L'ensemble ${\cal T}$ est un ensemble de repr\'esentants de classes de conjugaison par $H(F)$ dans un ensemble de tores $\underline{{\cal T}}$. De m\^eme, $\tilde{{\cal T}}$ est ensemble de repr\'esentants de classes de conjugaison par $\tilde{G}(F)$ dans un ensemble de tores $\tilde{\underline{{\cal T}}}$.  Montrons que
 
 (5) $\tilde{\underline{{\cal T}}}$ est l'ensemble des $T\in \underline{\cal T}$ qui sont inclus dans $\tilde{G}$.
 
 Soit $T$ un \'el\'ement de l'un ou l'autre de ces ensembles. On a une d\'ecomposition orthogonale $\tilde{V}=\tilde{V}'\oplus \tilde{V}''$ o\`u $d_{\tilde{V}'}$ est paire et $T$ est inclus dans le groupe sp\'ecial orthogonal $\tilde{G}'$ de $\tilde{V}'$. De plus $A_{T}=\{1\}$. Notons $W''$, resp. $V''$, l'orthogonal de $\tilde{V}'$ dans $W$, resp. $V$. Supposons $d_{V}$ pair. Alors $T$ appartient au deuxi\`eme ensemble si et seulement si $d_{an,W''}=1$. Puisque $d_{W}$ est impair, $T$ appartient au premier ensemble si et seulement si $d_{an,W''}=1$. Ces conditions  sont les m\^emes. Supposons $d_{V}$ impair. Alors $T$ appartient au deuxi\`eme ensemble si et seulement si $d_{an,V''}=1$. Il appartient au premier si et seulement si $d_{an,\tilde{V}''}=1$. Mais $V''$ est la somme orthogonale de $\tilde{V}''$ et de l'orthogonal de $\tilde{V}$ dans $V$. Ce dernier espace est hyperbolique. Donc $d_{an,V''}=d_{an,\tilde{V}''}$ et nos deux conditions sont encore \'equivalentes.
 
 Soit $T\in {\cal T}$. Montrons que
 
 (6) si $T$ n'est pas conjugu\'e \`a un \'el\'ement de $\tilde{\underline{\cal T}}$ par un \'el\'ement de $H(F)$, on a $m_{T}(\rho,\theta)=0$.
 
  Notons $W=W'\oplus W''$ la d\'ecomposition attach\'ee \`a $T$ et soit $t\in T(F)$ en position g\'en\'erale. Comme dans la premi\`ere partie de la preuve, $c_{\theta}(t)$ se calcule gr\^ace au lemme 2.3. Pour que ce terme soit non nul, il faut qu'il existe $g\in G(F)$ tel que $gtg^{-1}\in L(F)$. Comme plus haut, cette relation entra\^{\i}ne $gW'\subset \tilde{V}$, donc $gW'\subset W$. Notons $W''_{1}$ l'orthogonal de $gW'$ dans $W$. Remarquons qu'il n'est pas nul puisque $d_{\tilde{V}}<d_{W}$. Les deux sous-espaces quadratiques $gW'$ et $W'$ de $W$ sont isomorphes. Donc leurs orthogonaux $W''_{1}$ et $W''$ le sont aussi. Fixons un isomorphisme $\gamma$ de $W''$ sur $W''_{1}$. Notons $h$ l'automorphisme de $W$ qui agit par $g$ sur $W'$ et par $\gamma$ sur $W''$. C'est un \'el\'ement du groupe orthogonal de $W$. Quitte \`a multiplier $\gamma$ \`a gauche par un \'el\'ement du groupe orthogonal de $W''_{1}$ de d\'eterminant $-1$ (un tel \'el\'ement existe puisque $W''_{1}$ n'est pas nul), on peut supposer $h\in H(F)$. La d\'ecomposition associ\'ee \`a $hTh^{-1}$ est $gW'\oplus W''_{1}$. Puisque $gW'$ est inclus dans $\tilde{V}$,  $hTh^{-1}$ appartient \`a $\tilde{\underline{\cal T}}$ contrairement \`a l'hypoth\`ese. Cette contradiction montre que la fonction $c_{\theta}$ est nulle en un point g\'en\'eral de $T(F)$. D'o\`u la conclusion.
 
 Introduisons dans $\tilde{\underline{\cal T}}$ une nouvelle relation d'\'equivalence: deux \'el\'ements sont $H$-\'equivalents si et seulement s'ils sont conjugu\'es par un \'el\'ement de $H(F)$. Fixons un ensemble de repr\'esentants ${\cal T}^{\sharp}$ des classes de $H$-\'equivalence. On peut supposer ${\cal T}^{\sharp}\subset \tilde{{\cal T}}$. Les relations (5) et (6) montrent que, dans la formule exprimant $m_{geom}(\rho,\theta)$, on peut remplacer la somme sur $T\in {\cal T}$ par la somme sur $T\in {\cal T}^{\sharp}$. Fixons $T\in {\cal T}^{\sharp}$. Notons $\tilde{{\cal T}}_{T}$ l'ensemble des \'el\'ements de $\tilde{{\cal T}}$ qui sont $H$-\'equivalents \`a $T$. Il suffit de prouver que l'on a l'\'egalit\'e
 $$(7) \qquad m_{T}(\rho,\theta)=\sum_{\tilde{T}\in \tilde{{\cal T}}_{T}}m_{\tilde{T}}(\rho,\tilde{\theta})\prod_{j=1,...,s}m_{geom}(\theta_{i}).$$
On note $\tilde{V}=\tilde{V}'\oplus \tilde{V}''$ la d\'ecomposition attach\'ee \`a $T$, $W''$ l'orthogonal de $\tilde{V}'$ dans $W$, $\tilde{G}'$, $\tilde{G}''$ et $H''$ les groupes sp\'eciaux orthogonaux de $\tilde{V}'$, $\tilde{V}''$ et $W''$. Pour tout groupe sp\'ecial orthogonal, par exemple $G$, on note $G^+$ le groupe orthogonal correspondant. On s\'epare trois cas. Le cas (I) est celui o\`u $\tilde{V}'=\{0\}$ ou $\tilde{V}''\not=\{0\}$. Le cas (II) est celui o\`u $\tilde{V}'\not=\{0\}$, $\tilde{V}''=\{0\}$ et il existe un \'el\'ement $\gamma'\in Norm_{\tilde{G}^+(F)}(T)$ de d\'eterminant $-1$. Le cas (III) est celui o\`u $\tilde{V}'\not=\{0\}$, $\tilde{V}''=\{0\}$ et il n'existe pas de tel \'el\'ement $\gamma'$. Dans les cas (II) et (III),  on d\'efinit un \'el\'ement $\gamma$ de la fa\c{c}on suivante.  On fixe un \'el\'ement $\gamma'$  de $\tilde{G}^+(F)$ de d\'eterminant $-1$. Dans le cas (II), on suppose que $\gamma'$ appartient \`a $Norm_{\tilde{G}^+(F)}(T)$. On fixe un \'el\'ement $\gamma''$ de d\'eterminant $-1$ de $H^{_{''}+}(F)$ et on note $\gamma$ l'\'el\'ement de $H(F)$ qui agit par $\gamma'$ sur $\tilde{V}$ et par $\gamma''$ sur $W''$. Montrons que

(8)(i) dans le cas (I), $\tilde{{\cal T}}_{T}$ est r\'eduit \`a $\{T\}$ et on a $\vert W(H,T)\vert =\vert W(\tilde{G},T)\vert $;

(8)(ii) dans le cas (II), $\tilde{{\cal T}}_{T}$ est r\'eduit \`a $\{T\}$ et on a $\vert W(H,T)\vert =2\vert W(\tilde{G},T)\vert $;

(8)(iii) dans le cas (III), $\tilde{{\cal T}}_{T}$ est r\'eduit \`a deux \'el\'ements  et on peut supposer que $\tilde{{\cal T}}_{T}=\{T,\gamma T\gamma^{-1}\}$; on a $\vert W(H,T)\vert =\vert W(\tilde{G},T)\vert =\vert W(\tilde{G},\gamma T\gamma^{-1})\vert $.

Notons $W_{1}$ l'orthogonal de $\tilde{V}$ dans $W$, $H_{1}$ son groupe sp\'ecial orthogonal et posons $\Gamma=(\tilde{G}^+\times H_{1}^+)\cap H$. D\'ecrivons l'ensemble des \'el\'ements de $\tilde{\underline{\cal T}}$ qui sont $H$-\'equivalents \`a $T$. Ce sont les \'el\'ements de cet ensemble qui sont de la forme $hTh^{-1}$ pour un $h\in H(F)$. La m\^eme preuve qu'en (6) montre que l'on peut supposer que $h$ conserve $\tilde{V}$. Autrement dit $h\in  \Gamma(F)$. Inversement, si $h \in \Gamma(F)$, $hTh^{-1}$ appartient \`a notre ensemble de tores. L'application $h\mapsto hTh^{-1}$ se descend donc en une bijection de 
 $$ \Gamma(F)/(\Gamma(F)\cap Norm_{H(F)}(T))$$
 sur cet ensemble de tores. L'ensemble $\tilde{{\cal T}}_{T}$ est un ensemble de repr\'esentants des classes de conjugaison par $\tilde{G}(F)$ dans l'ensemble pr\'ec\'edent. Il est donc en bijection avec
 $$\tilde{G}(F)\backslash  \Gamma(F)/(\Gamma(F)\cap Norm_{H(F)}(T)).$$
 Si $\tilde{V}'=\{0\}$, $T=\{1\}$ et les assertions de (8)(i) sont \'evidentes. Supposons $\tilde{V}'\not=\{0\}$.
 L'ensemble $\Gamma$ a deux composantes connexes. Sa composante neutre est $\tilde{G}\times H_{1}$, qui est inclus dans $\tilde{G}Norm_{H}(T)$ car $H_{1}$ est inclus dans $Norm_{H}(T)$ (et m\^eme dans $Z_{H}(T)$). L'ensemble ci-dessus a donc un \'el\'ement si $Norm_{H(F)}(T)$ coupe les deux composantes de $\Gamma$, et deux \'el\'ements sinon. On a $\Gamma(F)\cap Norm_{H(F)}(T)=(Norm_{\tilde{G}^+(F)}(T)\times H_{1}^+(F))\cap H(F)$. Donc $Norm_{H(F)}(T)$ coupe les deux composantes de $\Gamma$ si et seulement si $Norm_{\tilde{G}^+(F)}(T)$ coupe les deux composantes de $\tilde{G}^+$. Si $\tilde{V}''\not=\{0\}$, cette propri\'et\'e est v\'erifi\'ee puisque $\tilde{G}^{_{''}+}(F)\subset Norm_{\tilde{G}^+(F)}(T)$. Elle l'est aussi dans le cas (II) par d\'efinition de l'\'el\'ement $\gamma'$ tandis qu'elle ne l'est pas dans le cas (III) par d\'efinition de ce cas. Remarquons que, dans ce cas (III), l'\'el\'ement $\gamma$ que l'on a d\'efini appartient \`a $\Gamma(F)$ et pas \`a sa composante neutre. On en d\'eduit les premi\`eres assertions de (8). 
 
   Soit $h\in Norm_{H(F)}(T)$. N\'ecessairement, $h$ conserve $\tilde{V}'$. Comme dans la preuve de (6), on montre qu'en le multipliant \`a droite par un \'el\'ement de $H''(F)$ (ce groupe est inclus dans $Z_{H(F)}(T)$), on peut supposer que $h$ conserve $\tilde{V}$, donc appartient \`a $\Gamma(F)\cap Norm_{H(F)}(T)$. Le quotient $W(H,T)$ est donc \'egal \`a
 $$(\Gamma(F)\cap Norm_{H(F)}(T))/(\Gamma(F)\cap Z_{H(F)}(T)).$$
 Posons 
 $$ \Delta=(\Gamma_{F}\cap Norm_{H(F)}(T))/\left((\Gamma(F)\cap Z_{H(F)}(T))(\Gamma^0(F)\cap Norm_{H(F)}(T))\right).$$
 Il y a un homomorphisme surjectif de $W(H,T)$ sur $\Delta$, de noyau  \'egal \`a 
 $$(\Gamma^0(F)\cap Norm_{H(F)}(T))/(\Gamma^0(F)\cap Z_{H(F)}(T)).$$
Mais $H_{1}$ est contenu dans $Z_{H}(T)$, donc ce quotient n'est autre que $Norm_{\tilde{G}(F)}(T)/Z_{\tilde{G}(F)}(T)$, autrement dit $W(\tilde{G},T)$. On obtient que $\vert W(H,T)\vert $ est \'egal \`a $\vert W(\tilde{G},T)\vert \vert \Delta\vert $  . Si $\tilde{V}''\not=\{0\}$, $Z_{H(F)}(T)$ coupe les deux composantes de $\Gamma$ car $\tilde{G}^{_{''}+}$ centralise $T$. Donc le groupe $\Delta$ n'a qu'un \'el\'ement.  Dans le cas (II), $Norm_{H(F)}(T)$ coupe les deux composantes de $\Gamma$ tandis que $Z_{H}(T)=T\times H''$ est inclus dans la composante neutre. Alors $\Delta$ a deux \'el\'ements. Dans le cas (III), $Norm_{H(F)}(T)$ est inclus dans la composante neutre de $\Gamma$ et $\Delta$ n'a qu'un \'el\'ement. D'o\`u les derni\`eres assertions de (8).
 
 Soit $t\in T(F)$ dont toutes les valeurs propres dans $\tilde{V}$ soient distinctes. Comme dans la premi\`ere partie de la preuve, on a $c_{\theta}(t)=c_{\theta,{\cal O}}(t)$ pour une certaine orbite nilpotente r\'eguli\`ere ${\cal O}$ de $\mathfrak{g}_{t}(F)$. On calcule ce terme en utilisant le lemme 2.3. Montrons que
 
 (9) on peut choisir pour ensemble ${\cal X}^L(t)$ l'ensemble $\{t\}$ dans le cas (I) et l'ensemble $\{t,\gamma t\gamma^{-1}\}$ dans les cas (II) et (III).
 
 C'est clair si $T=\{1\}$. Supposons $t\not=\{1\}$, ce qui \'equivaut \`a $\tilde{V}'\not=\{0\}$. Soit $g\in G(F)$ tel que $gtg^{-1}\in L(F)$. Comme dans la premi\`ere partie de la preuve, cela entra\^{\i}ne $g\tilde{V}'\subset \tilde{V}$. Quitte \`a multiplier $g$ \`a droite par un \'el\'ement de $\tilde{G}(F)$, on peut supposer que $g$ conserve $\tilde{V}'$. Notons $g'$ la restriction de $g$ \`a $\tilde{V}'$. Si $\tilde{V}''\not=\{0\}$, soit $g''$ un \'el\'ement de $\tilde{G}^{_{''}+}(F)$ tel que $det(g')det(g'')=1$. Soit $g_{1}$ l'\'el\'ement de $G(F)$ qui agit par $g'$ sur $\tilde{V}'$, par $g''$ sur $\tilde{V}''$ et par l'identit\'e sur l'orthogonal de $\tilde{V}$ dans $V$. Alors $gtg^{-1}=g_{1}tg_{1}^{-1}$ et $g_{1}$ appartient \`a $L(F)$. D'o\`u le r\'esultat. Supposons maintenant $\tilde{V}''=\{0\}$. Si $g'\in \tilde{G}'(F)$, la m\^eme construction s'applique et $gtg^{-1}$ est conjugu\'e \`a $t$ par un \'el\'ement de $L(F)$. Si $det(g')=-1$, on construit de m\^eme un \'el\'ement $g_{1}\in L(F)$ qui agit sur $\tilde{V}'$ par $\gamma'g^{_{'}-1}$. Alors $g_{1}gtg^{-1}g_{1}^{-1}=\gamma t\gamma^{-1}$. On peut donc supposer que ${\cal X}^{L}(t)$ est inclus dans $\{t,\gamma t\gamma^{-1}\}$. D'autre part, les deux \'el\'ements de cet ensemble ne sont pas conjugu\'es par un \'el\'ement de $L(F)$. Sinon, il existerait $\tilde{g}\in \tilde{G}(F)=\tilde{G}'(F)$ tel que $\tilde{g}t\tilde{g}^{-1}=\gamma' t\gamma^{_{'}-1}$ et $\tilde{g}^{-1}\gamma'$ serait un \'el\'ement de $\tilde{G}^{_{'}+}(F)$ de d\'eterminant  $-1$ et commutant \`a $T$. Un tel \'el\'ement n'existe pas, d'o\`u l'assertion et (9).

Achevons la preuve en supposant que l'on est dans le cas (II) et en laissant les autres cas au lecteur.  L'ensemble $\Gamma_{t}/G_{t}(F)$ du lemme 2.3 est r\'eduit \`a $\{1\}$ et l'ensemble $\Gamma_{\gamma t\gamma^{-1}}/G_{t}(F)$ est r\'eduit \`a la classe de $\gamma$. Les groupes $Z_{L}(t)$ et $Z_{L}(\gamma t\gamma^{-1})$ sont connexes. De plus $\gamma{\cal O}={\cal O}$.  Le lemme 2.3 entra\^{\i}ne
  $$c_{\theta}(t)=c_{\theta,{\cal O}}(t)=D^G(t)^{-1/2}(D^L(t)^{1/2}c_{\theta^L,{\cal O}^L}(t)+D^{L}(\gamma t\gamma^{-1})^{1/2}c_{\theta^L,{\cal O}^L}(\gamma t\gamma^{-1})),$$
  o\`u ${\cal O}^L$ est l'unique orbite nilpotente de $\mathfrak{l}_{t}(F)=\mathfrak{l}_{\gamma t\gamma^{-1}}(F)$ telle que $[{\cal O}:{\cal O}^L]=1$. D'o\`u
  $$m_{T}(\rho,\theta)=\vert W(H,T)\vert ^{-1}(\int_{T(F)}c_{\check{\rho}}(t)c_{\theta^L,{\cal O}^L}(t)D^G(t)^{-1/2}D^L(t)^{1/2}D^H(t)\Delta(t)^rdt$$
  $$+ \int_{T(F)}c_{\check{\rho}}(t)c_{\theta^L,{\cal O}^L}(\gamma t\gamma^{-1})D^G(t)^{-1/2}D^L(\gamma t\gamma^{-1})^{1/2}D^H(t)\Delta(t)^rdt).$$
  Dans la seconde int\'egrale, on effectue le changement de variables $t\mapsto \gamma^{-1}t\gamma$. Il laisse invariantes les fonctions $c_{\check{\rho}}$, $D^G$, $D^H$ et $\Delta$. La seconde int\'egrale est donc \'egale \`a la premi\`ere.  Gr\^ace \`a (8)(ii), on obtient
   $$(10) \qquad m_{T}(\rho,\theta)=\vert W(\tilde{G},T)\vert ^{-1}\int_{T(F)}c_{\check{\rho}}(t)c_{\theta^L,{\cal O}^L}(t)D^G(t)^{-1/2}D^L(t)^{1/2}D^H(t)\Delta(t)^rdt.$$
Soit $t$ comme pr\'ec\'edemment. On a encore l'\'egalit\'e  (4) et $c_{\theta_{j},{\cal O}_{GL_{k_{j}}}}(1)=m_{geom}(\theta_{j})$ pour tout $j$. On v\'erifie que $\tilde{{\cal O}}$ est pr\'ecis\'ement l'orbite qui sert \`a d\'efinir $c_{\tilde{\theta}}(t)$. Traitons seulement le cas le plus subtil o\`u $d_{V}$ est pair et $d_{\tilde{V}}\geq4$. D'apr\`es [W1] 7.3, ${\cal O}$ est l'orbite param\'etr\'ee par $\nu_{0}$, o\`u la forme $x\mapsto \nu_{0}x^2$ est le noyau anisotrope de l'orthogonal de $W$ dans $V$. L'orbite servant \`a d\'efinir $c_{\tilde{\theta}}(t)$ est param\'etr\'ee par $-\tilde{\nu}_{0}$, o\`u $x\mapsto \tilde{\nu}_{0}x^2$ est le noyau anisotrope de l'orthogonal de $\tilde{V}$ dans $W$. L'orbite $\tilde{{\cal O}}$ est, comme ${\cal O}$, param\'etr\'ee par $\nu_{0}$. Or $\nu_{0}=-\tilde{\nu}_{0}$ (dans $F^{\times}/F^{\times 2}$) car l'orthogonal de $\tilde{V}$ dans $V$ est hyperbolique. D'o\`u l'assertion. On obtient 
$$c_{\theta^L,{\cal O}^L}(t)=c_{\tilde{\theta}}(t)\prod_{j=1,...,s}m_{geom}(\theta_{j}).$$
 Un calcul analogue montre que $c_{\check{\rho}}(t)=\tilde{c}_{\check{\rho}}(t)$. On a encore l'\'egalit\'e (3). Comme dans la preuve de cette \'egalit\'e, on v\'erifie que
$$D^{\tilde{G}}(t)=D^{\tilde{G}'}(t)\Delta(t)^{d_{\tilde{V}}-d_{\tilde{V}'}}$$
et
$$D^H(t)=D^{\tilde{G}'}(t)\Delta(t)^{d_{W}-d_{\tilde{V}'}}.$$
D'o\`u
$$D^G(t)^{-1/2}D^L(t)^{1/2}D^H(t)\Delta(t)^r=D^{\tilde{G}}(t)\Delta(t)^{a},$$
o\`u $a=r-k+d_{W}-d_{\tilde{V}}$. Par d\'efinition, $d_{V}-d_{W}=2r+1$, $d_{W}-d_{\tilde{V}}=2\tilde{r}+1$ et $d_{V}-d_{\tilde{V}}=2k$. Alors $a=\tilde{r}$. Mais alors, le membre de droite de l'\'egalit\'e (10) n'est autre que $m_{T}(\rho,\tilde{\theta})\prod_{j=1,...,s}m_{geom}(\theta_{j})$. Gr\^ace \`a (8)(ii), on obtient l'\'egalit\'e (5), ce qui ach\`eve la d\'emonstration. $\square$

\bigskip

\subsection{Fonctions cuspidales sur les groupes sp\'eciaux orthogonaux}

Soit $(V,q_{V})$  un espace quadratique.     Rappelons qu'en 2.5, on a associ\'e un quasi-caract\`ere $I\theta_{f}$ \`a toute fonction cuspidale $f\in C_{c}^{\infty}(G(F))$.

\ass{Lemme}{ Soit $f\in C_{c}^{\infty}(G(F))$ une fonction cuspidale. On a l'\'egalit\'e
$$I\theta_{f}=\sum_{\pi\in \Pi_{ell}(G)}t(\check{\pi})^{-1}\theta_{\check{\pi}}(f)\theta_{\pi}.$$}

Preuve. Soit $f\in C_{c}^{\infty}(G(F))$ une fonction cuspidale. Reprenons la formule 2.5(1). Elle se simplifie puisque $A_{G}=\{1\}$ et devient
$$ I\theta_{f}=\sum_{\pi\in T_{ell}(G)}c(\pi)\theta_{\check{\pi}}(f)\theta_{\pi}.$$
D\'ecrivons l'ensemble $T_{ell}(G)$ et les constantes $c(\pi)$. Consid\'erons l'ensemble des couples $(L,\tau)$ tels que $L\in {\cal L}(M_{min})$, $\tau$ est une repr\'esentation admissible irr\'eductible de $L(F)$ de la s\'erie discr\`ete et $R(\tau)\cap W(L)_{reg}\not=\emptyset$. On d\'efinit de fa\c{c}on \'evidente la notion de conjugaison de tels couples et on fixe un ensemble de repr\'esentants $\tilde{T}_{ell}(G)$ des classes de conjugaison. Soit $(L,\tau)\in \tilde{T}_{ell}(G)$.  Comme on l'a dit en 4.1,  l'ensemble $R(\tau)\cap W(L)_{reg}$ a un unique \'el\'ement. Notons $t$  cet \'el\'ement. Pour tout $\zeta\in R(\tau)^{\vee}$, on d\'efinit la repr\'esentation elliptique $\pi(\zeta)=Ind_{Q}^G(\tau,\zeta) $ de $G(F)$, o\`u $Q$ est un \'el\'ement fix\'e de ${\cal P}(L)$. On d\'efinit la repr\'esentation virtuelle 
$$(1) \qquad \pi=\sum_{\zeta\in R(\tau)^{\vee}}\zeta(t)\pi(\zeta).$$
Alors $T_{ell}(G)$ est l'ensemble des ces repr\'esentations virtuelles quand $(L,\tau)$ d\'ecrit $\tilde{T}_{ell}(G)$. Le nombre $c(\pi)$ associ\'e \`a la repr\'esentation $\pi$ ci-dessus est $\vert R(\tau)\vert ^{-1}\vert det(t-1)_{\vert {\cal A}_{L}}\vert ^{-1}$. Rappelons que les repr\'esentations temp\'er\'ees elliptiques de $G(F)$  sont exactement les repr\'esentations $\pi(\zeta)$ introduites ci-dessus quand $(L,\tau)$ d\'ecrit $\tilde{T}_{ell}(G)$ et que l'on a $r(\pi(\zeta))=r(\pi(\zeta)\check{})=\vert R(\tau)\vert $ et $t(\pi(\zeta))=t(\pi(\zeta)\check{})=\vert det(t-1)_{{\cal A}_{L}}\vert $. Pour prouver l'\'egalit\'e de l'\'enonc\'e, on peut donc fixer $(L,\tau)\in \tilde{T}_{ell}(G)$ et prouver l'\'egalit\'e
$$(2) \qquad \sum_{\zeta\in R(\tau)^{\vee}}\theta_{\pi(\zeta)\check{}}(f)\theta_{\pi(\zeta)}=\vert R(\tau)\vert ^{-1}\theta_{\check{\pi}}(f)\theta_{\pi},$$
o\`u $\pi$ est d\'efinie par (1). Soit $r\in R(\tau)$, $r\not=t$. D'apr\`es [A5] proposition 2.1(b), la repr\'esentation virtuelle
$$\sum_{\zeta\in R(\tau)^{\vee}}\zeta(r)\pi(\zeta)\check{}$$
est une somme de repr\'esentations  proprement induites. Puisque $f$ est cuspidale, le caract\`ere de cette repr\'esentation annule $f$. Cela \'etant vrai pour tout $r\not=t$, on en d\'eduit que $\zeta(t)\theta_{\pi(\zeta)\check{}}(f)$ est ind\'ependant de $\zeta$. Ce nombre est donc \'egal \`a
$\vert R(\tau)\vert ^{-1}\theta_{\check{\pi}}(f)$. Le membre de gauche de la relation (2) est donc \'egal \`a
$$\vert R(\tau)\vert ^{-1}\sum_{\zeta\in R(\tau)^{\vee}}\theta_{\check{\pi}}(f)\zeta(t)\theta_{\pi(\zeta)}$$
et ceci n'est autre que le membre de droite de (2). Cela d\'emontre (2) et le lemme. $\square$

\bigskip

\subsection{Pseudo-coefficients}

Soit $(V,q_{V})$ un  espace quadratique. Soient $L\in {\cal L}(M_{min})$ et $\tau $ une repr\'esentation admissible irr\'eductible de la s\'erie discr\`ete de $L(F)$. Supposons $R(\tau)\cap W(L)_{reg}\not=\emptyset$, notons $t$ l'unique \'el\'ement de cet ensemble . Pour tout $\zeta\in R(\tau)^{\vee}$,  on introduit la repr\'esentation elliptique $\pi(\zeta)$ de $G(F)$ comme dans la preuve pr\'ec\'edente. 

\ass{Lemme}{Il existe une fonction cuspidale $f\in C_{c}^{\infty}(G(F))$ telle que

(i) $I\theta_{f}=\sum_{\zeta\in R(\tau)^{\vee}}\zeta(t)\theta_{\pi(\zeta)}$;

(ii) $\theta_{\pi(\zeta)\check{}}(f)=\zeta(t)t(\pi(\zeta)\check{})$ pour tout $\zeta\in R(\tau)\check{}$;

(iii) $\theta_{\sigma\check{}}(f)=0$ pour tout $\sigma\in Temp(G)$ qui n'est pas l'une des repr\'esentations $\pi(\zeta)$.}

Preuve. On d\'efinit  $T_{ell}(G)$ et la repr\'esentation virtuelle $\pi$ comme dans la preuve pr\'ec\'edente. D'apr\`es [A5], p.94, il existe une fonction cuspidale $f\in C_{c}^{\infty}(G(F))$ telle que $\theta_{\check{\pi}}(f)=\vert R(\tau)\vert \vert det(t-1)_{\vert {\cal A}_{L}}\vert $ et $\theta_{\check{\pi}'}(f)=0$ pour tout $\pi'\in T_{ell}(G)$, $\pi'\not=\pi$. Fixons une telle fonction. Le th\'eor\`eme 5.1 de [A5] affirme pr\'ecis\'ement que l'\'egalit\'e du (i) de l'\'enonc\'e est v\'erifi\'ee ( le terme $d(\tau)^{-1}$ d'Arthur est \'egal \`a $\vert det(t-1)_{\vert {\cal A}_{L}}\vert^{-1}$ et il y a encore un $\vert R(\tau)\vert ^{-1}$ cach\'e dans la d\'efinition de la mesure $d\tau$, cf. [A5] p. 96). Comme dans la preuve pr\'ec\'edente, on montre que, pour $\zeta\in R(\tau)^{\vee}$, on a
$$\theta_{\pi(\zeta)\check{}}(f)=\vert R(\tau)\vert ^{-1}\zeta(t)\theta_{\check{\pi}}(f).$$
D'o\`u (ii) puisque $t(\pi(\zeta)\check{})=\vert det(t-1)_{\vert {\cal A}_{L}}\vert $. Soit $\sigma\in Temp(G)$ qui n'est pas l'une des repr\'esentations $\pi(\zeta)$. Il existe alors $(L',\tau')\in \tilde{T}_{ell}(G)$ et $\zeta'\in R(\tau')$ de sorte que $\sigma=\pi'(\zeta')$, avec une notation \'evidente.  Comme ci-dessus, et avec des notations analogues, on a
$$\theta_{\check{\sigma}}=\vert R(\tau')\vert ^{-1}\zeta'(t')\theta_{\check{\pi}'}(f),$$
o\`u $\pi'$ est un \'el\'ement de $T_{ell}(G)$ diff\'erent de $\pi$. Donc $\theta_{\check{\sigma}}(f)=0$ et le (iii). $\square$

\bigskip

\subsection{Une cons\'equence du th\'eor\`eme}

Soient $(V,q_{V})$ et $(W,q_{W})$ deux espaces quadratiques compatibles.   Soient $\rho\in Temp(H)$ et $f\in C_{c}^{\infty}(G(F))$ une fonction cuspidale.  Posons
$$m_{spec}(\rho,f)=\sum_{\pi\in \Pi_{ell}(G); m(\check{\rho},\pi)=1}t(\pi)^{-1}\theta_{\pi}(f).$$

\ass{Corollaire}{On suppose v\'erifi\'e le th\'eor\`eme 7.1 pour nos deux espaces quadratiques. Alors, pour toute fonction cuspidale $f\in C_{c}^{\infty}(G(F))$ et toute repr\'esentation $\rho\in Temp(H)$, on a l'\'egalit\'e $m_{geom}(\rho,I\theta_{f})=m_{spec}(\rho,f)$.}

Preuve.  Supposons $d_{V}>d_{W}$, la preuve \'etant sym\'etrique dans le cas oppos\'e. En changeant $\pi$ en $\check{\pi}$ dans la d\'efinition de $m_{spec}(\rho,f)$ et en utilisant 7.1(1), on a
$$m_{spec}(\rho,f)=\sum_{\pi\in \Pi_{ell}(G)}m(\rho,\pi)t(\check{\pi})^{-1}\theta_{\check{\pi}}(f).$$
Gr\^ace au th\'eor\`eme, c'est aussi
$$m_{spec}(\rho,f) =\sum_{\pi\in \Pi_{ell}(G)}m_{geom}(\rho,\pi)t(\check{\pi})^{-1}\theta_{\check{\pi}}(f)$$
$$=\sum_{T\in {\cal T}}\vert W(H,T)\vert ^{-1}\int_{T(F)}c_{\check{\rho}}(t)Ic_{f}(t)D^H(t)\Delta(t)^rdt,$$
o\`u
$$Ic_{f}(t)=\sum_{\pi\in \Pi_{ell}(G)}t(\check{\pi})^{-1}\theta_{\check{\pi}}(f)c_{\pi}(t).$$
Il suffit de d\'emontrer que, pour tout $T\in {\cal T}$ et presque tout $t\in T(F)$, on a l'\'egalit\'e $c_{I\theta_{f}}(t)=Ic_{f}(t)$. Ces deux fonctions se d\'eduisent par la d\'efinition de [W1] 7.3 du quasi-caract\`ere $I\theta_{f}$ pour la premi\`ere, du quasi-caract\`ere
$$ \sum_{\pi\in \Pi_{ell}(G)}t(\check{\pi})^{-1}\theta_{\check{\pi}}(f)\theta_{\pi}$$
pour la seconde. Ces deux quasi-caract\`eres sont \'egaux d'apr\`es le lemme 7.3. $\square$ 

\bigskip

\subsection{Le cas du groupe lin\'eaire}

Soit $k\geq1$ un entier et $G=GL_{k}$. Soit $f\in C_{c}^{\infty}(G(F))$ une fonction cuspidale. Posons
 $$m_{spec}(f)=\sum_{{\cal O}\in \{\Pi_{ell}(G)\}}[i{\cal A}_{{\cal O}}^{\vee}:i{\cal A}_{G,F}^{\vee}]^{-1}\theta_{\pi}(f{\bf 1}_{H_{G}=0}),$$
o\`u, comme toujours, on a fix\'e un point-base $\pi$ dans chaque orbite.

\ass{Lemme}{Pour toute fonction cuspidale $f\in C_{c}^{\infty}(G(F))$, on a l'\'egalit\'e $m_{geom}(I\theta_{f})=m_{spec}(f)$.}

Preuve. On utilise encore la formule 2.5(1). Pour le groupe $GL_{k}$, les $R$-groupes sont triviaux et les coefficients $c({\cal O})$ sont \'egaux \`a $[i{\cal A}_{{\cal O}}^{\vee}:i{\cal A}_{G,F}^{\vee}]^{-1}$.  On obtient
$$I\theta_{f}=\sum_{{\cal O}\in \{\Pi_{ell}(G)\}}[i{\cal A}_{{\cal O}}^{\vee}:i{\cal A}_{G,F}^{\vee}]^{-1}\theta_{\check{\pi}}(f{\bf 1}_{H_{G}=0})\theta_{\pi}.$$
 D'o\`u
$$m_{geom}(I\theta_{f})=\sum_{{\cal O}\in \{\Pi_{ell}(G)\}}[i{\cal A}_{{\cal O}}^{\vee}:i{\cal A}_{G,F}^*]^{-1}\theta_{\check{\pi}}(f{\bf 1}_{H_{G}=0})c_{\theta_{\pi},{\cal O}_{GL_{k}}}.$$
D'apr\`es un r\'esultat de Rodier ([R] th\'eor\`eme p.161 et remarque 2, p.162), pour toute repr\'esentation admissible irr\'eductible $\pi$ de $G(F)$, on a $c_{\theta_{\pi},{\cal O}_{GL_{k}}}=1$ si $\pi$ admet un mod\`ele de Whittaker, $0$ sinon. Or les repr\'esentations temp\'er\'ees elliptiques de $GL_{k}(F)$, qui ne sont autres que les repr\'esentations de la s\'erie discr\`ete, poss\`edent toutes un mod\`ele de Whittaker. Les termes $c_{\theta_{\pi},{\cal O}_{GL_{k}}}$ intervenant dans le membre de droite de l'\'egalit\'e ci-dessus sont  tous \'egaux \`a $1$. On peut changer la somme sur $\pi$ en une somme sur $\check{\pi}$ et ce membre de droite   devient $m_{spec}(f)$. $\square$

\bigskip

\subsection{D\'ebut de la preuve; le cas o\`u $\pi$ est induite}

On d\'emontre le th\'eor\`eme par r\'ecurrence sur $sup(d_{V},d_{W})$. On suppose d\'esormais fix\'es deux espaces quadratiques compatibles $(V,q_{V})$ et $(W,q_{W})$, avec $d_{V}>d_{W}$ et on suppose le th\'eor\`eme v\'erifi\'e pour tout couple d'espaces quadratiques $(V',q_{V'})$, $(W',q_{W'})$ compatibles et tels que $sup(d_{V'},d_{W'})<d_{V}$. Le cas o\`u $V$ est de dimension $2$ et $q_{V}$ est hyperbolique est imm\'ediat. On exclut ce cas.

 Pour $\rho\in Temp(H)$, on prolonge par lin\'earit\'e les applications $\pi\mapsto m(\rho,\pi)$ et $\pi\mapsto m_{geom}(\rho,\pi)$ \`a l'espace des combinaisons lin\'eaires finies \`a coefficients complexes d'\'el\'ements de $Temp(G)$.

\ass{Lemme}{Soient $\pi$ une repr\'esentation temp\'er\'ee de $G(F)$ et $\rho\in Temp(H)$. Supposons $\pi$ proprement induite. Alors on a l'\'egalit\'e $m(\rho,\pi)=m_{geom}(\rho,\pi)$.}

Preuve. On peut trouver
 
 -un L\'evi $L=GL_{k}\times \tilde{G}$ de $G$, o\`u $k\geq1$ et $\tilde{G}$ est le groupe sp\'ecial orthogonal d'un sous-espace $\tilde{V}$ de $V$;
 
 - un \'el\'ement $Q\in {\cal P}(L)$;
 
 - des repr\'esentations admissibles irr\'eductibles et temp\'er\'ees $\mu$ de $GL_{k}(F)$ et $\tilde{\pi}$ de $\tilde{G}(F)$,
 
 \noindent de sorte que $\pi=Ind_{Q}^G(\mu\otimes\tilde{\pi})$.   Si $m(\rho,\tilde{\pi})=0$, on a $m(\rho,\pi')=0$ pour toute composante irr\'eductible $\pi'$ de $\pi$ d'apr\`es les propositions 5.3 et 5.8. Donc $m(\rho,\pi)=0$. Si $m(\rho,\tilde{\pi})=1$, les m\^emes propositions 5.3 et 5.8 et le lemme 5.5 montrent qu'il existe une unique composante irr\'eductible $\pi'$ de $\pi$ telle que $m(\rho,\pi')=1$. Donc $m(\rho,\pi)=1$. Dans les deux cas,   $m(\rho,\pi)=m(\rho,\tilde{\pi})$. D'autre part, d'apr\`es le lemme 7.2 et les d\'efinitions, on a l'\'egalit\'e 
 $$m_{geom}(\rho,\pi)=m_{geom}(\rho,\tilde{\pi})m_{geom}(\theta_{\mu}).$$
 D'apr\`es le m\^eme r\'esultat de Rodier que l'on a utilis\'e en 7.6, on a $m_{geom}(\theta_{\mu})=1$ puisque $\mu$ est temp\'er\'ee et irr\'eductible.  D'apr\`es l'hypoth\`ese de r\'ecurrence, on a l'\'egalit\'e $m(\rho,\tilde{\pi})=m_{geom}(\rho,\tilde{\pi})$.   La conclusion s'ensuit. $\square$

\bigskip

\subsection{Comparaison de deux limites}

On conserve la situation du paragraphe pr\'ec\'edent et on reprend les notations du paragraphe 7.5.

\ass{Proposition}{Pour toute fonction cuspidale $f\in C_{c}^{\infty}(G(F))$ et toute repr\'esentation $\rho\in Temp(H)$, on a l'\'egalit\'e $m_{geom}(\rho,I\theta_{f})=m_{spec}(\rho,f).$}

Preuve. Les deux membres ne d\'ependent de $f$ qu'\`a \'equivalence pr\`es. D'apr\`es le lemme 2.7, on peut supposer $f$ tr\`es cuspidale. Le th\'eor\`eme 6.1, appliqu\'e \`a $\check{\rho}$, calcule la limite quand $N$ tend vers l'infini de $I_{N}(\theta_{\check{\rho}},f)$. On a not\'e cette limite $I_{spec}(\theta_{\check{\rho}},f)$. Mais  le th\'eor\`eme 7.8 de [W1], appliqu\'e \`a $\theta=\theta_{\check{\rho}}$, calcule la m\^eme limite. On l'a not\'e $I(\theta_{\check{\rho}},f)$. Avec les d\'efinitions de 7.2, c'est simplement $m_{geom}(\rho,\theta_{f})$.   On a donc
$$I_{spec}(\theta_{\check{\rho}},f)=m_{geom}(\rho,\theta_{f}).$$
 
{\bf Remarque.} Dans [W1], on avait utilis\'e d'autres mesures. On utilise ici celles que l'on a d\'efinies en 1.2.  Le quasi-caract\`ere $\theta_{f}$  est normalis\'e comme en 2.6.

Reprenons la d\'efinition de 6.1. On a
$$I_{spec}(\theta_{\check{\rho}},f)=\sum_{L\in {\cal L}(M_{min})}\vert W^L\vert \vert W^G\vert ^{-1}(-1)^{a_{L}}I_{spec,L}(\theta_{\check{\rho}},f),$$
o\`u
$$I_{spec,L}(\theta_{\check{\rho}},f)=\sum_{{\cal O}\in\{\Pi_{ell}(L)\}; m({\cal O},\check{\rho})=1}[i{\cal A}_{{\cal O}}^{\vee}:i{\cal A}_{L,F}^{\vee}]^{-1}t(\pi)^{-1}\int_{i{\cal A}_{L,F}^*}J_{L}^G(\pi_{\lambda},f)d\lambda.$$
D'apr\`es la d\'efinition de 7.5, on a $I_{spec,G}(\theta_{\check{\rho}},f)=m_{spec}(\rho,f)$. Soit $L\in {\cal L}(M_{min})$, $L\not=G$. Introduisons la fonction $f_{L}=\phi_{L}(f){\bf 1}_{H_{L}=0}$. Alors
$$I_{spec,L}(\theta_{\check{\rho}},f)=\sum_{{\cal O}\in\{\Pi_{ell}(L)\}; m({\cal O},\check{\rho})=1}[i{\cal A}_{{\cal O}}^{\vee}:i{\cal A}_{L,F}^{\vee}]^{-1}t(\pi)^{-1}\theta_{\pi}(f_{L}).$$
Ecrivons 
$$L=GL_{k_{1}}\times...\times GL_{k_{s}}\times \tilde{G}.$$
La fonction $f_{L}$ est cuspidale d'apr\`es le lemme 2.6(i) et l'espace des fonctions cuspidales sur $L(F)$ est le produit tensoriel des espaces des fonctions cuspidales sur chacun des facteurs de $L(F)$. Sur l'espace des fonctions cuspidales sur $\tilde{G}(F)$, on a d\'efini en 7.5 une forme lin\'eaire $f'\mapsto m_{spec}(\rho,f')$. Sur l'espace des fonctions cuspidales sur un facteur $GL_{k_{j}}(F)$, on a d\'efini en 7.6 une forme lin\'eaire $f'\mapsto m_{spec}(f')$. Notons $f'\mapsto m_{spec}(\rho,f')$ le produit tensoriel des ces formes lin\'eaires. Montrons que
$$(1) \qquad I_{spec,L}(\theta_{\check{\rho}},f)=m_{spec}(\rho,f_{L}).$$
 On a
 $$\{\Pi_{ell}(L)\}=\{\Pi_{ell}(GL_{k_{1}})\}\times...\times \{\Pi_{ell}(GL_{k_{s}})\}\times \Pi_{ell}(\tilde{G}).$$
 Pour ${\cal O}={\cal O}_{1}\times...\times{\cal O}_{s}\times \tilde{\pi}\in \{\Pi_{ell}(L)\}$, on a
 $$i{\cal A}_{{\cal O}}^{\vee}=i{\cal A}_{{\cal O}_{1}}^{\vee}\oplus...\oplus i{\cal A}_{{\cal O}_{s}}^{\vee},$$
 $$t(\pi)=t(\tilde{\pi}) \text{ et } m({\cal O},\check{\rho})=m(\tilde{\pi},\check{\rho}).$$
 Si $f_{L}$ est produit tensoriel de fonctions sur chaque facteur, l'\'egalit\'e (1) est donc quasiment tautologique et le cas g\'en\'eral s'en d\'eduit par lin\'earit\'e.
 
 On d\'efinit $m_{geom}(\rho,I\theta_{f_{L}})$ comme en 7.2. Ce  terme est \'egal \`a $m_{spec}(\rho,f_{L})$: il suffit d'appliquer le lemme 7.6 \`a chaque facteur $GL$ et le lemme 7.5 aux espaces $\tilde{V}$ et $W$. C'est loisible d'apr\`es l'hypoth\`ese de r\'ecurrence. A ce point, nous avons d\'emontr\'e l'\'egalit\'e
 $$(1) \qquad m_{spec}(\rho,f)= m_{geom}(\rho,\theta_{f})-\sum_{L\in {\cal L}(M_{min}); L\not=G}\vert W^L\vert \vert W^G\vert ^{-1}(-1)^{a_{L}}m_{geom}(\rho,I\theta_{f_{L}}).$$
 Gr\^ace au lemme 2.6(ii), on a l'\'egalit\'e
 $$m_{geom}(\rho,\theta_{f})=\sum_{L\in {\cal L}(M_{min})}\vert W^L\vert \vert W^G\vert ^{-1}(-1)^{a_{L}}m_{geom}(\rho,Ind_{L}^G(I\theta_{f_{L}})).$$
 Le terme index\'e par $L=G$ n'est autre que $m_{geom}(\rho,I\theta_{f})$. Pour $L\not=G$, le terme $m_{geom}(\rho,Ind_{L}^G(I\theta_{f_{L}}))$ est \'egal d'apr\`es le lemme 7.2 au terme $m_{geom}(\rho,I_{\theta_{f_{L}}})$ qui intervient dans (1). A cause du signe n\'egatif pr\'esent dans (1), ces termes disparaissent et  (1) devient l'\'egalit\'e de l'\'enonc\'e. $\square$
 
 \bigskip
 
 \subsection{Fin de la preuve}
 
 On veut prouver que $m(\rho,\pi)=m_{geom}(\rho,\pi)$ pour tout $\rho\in Temp(H)$ et tout $\pi\in Temp(G)$. Fixons $\rho$. On peut aussi bien d\'emontrer l'\'egalit\'e pr\'ec\'edente pour $\pi$ parcourant un ensemble de repr\'esentations virtuelles tel que tout \'el\'ement de $ Temp(G)$ soit combinaison lin\'eaire d'\'el\'ements de cet ensemble. La r\'eunion de $T_{ell}(G)$ et de l'ensemble des repr\'esentations temp\'er\'ees qui sont des induites propres convient. Le cas d'une telle induite est r\'egl\'e par le lemme 7.7.  Reste le cas d'un \'el\'ement $\pi$ de $T_{ell}(G)$ pour lequel on reprend les notations de 7.4. Soit $f$ v\'erifiant les conditions du lemme de ce paragraphe. D'apr\`es le (i) de ce lemme, on a
 $$m_{geom}(\rho,I\theta_{f})=m_{geom}(\rho,\pi).$$
 D'apr\`es les (ii) et (iii) du lemme, on a aussi
 $$m_{spec}(\rho,f)=\sum_{\zeta\in R(\tau)^{\vee};m(\rho,\pi(\zeta))=1}t(\pi(\zeta)\check{})^{-1}\theta_{\pi(\zeta)\check{}}(f)$$
 $$=\sum_{\zeta\in R(\tau)^{\vee}}m(\rho,\pi(\zeta))\zeta(t)=m(\rho,\pi).$$
 Alors l'\'egalit\'e voulue r\'esulte du lemme 7.8 $\square$
 
 \bigskip
 
 \subsection{Cons\'equence pour la conjecture locale de Gross-Prasad}
 
 Soient $(V_{i},q_{V_{i}}$ et $(W_{i},q_{W_{i}})$ deux espaces quadratiques compatibles tels que $d_{V_{i}}>d_{W_{i}}$. On note $G_{i}$ et $H_{i}$ leurs groupes sp\'eciaux orthogonaux et on suppose ces groupes quasi-d\'eploy\'es sur $F$. A \'equivalence pr\`es, il existe au plus un espace quadratique $(V',q_{V'})$ tel que $d_{V'}=d_{V_{i}}$ et que les discriminants de $q_{V'}$ et $q_{V_{i}}$ soient \'egaux mais leurs indices de Witt soient distincts. Si cet espace existe (ce qui est toujours le cas si $d_{V_{i}}\geq3$) on le note $(V_{a},q_{V_{a}})$. On introduit de m\^eme un \'eventuel espace quadratique $(W_{a},q_{W_{a}})$. On note $G_{a}$ et $H_{a}$ leurs groupes sp\'eciaux orthogonaux. Le groupe $G_{a}$, resp. $H_{a}$, est une forme int\'erieure de $G_{i}$, resp. $H_{i}$. Les espaces quadratiques $(V_{a},q_{V_{a}})$ et $(W_{a},q_{W_{a}})$ sont compatibles.
 
 Nous admettons que les ensembles de repr\'esentations $Temp(G_{i})$, $Temp(G_{a})$, $Temp(H_{i})$ et $Temp(H_{a})$ se d\'ecomposent en unions disjointes de $L$-paquets de sorte que les propri\'et\'es (1), (2) et (3) de [W1] 13.2 soient v\'erifi\'ees. Soient $\Pi_{i}$ un $L$-paquet dans $Temp(G_{i})$ et $\Sigma_{i}$ un $L$-paquet dans $Temp(H_{i})$. Si l'espace $(V_{a},q_{V_{a}})$ existe, il peut correspondre \`a $\Pi_{i}$ un $L$-paquet  de $Temp(G_{a})$. On le note $\Pi_{a}$. Dans les autres cas, c'est-\`a-dire ou bien l'espace $(V_{a},q_{V_{a}})$ existe et aucun $L$-paquet de $Temp(G_{a})$ ne correspond \`a $\Pi_{i}$, ou bien l'espace $(V_{a},q_{V_{a}})$ n'existe pas, on pose $\Pi_{a}=\emptyset$. On d\'efinit de fa\c{c}on similaire $\Sigma_{a}$.
 
 \ass{Th\'eor\`eme}{Il existe un unique couple $(\rho,\pi)\in (\Sigma_{i}\times \Pi_{i})\cup (\Sigma_{a}\times \Pi_{a})$ tel que $m(\rho,\pi)=1$.}
 
 La preuve est la m\^eme que celle du th\'eor\`eme 13.3 de [W1]. Dans cette r\'ef\'erence, l'hypoth\`ese de cuspidalit\'e des \'el\'ements de $\Pi_{i}\cup \Pi_{a}$ ne servait qu'\`a utiliser l'\'egalit\'e $m(\rho,\pi)=m_{geom}(\rho,\pi)$ qui n'\'etait alors d\'emontr\'ee que pour $\pi$ cuspidale. Maintenant que l'on dispose de cette \'egalit\'e pour toutes les repr\'esentations $\rho$, $\pi$ temp\'er\'ees, cette hypoth\`ese ne sert plus et on obtient le r\'esultat pour tous les paquets temp\'er\'es. $\square$

\bigskip

{\bf Bibliographie}

\bigskip

[AGRS] A. Aizenbud, D. Gourevitch, S. Rallis, G. Schiffmann: {\it Multiplicity one theorems}, pr\'ep. 2007

[A1] J. Arthur: {\it The trace formula in invariant form}, Annals of Math. 114 (1981), p.1-74

[A2] -----------: {\it The invariant trace formula I. Local theory}, J. AMS 1 (1988), p. 323-383

[A3] -----------: {\it A local trace formula}, Publ. IHES 73 (1991), p.5-96

[A4] -----------: {\it Intertwining operators and residues I. Weighted characters}, J. Funct. Analysis 84 (1989), p. 19-84

[A5] -----------: {\it On elliptic tempered characters}, Acta Math. 171 (1993), p.73-138

[A6] -----------: {\it On a family of distributions obtained from Eisenstein series II: explicit formulas}, Amer. J. Math. 104 (1982), p.1289-1336

 [GGP]  W. T. Gan, B. Gross, D. Prasad:{\it Symplectic local root numbers, central critical $L$-values and restriction problems in the representation theory of classical groups}, pr\'ep. 2008
 
 [GP] B. Gross, D. Prasad: {\it On irreducible repr\'esentations of $SO_{2n+1}\times SO_{2m}$}, Can. J. Math. 46 (1994), p.930-950

[HCDeBS]: Harish-Chandra, S. DeBacker, P. Sally: {\it Admissible invariant distributions on reductive $p$-adic groups}, AMS Univ. lecture series 16 (1999)

[II] A. Ichino, T. Ikeda: {\it On the periods of automorphic forms on special orthogonal groups and the Gross-Prasad conjecture} , pr\'ep. 2008

[R] F. Rodier: {\it Mod\`ele de Whittaker et caract\`eres de repr\'esentations}, in Non commutative harmonic analysis, J. Carmona, J. Dixmier, M. Vergne \'ed. Springer LN 466 (1981), p.151-171

[W1] J.-L. Waldspurger: {\it Une formule int\'egrale reli\'ee \`a la conjecture de Gross-Prasad}, pr\'epublication 2008

[W2] -----------------------: {\it La formule de Plancherel pour les groupes $p$-adiques, d'apr\`es Harish-Chandra}, J. of the Inst. of Math. Jussieu 2 (2003), p.235-333

\bigskip
 Institut de Math\'ematiques de Jussieu-CNRS

175, rue du Chevaleret

75013 Paris

e-mail: waldspur@math.jussieu.fr

\end{document}